%% file: 1006.2597.Russian.tex
\documentclass{amsart}
\def\UseRussian{}
\input{1006.2597}

\begin{document}
\title{Производная Гато и интеграл над банаховой алгеброй}

\begin{abstract}
\ShowEq{Abstract}
\end{abstract}
\shorthandoff{"}
\ShowEq{contents}
\shorthandon{"}%
\end{document}

%% file: 1006.2597.tex
\def\BookNumber{1006.2597}
\input{Commands}
\xRefDef{0812.4763}
\xRefDef{0912.3315}
\xRefDef{1003.1544}

\DefEq
{
\maketitle
\tableofcontents
\input{\FilePrefix Preface.1006.2597.English}

\input{\FilePrefix Convention.English}

\input{\FilePrefix Linear.Mapping.English}
\input{\FilePrefix Derivative.English}
\input{\FilePrefix Second.Derivative.English}
\input{\FilePrefix Biblio.English}
\input{\FilePrefix Index.English}
\input{\FilePrefix Symbol.English}
}
{contents}

\DefEq
{
\input{\FilePrefix Abstract.1006.2597.English}
}
{Abstract}

%% file: Commands.tex
\def\Defined{}
\scrollmode
\ifx\FilePrefix\undefined
\newcommand{\FilePrefix}{}
\fi
\ifx\UseRussian\Defined
	\usepackage{cmap}
		\usepackage[T2A,T2B]{fontenc}
	\usepackage[cp1251]{inputenc}
	\usepackage[english,russian]{babel}
\fi
\ifx\Presentation\Defined
	\usepackage[margin=1cm]{geometry}
\fi
\ifx\PublishBook\Defined
	\def\PrintPaper{}
	\usepackage{setspace}
	\ifx\UseRussian\undefined
		\usepackage{pslatex}
	\fi
	\usepackage[margin=2cm]{geometry}
\fi
\usepackage{footmisc}
\usepackage[all]{xy}
\usepackage{color}
\ifx\PrintPaper\undefined
	\definecolor{UrlColor}{rgb}{.9,0,.3}
	\definecolor{SymbColor}{rgb}{.4,0,.9}
	\definecolor{IndexColor}{rgb}{1,.3,.6}
	\newcommand\BlueText[1]{\textcolor{blue}{#1}}
	
\else
	\definecolor{UrlColor}{rgb}{.1,.1,.1}
	\definecolor{SymbColor}{rgb}{.1,.1,.1}
	\definecolor{IndexColor}{rgb}{.1,.1,.1}
	\newcommand\BlueText[1]{#1}
	
\fi

\usepackage{xr-hyper}
	\usepackage[unicode]{hyperref}
	\hypersetup{pdfdisplaydoctitle=true}
	\hypersetup{colorlinks}
	\hypersetup{citecolor=UrlColor}
	\hypersetup{urlcolor=UrlColor}
	\hypersetup{linkcolor=UrlColor}
	\hypersetup{pdffitwindow=true}
	\hypersetup{pdfnewwindow=true}
	\hypersetup{pdfstartview={FitH}}
\def\hyph{\penalty0\hskip0pt\relax-\penalty0\hskip0pt\relax}
\def\Hyph{-\penalty0\hskip0pt\relax}

\def\ValueOff{off}
\def\ValueOn{on}
\def\Items#1{\ItemList#1,LastItem,}%
\def\LastItem{LastItem}%
\def\ItemList#1,{\def\ViewBook{#1}%
	\ifx\ViewBook\LastItem%
	\else%
		\ifx\ViewBook\BookNumber%
			\def\Semafor{on}%
		\fi%
		\expandafter\ItemList%
	\fi%
}%

\newcommand{\ePrints}[1]
{%
	\def\Semafor{off}%
	\Items{#1}%
}%

\newcommand{\Basis}[1]{\overline{\overline{#1}}{}}
\newcommand{\Vector}[1]{\overline{#1}{}}
\ifx\PrintPaper\undefined
	\newcommand{\gi}[1]{\boldsymbol{\textcolor{IndexColor}{#1}}}
\else
	\newcommand{\gi}[1]{\boldsymbol{#1}}
\fi

\newcommand{\VX}[1]{\Vector{#1}_{[1]}}
\makeatletter
\newcommand{\NameDef}[1]{%
	\expandafter\gdef\csname #1\endcsname%
}%
\newcommand{\ShowSymbol}[1]{%
	\@nameuse{ViewSymbol#1}%
}%
\ifx\PrintPaper\undefined
	\newcommand{\symb}[3]{%
		\@ifundefined{ViewSymbol#3}{%
			\NameDef{ViewSymbol#3}{\textcolor{SymbColor}{#1}}%
			\NameDef{RefSymbol#3}{\pageref{symbol: #3}}%
			\@namedef{LabelSymbol#3}{\label{symbol: #3}}%
		}{%
			\NameDef{RefSymbol#3}{\pageref{symbol: #3}, \pageref{symbol 1: #3}}%
			\@namedef{LabelSymbol#3}{\label{symbol 1: #3}}%
		}%
		\ifcase#2%0
		\or%1
			$\@nameuse{ViewSymbol#3}$%
		\or%2
			\[\@nameuse{ViewSymbol#3}\]%
		\else%
		\fi%
		\@nameuse{LabelSymbol#3}%
	}%
\else
	\newcommand{\symb}[3]{%
		\@ifundefined{ViewSymbol#3}{%
			\NameDef{ViewSymbol#3}{#1}%
			\NameDef{RefSymbol#3}{\pageref{symbol: #3}}%
			\@namedef{LabelSymbol#3}{\label{symbol: #3}}%
		}{%
			\message {error: extra symb #3}
		}%
		\ifcase#2%0
		\or%1
			$\@nameuse{ViewSymbol#3}$%
		\or%2
			\[\@nameuse{ViewSymbol#3}\]%
		\else%
		\fi%
		\@nameuse{LabelSymbol#3}%
	}%
\fi
\newcommand{\DefEq}[2]{%
	\@ifundefined{ViewEq#2}{%
		\NameDef{ViewEq#2}{#1}%
	}{%
	}%
}%
\newcommand{\DefEquation}[2]{%
	\DefEq
	{
	\begin{equation}
	#1
	\EqLabel{#2}
	\end{equation}
	}
	{#2}
}%
\newcommand\EqRef[1]{\eqref{eq: #1}}%
\newcommand\EqLabel[1]{\label{eq: #1}}%
\newcommand\ShowEq[1]{%
	\@ifundefined{ViewEq#1}{%
		\message {error: missed ShowEq #1}
  }{%
	\@nameuse{ViewEq#1}%
	}%
}%
\newcommand\DrawEq[2]{%
	\@ifundefined{ViewEq#1}{%
		\message {error: missed ShowEq #1}
  }{%
	\def\Temp{}
	\def\Tempa{#2}
	\ifx\Tempa\Temp
		\[%
		\@nameuse{ViewEq#1}%
		\]%
	\else
		\begin{equation}%
		\@nameuse{ViewEq#1}%
		\EqLabel{#1, #2}%
		\end{equation}%
	\fi
	}%
}%
\makeatother
\DeclareMathOperator{\rank}{\mathrm{rank}}

\newcommand{\subs}{${}_*$\Hyph}
\newcommand{\sups}{${}^*$\Hyph}

\newcommand{\CRstar}{{}^*{}_*}
\newcommand{\RCstar}{{}_*{}^*}
\newcommand{\CRcirc}{{}^{\circ}{}_{\circ}}
\newcommand{\RCcirc}{{}_{\circ}{}^{\circ}}

\newcommand{\RC}{$\RCstar$\Hyph}
\newcommand{\CR}{$\CRstar$\Hyph}
\newcommand{\drc}{$D\RCstar$\Hyph}
\newcommand{\Drc}{$\mathcal D\RCstar$\Hyph}
\newcommand{\dcr}{$D\CRstar$\hyph}
\newcommand{\rcd}{$\RCstar D$\Hyph}
\newcommand{\crd}{$\CRstar D$\Hyph}

\newcommand{\Acr}{$A\CRcirc$\Hyph}

\newcommand\sT{$\star T$\Hyph}%
\newcommand\Ts{$T\star$\Hyph}%
\newcommand\sD{$\star D$\Hyph}%
\newcommand\Ds{$D\star$\Hyph}%
\newcommand\VirtFrac{\vphantom{\overset{\rightarrow}{\frac 11}^{\frac 11}}}
\newcommand\VirtVar{\vphantom{\overset{\rightarrow}{1}^1}}
\newcommand\pC[2]{{}_{#1\cdot #2}}%
\newcommand\DrcPartial[1]%
{%
	\def\tempa{}%
	\def\tempb{#1}%
	\ifx\tempa\tempc%
		(\partial\RCstar)%
	\else%
		({}_{#1}\partial\RCstar)%
	\fi%
}%
\newcommand\crDPartial[1]%
{%
	\def\tempa{}%
	\def\tempb{#1}%
	\ifx\tempa\tempc%
		(\CRstar\partial)%
	\else%
		(\CRstar\partial_{#1})%
	\fi%
}%
\newcommand\StandPartial[3]%
{%
	\frac{\partial^{\gi{#3}} #1}{\partial #2}%
}%

	\renewcommand{\uppercasenonmath}[1]{}
\makeatletter
\newcommand\@dotsep{4.5}
\def\@tocline#1#2#3#4#5#6#7
{\relax
		\par \addpenalty\@secpenalty\addvspace{#2}%
		\begingroup %\hyphenpenalty\@M
		\@ifempty{#4}{%
			\@tempdima\csname r@tocindent\number#1\endcsname\relax
		}{%
			\@tempdima#4\relax
		}%
		\parindent\z@ \leftskip#3\relax \advance\leftskip\@tempdima\relax
		\rightskip\@pnumwidth plus1em \parfillskip-\@pnumwidth
		#5\leavevmode\hskip-\@tempdima #6\relax
		\leaders\hbox{$\m@th
		\mkern \@dotsep mu\hbox{.}\mkern \@dotsep mu$}\hfill
		\hbox to\@pnumwidth{\@tocpagenum{#7}}\par
		\nobreak
		\endgroup
}
\makeatother 

\ifx\PrintBook\undefined
	\def\Chapter{\section}
	\def\Section{\subsection}
	\makeatletter
		\renewcommand{\@indextitlestyle}{%
			\twocolumn[\section{\indexname}]%
			\def\IndexSpace{off}%
		}
	\makeatother 
		\ifx\PrintPaper\undefined
			\thanks{\href{mailto:Aleks\_Kleyn@MailAPS.org}{Aleks\_Kleyn@MailAPS.org}}
			\ePrints{1102.1776,Finsler}
			\ifx\Semafor\ValueOff
					\thanks{\ \ \ \url{http://sites.google.com/site/AleksKleyn/}}
					\thanks{\ \ \ \url{http://arxiv.org/a/kleyn\_a\_1}}
					\thanks{\ \ \ \url{http://AleksKleyn.blogspot.com/}}
			\fi
		\fi
\else

	\def\Chapter{\chapter}
	\def\Section{\section}
\makeatletter
\def\@maketitle{%
  \cleardoublepage \thispagestyle{empty}%
  \begingroup \topskip\z@skip
  \null\vfil
  \begingroup
  \LARGE\bfseries \centering
  \openup\medskipamount
  \@title\par\vspace{24pt}%
  \def\and{\par\medskip}\centering
  \mdseries\authors\par\bigskip
  \endgroup
  \vfil\vspace{24pt}
  \ifx\@empty\addresses \else \@setaddresses \fi
  \vfil
  \ifx\@empty\@dedicatory
  \else \begingroup
    \centering{\footnotesize\itshape\@dedicatory\@@par}%
    \endgroup
  \fi
  \vfill
  \newpage\thispagestyle{empty}
  \@setabstract
  \begin{center}
    \ifx\@empty\@subjclass\else\@setsubjclass\fi
    \ifx\@empty\@keywords\else\@setkeywords\fi
    \ifx\@empty\@translators\else\vfil\@settranslators\fi
    \ifx\@empty\thankses\else\vfil\@setthanks\fi
  \end{center}
  \vfil
  \endgroup}
\makeatother 

	\makeatletter
	\renewcommand{\@indextitlestyle}{%
		\twocolumn[\chapter{\indexname}]%
		\def\IndexSpace{off}%
		\let\@secnumber\@empty
		\chaptermark{\indexname}%
	}
	\makeatother 
	\email{\href{mailto:Aleks\_Kleyn@MailAPS.org}{Aleks\_Kleyn@MailAPS.org}}
	\ePrints{1102.1776,Finsler}
	\ifx\Semafor\ValueOff
		\urladdr{\url{http://sites.google.com/site/alekskleyn/}}
		\urladdr{\url{http://arxiv.org/a/kleyn\_a\_1}}
		\urladdr{\url{http://AleksKleyn.blogspot.com/}}
	\fi
\fi

\ifx\SelectlEnglish\undefined
	\ifx\UseRussian\undefined
		\def\SelectlEnglish{}
	\fi
\fi

\newcommand\arXivRef{http://arxiv.org/PS_cache/}
\newcommand\wRefDef[2]
	{
	\def\Tempa{#1}
	\def\Tempb{0405.027}
	\ifx\Tempa\Tempb
	\def\wRef{\arXivRef gr-qc/pdf/0405/0405027v3.pdf}
	\fi
	\def\Tempb{0405.028}
	\ifx\Tempa\Tempb
	\def\wRef{\arXivRef gr-qc/pdf/0405/0405028v5.pdf}
	\fi
	\def\Tempb{0412.391}
	\ifx\Tempa\Tempb
	\def\wRef{\arXivRef math/pdf/0412/0412391v4.pdf}
	\fi
	\def\Tempb{0612.111}
	\ifx\Tempa\Tempb
	\def\wRef{\arXivRef math/pdf/0612/0612111v2.pdf}
	\fi
	\def\Tempb{0701.238}
	\ifx\Tempa\Tempb
	\def\wRef{\arXivRef math/pdf/0701/0701238v4.pdf}
	\fi
	\def\Tempb{0702.561}
	\ifx\Tempa\Tempb
	\def\wRef{\arXivRef math/pdf/0702/0702561v3.pdf}
	\fi
	\def\Tempb{0707.2246}
	\ifx\Tempa\Tempb
	\def\wRef{\arXivRef arxiv/pdf/0707/0707.2246v2.pdf}
	\fi
	\def\Tempb{0803.3276}
	\ifx\Tempa\Tempb
	\def\wRef{\arXivRef arxiv/pdf/0803/0803.3276v3.pdf}
	\fi
	\def\Tempb{0812.4763}
	\ifx\Tempa\Tempb
	\def\wRef{\arXivRef arxiv/pdf/0812/0812.4763v6.pdf}
	\fi
 	\def\Tempb{0906.0135}
	\ifx\Tempa\Tempb
	\def\wRef{\arXivRef arxiv/pdf/0906/0906.0135v3.pdf}
	\fi
 	\def\Tempb{0909.0855}
	\ifx\Tempa\Tempb
	\def\wRef{\arXivRef arxiv/pdf/0909/0909.0855v5.pdf}
	\fi
 	\def\Tempb{0912.3315}
	\ifx\Tempa\Tempb
	\def\wRef{\arXivRef arxiv/pdf/0912/0912.3315v2.pdf}
	\fi
 	\def\Tempb{0912.4061}
	\ifx\Tempa\Tempb
	\def\wRef{\arXivRef arxiv/pdf/0912/0912.4061v2.pdf}
	\fi
 	\def\Tempb{1003.1544}
	\ifx\Tempa\Tempb
	\def\wRef{\arXivRef arxiv/pdf/1003/1003.1544v2.pdf}
	\fi
 	\def\Tempb{1104.5197}
	\ifx\Tempa\Tempb
	\def\wRef{\arXivRef arxiv/pdf/1104/1104.5197.pdf}
	\fi
 	\def\Tempb{1107.1139}
	\ifx\Tempa\Tempb
	\def\wRef{\arXivRef arxiv/pdf/1104/1107.1139.pdf}
	\fi
 	\def\Tempb{1107.5037}
	\ifx\Tempa\Tempb
	\def\wRef{\arXivRef arxiv/pdf/1107/1107.5037.pdf}
	\fi
 	\def\Tempb{8433-5163}
	\ifx\Tempa\Tempb
	\def\wRef{http://www.amazon.com/}
	\fi
 	\def\Tempb{8443-0072}
	\ifx\Tempa\Tempb
	\def\wRef{http://www.amazon.com/}
	\fi
	\externaldocument[#1-#2-]{\FilePrefix #1.#2}[\wRef]
	}
\newcommand\LanguagePrefix{}%
\makeatletter
\newcommand\StartLabelItem[1]
{
\def\theenumi{\ref*{#1}.\@arabic\c@enumi}
\def\labelenumi{\theenumi:}
}
\newcommand\StopLabelItem
{
\def\theenumi{\@arabic\c@enumi}
\def\labelenumi{(\theenumi)}
}
\makeatother

\ifx\SelectlEnglish\undefined
	\newcommand\CurrentLanguage{Russian.}%
	\author{Александр Клейн}
		\newtheorem{theorem}{Теорема}[section]
		\newtheorem{corollary}[theorem]{Следствие}
		\newtheorem{convention}[theorem]{Соглашение}
		\theoremstyle{definition}
		\newtheorem{definition}[theorem]{Определение}
		\newtheorem{example}[theorem]{Пример}
		\newtheorem{xca}[theorem]{Exercise}
		\theoremstyle{remark}
		\newtheorem{remark}[theorem]{Замечание}
		\newtheorem{lemma}[theorem]{Лемма}
		
	\newcommand\xRefDef[1]
		{
		\wRefDef{#1}{Russian}
		\NameDef{xRefDef#1}{}%
		}
	\makeatletter
	\newcommand\xRef[2]%
	{%
		\@ifundefined{xRefDef#1}{%
		\ref{#2}%
		}{% 
		\citeBib{#1}-\ref{#1-Russian-#2}%
		}%
	}%
	\newcommand\xEqRef[2]%
	{%
		\@ifundefined{xRefDef#1}{%
		\eqref{eq: #2}%
		}{%
		\citeBib{#1}-\eqref{#1-Russian-eq: #2}%
		}%
	}%
	\makeatother
	\ifx\PrintBook\undefined
		\newcommand{\BibTitle}{%
			\section{Список литературы}%
		}
	\else
		\newcommand{\BibTitle}{%
			\chapter{Список литературы}%
		}
	\fi
\else
	\newcommand\CurrentLanguage{English.}%
	\author{Aleks Kleyn}
		\newtheorem{theorem}{Theorem}[section]
		\newtheorem{corollary}[theorem]{Corollary}
		\newtheorem{convention}[theorem]{Convention}
		\theoremstyle{definition}
		\newtheorem{definition}[theorem]{Definition}
		\newtheorem{example}[theorem]{Example}
		
		\theoremstyle{remark}
		\newtheorem{remark}[theorem]{Remark}
	\newcommand\xRefDef[1]
		{
		\wRefDef{#1}{English}
		\NameDef{xRefDef#1}{}%
		}
	\makeatletter
	\newcommand\xRef[2]%
	{%
		\@ifundefined{xRefDef#1}{%
		\ref{#2}%
		}{% 
		\citeBib{#1}-\ref{#1-English-#2}%
		}%
	}%
	\newcommand\xEqRef[2]%
	{%
		\@ifundefined{xRefDef#1}{%
		\eqref{eq: #2}%
		}{%
		\citeBib{#1}-\eqref{#1-English-eq: #2}%
		}%
	}%
	\makeatother
	\ifx\PrintBook\undefined
		\newcommand{\BibTitle}{%
			\section{References}%
		}
	\else
		\newcommand{\BibTitle}{%
			\chapter{References}%
		}
	\fi
\fi

\ifx\PrintBook\undefined
		\numberwithin{Hfootnote}{section}
\else
	\numberwithin{section}{chapter}
	\numberwithin{footnote}{chapter}
	\numberwithin{Hfootnote}{chapter}
\fi

\ifx\Presentation\undefined
	\numberwithin{equation}{section}
	\numberwithin{figure}{section}
	\numberwithin{table}{section}
	\numberwithin{Item}{section}
\fi

\makeatletter
\newcommand\org@maketitle{}
\let\org@maketitle\maketitle
\def\maketitle{%
	\hypersetup{pdftitle={\@title}}%
	\hypersetup{pdfauthor={\authors}}%
	\hypersetup{pdfsubject=\@keywords}%
	\ifx\UseRussian\Defined
		\pdfbookmark[1]{\@title}{TitleRussian}
	\else
		\pdfbookmark[1]{\@title}{TitleEnglish}
	\fi
	\org@maketitle
}
\def\make@stripped@name#1{%
	\begingroup
		\escapechar\m@ne
		\global\let\newname\@empty
		\protected@edef\Hy@tempa{\CurrentLanguage #1}%
		\edef\@tempb{%
			\noexpand\@tfor\noexpand\Hy@tempa:=%
			\expandafter\strip@prefix\meaning\Hy@tempa
		}%
		\@tempb\do{%
			\if\Hy@tempa\else
				\if\Hy@tempa\else
					\xdef\newname{\newname\Hy@tempa}%
				\fi
			\fi
		}%
	\endgroup
}%
\newenvironment{enumBib}{%
	\BibTitle
	\advance\@enumdepth \@ne
	\edef\@enumctr{enum\romannumeral\the\@enumdepth}\list
	{\csname biblabel\@enumctr\endcsname}{\usecounter
	{\@enumctr}\def\makelabel##1{\hss\llap{\upshape##1}}}
}{%
	\endlist
}
\makeatletter

\newcommand{\BiblioItem}[2]
{
	\def\Semafor{off}
	\@ifundefined{\LanguagePrefix ViewCite#1}{}{%
		\def\Semafor{on}%
	}%
	\ifx\Semafor\ValueOff
		\@ifundefined{xRefDef#1}{}{% 
		\def\Semafor{on}%
		}%
	\fi
	\ifx\Semafor\ValueOn
		\ifx\IndexState\ValueOff
			\begin{enumBib}
			\def\IndexState{on}
		\fi
		\item \label{\LanguagePrefix bibitem: #1}#2%
	\fi
}
\makeatother
\newcommand{\OpenBiblio}
{
	\def\IndexState{off}
}
\newcommand{\CloseBiblio}
{
	\ifx\IndexState\ValueOn
		\end{enumBib}
		\def\IndexState{off}
	\fi
}

\makeatletter
\def\StartCite{[}%
\def\citeBib#1{\protect\showCiteBib#1,endCite,}%
\def\endCite{endCite}%
\def\showCiteBib#1,{\def\temp{#1}%
\ifx\temp\endCite
]%
\def\StartCite{[}%
\else
	\StartCite\LanguagePrefix \ref{\LanguagePrefix bibitem: #1}%
	\@ifundefined{\LanguagePrefix ViewCite#1}{%
		\NameDef{\LanguagePrefix ViewCite#1}{}%
	}{%
	}%
	\def\StartCite{, }%
\expandafter\showCiteBib%
\fi}%
\makeatother

\newcommand{\arp}{\ar @{-->}}
\newcommand{\ars}{\ar @{.>}}
\newcommand\Bundle[1]{{\mathbb #1}}%{{\bf #1}}
\newcommand{\bundle}[4]%
{%
	\def\tempa{}%
	\def\tempb{#3}%
	\def\tempc{#1}%
	\ifx\tempa\tempb%
		\ifx\tempa\tempc%
			#2%
		\else%
			\xymatrix{#2:#1\arp[r]&#4}%
		\fi%
	\else%
		\ifx\tempa\tempc%
			#2[#3]%
		\else%
			\xymatrix{#2[#3]:#1\arp[r]&#4}%
		\fi%
	\fi%
}%
\makeatletter
\newcommand{\AddIndex}[2]%
{%
	\@ifundefined{RefIndex#2}{%
		\NameDef{RefIndex#2}{\pageref{index: #2}}%
		\label{index: #2}%
	}{%
		\NameDef{RefIndex#2}{\pageref{index: #2}, \pageref{index 1: #2}}%
		\label{index 1: #2}%
	}%
	{\bf #1}%
}%
\newcommand{\Index}[2]%
{%
	\def\Semafor{off}%
	\@ifundefined{RefIndex#2}{%
	}{%
		\def\Semafor{on}
	}%
	\ifx\Semafor\ValueOn%
		\def\tempa{}%
		\def\tempb{#2}%
		\ifx\IndexState\ValueOff%
			\begin{theindex}%
			\def\IndexState{on}%
		\fi%
		\ifx\IndexSpace\ValueOn%
			\indexspace%
			\def\IndexSpace{off}%
		\fi%
		\item #1%
		\ifx\tempa\tempb%
		\else%
			\ \@nameuse{RefIndex#2}%
		\fi%
	\fi%
}%
\newcommand{\Symb}[2]
{
	\def\Semafor{off}
	\@ifundefined{ViewSymbol#2}{%
	}{%
		\def\Semafor{on}
	}%
	\ifx\Semafor\ValueOn
		\ifx\IndexState\ValueOff
			\begin{theindex}
			\def\IndexState{on}
		\fi
		\ifx\IndexSpace\ValueOn
			\indexspace
			\def\IndexSpace{off}
		\fi
		\item $\displaystyle\@nameuse{ViewSymbol#2}$\ \ #1
		\@nameuse{RefSymbol#2}%
	\fi
}

\makeatother
\newcommand{\SetIndexSpace}%
{%
	\def\IndexSpace{on}%
}%

\newcommand{\OpenIndex}
{
	\def\IndexState{off}
}
\newcommand{\CloseIndex}
{
	\ifx\IndexState\ValueOn
		\end{theindex}
		\def\IndexState{off}
	\fi
}

\def\LastMemo{LastMemo}%
\def\MemoList#1//{\def\temp{#1}%
	\ifx\temp\LastMemo
	\else%
		\par
		\BlueText{#1}%
		\expandafter\MemoList%
	\fi%
}%

%% file: Preface.1006.2597.English.tex
%auto-ignore
\input{Preface.1006.2597.Eq}

%\chapter{Preface}

\ePrints{2010.05.23}
\ifx\Semafor\ValueOn
\newpage
\fi
\section{Preface}

The possibility of linear approximation of mapping
is at the heart of calculus and
main constructions of calculus have their roots in linear
algebra.
Therefore, before we give the definition of the differentiable function,
we need to have an idea of the mappings that
we are going to use to approximate the behavior of the original function.

Since the product in the field is
commutative, then linear algebra over a field is relatively simple.
When we explore division ring where product is not commutative,
we still see some familiar statements of linear
algebra; however, we meet new statements, which
change the landscape of linear algebra.

\ePrints{2010.05.23}
\ifx\Semafor\ValueOn
\newpage
\fi

Here I want to draw attention to the evolution of the concept
of the derivative from the time of Newton.
When we study functions of single variable, the derivative in selected point
is a number.
\ePrints{2010.05.23}
\ifx\Semafor\ValueOn
\ShowEq{derivative single variable}
\fi
When we study function of multiple variables, we realize
that it is not enough to use number. The derivative becomes vector or gradient.
\ePrints{2010.05.23}
\ifx\Semafor\ValueOn
\ShowEq{derivative multiple variables}
\newpage
\fi
When we study maps of vector spaces, this is a first time that we tell about derivative
as operator.
\ePrints{2010.05.23}
\ifx\Semafor\ValueOn
\ShowEq{derivative vector variables}
\fi
However since this operator is linear, then we can represent
derivative as matrix.
Again we express a vector of increment of function
as product of a matrix of derivative (Jacobian matrix) over
vector of increment of argument.
\ePrints{2010.05.23}
\ifx\Semafor\ValueOn
\ShowEq{derivative vector variables 1}
\newpage
\else

Surely, such behavior of derivative weakens our attention.
When we consider objects which are more complex than fields or vector spaces,
we still try to see an object which can be written as a factor before an increment
and which does not depend on the increment.
\fi

The assumption that the derivative of the mapping $f$
over algebra $A$ is defined by equation
\ShowEq{Frechet derivative, algebra}
initially looks attractive.
\ePrints{2010.05.23}
\ifx\Semafor\ValueOn
However, let $f(x)=f_1(x)f_2(x)$. According to the definition
\EqRef{Frechet derivative, algebra},\footnote{In the last chain of the equation
\EqRef{Frechet derivative, algebra 1},
we used the statement
\ShowEq{Frechet derivative, algebra 2}
}
\ShowEq{Frechet derivative, algebra 1}
From the equation
\EqRef{Frechet derivative, algebra},
it also follows that
\ShowEq{Frechet derivative, algebra 3}
From the equations
\EqRef{Frechet derivative, algebra 1}
and
\EqRef{Frechet derivative, algebra 3},
it follows that
\ShowEq{Frechet derivative, algebra 4}
We cannot get rid of the value of $h$ in the equation
\EqRef{Frechet derivative, algebra 4}.
Therefore,
\else
At first glance, such a definition satisfies the classical properties
of derivative of mapping over field.
However, in general,
\fi
the product of differentiable functions
is not differentiable function.
It causes that the set of differentiable functions is
very small, and such theory of differentiation is not
interesting.

\ePrints{2010.05.23}
\ifx\Semafor\ValueOn
\newpage
Is there a wider class of mappings that could
contain derivatives of various functions?
To answer this question, consider the mapping $y=x^2$
\ShowEq{mapping y=x2}
The mapping
\ShowEq{mapping y=x2 1}
is not linear mapping over algebra $A$,
but it is linear mapping over center of algebra $A$.
\else
Since the algebra is a module over some commutative
ring, there exist two ways to explore structures
generated over the algebra.\footnote{To explore the differentiation
in the algebra, there is one more method of studying.
Considering a certain set of functions,
we can define differential operators acting
on this set (\citeBib{Serge Lang}, p. 368).

For instance, in the paper \citeBib{math.CV-0405471},
Ludkovsky considers differential operators $\Pza$ and $\Pzb$
in Cayley - Dickson algebra such that
\ShowEq{Ludkovsky differential operators}
Considering given properties of differential operator,
Ludkovsky study its structure.

The exploration of differential operator from different points of view
gives more deep knowledge and in the future I suppose to consider
relation between different way of study theory of differentiation.
}

If the algebra is a free module, then we can choose a basis and
consider all operations in the coordinates relative a given basis.
Although the base can be arbitrary, we can choose the most simple
basis in terms of algebraic operations.
Beyond doubt, this approach has the advantage
that we are working in commutative ring where all
operations are well studied.

Exploration of operations in algebra regardless of
the chosen basis gives an opportunity to consider elements
of algebra as independent objects.
However, noncommutativity of the product in the algebra is a source
of a lot of difficulties on this way.

The question arises as to whether there exists an alternative method
if the definition of derivative \EqRef{Frechet derivative, algebra}
restricts our ability to study infinitesimal behavior of map?
The answer on this question is affirmative.
We explore the calculus in algebra
that is normed module.
We know two types of derivative in normed space. Strong derivative or
the Fr\'echet derivative is analogue of the derivative
\EqRef{Frechet derivative, algebra}.
Besides strong derivative there exits weak derivative or the G\^ateaux derivative.
The main idea is that differential may depend on direction.

The algebra $A$ is the
module over commutative ring $D$. If we relax the definition of
derivative and require that derivative of the mapping $f$ is a linear
mapping of the module $A$, then we see that at least
polynomial in algebra $A$ is differentiable mapping.
Defined in this way the derivative has many properties
of derivative of the mapping over field.
\fi
Therefore, in algebra, because of the noncommutativity of multiplication
differential of the function has terms of the form
\ShowEq{adxb}
and we cannot write the differential of the mapping
as a product of the derivative and the differential of the argument.

Hamilton was first to explore
the differential of mapping of quaternion algebra
(\citeBib{Hamilton Elements of Quaternions 1}).
Apparently, his results appeared so out of the ordinary that
his contemporaries found it difficult
to embrace this Hamilton's idea.
\ePrints{2010.05.23}
\ifx\Semafor\ValueOn
Next generations forgot this research.
\else
Next generations forgot this research.\footnote{In the paper
\citeBib{W.Bertram H.Glockner K.Neeb},
Bertram explored the G\^ateaux derivative of mappings over the
commutative ring.
Since the product is commutative, then statements in the
paper are close to statements of classical calculus.
The derivative is defined as the mapping of the differential of the argument
into the differential of the function.}
\fi

The problem of inseparability of
derivative and differential is so serious that when G\^ateaux defined
weak differentiation, he considered a derivative only in
case when he was able to separate increment of an argument as factor.

\ePrints{2010.05.23}
\ifx\Semafor\ValueOn
\newpage
\fi
However, how is this obstacle serious?
Derivative is a mapping of differential of an argument
to differential of function.
In other words, derivative is some algorithm
whose input is differential of argument and
the output is the differential of function.
If we consider functional notation for linear mapping, namely
\ShowEq{functional notation for linear mapping}
then we can formally separate the differential of the argument
from the derivative and
make a notation of the differential of function
\ShowEq{Gateaux differential}
more familiar, namely
\ShowEq{Gateaux differential, 1}
\ePrints{2010.05.23}
\ifx\Semafor\ValueOff
Tensor product of algebras allows us to write the
structure of derivative as operator.

In addition, new notation allows us to simplify many expressions.
For instance, an expression
\ShowEq{differential vector function}
gets form
\ShowEq{differential vector function, 1}

As in the case of mappings over the field, the differential
of mapping over the algebra is polynomial of first power
with respect to increment of argument.
The structure of polynomial over division ring is different from structure of polynomial over field.
I consider some properties of polynomial in section
\ref{section: Taylor Series}.
Using obtained theorems I explore
Taylor series expansion of map and method
to find solution of differential equation.
\fi

%% file: Preface.1006.2597.Eq.tex
%auto-ignore

\def\zb{\overline z}
\def\Pza{\partial_z}
\def\Pzb{\partial_{\zb}}

\DefEq
{
\ePrints{2010.05.23}
\ifx\Semafor\ValueOff
\[
f\circ x=f(x)
\]
\else
\begin{equation}
f\circ x=f(x)
\end{equation}
\fi
}
{functional notation for linear mapping}

\DefEq
{
\begin{align*}
z&=x^2+y^3
\\
dz&=\textcolor{red}{2x}\ dx+\textcolor{red}{3y^2}\ dy
\end{align*}
}
{derivative multiple variables}

\DefEq
{
\[
\begin{matrix}
x=u\sin v&y=u\cos v&z=u
\\
dx=\textcolor{red}{\sin v}\ du+\textcolor{red}{u\cos v}\ dv&
dy=\textcolor{red}{\cos v}\ du-\textcolor{red}{u\sin v}\ dv&
dz=\textcolor{red}{1}\ du+\textcolor{red}{0}\ dv
\end{matrix}
\]
}
{derivative vector variables}

\DefEq
{
$h\partial f_2(x)h=o(h)$
}
{Frechet derivative, algebra 2}

\DefEquation
{
\begin{array}{r@{\,}l}
&f_1(x+h)f_2(x+h)-f_1(x)f_2(x)
\\
=&f_1(x+h)f_2(x+h)-f_1(x)f_2(x+h)
+f_1(x)f_2(x+h)-f_1(x)f_2(x)
\\
=&(f_1(x+h)-f_1(x))f_2(x+h)+f_1(x)(f_2(x+h)-f_2(x))
\\
=&\partial f_1(x)hf_2(x)+f_1(x)\partial f_2(x)h
\end{array}
}
{Frechet derivative, algebra 1}

\DefEquation
{
\partial(f_1(x)f_2(x))h
=
\partial f_1(x)hf_2(x)+f_1(x)\partial f_2(x)h
}
{Frechet derivative, algebra 4}

\DefEquation
{
a\,dx\,b
}
{adxb}

\DefEquation
{
f(h)=xh+hx
}
{mapping y=x2 1}

\DefEquation
{
(x+h)^2-x^2=x^2+xh+hx+h^2-x^2=xh+hx+o(h)
}
{mapping y=x2}

\DefEquation
{
f_1(x+h)f_2(x+h)-f_1(x)f_2(x)=
\partial(f_1(x)f_2(x))h
}
{Frechet derivative, algebra 3}

\DefEq
{
\[
\begin{pmatrix}
dx\\dy\\dz
\end{pmatrix}
=
\begin{pmatrix}
\sin v&u\cos v\\
\cos v&-u\sin v\\
1&0
\end{pmatrix}
\begin{pmatrix}
du\\dv
\end{pmatrix}
\]
}
{derivative vector variables 1}

\DefEq
{
\[
d(x^2)=\textcolor{red}{2x}\ dx
\]
}
{derivative single variable}

\DefEq
{
\[
\begin{matrix}
\Pza z=1&\Pza\zb=0
\\
\Pzb z=0&\Pzb\zb=1
\end{matrix}
\]
}
{Ludkovsky differential operators}

\DefEq
{
\[
\partial \Vector f(\Vector x)(d\Vector x)
=
\Vector r_j
\frac{\partial\pC{s}{0} f^j(\Vector x)}{\partial x^i}
dx^i
\frac{\partial\pC{s}{1} f^j(\Vector x)}{\partial x^i}
\]
}
{differential vector function}

\DefEq
{
\begin{align*}
&
\partial \Vector f(\Vector x)\RCcirc d\Vector x
=
\begin{pmatrix}
\partial f^1(\Vector x)\RCcirc d\Vector x
\\...\\
\partial f^n(\Vector x)\RCcirc d\Vector x
\end{pmatrix}
\\
=&
\begin{pmatrix}
\displaystyle
\frac{\partial\pC{s}{0} f^1(\Vector x)}{\partial x^1}
\otimes
\frac{\partial\pC{s}{1} f^1(\Vector x)}{\partial x^1}
&...&
\displaystyle
\frac{\partial\pC{s}{0} f^1(\Vector x)}{\partial x^m}
\otimes
\frac{\partial\pC{s}{1} f^1(\Vector x)}{\partial x^m}
\\...\\
\displaystyle
\frac{\partial\pC{s}{0} f^n(\Vector x)}{\partial x^1}
\otimes
\frac{\partial\pC{s}{1} f^n(\Vector x)}{\partial x^1}
&...&
\displaystyle
\frac{\partial\pC{s}{0} f^n(\Vector x)}{\partial x^m}
\otimes
\frac{\partial\pC{s}{1} f^n(\Vector x)}{\partial x^m}
\end{pmatrix}
\RCcirc
\begin{pmatrix}
dx^1\\...\\dx^m
\end{pmatrix}
\end{align*}
}
{differential vector function, 1}

\DefEq
{
\[
\partial f(x)(dx)
\]
}
{Gateaux differential}

\DefEq
{
\[
\partial f(x)\circ dx
\]
}
{Gateaux differential, 1}

\DefEquation
{
f(x+h)-f(x)=\partial f(x)h+o(h)
}
{Frechet derivative, algebra}

%% file: Convention.English.tex
%auto-ignore
\input{\FilePrefix Convention.Eq}

\section{Conventions}

\ePrints{0812.4763,0906.0135,0908.3307,0909.0855,0912.3315,1003.1544}
\Items{1006.2597}
\ifx\Semafor\ValueOn
\begin{convention}
Function and mapping are synonyms. However according to
tradition, correspondence between either rings or vector
spaces is called mapping and a mapping of
either real field or quaternion algebra is called function.
\qed
\end{convention}
\fi

\ePrints{0701.238,0812.4763,0908.3307,0912.4061,1001.4852}
\Items{1003.1544}
\ifx\Semafor\ValueOn
\begin{convention}
In any expression where we use index I assume
that this index may have internal structure.
For instance, considering the algebra $A$ we enumerate coordinates of
$a\in A$ relative to basis $\Basis e$ by an index $i$.
This means that $a$ is a vector. However, if $a$
is matrix, then we need two indexes, one enumerates
rows, another enumerates columns. In the case, when index has
structure, we begin the index from symbol $\cdot$ in
the corresponding position. 
For instance, if I consider the matrix $a^i_j$ as an element of a vector
space, then I can write the element of matrix as $a^{\cdot}{}^i_j$.
\qed
\end{convention}
\fi

\ePrints{0701.238,0812.4763,0908.3307,0912.4061,1006.2597,1011.3102}
\Items{Calculus.Paper}
\ifx\Semafor\ValueOn
\begin{convention}
I assume sum over index $s$ in expression like
\ShowEq{Sum over repeated index}
\qed
\end{convention}
\fi

\ePrints{0701.238,0812.4763,0906.0135,0908.3307,0909.0855}
\ifx\Semafor\ValueOn
\begin{convention}
We can consider division ring $D$ as $D$\Hyph vector space
of dimension $1$. According to this statement, we can explore not only
homomorphisms of division ring $D_1$ into division ring $D_2$,
but also linear maps of division rings.
This means that map is multiplicative over
maximum possible field. In particular, linear map
of division ring $D$ is multiplicative over center $Z(D)$. This statement
does not contradict with
definition of linear map of field because for field $F$ is true
$Z(F)=F$.
When field $F$ is different from
maximum possible, I explicit tell about this in text.
\qed
\end{convention}
\fi

\ePrints{0912.4061}
\ifx\Semafor\ValueOn
\begin{convention}
For given field $F$, unless otherwise stated,
we consider finite dimensional $F$\Hyph algebra.
%without zero divisors.
\qed
\end{convention}
\fi

\ePrints{0701.238,0812.4763,0906.0135,0908.3307}
\ifx\Semafor\ValueOn
\begin{convention}
In spite of noncommutativity of product a lot of statements
remain to be true if we substitute, for instance, right representation by
left representation or right vector space by left
vector space.
To keep this symmetry in statements of theorems
I use symmetric notation.
For instance, I consider \Ds vector space
and \sD vector space.
We can read notation \Ds vector space
as either D\Hyph star\Hyph vector space or
left vector space.
We can read notation \Ds linear dependent vectors
as either D\Hyph star\Hyph linear dependent vectors or
vectors that are linearly dependent from left.
\qed
\end{convention}
\fi

\ePrints{0701.238,0812.4763,0906.0135,0908.3307,0909.0855,0912.4061}
\Items{1001.4852,1003.1544,1006.2597,1104.5197,1105.4307,1107.1139}
\ifx\Semafor\ValueOn
\begin{convention}
Let $A$ be free finite
dimensional algebra.
Considering expansion of element of algebra $A$ relative basis $\Basis e$
we use the same root letter to denote this element and its coordinates.
However we do not use vector notation in algebra.
In expression $a^2$, it is not clear whether this is component
of expansion of element
$a$ relative basis, or this is operation $a^2=aa$.
To make text clearer we use separate color for index of element
of algebra. For instance,
\ShowEq{Expansion relative basis in algebra}
\qed
\end{convention}

\begin{convention}
If free finite dimensional algebra has unit, then we identify
the vector of basis $\Vector e_{\gi 0}$ with unit of algebra.
\qed
\end{convention}
\fi

\ePrints{1104.5197,1105.4307}
\ifx\Semafor\ValueOn
\begin{convention}
Although the algebra is a free module over some
ring, we do not use the vector notation
to write elements of algebra. In the case when I consider the
matrix of coordinates of element of algebra, I will use vector
notation to write corresponding element.
In order to avoid ambiguity when I use conjugation,
I denote $a^*$ element conjugated to element $a$.
\qed
\end{convention}
\fi

\ePrints{0906.0135,0912.3315,8443-0072,1111.6035,1102.5168}
\ifx\Semafor\ValueOn
\begin{convention}
In \citeBib{Cohn: Universal Algebra},
an arbitrary operation of algebra is denoted by letter $\omega$,
and $\Omega$ is the set of operations of some universal algebra.
Correspondingly, the universal algebra with the set of operations
$\Omega$ is denoted as $\Omega$\Hyph algebra.
Similar notations we see in
\citeBib{Burris Sankappanavar} with small difference
that an operation in the algebra is denoted by letter $f$
and $\mathcal F$ is the set of operations.
I preferred first case of notations because in this case it is
easier to see where I use operation.
\qed
\end{convention}
\fi

\ePrints{0906.0135,0912.3315,8443-0072}
\ifx\Semafor\ValueOn
\begin{convention}
Since the number of universal algebras
in the tower of representations is varying,
then we use vector notation for a tower of
representations. We denote the set
$(A_1,...,A_n)$ of $\Omega_i$\Hyph algebras $A_i$, $i=1$, ..., $n$
as $\Vector A$. We denote the set of representations
$(f_{1,2},...,f_{n-1,n})$ of these algebras as $\Vector f$.
Since different algebras have different type, we also
talk about the set of $\Vector{\Omega}$\Hyph algebras.
\ePrints{8443-0072}
\ifx\Semafor\ValueOn
We
\else
In relation to the set $\Vector A$,
we also use matrix notations 
that we discussed
in section \xRef{0701.238}{section: Concept of Generalized Index}.
For instance, we
\fi
use the symbol $\Vector A_{[1]}$ to denote the
set of $\Vector{\Omega}$\Hyph algebras $(A_2,...,A_n)$.
In the corresponding notation $(\VX A,\Vector f)$ of tower
of representation, we assume that $\Vector f=(f_{2,3},...,f_{n-1,n})$.
\qed
\end{convention}

\begin{convention}
Since we use vector notation for elements of the
tower of representations, we need convention about notation of operation.
We assume that we get result of operation componentwise. For instance,
\ShowEq{vector notation in tower of representations}
\qed
\end{convention}
\fi

\ePrints{8443-0072,1111.6035,0906.0135,NewAffine,1102.5168}
\ifx\Semafor\ValueOn
\begin{convention}
Let $A$ be $\Omega_1$\Hyph algebra.
Let $B$ be $\Omega_2$\Hyph algebra.
Notation
\ShowEq{A->*B}
means that there is representation of $\Omega_1$\Hyph algebra $A$
in $\Omega_2$\Hyph algebra $B$.
\qed
\end{convention}
\fi

\ePrints{0702.561,0707.2246,0803.2620}
\ifx\Semafor\ValueOn
\begin{convention}
I use arrow $\xymatrix{\arp[r]&}$ to represent
projection of bundle on diagram.
I use arrow $\xymatrix{\ars[r]&}$ to represent
section of bundle on diagram.
\qed
\end{convention}
\fi

\ePrints{0912.3315}
\ifx\Semafor\ValueOn
\begin{remark}
I believe that diagrams of maps are an important tool.
However, sometimes I want
to see the diagram as three dimensional figure
and I expect that this would increase its expressive
power. Who knows what surprises the future holds.
In 1992, at a conference in Kazan, I have described to my colleagues
what advantages the computer preparation of papers has.
8 years later I learned from the letter from Kazan that now we can
prepare paper using LaTeX.
\qed
\end{remark}
\fi

\ePrints{1001.4852,1003.1544,1006.2597,1011.3102}
\Items{Calculus.Paper}
\ifx\Semafor\ValueOn
\begin{convention}
If, in a certain expression, we use several operations
which include the operation $\circ$, then
it is assumed that the operation $\circ$ is executed first.
Below is an example of equivalent expressions.
\ShowEq{list circ expressions}
\qed
\end{convention}
\fi

%\ePrints{0906.0135,NewAffine}
%\ifx\Semafor\ValueOn
%\item
%Let $\VX X$ be the basis of the tower of representations $(\Vector A,\Vector f)$.
%A $\Omega_n$\Hyph word
%$w_n(\Vector f,\VX X,x_n)$ is called coordinates $x_n\in A_n$
%relative to basis $\VX X$.
%\fi

\ePrints{1107.1139}
\ifx\Semafor\ValueOn
\begin{convention}
For given $D$\Hyph algebra $A$
we define left shift
\ShowEq{left shift, D algebra}
by the equation
\ShowEq{left shift 1, D algebra}
and right shift
\ShowEq{right shift, D algebra}
by the equation
\ShowEq{right shift 1, D algebra}
\qed
\end{convention}
\fi

\ifx\PrintPaper\undefined
Without a doubt, the reader may have questions,
comments, objections. I will appreciate any response.
\fi

%% file: Convention.Eq.tex
%auto-ignore

\DefEq
{
\[
\Vector r(\Vector a)=(r_1(a_1),...,r_n(a_n))
\]
}
{vector notation in tower of representations}

\DefEq
{
\[
\xymatrix
{
A\ar[r]|{*}&B
}
\]
}
{A->*B}

\DefEq
{
\[
\begin{array}{r@{\ }lr@{\ }l}
f\circ xy&\equiv f(x)y
&
f\circ(xy)&\equiv f(xy)
\\
f\circ x+y&\equiv f(x)+y
&
f\circ (x+y)&\equiv f(x+y)
\end{array}
\]
}
{list circ expressions}

\DefEq
{
\[
a\pC s0xa\pC s1
\]
}
{Sum over repeated index}

\DefEq
{
\[
a=a^{\gi i}\Vector e_{\gi i}
\]
}
{Expansion relative basis in algebra}

\DefEq
{
\symb{a\circ}1{left shift, D algebra}
}
{left shift, D algebra}

\DefEq
{
\symb{a\star}1{right shift, D algebra}
}
{right shift, D algebra}

\DefEq
{
\[
\ShowSymbol{left shift, D algebra}x=ax
\]
}
{left shift 1, D algebra}

\DefEq
{
\[
\ShowSymbol{right shift, D algebra}x=xa
\]
}
{right shift 1, D algebra}

%% file: Linear.Mapping.English.tex
%auto-ignore
\input{\FilePrefix Linear.Mapping.Eq}

\ePrints{1011.3102}
\ifx\Semafor\ValueOff
\ePrints{2010.05.23}
\ifx\Semafor\ValueOn
\newpage
\fi
\Chapter{Linear Mapping of Algebra}

\ePrints{2010.05.23}
\ifx\Semafor\ValueOff
\Section{Module}
\label{Section: Module}

\begin{theorem}
\label{theorem: effective representation of the ring}
Let ring $D$ has unit $e$.
Representation
\ShowEq{representation of the ring}
of the ring $D$
in an Abelian group $A$ is
\AddIndex{effective}{effective representation of ring}
iff $a=0$ follows from equation $f(a)=0$.
\end{theorem}
\begin{proof}
We define the sum of transformations $f$ and $g$ of an Abelian group
according to rule
\ShowEq{sum of transformations of Abelian group}
Therefore, considering the representation of the ring $D$ in
the Abelian group $A$, we assume
\ShowEq{sum of transformations of Abelian group, 1}
We define the product of transformation of representation
according to rule
\ShowEq{product of transformations of representation}

Suppose $a$, $b\in R$
cause the same transformation. Then
\ShowEq{representation of ring, 1}
for any $m\in A$.
From the equation
\EqRef{representation of ring, 1}
it follows that $a-b$ generates zero transformation
\ShowEq{representation of ring, 2}
Element $e+a-b$ generates an identity transformation.
Therefore, the representation $f$ is effective iff $a=b$.
\end{proof}

\begin{definition}
\label{definition: module over ring}
Let $D$ be commutative ring.
$A$ is a \AddIndex{module over ring}{module over ring} $D$
if $A$ is an Abelian group and
there exists effective representation of ring $D$
in an Abelian group $A$.
\qed
\end{definition}

\begin{definition}
\label{definition: free module over ring}
We call set of vectors
\ShowEq{basis, module}
a \AddIndex{$D\star$\Hyph basis for module}{D basis, module}
if vectors $\Veb$ are
$D\star$\Hyph linearly independent and adding to this system any other vector
we get a new system which is $D\star$\Hyph linearly dependent.
$A$ is \AddIndex{free module over ring}{free module over ring} $D$,
if $A$ has basis
over ring $D$.\footnote{I follow to the
definition in \citeBib{Serge Lang}, p. 135.}
\qed
\end{definition}

\begin{theorem}
\label{theorem: free module over ring, change basis}
Let $A$ be free module over ring $D$.
Coordinates $a^{\gi j}$ of vector $a\in A$
are coordinates of $D$-valued contravariant tensor
\ShowEq{free module over ring, change basis, 3}
\end{theorem}
\begin{proof}
Let $\Basis e'$ be another basis. Let
\ShowEq{free module over ring, change basis}
be transformation, mapping basis $\Basis e$ into
basis $\Basis e'$.
Because vector $a$ does not change, then
\ShowEq{free module over ring, change basis, 1}
From equations \EqRef{free module over ring, change basis}
and \EqRef{free module over ring, change basis, 1}
it follows that
\ShowEq{free module over ring, change basis, 2}
Because vectors $\Vector e_{\gi j}$ are linear independent, then equation
\EqRef{free module over ring, change basis, 3}
follows
from equation
\EqRef{free module over ring, change basis, 2}.
Therefore, coordinates of vector are tensor.
\end{proof}

\ePrints{1011.3102}
\ifx\Semafor\ValueOff
Following definition is consequence of definitions
\ref{definition: module over ring}
and \xRef{0912.3315}{definition: morphism of representations of F algebra}.

\begin{definition}
\label{definition: linear map from A1 R1 to A2 R2}
Let $A_1$ be module over ring $R_1$.
Let $A_2$ be module over ring $R_2$.
Morphism
\ShowEq{linear map from A1 R1 to A2 R2}
of representation of ring $R_1$ in the Abelian group $A_1$
into representation of ring $R_2$ in the Abelian group $A_2$
is called
\AddIndex{linear mapping of $R_1$\Hyph module $A_1$
into $R_2$\Hyph module $A_2$}
{linear mapping of R_1 module into R_2 module}.
\qed
\end{definition}

\begin{theorem}
Linear mapping
\ShowEq{linear map from A1 R1 to A2 R2}
of $R_1$\Hyph module $A_1$
into $R_2$\Hyph module $A_2$
satisfies to equations\footnote{In classical notation, proposed equations
have quite familiar form
\ShowEq{linear map from A1 R1 to A2 R2, 1, old}}
\ShowEq{linear map from A1 R1 to A2 R2, 1}
\end{theorem}
\begin{proof}
From definitions
\ref{definition: linear map from A1 R1 to A2 R2}
and \xRef{0912.3315}{definition: morphism of representations of F algebra}
it follows that
\begin{itemize}
\item the mapping $f$ is a homomorphism of the ring $R_1$
into the ring $R_2$ (the equation
\EqRef{linear map from A1 R1 to A2 R2, 1 3})
\item the mapping $g$ is a homomorphism of the Abelian group $A_1$
into the Abelian group $A_2$ (the equation
\EqRef{linear map from A1 R1 to A2 R2, 1 1})
\end{itemize}
The equation
\EqRef{linear map from A1 R1 to A2 R2, 1 2}
follows from the equation
\xEqRef{0912.3315}{morphism of representations of F algebra, definition, 2}.
\end{proof}

According to the theorem
\xRef{0912.3315}{theorem: morphism of representations of algebra, reduce},
in the study of linear mappings, without loss of generality,
we can assume $R_1=R_2$.
\fi

\begin{definition}
\label{definition: linear map from A1 to A2, module}
Let $A_1$ and
$A_2$ be modules over the ring $R$.
Morphism
\ShowEq{linear map from A1 to A2}
of representation of the ring $D$ in the Abelian group $A_1$
into representation of the ring $D$ in the Abelian group $A_2$
is called
\AddIndex{linear mapping of $D$\Hyph module $A_1$
into $D$\Hyph module $A_2$}
{linear mapping of R modules}.
\qed
\end{definition}

\begin{theorem}
\label{theorem: linear map from A1 to A2}
Linear mapping
\ShowEq{linear map from A1 to A2}
of $D$\Hyph module $A_1$
into $D$\Hyph module $A_2$
satisfies to equations\footnote{In classical notation, proposed equations
have form
\ShowEq{linear map from A1 to A2, 1 old}
In some books
(for instance, \citeBib{Serge Lang}, p. 119) the theorem \ref{theorem: linear map from A1 to A2}
is considered as a definition.}
\ShowEq{linear map from A1 to A2, 1}
\end{theorem}
\begin{proof}
From definition
\ref{definition: linear map from A1 to A2, module}
and theorem
\xRef{0912.3315}{theorem: morphism of representations of algebra, reduce}
it follows that
the mapping $g$ is a homomorphism of the Abelian group $A_1$
into the Abelian group $A_2$ (the equation
\EqRef{linear map from A1 to A2, 1 1})
The equation
\EqRef{linear map from A1 to A2, 1 2}
follows from the equation
\xEqRef{0912.3315}{morphism of representations of F algebra}.
\end{proof}

\begin{definition}
\label{definition: polylinear map of modules}
Let $D$ be the commutative ring.
Let $A_1$, ..., $A_n$, $S$ be $D$\Hyph modules.
We call map
\ShowEq{polylinear map of algebras}
\AddIndex{polylinear mapping of modules}{polylinear map of modules}
$A_1$, ..., $A_n$
into module
$S$,
if
\ShowEq{polylinear map of algebras, 1}
\qed
\end{definition}

\Section{Algebra over Ring}
\label{section: Algebra over Ring}

\fi%2010.05.23
\else%1011.3102
\section{Algebra over Ring}
\fi%1011.3102
\begin{definition}
\label{definition: algebra over ring}
Let $D$ be commutative ring.
\ePrints{1003.1544,8433-5163}
\ifx\Semafor\ValueOn
Let $A$ be module over ring $D$.\footnote{There are
several equivalent definitions of algebra.
Initially I supposed to consider a representation of the ring $D$
in the Abelian group of the ring $A$.
But I had to explain why the product of elements of the ring $D$
and of algebra $A$ is commutative.
This required a definition of the center of the algebra $A$.
After careful analysis I have chosen the definition given in
\citeBib{Richard D. Schafer}, p. 1,
\citeBib{0105.155}, p. 4.
%Эта точка зрения имеет ещё два преимущества.
%Во-первых, как следует из теоремы
%если алгебра не имеет единицы, то не существует
%изоморфного образа кольца $D$ в алгебре $A$.
%Во-вторых, я не исключаю возможность рассматривать
%векторное пространство над алгеброй как алгебру.
}
\else
\ePrints{1011.3102}
\ifx\Semafor\ValueOn
Let $A$ be module over ring $D$.\footnote{This subsection
is written on the base of the section
\xRef{8433-5163}{section: Algebra over Ring}.}
\else
Let $A$ be module over ring $D$.
\fi
\fi
For given bilinear mapping
\ShowEq{product in algebra, definition 1}
we define product in $A$
\ShowEq{product in algebra, definition 2}
$A$ is a \AddIndex{algebra over ring}{algebra over ring} $D$
if $A$ is $D$\Hyph module and
we defined product
\EqRef{product in algebra, definition 2}
in $A$.
\ePrints{2010.05.23,1011.3102}
\ifx\Semafor\ValueOff
Algebra
\ShowEq{opposite algebra}
is called
\AddIndex{the opposite algebra to algebra $A$}{opposite algebra}
if we define a product in the module $A$ according to
rule\footnote{I made the definition by analogy
with the definition \citeBib{Bourbaki: Algebra 1}-2, p. 2.}
\ShowEq{opposite algebra, product}
\fi
If $A$ is free
$D$\Hyph module, then $A$ is called
\AddIndex{free algebra over ring}{free algebra over ring} $D$.
\qed
\end{definition}

\ePrints{2010.05.23,1011.3102}
\ifx\Semafor\ValueOff
\begin{remark}
\label{remark: opposite algebra}
Algebra $A$ and opposite algebra coincide as
modules.
\qed
\end{remark}

\begin{theorem}
The multiplication in the algebra $A$ is distributive over addition.
\end{theorem}
\begin{proof}
The statement of the theorem follows from the chain of equations
\ShowEq{product distributive in algebra}
\end{proof}

The multiplication in algebra can be neither commutative
nor associative. Following definitions are based
on definitions given in \citeBib{Richard D. Schafer}, p. 13.

\begin{definition}
\label{definition: commutator of algebra}
The \AddIndex{commutator}{commutator of algebra}
\ShowEq{commutator of algebra}
measures commutativity in $D$\Hyph algebra $A$.
$D$\Hyph algebra $A$ is called
\AddIndex{commutative}{commutative D algebra},
if
\ShowEq{commutative D algebra}
\qed
\end{definition}

\begin{definition}
\label{definition: associator of algebra}
The \AddIndex{associator}{associator of algebra}
\ShowEq{associator of algebra}
measures associativity in $D$\Hyph algebra $A$.
$D$\Hyph algebra $A$ is called
\AddIndex{associative}{associative D algebra},
if
\ShowEq{associative D algebra}
\qed
\end{definition}

\begin{theorem}
Let $A$ be algebra over commutative ring $D$.\footnote{The statement of the
theorem is based on the equation
\citeBib{Richard D. Schafer}-(2.4).}
\ShowEq{associator of algebra, 1}
for any $a$, $b$, $c$, $d\in A$.
\end{theorem}
\begin{proof}
The equation \EqRef{associator of algebra, 1}
follows from the chain of equations
\ShowEq{associator of algebra, 2}
\end{proof}

\begin{definition}
\label{definition: nucleus of algebra}
The set\footnote{The definition is based on
the similar definition in
\citeBib{Richard D. Schafer}, p. 13}
\ShowEq{nucleus of algebra}
is called the
\AddIndex{nucleus of an $D$\Hyph algebra $A$}{nucleus of algebra}.
\qed
\end{definition}

\begin{definition}
\label{definition: center of algebra}
The set\footnote{The definition is based on
the similar definition in
\citeBib{Richard D. Schafer}, p. 14}
\ShowEq{center of algebra}
is called the
\AddIndex{center of an $D$\Hyph algebra $A$}{center of algebra}.
\qed
\end{definition}

\begin{theorem}
\label{theorem: unit of algebra and ring}
Let $D$ be commutative ring.
If $D$\Hyph algebra $A$ has unit, then there exits
an isomorphism $f$ of the ring $D$ into the center of the algebra $A$.
\end{theorem}
\begin{proof}
Let $e\in A$ be the unit of the algebra $A$.
Then $f\circ a=ae$.
\end{proof}
\else
\newpage
\fi

Let $\Basis e$ be the basis of free algebra $A$ over ring $D$.
If algebra $A$ has unit,
then we assume that $\Vector e_{\gi 0}$ is the unit of algebra $A$.

\begin{theorem}
Let $\Basis e$ be the basis of free algebra $A$ over ring $D$.
Let
\ShowEq{a b in basis of algebra}
We can get the product of $a$, $b$ according to rule
\ShowEq{product in algebra}
where
\ShowEq{structural constants of algebra}
are \AddIndex{structural constants of algebra $A$ over ring $D$}
{structural constants of algebra}.
The product of basis vectors in the algebra $A$ is defined according to rule
\ShowEq{product of basis vectors, algebra}
\ePrints{2010.05.23}
\ifx\Semafor\ValueOn
\qed
\end{theorem}
\newpage
\else
\end{theorem}
\begin{proof}
The equation
\EqRef{product of basis vectors, algebra}
is corollary of the statement that $\Basis e$
is the basis of the algebra $A$.
Since the product in the algebra is a bilinear mapping,
then we can write the product of $a$ and $b$ as
\ShowEq{product in algebra, 1}
From equations
\EqRef{product of basis vectors, algebra},
\EqRef{product in algebra, 1},
it follows that
\ShowEq{product in algebra, 2}
Since $\Basis e$ is a basis of the algebra $A$, then the equation
\EqRef{product in algebra}
follows from the equation
\EqRef{product in algebra, 2}.
\end{proof}

\ePrints{1011.3102}
\ifx\Semafor\ValueOff
\begin{theorem}
Since the algebra $A$ is commutative, then
\ShowEq{commutative product in algebra, 1}
Since the algebra $A$ is associative, then
\ShowEq{associative product in algebra, 1}
\end{theorem}
\begin{proof}
For commutative algebra,
the equation
\EqRef{commutative product in algebra, 1}
follows from equation
\ShowEq{commutative product in algebra}
For associative algebra,
the equation
\EqRef{associative product in algebra, 1}
follows from equation
\ShowEq{associative product in algebra}
\end{proof}
\fi
\fi

\ePrints{1011.3102}
\ifx\Semafor\ValueOff
\Section{Linear Mapping of Algebra}
\label{Section: Linear Mapping of Algebra}
\fi

\ePrints{2010.05.23,1011.3102}
\ifx\Semafor\ValueOff
Algebra is a ring. A mapping, preserving
the structure of algebra as a ring, is called homomorphism of
algebra. However, the statement that algebra
is a module over a commutative ring is more important for us.
A mapping, preserving the
structure of algebra as module, is called a linear mapping of
algebra.
Thus,
the following definition
is based on the definition
\ref{definition: linear map from A1 to A2, module}.
\fi

\begin{definition}
\label{definition: linear map from A1 to A2, algebra}
Let $A_1$ and
$A_2$ be algebras over ring $D$.
\ePrints{1011.3102}
\ifx\Semafor\ValueOff
Morphism
\ShowEq{linear map from A1 to A2}
of the representation of the ring $D$ in the Abelian group $A_1$
into the representation of the ring $D$ in the Abelian group $A_2$
\else
The linear mapping\footnote{This subsection
is written on the base of the section
\xRef{8433-5163}{Section: Linear Mapping of Algebra}.}
\ShowEq{linear map from A1 to A2}
of the $D$\hyph module $A_1$
into the $D$\hyph module $A_2$
\fi
is called
\AddIndex{linear mapping of $D$\Hyph algebra $A_1$
into $D$\Hyph algebra $A_2$}
{linear mapping of R algebras}.
Let us denote
\ShowEq{set linear mappings, algebra}
set of linear mappings
of algebra
$A_1$
into algebra
$A_2$.
\qed
\end{definition}
\ePrints{2010.05.23}
\ifx\Semafor\ValueOn
\newpage
\fi

\ePrints{1011.3102}
\ifx\Semafor\ValueOff
\begin{theorem}
Linear mapping
\ShowEq{linear map from A1 to A2}
of $D$\Hyph algebra $A_1$
into $D$\Hyph algebra $A_2$
satisfies to equations
\ShowEq{linear map from A1 to A2, algebra}
\ePrints{2010.05.23}
\ifx\Semafor\ValueOn
\qed
\end{theorem}
\newpage
\else
\end{theorem}
\begin{proof}
The statement of theorem is a corollary of the theorem
\ref{theorem: linear map from A1 to A2}.
\end{proof}

\begin{theorem}
\label{theorem: linear map times constant, algebra}
Consider $D$\Hyph algebra $A_1$ and $D$\Hyph algebra $A_2$.
Let map
\ShowEq{linear map from A1 to A2}
be linear map.
Then maps
\ShowEq{linear map times constant, algebra}
defined by equations
\ShowEq{linear map times constant, 0, algebra}
are linear.
\end{theorem}
\begin{proof}
Statement of theorem follows from chains of equations
\ShowEq{linear map times constant, 1, algebra}
\end{proof}
\fi
\fi

\begin{theorem}
\label{theorem: linear map, 0, D algebra}
Consider $D$\Hyph algebra $A_1$ and $D$\Hyph algebra $A_2$.
Let map
\ShowEq{linear map from A1 to A2}
be linear map.
Then
\ShowEq{linear map, 0, D algebra}
\end{theorem}
\begin{proof}
Corollary of equation
\ShowEq{linear map, 0, D algebra, 1}
\end{proof}

\ePrints{2010.05.23}
\ifx\Semafor\ValueOff
\ePrints{8433-5163,1011.3102}
\ifx\Semafor\ValueOff
\Section{Polylinear Mapping of Algebra}
\fi

\begin{definition}
\label{definition: polylinear map of algebras}
Let $D$ be the commutative associative ring.
Let $A_1$, ..., $A_n$ be
$D$\Hyph algebras and $S$ be $D$\Hyph module.
We call map
\ShowEq{polylinear map of algebras}
\AddIndex{polylinear mapping of algebras}{polylinear map of algebras}
$A_1$, ..., $A_n$
into module
$S$,
if
\ShowEq{polylinear map of algebras, 1}
Let us denote
\ShowEq{set polylinear mappings, algebra}
set of polylinear mappings
of algebras
$A_1$, ..., $A_n$
into module
$S$.
Let us denote
\ShowEq{set polylinear mappings An, algebra}
set of $n$\hyph linear mappings
of algebra $A$ ($A_1=...=A_n=A$)
into module
$S$.
\qed
\end{definition}

\ePrints{8433-5163}
\ifx\Semafor\ValueOff
\begin{theorem}
\label{theorem: sum of polylinear maps, algebra}
Let $D$ be the commutative associative ring.
Let $A_1$, ..., $A_n$ be
$D$\Hyph algebras and $S$ be $D$\Hyph module.
Let mappings
\ShowEq{sum of polylinear maps, algebra}
be polylinear mappings.
Then mapping $f+g$
defined by equation
\ShowEq{sum of polylinear maps, 0, algebra}
is polylinear.
\end{theorem}
\begin{proof}
Statement of theorem follows from chains of equations
\ShowEq{sum of polylinear maps, 1, algebra}
\end{proof}

\begin{corollary}
Consider algebra $A_1$ and algebra $A_2$.
Let mappings
\ShowEq{sum of linear maps, algebra}
be linear mappings.
Then mapping $f+g$
defined by equation
\ShowEq{sum of linear maps, 0, algebra}
is linear.
\qed
\end{corollary}

\begin{theorem}
\label{theorem: polylinear map times scalar, algebra}
Let $D$ be the commutative associative ring.
Let $A_1$, ..., $A_n$ be
$D$\Hyph algebras and $S$ be $D$\Hyph module.
Let mapping
\ShowEq{polylinear map of algebras}
be polylinear mapping.
Then mapping
\ShowEq{linear map times scalar, algebra}
defined by equation
\ShowEq{linear map times scalar, 0, algebra}
is polylinear.
This holds
\ShowEq{linear map times scalar, pf, algebra}
\end{theorem}
\begin{proof}
Statement of theorem follows from chains of equations
\ShowEq{polylinear map times scalar, 1, algebra}
\end{proof}

\begin{corollary}
Consider algebra $A_1$ and algebra $A_2$.
Let mapping
\ShowEq{linear map from A1 to A2}
be linear mapping.
Then mapping
\ShowEq{linear map times scalar, algebra}
defined by equation
\ShowEq{linear map times scalar, 0, algebra}
is linear.
This holds
\ShowEq{linear map times scalar, pf, algebra}
\qed
\end{corollary}

\begin{theorem}
\label{theorem: module of polylinear mappings}
Let $D$ be the commutative associative ring.
Let $A_1$, ..., $A_n$ be
$D$\Hyph algebras and $S$ be $D$\Hyph module.
The set \ShowEq{module of polylinear mappings}
is a $D$\Hyph module.
\end{theorem}
\begin{proof}
The theorem \ref{theorem: sum of polylinear maps, algebra}
determines the sum of polylinear mappings
into $D$\Hyph module $S$.
Let
\ShowEq{module of polylinear mappings, 1}
For any
\ShowEq{module of polylinear mappings, 2}
Therefore, sum of polylinear mappings is commutative and associative.

The mapping $z$ defined by equation
\ShowEq{module of polylinear mappings, 3}
is zero of addition, because
\ShowEq{module of polylinear mappings, 4}
For a given mapping $f$ a mapping $g$
defined by equation
\ShowEq{module of polylinear mappings, 5}
satisfies to equation
\ShowEq{module of polylinear mappings, 6}
because
\ShowEq{module of polylinear mappings, 7}
Therefore, the set $\mathcal L(A_1;A_2)$
is an Abelian group.

From the theorem
\ref{theorem: polylinear map times scalar, algebra},
it follows that the representation
of the ring $D$ in the Abelian group $\mathcal L(A_1,...,A_n;S)$ is defined.
Since the ring $D$ has unit,
\ePrints{1011.3102}
\ifx\Semafor\ValueOff
then, according to the theorem
\ref{theorem: effective representation of the ring},
\else
then, according to the theorem
\xRef{8433-5163}{theorem: effective representation of the ring},
\fi
specified representation is effective.
\end{proof}

\begin{corollary}
Let $D$ be commutative ring with unit.
Consider $D$\Hyph algebra $A_1$ and $D$\Hyph algebra $A_2$.
The set $\mathcal L(A_1;A_2)$
is an $D$\Hyph module.
\qed
\end{corollary}
\fi
\fi

\ePrints{1011.3102}
\ifx\Semafor\ValueOff
\Section{Algebra \texorpdfstring{$\mathcal L(A;A)$}{L(A;A)}}

\begin{theorem}
\label{theorem: product of linear mapping, algebra}
Let $A$, $B$, $C$ be algebras over commutative ring $D$.
Let $f$ be linear mapping from $D$\Hyph algebra $A$ into $D$\Hyph algebra $B$.
Let $g$ be linear mapping from $D$\Hyph algebra $B$ into $D$\Hyph algebra $C$.
The mapping $g\circ f$ defined by diagram
\ShowEq{product of linear mapping, algebra}
is linear mapping from $D$\Hyph algebra $A$ into $D$\Hyph algebra $C$.
\end{theorem}
\begin{proof}
The proof of the theorem follows from chains of equations
\ShowEq{product of linear mapping, algebra, 1}
\end{proof}
\ePrints{2010.05.23}
\ifx\Semafor\ValueOn
\newpage
\fi

\begin{theorem}
\label{theorem: product of linear mapping, algebra, f}
Let $A$, $B$, $C$ be algebras over the commutative ring $D$.
Let $f$ be a linear mapping from $D$\Hyph algebra $A$ into $D$\Hyph algebra $B$.
The mapping $f$ generates a linear mapping
\ShowEq{product of linear mapping, algebra, f}
\ShowEq{product of linear mapping, algebra, f diagram}
\end{theorem}
\begin{proof}
The proof of the theorem follows from chains of equations\footnote{We use
following definitions of operations over mappings
\ShowEq{operations over linear map}
}
\ShowEq{product of linear mapping, algebra, f1}
\end{proof}
\ePrints{2010.05.23}
\ifx\Semafor\ValueOn
\newpage
\fi

\begin{theorem}
\label{theorem: product of linear mapping, algebra, g}
Let $A$, $B$, $C$ be algebras over the commutative ring $D$.
Let $g$ be a linear mapping from $D$\Hyph algebra $B$ into $D$\Hyph algebra $C$.
The mapping $g$ generates a linear mapping
\ShowEq{product of linear mapping, algebra, g}
\ShowEq{product of linear mapping, algebra, g diagram}
\end{theorem}
\begin{proof}
The proof of the theorem follows from chains of equations\footnote{We use
following definitions of operations over mappings
\ShowEq{operations over linear map}
}
\ShowEq{product of linear mapping, algebra, g1}
\end{proof}
%\fi
\ePrints{2010.05.23}
\ifx\Semafor\ValueOn
\newpage
\fi

\begin{theorem}
\label{theorem: product of linear mapping, algebra, gf}
Let $A$, $B$, $C$ be algebras over the commutative ring $D$.
The mapping
\ShowEq{product of linear mapping, algebra, gf}
is bilinear mapping.
\ePrints{2010.05.23}
\ifx\Semafor\ValueOn
\qed
\end{theorem}
\else
\end{theorem}
\begin{proof}
The theorem follows from theorems
\ref{theorem: product of linear mapping, algebra, f},
\ref{theorem: product of linear mapping, algebra, g}.
\end{proof}
\fi%2010.05.23
\fi

\begin{theorem}
\label{theorem: module L(A;A) is algebra}
Let $A$ be algebra over commutative ring $D$.
$D$\Hyph module $\mathcal L(A;A)$ equiped by product
\ShowEq{module L(A;A) is algebra}
\ePrints{1011.3102}
\ifx\Semafor\ValueOn
\ShowEq{product of linear mapping, algebra 1}
\fi
is algebra over $D$.
\end{theorem}
\begin{proof}
\ePrints{1011.3102}
\ifx\Semafor\ValueOff
The theorem follows from definition
\ref{definition: algebra over ring}
and theorem
\ref{theorem: product of linear mapping, algebra, gf}.
\else
See the proof of the theorem
\xRef{8433-5163}{theorem: module L(A;A) is algebra}.
\fi
\end{proof}
\ePrints{2010.05.23}
\ifx\Semafor\ValueOn
\newpage
\fi

\input{\FilePrefix Tensor.Algebra.English}

\ePrints{1011.3102}
\ifx\Semafor\ValueOff
\Section{Linear Mapping into Associative Algebra}
\else
\section{Linear Mapping into Associative Algebra}
\fi

\ePrints{2010.05.23}
\ifx\Semafor\ValueOff
\ePrints{1011.3102}
\ifx\Semafor\ValueOff
\begin{theorem}
\label{theorem: bilinear mapping A LA}
Consider $D$\Hyph algebras $A_1$ and $A_2$.
For given mapping $f\in\mathcal L(A_1;A_2)$,
the mapping
\ShowEq{bilinear mapping A LA}
is bilinear mapping.
\end{theorem}
\begin{proof}
The statement of theorem follows from chains of equations
\ShowEq{bilinear mapping A LA, 1}
\end{proof}
\fi

\begin{theorem}
\label{theorem: linear map AA LAA}
Consider $D$\Hyph algebras $A_1$ and $A_2$.
For given mapping $f\in\mathcal L(A_1;A_2)$,
there exists linear mapping
\ShowEq{linear map AA LAA}
defined by the equation
\ShowEq{linear map AA LAA, 1}
\end{theorem}
\begin{proof}
\ePrints{1011.3102}
\ifx\Semafor\ValueOn
See the proof of the theorems
\xRef{8433-5163}{theorem: bilinear mapping A LA}
and
\xRef{8433-5163}{theorem: linear map AA LAA}.
\else
The statement of the theorem is corollary of theorems
\ref{theorem: tensor product and polylinear mapping},
\ref{theorem: bilinear mapping A LA}.
\fi
\end{proof}

\begin{theorem}
\label{theorem: algebra A2 representation in LA}
Consider $D$\Hyph algebras $A_1$ and $A_2$.
Let us define product in algebra $\ATwo$
according to rule
\ShowEq{product in algebra AA}
A linear mapping
\ShowEq{representation A2 in LA}
defined by the equation
\ShowEq{representation A2 in LA, 1}
is representation\footnote{See the definition of representation
of $\Omega$\Hyph algebra in the definition
\xRef{0912.3315}{definition: Tstar representation of algebra}.
}
of algebra $\ATwo$
in module $\mathcal L(A_1;A_2)$.
\end{theorem}
\begin{proof}
\ePrints{1011.3102}
\ifx\Semafor\ValueOn
See the proof of the theorem
\xRef{8433-5163}{theorem: algebra A2 representation in LA}.
\else
According to theorem
\ref{theorem: linear map times constant, algebra},
mapping \EqRef{representation A2 in LA, 1}
is transformation of module $\mathcal L(A_1;A_2)$.
For a given tensor $c\in \ATwo$,
a transformation $h(c)$ is a linear transformation
of module $\mathcal L(A_1;A_2)$, because
\ShowEq{representation A2 in LA, 2}
According to theorem
\ref{theorem: linear map AA LAA},
mapping \EqRef{representation A2 in LA, 1}
is linear mapping.

Let
\ShowEq{algebra A2 representation in LA, 1}
According to the theorem
\ref{theorem: linear map AA LAA}
\ShowEq{algebra A2 representation in LA, 2}
Therefore, according to the theorem
\ref{theorem: linear map AA LAA}
\ShowEq{algebra A2 representation in LA, 3}
Since the product in algebra $A_2$ is associative, then
\ShowEq{algebra A2 representation in LA, 4}
Therefore, since we define the product in algebra $\ATwo$
according to equation
\EqRef{product in algebra AA},
then the mapping \EqRef{representation A2 in LA}
is morphism of algebras.
According to the definition
\xRef{0912.3315}{definition: Tstar representation of algebra}
mapping \EqRef{representation A2 in LA, 1}
is a representation of the algebra $\ATwo$
in the module $\mathcal L(A_1;A_2)$.
\fi
\end{proof}

\begin{theorem}
\label{theorem: representation of algebra A2 in LA}
Consider $D$\Hyph algebra $A$.
Let us define product in algebra $A\otimes A$
according to rule
\EqRef{product in algebra AA}.
A representation
of algebra $A\otimes A$
\ShowEq{representation AA in LA}
in module $\mathcal L(A;A)$
defined by the equation
\ShowEq{representation AA in LA, 1}
allows us to identify tensor
\ShowEq{product in algebra AA 1}
and mapping \ShowEq{product in algebra AA 2} where
\ShowEq{product in algebra AA 3}
is identity mapping.
\end{theorem}
\begin{proof}
\ePrints{1011.3102}
\ifx\Semafor\ValueOn
See the proof of the theorem
\xRef{8433-5163}{theorem: representation of algebra A2 in LA}.
\else
According to the theorem
\ref{theorem: linear map AA LAA},
the mapping $f\in\mathcal L(A;A)$ and the tensor $d\in A\otimes A$
generate the mapping
\ShowEq{tensor and mapping in A, algebra}
If we assume $f=\delta$, $d=a\otimes b$,
then the equation \EqRef{tensor and mapping in A, algebra}
gets form
\ShowEq{tensor and mapping in A, 1, algebra}
If we assume
\ShowEq{tensor and mapping in A, 2, algebra}
then comparison of equations
\EqRef{tensor and mapping in A, 1, algebra} and
\EqRef{tensor and mapping in A, 2, algebra}
gives a basis to identify the action of the tensor
$a\otimes b$ and
transformation $(a\otimes b)\circ\delta$.
\fi
\end{proof}

From the theorem \ref{theorem: representation of algebra A2 in LA},
it follows that we can consider the mapping
\EqRef{representation A2 in LA, 1}
as the product of mappings
$a\otimes b$ and $f$.
\ePrints{1011.3102}
\ifx\Semafor\ValueOff
The tensor $a\in \ATwo$ is
\AddIndex{nonsingular}{nonsingular tensor, algebra},
if there exists
the tensor $b\in \ATwo$ such that $a\circ b=1\otimes 1$.

\begin{definition}
Consider the representation of algebra $\ATwo$
in the module $\mathcal L(A_1;A_2)$.\footnote{The definition is made by
analogy with the definition
\xRef{0912.3315}{definition: orbit of Tstar representation of group}.}
The set
\ShowEq{orbit of linear mapping}
is called \AddIndex{orbit of linear mapping}
{orbit of linear mapping}
$f\in\mathcal L(A_1;A_2)$
\qed
\end{definition}
\fi
\fi

\begin{theorem}
\label{theorem: h generated by f, associative algebra}
Consider $D$\Hyph algebra $A_1$
and associative $D$\Hyph algebra $A_2$.
Consider the representation of algebra $\ATwo$
in the module $\mathcal L(A_1;A_2)$.
The mapping
\ShowEq{linear map h A1 A2}
generated by the mapping
\ShowEq{linear map f A1 A2}
has form
\ShowEq{h generated by f, associative algebra}
\ePrints{2010.05.23}
\ifx\Semafor\ValueOn
\qed
\end{theorem}
\newpage
\else
\end{theorem}
\begin{proof}
\ePrints{1011.3102}
\ifx\Semafor\ValueOn
See the proof of the theorem
\xRef{8433-5163}{theorem: h generated by f, associative algebra}
\else%
We can represent any tensor $a\in\ATwo$ in the form
\ShowEq{h generated by f 1, associative algebra}
According to the theorem \ref{theorem: algebra A2 representation in LA},
the mapping \EqRef{representation A2 in LA, 1}
is linear.
This proofs the statement of the theorem.
\fi%1011.3102
\end{proof}

\ePrints{1011.3102}
\ifx\Semafor\ValueOff
\begin{theorem}
\label{theorem: orbit of linear mapping}
Let $A_2$ be algebra with unit $e$.
Let
\ShowEq{orbit of linear mapping, 02}
be a nonsingular tensor.
Orbits of linear mappings
$f\in\mathcal L(A_1;A_2)$
and \ShowEq{orbit of linear mapping, 01}
coincide
\ShowEq{orbit of linear mapping, f = g}
\end{theorem}
\begin{proof}
If
\ShowEq{orbit of linear mapping, 1}
then there exists
\ShowEq{orbit of linear mapping, 2}
such that
\ShowEq{orbit of linear mapping, 3}
In that case
\ShowEq{orbit of linear mapping, 4}
Therefore,
\ShowEq{orbit of linear mapping, 5}
\ShowEq{orbit of linear mapping, g in f}

Since $a$ is nonsingular tensor, then
\ShowEq{orbit of linear mapping, 6}
If
\ShowEq{orbit of linear mapping, 5}
then there exists
\ShowEq{orbit of linear mapping, 2}
such that
\ShowEq{orbit of linear mapping, 7}
From equations
\EqRef{orbit of linear mapping, 6},
\EqRef{orbit of linear mapping, 7},
it follows that
\ShowEq{orbit of linear mapping, 8}
Therefore,
\ShowEq{orbit of linear mapping, 1}
\ShowEq{orbit of linear mapping, f in g}

\EqRef{orbit of linear mapping, f = g}
follows from equations
\EqRef{orbit of linear mapping, g in f},
\EqRef{orbit of linear mapping, f in g}.
\end{proof}

From the theorem
\ref{theorem: orbit of linear mapping}, it also follows that if
\ShowEq{orbit of linear mapping, 01}
and \ShowEq{orbit of linear mapping, 02}
is a singular tensor, then
relationship
\EqRef{orbit of linear mapping, g in f} is true.
However, the main result of the theorem
\ref{theorem: orbit of linear mapping}
is that the representations of the algebra
$\ATwo$ in module $\mathcal L(A_1;A_2)$
generates an equivalence in the module $\mathcal L(A_1;A_2)$.
If we successfully choose the representatives of each equivalence class,
then the resulting set will be generating set
of considered representation.\footnote{Generating set
of representation is defined in definition
\xRef{0912.3315}{definition: generating set of representation}.}
\fi
\fi

\ePrints{1011.3102}
\ifx\Semafor\ValueOff
\Section{Linear Mapping into Free
Finite Dimensional Associative Algebra}

\begin{theorem}
Let $A_1$ be algebra
over the ring $D$.
Let $A_2$ be free finite dimensional associative algebra
over the ring $D$.
Let $\Basis e$ be basis of the algebra $A_2$ over the ring $D$.
The mapping
\ShowEq{standard representation of mapping A1 A2, 1, associative algebra}
generated by the mapping $f\in\mathcal L(A_1;A_2)$
through the tensor $a\in\ATwo$, has the standard representation
\ShowEq{standard representation of mapping A1 A2, 2, associative algebra}
\ePrints{2010.05.23}
\ifx\Semafor\ValueOn
\qed
\end{theorem}
\newpage
\else%2010.05.23
\end{theorem}
\begin{proof}
According to theorem
\ref{theorem: standard component of tensor, algebra},
the standard representation of the tensor $a$ has form
\ShowEq{standard representation of mapping A1 A2, 3, associative algebra}
The equation
\EqRef{standard representation of mapping A1 A2, 2, associative algebra}
follows from equations
\EqRef{standard representation of mapping A1 A2, 1, associative algebra},
\EqRef{standard representation of mapping A1 A2, 3, associative algebra}.
\end{proof}
\fi%2010.05.23

\begin{theorem}
\label{theorem: coordinates of mapping A1 A2, algebra}
Let $\Basis e_1$ be basis of the free finite dimensional
$D$\Hyph algebra $A_1$.
Let $\Basis e_2$ be basis of the free finite dimensional associative
$D$\Hyph algebra $A_2$.
Let $C_{2\cdot}{}_{\gi{kl}}^{\gi p}$ be structural constants of algebra $A_2$.
Coordinates of the mapping
\ShowEq{coordinates of mapping A1 A2, 1}
generated by the mapping $f\in\mathcal(A_1;A_2)$
through the tensor $a\in\ATwo$ and its standard components
are connected by the equation
\ShowEq{coordinates of mapping A1 A2, 2, associative algebra}
\ePrints{2010.05.23}
\ifx\Semafor\ValueOn
\qed
\end{theorem}
\newpage
\else%2010.05.23
\end{theorem}
\begin{proof}
Relative to bases
$\Basis e_1$ and $\Basis e_2$, linear mappings $f$ and $g$ have form
\ShowEq{coordinates of mappings f g, associative algebra}
From equations
\EqRef{coordinates of mapping f, associative algebra},
\EqRef{coordinates of mapping g, associative algebra},
\EqRef{standard representation of mapping A1 A2, 2, associative algebra}
it follows that
\ShowEq{coordinates of mapping A1 A2, 3, associative algebra}
Since vectors $\Vector e_{2\cdot\gi k}$
are linear independent and $x^{\gi i}$ are arbitrary,
then the equation
\EqRef{coordinates of mapping A1 A2, 2, associative algebra}
follows from the equation
\EqRef{coordinates of mapping A1 A2, 3, associative algebra}.
\end{proof}

\begin{theorem}
\label{theorem: linear map over ring, matrix}
Let $D$ be field.
Let $\Basis e_1$ be basis of the free finite dimensional
$D$\Hyph algebra $A_1$.
Let $\Basis e_2$ be basis of the free finite dimensional associative
$D$\Hyph algebra $A_2$.
Let $C_{2\cdot}{}_{\gi{kl}}^{\gi p}$ be structural constants of algebra $A_2$.
Consider matrix
\ShowEq{linear map over ring, matrix}
whose rows and columns are indexed by 
\ShowEq{linear map over ring, row of matrix}
and
\ShowEq{linear map over ring, column of matrix},
respectively.
If
matrix $\mathcal B$ is nonsingular,
then, for given coordinates of linear transformation
\ShowEq{linear map over ring, left side}
and for mapping
\ShowEq{linear map over ring, f=1},
the system of linear equations
\EqRef{coordinates of mapping A1 A2, 2, associative algebra}
with standard components of this transformation
$g^{\gi{kr}}$ has the unique
solution.

If
matrix $\mathcal B$ is singular,
then
the equation
\ShowEq{linear map over ring, determinant=0, 1}
is the condition for the existence of solutions
of the system of linear equations
\EqRef{coordinates of mapping A1 A2, 2, associative algebra}.
In such case the system of linear equations
\EqRef{coordinates of mapping A1 A2, 2, associative algebra}
has infinitely many solutions and there exists linear
dependence between values $g_{\gi m}^{\gi k}$.
\end{theorem}
\begin{proof}
The statement of the theorem
is corollary of
the theory of linear equations over field.
\end{proof}
\fi%2010.05.23

\begin{theorem}
\label{theorem: linear map over ring, determinant=0}
Let $A$ be free finite dimensional associative algebra
over the field $D$.
Let $\Basis e$ be the basis of the algebra $A$ over the field $D$.
Let $C_{\gi{kl}}^{\gi p}$
be structural constants of algebra $A$.
Let matrix \EqRef{linear map over ring, matrix} be singular.
Let the linear mapping $f\in\mathcal L(A;A)$ be nonsingular.
If coordinates of linear transformations $f$ and $g$
satisfy to the equation
\ShowEq{linear map over ring, determinant=0, 2}
then the system of linear equations
\ShowEq{coordinates of mapping A}
has infinitely many solutions.
\ePrints{2010.05.23}
\ifx\Semafor\ValueOn
\qed
\end{theorem}
\newpage
\else%2010.05.23
\end{theorem}
\begin{proof}
According to the equation
\EqRef{linear map over ring, determinant=0, 2}
and the theorem
\ref{theorem: linear map over ring, matrix},
the system of linear equations
\ShowEq{coordinates of mapping A1 A2, 4f}
has infinitely many solutions
corresponding to linear mapping
\ShowEq{coordinates of mapping A1 A2, 5f}
According to the equation
\EqRef{linear map over ring, determinant=0, 2}
and the theorem
\ref{theorem: linear map over ring, matrix},
the system of linear equations
\ShowEq{coordinates of mapping A1 A2, 4g}
has infinitely many solutions
corresponding to linear mapping
\ShowEq{coordinates of mapping A1 A2, 5g}
Mappings $f$ and $g$ are generated by the mapping $\delta$.
According to the theorem
\ref{theorem: orbit of linear mapping},
the mapping $f$ generates the mapping $g$.
This proves the statement of the theorem.
\end{proof}
\fi%2010.05.23
\fi%1011.3102

\begin{theorem}
\label{theorem: linear mapping in L(A,A), associative algebra}
Let $A$ be free finite dimensional associative algebra
over the field $D$.
The representation of algebra $A\otimes A$ in algebra $\mathcal L(A;A)$ has
finite \AddIndex{basis}{basis of algebra L(A,A)} $\Basis I$.
\begin{enumerate}
\item \label{f in L(A,A), 1, associative algebra}
The linear mapping
\ShowEq{f in L(A,A)}
has form
\ShowEq{f in L(A,A), 1, associative algebra}
\item \label{f in L(A,A), 2, associative algebra}
Its standard representation has form
\ShowEq{f in L(A,A), 2, associative algebra}
\end{enumerate}
\ePrints{2010.05.23}
\ifx\Semafor\ValueOn
\qed
\end{theorem}
\newpage
\else%2010.05.23
\end{theorem}
\begin{proof}
\ePrints{1011.3102}
\ifx\Semafor\ValueOn
See the proof of the theorem
\xRef{8433-5163}{theorem: linear mapping in L(A,A), associative algebra}
\else%1011.3102
From the theorem \ref{theorem: linear map over ring, determinant=0},
it follows that if matrix $\mathcal B$ is singular
and the mapping $f$ satisfies to the equation
\ShowEq{linear map over ring, determinant=0, 3}
then the mapping $f$ generates the same set of mappings
that is generated by the mapping $\delta$.
Therefore, to build the basis of representation
of the algebra $A\otimes A$ in the module $\mathcal L(A;A)$,
we must perform the following construction.

The set of solutions of system of equations
\EqRef{coordinates of mapping A}
generates a free submodule $\mathcal L$ of the module $\mathcal L(A;A)$.
We build the basis $(\Vector h_1, ..., \Vector h_k)$
of the submodule $\mathcal L$.
Then we supplement this basis by linearly independent vectors
$\Vector h_{k+1}$, ..., $\Vector h_m$, that do not belong
to the submodule $\mathcal L$ so that the set of vectors
$\Vector h_1$, ..., $\Vector h_m$
forms a basis of the module $\mathcal L(A;A)$.
The set of orbits $(A\otimes A)\circ\delta$,
$(A\otimes A)\circ\Vector h_{k+1}$, ..., $(A\otimes A)\circ\Vector h_m$
generates the module $\mathcal L(A;A)$.
Since the set of orbits is finite, we can choose the orbits
so that they do not intersect. For each orbit we can choose a
representative which generates the orbit.
\fi%1011.3102
\end{proof}
\fi%2010.05.23

\ePrints{1011.3102}
\ifx\Semafor\ValueOff
\begin{example}
For complex field, the algebra $\mathcal L(C;C)$ has basis
\ShowEq{basis L(C,C)}
For quaternion algebra, the algebra $\mathcal L(H;H)$ has basis
\ShowEq{basis L(H,H)}
\qed
\end{example}
\fi%1011.3102

\ePrints{1011.3102}
\ifx\Semafor\ValueOff
\Section{Linear Mapping into Nonassociative Algebra}
\else
\section{Linear Mapping into Nonassociative Algebra}
\fi
\label{Section: Linear Mapping into Nonassociative Algebra}

Since the product is nonassociative, we may assume
that action of $a$, $b\in A$ over the mapping $f$ may have
form either $a(fb)$, or $(af)b$.
\ePrints{2010.05.23,1011.3102}
\ifx\Semafor\ValueOff
However this assumption leads us to a rather complex structure
of the linear mapping.
To better understand how complex the structure of the linear
mapping, we begin by considering the left and right shifts
in nonassociative algebra.

\begin{theorem}
Let
\ShowEq{left shift algebra}
be mapping of left shift.
Then
\ShowEq{left shift algebra, 1}
where we introduced linear mapping
\ShowEq{left shift algebra, 2}
\end{theorem}
\begin{proof}
From the equations
\EqRef{associator of algebra},
\EqRef{left shift algebra},
it follows that
\ShowEq{left shift algebra, 3}
The equation \EqRef{left shift algebra, 1}
follows from equation \EqRef{left shift algebra, 3}.
\end{proof}
%\ShowEq{left shift algebra, 2a}

\begin{theorem}
Let
\ShowEq{right shift algebra}
be mapping of right shift.
Then
\ShowEq{right shift algebra, 1}
where we introduced linear mapping
\ShowEq{right shift algebra, 2}
\end{theorem}
\begin{proof}
From the equations
\EqRef{associator of algebra},
\EqRef{right shift algebra}
it follows that
\ShowEq{right shift algebra, 3}
The equation \EqRef{right shift algebra, 1}
follows from equation \EqRef{right shift algebra, 3}.
\end{proof}

Let
\ShowEq{f over A}
be linear mapping of the algebra $A$.
According to the theorem \ref{theorem: linear map times constant, algebra},
the mapping
\ShowEq{g over A}
is also a linear mapping.
However, it is not obvious whether we can write the mapping $g$ as
a sum of terms of type $(ax)b$ and $a(xb)$.

If $A$ is free finite dimensional algebra, then we can assume
that the linear mapping has the standard representation like\footnote{The choice
is arbitrary. We may consider the standard representation like
\ShowEq{linear map, standard representation, nonassociative algebra, 1}
Then the equation
\EqRef{coordinates of mapping A1 A2, 2, nonassociative algebra}
has form
\ShowEq{coordinates of mapping A1 A2, 21, nonassociative algebra}
I chose the expression
\EqRef{linear map, standard representation, nonassociative algebra}
because order of the factors corresponds to the order chosen in the theorem
\ref{theorem: linear mapping in L(A,A), associative algebra}.}
\ShowEq{linear map, standard representation, nonassociative algebra}
In this case we can use
the theorem \ref{theorem: linear mapping in L(A,A), associative algebra}
for mappings into nonassociative algebra.
\else%2010.05.23,1011.3102
\ePrints{2010.05.23}
\ifx\Semafor\ValueOn
\newpage
\fi%2010.05.23
\fi%2010.05.23,1011.3102

\begin{theorem}
\label{theorem: coordinates of mapping A1 A2, nonassociative algebra}
Let $\Basis e_1$ be basis of the free finite dimensional
$D$\Hyph algebra $A_1$.
Let $\Basis e_2$ be basis of the free finite dimensional nonassociative
$D$\Hyph algebra $A_2$.
Let $C_{2\cdot}{}_{\gi{kl}}^{\gi p}$ be structural constants of algebra $A_2$.
Let the mapping
\ShowEq{standard representation of mapping A1 A2, 1, nonassociative algebra}
generated by the mapping $f\in\mathcal(A_1;A_2)$
through the tensor $a\in\ATwo$, has the standard representation
\ShowEq{standard representation of mapping A1 A2, 2, nonassociative algebra}
Coordinates of the mapping
\EqRef{standard representation of mapping A1 A2, 1, nonassociative algebra}
and its standard components
are connected by the equation
\ShowEq{coordinates of mapping A1 A2, 2, nonassociative algebra}
\ePrints{2010.05.23}
\ifx\Semafor\ValueOn
\qed
\end{theorem}
\newpage
\else%2010.05.23
\end{theorem}
\begin{proof}
\ePrints{1011.3102}
\ifx\Semafor\ValueOn
See the proof of the theorem
\xRef{8433-5163}{theorem: coordinates of mapping A1 A2, nonassociative algebra}
\else%1011.3102
Relative to bases
$\Basis e_1$ and $\Basis e_2$, linear mappings $f$ and $g$ have form
\ShowEq{coordinates of mappings f g, nonassociative algebra}
From equations
\EqRef{coordinates of mapping f, nonassociative algebra},
\EqRef{coordinates of mapping g, nonassociative algebra},
\EqRef{standard representation of mapping A1 A2, 2, nonassociative algebra}
it follows that
\ShowEq{coordinates of mapping A1 A2, 3, nonassociative algebra}
Since vectors $\Vector e_{2\cdot\gi k}$
are linear independent and $x^{\gi i}$ are arbitrary,
then the equation
\EqRef{coordinates of mapping A1 A2, 2, nonassociative algebra}
follows from the equation
\EqRef{coordinates of mapping A1 A2, 3, nonassociative algebra}.
\fi%1011.3102
\end{proof}
\fi%2010.05.23

\begin{theorem}
\label{theorem: linear mapping in L(A,A), nonassociative algebra}
Let $A$ be free finite dimensional nonassociative algebra
over the ring $D$.
The representation of algebra $A\otimes A$ in algebra $\mathcal L(A;A)$ has
finite basis $\Basis I$.
\begin{enumerate}
\item \label{f in L(A,A), 1, nonassociative algebra}
The linear mapping
\ShowEq{f in L(A,A)}
has form
\ShowEq{f in L(A,A), 1, nonassociative algebra}
\item \label{f in L(A,A), 2, nonassociative algebra}
Its standard representation has form
\ShowEq{f in L(A,A), 2, nonassociative algebra}
\end{enumerate}
\ePrints{2010.05.23}
\ifx\Semafor\ValueOn
\qed
\end{theorem}
\newpage
\else%2010.05.23
\end{theorem}
\begin{proof}
\ePrints{1011.3102}
\ifx\Semafor\ValueOn
See the proof of the theorem
\xRef{8433-5163}{theorem: linear mapping in L(A,A), nonassociative algebra}
\else%1011.3102
Consider matrix
\EqRef{linear map over ring, matrix}.
If
matrix $\mathcal B$ is nonsingular,
then, for given coordinates of linear transformation
\ShowEq{linear map over ring, left side}
and for mapping
\ShowEq{linear map over ring, f=1},
the system of linear equations
\EqRef{coordinates of mapping A1 A2, 2, nonassociative algebra}
with standard components of this transformation
$g^{\gi{kr}}$ has the unique
solution.
If
matrix $\mathcal B$ is singular,
then according to the theorem
\ref{theorem: linear mapping in L(A,A), associative algebra}
there exists finite basis $\Basis I$
generating the set of linear mappings.
\fi%1011.3102
\end{proof}
\fi%2010.05.23

\ePrints{1011.3102}
\ifx\Semafor\ValueOff
Unlike the case of an associative algebra,
the set of generators $I$ in the theorem
\ref{theorem: linear mapping in L(A,A), nonassociative algebra}
is not minimal.
From the equation \EqRef{left shift algebra, 1}
it follows that the equation
\EqRef{orbit of linear mapping, 4} does not hold.
Therefore, orbits of mappings
$I_k$ do not generate an equivalence relation in the algebra $L(A;A)$.
Since we consider only mappings like $(aI_k)b$, then
it is possible that for $k\ne l$ the mapping $I_k$ generates the mapping $I_l$,
if we consider all possible operations in the algebra $A$.
Therefore, the set of generators $I_k$ of nonassociative algebra $A$
does not play such a critical role as conjugation in
complex field.
The answer to the question of how important it is the mapping $I_k$ in
nonassociative algebra requires additional research.
\else%1011.3102

\ePrints{1011.3102}
\ifx\Semafor\ValueOff
\Section{Polylinear Mapping into Associative Algebra}
\else
\section{Polylinear Mapping into Associative Algebra}
\fi

\begin{theorem}
\label{theorem: polylinear mapping A 1n A}
Let $A_1$, ..., $A_n$, $A$ be associative $D$\Hyph algebras.
Let
\ShowEq{set of mappings, L(A;A)}
\ShowEq{set of elements A}
For given transposition $\sigma$ of $n$ variables,
the mapping
\ShowEq{n linear mapping A LA}
is $n$\hyph linear mapping into algebra $A$.
\end{theorem}
\begin{proof}
The statement of theorem follows from chains of equations
\ShowEq{polylinear mapping A 1n A, 1}
\end{proof}

In the equation \EqRef{n linear mapping A LA},
as well as in other expressions of polylinear mapping,
we have convention that mapping $f_i$ has variable $x_i$ as argument.

\begin{theorem}
\label{theorem: polylinear mapping A LA}
Let $A_1$, ..., $A_n$, $A$ be associative $D$\Hyph algebras.
For given set of mappings
\ShowEq{set of mappings, L(A;A)}
the mapping
\ShowEq{n+1 linear mapping A LA}
defined by equation
\ShowEq{n+1 linear mapping A LA 1}
is $n+1$\hyph linear mapping into $D$\Hyph module $\LAAnA$.
\end{theorem}
\begin{proof}
The statement of theorem follows from chains of equations
\ShowEq{n+1 linear mapping A LA, 1}
\end{proof}

\begin{theorem}
\label{theorem: linear map An LAAnA}
Let $A_1$, ..., $A_n$, $A$ be associative $D$\Hyph algebras.
For given set of mappings
\ShowEq{set of mappings, L(A;A)}
there exists linear mapping
\ShowEq{linear map An LAAnA}
defined by the equation
\ShowEq{linear map An LAAnA, 1}
\end{theorem}
\begin{proof}
The statement of the theorem is corollary of theorems
\ref{theorem: tensor product and polylinear mapping},
\ref{theorem: polylinear mapping A LA}.
\end{proof}

\begin{theorem}
\label{theorem: polylinear mapping A LA, 1}
Let $A_1$, ..., $A_n$, $A$ be associative $D$\Hyph algebras.
For given tensor
\ShowEq{tensor over A} and given transposition $\sigma\in S_n$
the mapping
\ShowEq{n linear mapping LA LA}
defined by equation
\ShowEq{n linear mapping LA LA 1}
is $n$\hyph linear mapping into $D$\Hyph module $\LAAnA$.
\end{theorem}
\begin{proof}
The statement of theorem follows from chains of equations
\ShowEq{representation An in LAAnA, 3}
\end{proof}

\begin{theorem}
\label{theorem: linear mapping LA LA, 1}
Let $A_1$, ..., $A_n$, $A$ be associative $D$\Hyph algebras.
For given tensor
\ShowEq{tensor over A} and given transposition $\sigma\in S_n$
there exists linear mapping
\ShowEq{linear mapping LA LA}
defined by the equation
\ShowEq{linear mapping LA LA, 1}
\end{theorem}
\begin{proof}
The statement of the theorem is corollary of theorems
\ref{theorem: tensor product and polylinear mapping},
\ref{theorem: polylinear mapping A LA, 1}.
\end{proof}

\begin{theorem}
\label{theorem: polylinear map, algebra} 
Let $A$ be associative $D$\Hyph algebra.
Polylinear map
\ShowEq{polylinear map, algebra}
has form
\ShowEq{polylinear map, algebra, canonical morphism}
where $\sigma_s$ is a transposition of set of variables
\ShowEq{transposition of set of variables, algebra}
\end{theorem}
\begin{proof}
We prove statement by induction on $n$.

When $n=1$ the statement of theorem is corollary of
the statement \eqref{f in L(A,A), 1, associative algebra} of the theorem
\ref{theorem: linear mapping in L(A,A), associative algebra}.
In such case we may identify\footnote{In representation
\EqRef{polylinear map, algebra, canonical morphism}
we will use following rules.
\begin{itemize}
\item If range of any index is set
consisting of one element, then we will omit corresponding
index.
\item If $n=1$, then $\sigma_s$ is identical transformation.
We will not show such transformation in the expression.
\end{itemize}}
\ShowEq{polylinear map, associative algebra, 1, canonical morphism}

Let statement of theorem be true for $n=k-1$.
Then it is possible to represent mapping
\EqRef{polylinear map, algebra}
as
\ShowEq{polylinear map, induction on n, 1, associative algebra}
According to statement of induction polylinear mapping
$h$ has form
\ShowEq{polylinear map, induction on n, 2, associative algebra}
According to construction $h=g\circ a_k$.
Therefore, expressions $h_{t\cdot p}$
are functions of $a_k$.
Since $g\circ a_k$ is linear mapping of $a_k$,
then only one expression $h_{t\cdot p}$
is linear mapping of $a_k$, and rest expressions
$h_{t\cdot q}$
do not depend on $a_k$.

Without loss of generality, assume $p=0$.
According to the equation
\EqRef{h generated by f, associative algebra}
for given $t$
\ShowEq{polylinear map, induction on n, 3, associative algebra}
Assume $s=tr$. Let us define transposition $\sigma_s$ according to rule
\ShowEq{polylinear map, induction on n, 4, associative algebra}
Suppose
\ShowEq{polylinear map, induction on n, 5, associative algebra}
We proved step of induction.
\end{proof}

\begin{definition}
\begin{sloppypar}
Expression
\ShowEq{component of polylinear map, associative algebra}
in equation \EqRef{polylinear map, algebra, canonical morphism}
is called
\AddIndex{component of polylinear map $f$}
{component of polylinear map, associative algebra}.
\qed
\end{sloppypar}
\end{definition}

\begin{theorem}
\label{theorem: representation module An in LAAnA}
Consider $D$\Hyph algebras $A$.
\ifx\texFuture\Defined
Let us define product in algebra $\AOn$
according to rule
\ShowEq{product in algebra An}
\fi
A linear mapping
\ShowEq{representation An in LAnA}
defined by the equation
\ShowEq{representation An in LAnA, 1}
is representation\footnote{See the definition of representation
of $\Omega$\Hyph algebra in the definition
\xRef{0912.3315}{definition: Tstar representation of algebra}.}
of algebra
\ShowEq{representation An in LAnA, 2}
in $D$\Hyph module $\LAnA$.
\end{theorem}
\begin{proof}
According to the theorems
\ref{theorem: linear mapping in L(A,A), associative algebra},
\ref{theorem: polylinear map, algebra},
we can represent $n$\hyph linear mapping as sum of terms
\EqRef{n linear mapping A LA},
where $f_i$, $i=1$, ..., $n$, are generators of representation
\EqRef{representation A2 in LA}.
Let us write the term $s$ of the expression
\EqRef{polylinear map, algebra, canonical morphism}
as
\ShowEq{n linear mapping A LA, 1}
where
\ShowEq{n linear mapping A LA, 2}
Let us assume
\ShowEq{n linear mapping A LA, 3}
in equation
\EqRef{n linear mapping A LA, 1}.
Therefore, according to theorem
\ref{theorem: linear map An LAAnA},
mapping \EqRef{representation An in LAnA, 1}
is transformation of module $\LAnA$.
For a given tensor $c\in \AOn$
and given transposition $\sigma\in S_n$,
a transformation $h(c,\sigma)$ is a linear transformation
of module $\LAnA$
according to the theorem \ref{theorem: linear mapping LA LA, 1}.
According to the theorem
\ref{theorem: linear map An LAAnA},
mapping \EqRef{representation An in LAnA, 1}
is linear mapping.
\ifx\texFuture\Defined
Let
\ShowEq{algebra An representation in LAnA, 1}
According to the theorem
\ref{theorem: linear mapping LA LA, 1}
\ShowEq{algebra An representation in LAnA, 2}
Therefore, according to the theorem
\ref{theorem: linear mapping LA LA, 1}
\ShowEq{algebra An representation in LAnA, 3}
Since the product in algebra $A$ is associative, then
\ShowEq{algebra An representation in LAnA, 4}
Therefore, if we define product in algebra $\AOn$
according to equation
\EqRef{product in algebra An},
then the mapping \EqRef{representation An in LAnA, 1}
is morphism of algebras.
\fi
According to the definition
\xRef{0912.3315}{definition: Tstar representation of algebra}
mapping \EqRef{representation An in LAnA, 1}
is a representation of the algebra
\ShowEq{representation An in LAnA, 2}
in the module $\LAnA$.
\end{proof}

\begin{theorem}
\label{theorem: representation of algebra An in LAnA}
Consider $D$\Hyph algebra $A$.
A representation
\ShowEq{representation An in LAnA}
of algebra $\AOn$
in module $\LAnA$
defined by the equation
\EqRef{representation An in LAnA, 1}
allows us to identify tensor
\ShowEq{product in algebra An 1}
and transposition $\sigma\in S^n$
with mapping \ShowEq{product in algebra An 2} where
\ShowEq{product in algebra AA 3}
is identity mapping.
\end{theorem}
\begin{proof}
If we assume $f_i=\delta$,
\ShowEq{n tensor and mapping in A, 3, algebra}
in the equation \EqRef{linear map An LAAnA, 1},
then the equation \EqRef{linear map An LAAnA, 1}
gets form
\ShowEq{n tensor and mapping in A, 1, algebra}
If we assume
\ShowEq{n tensor and mapping in A, 2, algebra}
then comparison of equations
\EqRef{n tensor and mapping in A, 1, algebra} and
\EqRef{n tensor and mapping in A, 2, algebra}
gives a basis to identify the action of the tensor
\ShowEq{n tensor and mapping in A, 3, algebra}
and transposition $\sigma\in S^n$
with mapping \EqRef{product in algebra An 2}.
\end{proof}

Instead of notation $(a_0\otimes...a_n,\sigma)$, we also use
notation
\[
a_0\otimes_{\sigma(1)}...\otimes_{\sigma(n)}a_n
\]
when we want to show order of arguments in expression.
For instance, the following expressions are equivalent
\ShowEq{example of tensor notation}

\ePrints{1011.3102}
\ifx\Semafor\ValueOff
\Section{Polylinear Map into Free
Finite Dimensional Associative Algebra}
\else
\section{Polylinear Map into Free
Finite Dimensional Associative Algebra}
\fi

\begin{theorem}
Let $A$ be free finite dimensional associative algebra
over the ring $D$.
Let $\Basis I$ be basis of algebra $\mathcal L(A;A)$.
Let $\Basis e$ be the basis of the algebra $A$ over the ring $D$.
\AddIndex{Standard representation of polylinear mapping into associative algebra}
{polylinear map, standard representation, associative algebra}
has form
\ShowEq{polylinear map, associative algebra, standard representation}
Index $t$ enumerates every possible transpositions
$\sigma_t$ of the set of variables
\VarSet.
Expression
$\ShowSymbol{standard component of polylinear map, associative algebra}$
in equation \EqRef{polylinear map, associative algebra, standard representation}
is called
\AddIndex{standard component of polylinear mapping $f$}
{standard component of polylinear map, division ring}.
\end{theorem}
\begin{proof}
We change index $s$ in the equation
\EqRef{polylinear map, algebra, canonical morphism}
so as to group the terms with the same set of generators $I_k$.
Expression \EqRef{polylinear map, algebra, canonical morphism}
gets form
\ShowEq{polylinear map, algebra, canonical morphism, 1}
We assume that the index $s$ takes values depending on $k_1$, ..., $k_n$.
Components of polylinear map $f$
have expansion
\ShowEq{polylinear map, associative algebra, components extention}
relative to basis $\Basis e$.
If we substitute \EqRef{polylinear map, associative algebra, components extention}
into \EqRef{polylinear map, algebra, canonical morphism},
we get
\ShowEq{polylinear map, associative algebra, standard representation, 1}
Let us consider expression
\ShowEq{polylinear map, associative algebra, standard representation, 2}
The right-hand side is supposed to be the sum of the terms
with the index $s$, for which the transposition $\sigma_s$ is the same.
Each such sum
has a unique index $t$.
If we substitute expression
\EqRef{polylinear map, associative algebra, standard representation, 2}
into equation \EqRef{polylinear map, associative algebra, standard representation, 1}
we get equation \EqRef{polylinear map, associative algebra, standard representation}.
\end{proof}

\begin{theorem}
\label{theorem: polylinear map over field, associative algebra}
Let $A$ be free finite dimensional associative algebra
over the ring $D$.
Let $\Basis e$ be the basis of the algebra $A$ over the ring $D$.
Polylinear map
\EqRef{polylinear map, algebra}
can be represented as $D$-valued form of degree $n$
over ring $D$\footnote{We proved the theorem by analogy with
theorem in \citeBib{Rashevsky}, p. 107, 108}
\ShowEq{polylinear map over ring, associative algebra}
where
\ShowEq{polylinear map over ring, coordinates, associative algebra}
and values $f_{\gi{i_1...i_n}}$
are coordinates of $D$-valued covariant tensor.
\end{theorem}
\begin{proof}
According to the definition
\ref{definition: polylinear map of algebras},
the equation
\EqRef{polylinear map over ring, associative algebra}
follows from the chain of equations
\ShowEq{polylinear map over ring, 1, associative algebra}
Let $\Basis e'$ be another basis. Let
\ShowEq{polylinear map over ring, associative algebra, change basis}
be transformation, mapping basis $\Basis e$ into
basis $\Basis e'$.
From equations \EqRef{polylinear map over ring, associative algebra, change basis}
and \EqRef{polylinear map over ring, coordinates, associative algebra}
it follows
\ShowEq{polylinear map over ring, associative algebra, change coordinates}
From equation \EqRef{polylinear map over ring, associative algebra, change coordinates}
the tensor law of transformation of coordinates of polylinear map follows.
From equation \EqRef{polylinear map over ring, associative algebra, change coordinates}
and theorem
\ePrints{1011.3102}
\ifx\Semafor\ValueOff
\ref{theorem: free module over ring, change basis}
\else
\xRef{8433-5163}{theorem: free module over ring, change basis}
\fi
it follows that value of the mapping $f\circ(\Vector a_1,...,\Vector a_n)$ does not depend from
choice of basis.
\end{proof}

Polylinear mapping
\EqRef{polylinear map, algebra}
is \AddIndex{symmetric}
{polylinear map symmetric, associative algebra}, if
\ShowEq{polylinear map symmetric, associative algebra, definition}
for any transposition $\sigma$ of set \VarSet.

\begin{theorem}
\label{theorem: polylinear map symmetric, associative algebra}
If polyadditive map $f$ is symmetric,
then
\ShowEq{polylinear map symmetric, associative algebra}
\end{theorem}
\begin{proof}
Equation
\EqRef{polylinear map symmetric, associative algebra}
follows from equation
\ShowEq{polylinear map symmetric, 1, associative algebra}
\end{proof}

Polylinear mapping
\EqRef{polylinear map, algebra}
is \AddIndex{skew symmetric}
{polylinear map skew symmetric, associative algebra}, if
\ShowEq{polylinear map skew symmetric, associative algebra, definition}
for any transposition $\sigma$ of set \VarSet.
Here
\[
|\sigma|=
\left\{
\begin{matrix}
1&\textrm{transposition }\sigma\textrm{ even}
\\
-1&\textrm{transposition }\sigma\textrm{ odd}
\end{matrix}
\right.
\]

\begin{theorem}
\label{theorem: polylinear map skew symmetric, associative algebra}
If polylinear map $f$ is skew symmetric,
then
\ShowEq{polylinear map skew symmetric, associative algebra}
\end{theorem}
\begin{proof}
Equation
\EqRef{polylinear map skew symmetric, associative algebra}
follows from equation
\ShowEq{polylinear map skew symmetric, 1, associative algebra}
\end{proof}

\begin{theorem}
\label{theorem: coordinates of polylinear map, associative algebra over ring}
Let $A$ be free finite dimensional associative algebra
over the ring $D$.
Let $\Basis I$ be basis of algebra $\mathcal L(A;A)$.
Let $\Basis e$ be the basis of the algebra $A$ over the ring $D$.
Let polylinear over ring $D$ mapping $f$
be generated by set of mappings $(I_{k_1},...,I_{k_n})$.
Coordinates of the mapping $f$
and its components relative basis $\Basis e$
satisfy to the equation
\ShowEq{coordinates of polylinear map, associative algebra over ring}
\end{theorem}
\begin{proof}
In equation \EqRef{polylinear map, associative algebra, standard representation},
we assume
\ShowEq{coordinates of polylinear map, 1, associative algebra over ring}
Then equation \EqRef{polylinear map, associative algebra, standard representation}
gets form
\ShowEq{polylinear map, associative algebra, standard representation, 1a}
From equation \EqRef{polylinear map over ring, associative algebra}
it follows that
\ShowEq{polylinear map over field, associative algebra, 1a}
Equation
\EqRef{coordinates of polylinear map, associative algebra over ring}
follows from comparison of equations
\EqRef{polylinear map, associative algebra, standard representation, 1a}
and
\EqRef{polylinear map over ring, associative algebra}.
Equation
\EqRef{coordinates of polylinear map, associative algebra over ring, 1}
follows from comparison of equations
\EqRef{polylinear map, associative algebra, standard representation, 1a}
and
\EqRef{polylinear map over field, associative algebra, 1a}.
\end{proof}
\fi%1011.3102

%% file: Linear.Mapping.Eq.tex
%auto-ignore

\def\ATwo{A_2\otimes A_2}
\def\AOn{A^{\otimes n+1}}
\def\Eij{e_{\gi i}e_{\gi j}}
\def\VarSet{$\{a_1,...,a_n\}$}
\def\LAnA{\mathcal L(A^n;A)}
\def\LAAnA{\mathcal L(A_1,...,A_n;A)}
\def\Veb{\Vector e_i}
\def\Ii{i\in I}

\DefEq
{
\symb{\Basis{e}=(\Veb,\Ii)}1
{basis, module}
}
{basis, module}

\DefEq
{
\symb{\mathcal L(A_1,...,A_n;S)}
1{set polylinear mappings, algebra}
}
{set polylinear mappings, algebra}

\DefEq
{
$f\in \mathcal L(A_1;A_2)$, $a\otimes b$, $c\otimes d\in\ATwo$.
}
{algebra A2 representation in LA, 1}

\DefEq
{
\[
\begin{matrix}
f_1,...,f_n\in \mathcal L(A;A)
\\
a_0\otimes a_1\otimes...\otimes a_n, b_0\otimes b_1\otimes...\otimes b_n
\in \AOn
\end{matrix}
\]
}
{algebra An representation in LAnA, 1}

\DefEq
{
\[
(a\otimes b)\circ f=afb\in\mathcal L(A_1;A_2)
\]
}
{algebra A2 representation in LA, 2}

\DefEq
{
\[
\begin{matrix}
(a_0\otimes a_1\otimes...\otimes a_n,\sigma_1)\circ (f_1\otimes...\otimes f_n)
=a_0\sigma_1(f_1)a_1...\sigma_1(f_n)a_n\in\LAnA
\\
a'_1\sigma_1(f_1)a_1,a'_2\sigma_1(f_2)a_2,...,a'_n\sigma_1(f_n)a_n\in\mathcal L(A;A)
\\
\begin{matrix}
a'_1=a_0&a'_2=...=a'_n=e
\end{matrix}
\end{matrix}
\]
}
{algebra An representation in LAnA, 2}

\DefEq
{
\[
(c\otimes d)\circ((a\otimes b)\circ f)=c(afb)d
\]
}
{algebra A2 representation in LA, 3}

\DefEq
{
\begin{align*}
&(b_0\otimes b_1\otimes...\otimes b_n,\sigma_2)\circ
((a_0\otimes a_1\otimes...\otimes a_n,\sigma_1)\circ (f_1\otimes...\otimes f_n))
\\
=&(b_0\otimes b_1\otimes...\otimes b_n,\sigma_2)\circ
((a'_1\sigma_1(f_1) a_1)\otimes (a'_2\sigma_1(f_2)a_2)\otimes...\otimes (a'_n\sigma_1(f_n)a_n))
\\
=&b_0\sigma_2(a'_1\sigma_1(f_1) a_1)b_1...\sigma_2(a'_n\sigma_1(f_n)a_n)b_n
\end{align*}
}
{algebra An representation in LAnA, 3}

\DefEq
{
\[
(c\otimes d)\circ((a\otimes b)\circ f)=c(afb)d=(ca)f(bd)=
(ca\otimes bd)\circ f
\]
}
{algebra A2 representation in LA, 4}

\DefEq
{
\begin{align*}
&(b_0\otimes b_1\otimes...\otimes b_n,\sigma_2)\circ
((a_0\otimes a_1\otimes...\otimes a_n,\sigma_1)\circ (f_1\otimes...\otimes f_n))
\\
=&(b_0a_0)f_1(a_1b_1)...f_n(a_nb_n)
\\
=&((b_0a_0)\otimes(a_1b_1)\otimes...\otimes(a_nb_n))\circ(f_1\otimes...\otimes f_n)
\end{align*}
}
{algebra An representation in LAnA, 4}

\DefEq
{
\[
(a_0\otimes a_1\otimes a_2\otimes a_3,(2,1,3))\circ(x_1,x_2,x_3)=
a_0x_2a_1x_1a_2x_3a_3
\]
\[
(a_0\otimes_2 a_1\otimes_1 a_2\otimes_3 a_3)\circ(x_1,x_2,x_3)=
a_0x_2a_1x_1a_2x_3a_3
\]
}
{example of tensor notation}

\DefEq
{
\symb{\mathcal L(A^n;S)}
1{set polylinear mappings An, algebra}
}
{set polylinear mappings An, algebra}

\DefEq
{
$\delta\in\mathcal L(A;A)$
}
{product in algebra AA 3}

\DefEq
{
\begin{equation}
(ab)^{\gi k}=C^{\gi k}_{\gi{ij}}a^{\gi i}b^{\gi j}
\EqLabel{product in algebra}
\end{equation}
}
{product in algebra}

\DefEq
{
$d\circ\delta\in\mathcal L(A;A)$%
}
{product in algebra AA 2}%

\DefEquation
{
\begin{matrix}
(d,\sigma)\circ(f_1,...,f_n)&f_i=\delta\in\mathcal L(A;A)
\end{matrix}
}
{product in algebra An 2}%

\DefEq
{
$d\in A\otimes A$
}
{product in algebra AA 1}

\DefEq
{
$d\in\AOn$
}
{product in algebra An 1}

\DefEq
{
$\AOn\times S^n$
}
{representation An in LAnA, 2}

\DefEq
{
$\{a_1,...,a_n\}$
\[
\sigma_s=
\begin{pmatrix}
a_1&...&a_n
\\
\sigma_s(a_1)&...&\sigma_s(a_n)
\end{pmatrix}
\]
}
{transposition of set of variables, algebra}

\DefEquation
{
a=f\pC{s}{0}^n\ \sigma_s(I_{1\cdot s}\circ a_1)
\ f\pC{s}{1}^n\ ...\ \sigma_s(I_{n\cdot s}\circ a_n)\ f\pC{s}{n}^n
}
{polylinear map, algebra, canonical morphism}

\DefEquation
{
a=f\pC{k_1...k_n\cdot s}{0}^n\ \sigma_s(I_{k_1\cdot s}\circ a_1)
\ f\pC{k_1...k_n\cdot s}{1}^n\ ...\ \sigma_s(I_{k_n\cdot s}\circ a_n)\ f\pC{k_1...k_n\cdot s}{n}^n
}
{polylinear map, algebra, canonical morphism, 1}

\DefEquation
{
f:A^n\rightarrow A,
a=f\circ(a_1,...,a_n)
}
{polylinear map, algebra}

\DefEq
{
\begin{align}
&a^{\gi i\gi j}_1
\Vector e_{\gi i}(x\Vector e_{\gi j})
+a^{\gi i\gi j}_2
(\Vector e_{\gi i}x)\Vector e_{\gi j}
=b
\EqLabel{linear equation in nonassociative algebra, standard form}
\\
&\begin{matrix}
a^{\gi i\gi j}_1=\pC{s}{0}a^{\gi i}_1\ \pC{s}{1}a^{\gi j}_1
&\pC{s}{0}a_1=\pC{s}{0}a^{\gi i}_1e_{\gi i}
&\pC{s}{1}a_1=\pC{s}{1}a^{\gi i}_1e_{\gi i}
\\
\VirtVar
a^{\gi i\gi j}_2=\pC{s}{0}b^{\gi i}_2\ \pC{s}{1}b^{\gi j}_2
&\pC{s}{0}a_2=\pC{s}{0}a^{\gi i}_2e_{\gi i}
&\pC{s}{1}a_2=\pC{s}{1}a^{\gi i}_2e_{\gi i}
\end{matrix}
\nonumber
\end{align}
}
{linear equation in nonassociative algebra, standard form}

\DefEquation
{
\Vector e_{\gi i}\Vector e_{\gi j}=
C^{\gi k}_{\gi{ij}}\Vector e_{\gi k}
}
{product of basis vectors, algebra}

\DefEquation
{
C^{\gi p}_{\gi{ij}}
=
C^{\gi p}_{\gi{ji}}
}
{commutative product in algebra, 1}

\DefEquation
{
C^{\gi p}_{\gi{ij}}C^{\gi q}_{\gi{pk}}
=
C^{\gi q}_{\gi{ip}}C^{\gi p}_{\gi{jk}}
}
{associative product in algebra, 1}

\DefEq
{
\[
\Vector e_{\gi i}\Vector e_{\gi j}
=
\Vector e_{\gi j}\Vector e_{\gi i}
\]
}
{commutative product in algebra}

\DefEq
{
\begin{align*}
I_0\circ z&= z
\\
I_1\circ z&= \overline z
\end{align*}
}
{basis L(C,C)}

\DefEq
{
\begin{align*}
I_0\circ z&= z
\end{align*}
}
{basis L(H,H)}

\DefEquation
{
f=
(
a_{k\cdot s_k\cdot 0}
\otimes
a_{k\cdot s_k\cdot 1}
)
\circ I_k
=
\sum_{\gi k}
a_{k\cdot s_k\cdot 0}
I_k
a_{k\cdot s_k\cdot 1}
}
{f in L(A,A), 1, associative algebra}

\DefEq
{
\begin{align*}
f\circ(
a_1, ...,
a_i+ b_i, ...,
a_n)
&=
f\circ(
a_1, ...,
a_i, ...,
a_n)
+
f\circ(
a_1, ...,
b_i, ...,
a_n)
\\
f\circ(
a_1, ...,
pa_i, ...,
a_n)
&=
pf\circ(
a_1, ...,
a_i, ...,
a_n)
\end{align*}
\[
\begin{matrix}
1\le i\le n
&
a_i, b_i \in A_i
&
p\in D
\end{matrix}
\]
}
{polylinear map of algebras, 1}

\DefEq
{
\[
f:A_1\times...\times A_n\rightarrow
S
\]
}
{polylinear map of algebras}

\DefEq
{
\[
\xymatrix{
A^k\ar[rr]^f\ar@{=>}[drr]_{g\circ a_k}
& & A
\\
& &
\\
&&A^{k-1}\ar[uu]_h &
}
\]
\[
a=f\circ(a_1,...,a_k)=(g\circ a_k)\circ(a_1,...,a_{k-1})
\]
}
{polylinear map, induction on n, 1, associative algebra}

\DefEq
{
\[
a=
h\pC{t}{0}^{k-1}\ \sigma_t(I_{1\cdot t}\circ a_1)
\ h\pC{t}{1}^{k-1}\ ...
\ \sigma_t(I_{k-1\cdot t}\circ a_{k-1})\ h\pC{t}{k-1}^{k-1}
\]
}
{polylinear map, induction on n, 2, associative algebra}

\DefEq
{
\[
\begin{matrix}
f\pC{s}{p}^1=f\pC{s}{p}
&p=0, 1
\end{matrix}
\]
}
{polylinear map, associative algebra, 1, canonical morphism}

\DefEquation
{
f=
a^{k\cdot \gi{ij}}(\Vector e_{\gi i}\otimes\Vector e_{\gi j})\circ I_k
=
a^{k\cdot \gi{ij}}\Vector e_{\gi i}I_k\Vector e_{\gi j}
}
{f in L(A,A), 2, associative algebra}

\DefEquation
{
f
=
(
a_{k\cdot s_k\cdot 0}
\otimes
a_{k\cdot s_k\cdot 1}
)
\circ I_k
=
(a_{k\cdot s_k\cdot 0}
I_k)
a_{k\cdot s_k\cdot 1}
}
{f in L(A,A), 1, nonassociative algebra}

\DefEquation
{
f
=a^{k\cdot \gi{ij}}(\Vector e_{\gi i}\otimes\Vector e_{\gi j})\circ I_k
=
a^{k\cdot \gi{ij}}(\Vector e_{\gi i}I_k)\Vector e_{\gi j}
}
{f in L(A,A), 2, nonassociative algebra}

\DefEq
{
$f\in\mathcal L(A;A)$
}
{f in L(A,A)}

\DefEq
{
\symb{N(A)}0
{nucleus of algebra}
\[
\ShowSymbol{nucleus of algebra}=
\{
a\in A:
\forall b, c\in A,
(a,b,c)=(b,a,c)=(b,c,a)=0
\}
\]
}
{nucleus of algebra}

\DefEq
{
\symb{Z(A)}0
{center of algebra}
\[
\ShowSymbol{center of algebra}=
\{
a\in A:
a\in N(A),
\forall b\in A,
ab=ba
\}
\]
}
{center of algebra}

\DefEq
{
\begin{align*}
&(f+g)\circ(x_1,...,x_i+y_i,...,x_n)
\\
=&f\circ(x_1,...,x_i+y_i,...,x_n)
+g\circ(x_1,...,x_i+y_i,...,x_n)
\\
=&f\circ (x_1,...,x_i,...,x_n)+f\circ (x_1,...,y_i,...,x_n)
\\
+&g\circ (x_1,...,x_i,...,x_n)+g\circ (x_1,...,y_i,...,x_n)
\\
=&(f+g)\circ (x_1,...,x_i,...,x_n)
+(f+g)\circ (x_1,...,y_i,...,x_n)
\\[4pt]
&(f+g)\circ(x_1,...,px_i,...,x_n)
\\
=&f\circ(x_1,...,px_i,...,x_n)+g\circ(x_1,...,px_i,...,x_n)
\\
=&pf\circ (x_1,...,x_i,...,x_n)+pg\circ (x_1,...,x_i,...,x_n)
\\
=&p(f\circ (x_1,...,x_i,...,x_n)+g\circ (x_1,...,x_i,...,x_n))
\\
=&p(f+g)\circ (x_1,...,x_i,...,x_n)
\end{align*}
}
{sum of polylinear maps, 1, algebra}

\DefEq
{
\[
(f+g)\circ a=f\circ a+g\circ a
\]
}
{sum of linear maps, 0, algebra}

\DefEq
{
\[
(f+g)\circ (a_1,...,a_n)=f\circ (a_1,...,a_n)+g\circ (a_1,...,a_n)
\]
}
{sum of polylinear maps, 0, algebra}

\DefEq
{
$af$, $fb$, $a$, $b\in A_2$,
}
{linear map times constant, algebra}

\DefEq
{
$pf$, $p\in D$,
}
{linear map times scalar, algebra}

\DefEq
{
\begin{align*}
(af)\circ(x+y)
=&a\ f\circ(x+y)
=a\ (f\circ x+f\circ y)
=a\ f\circ x+a\ f\circ y
\\
=&(af)\circ x+(af)\circ y
\\
(af)\circ(px)
=&a\ f\circ(px)
=ap\ f\circ x
=pa\ f\circ x
\\
=&p\ (af)\circ x
\\
(fb)\circ(x+y)
=&f\circ(x+y)\ b
=(f\circ x+f\circ y)\ b
=f\circ x\ b+f\circ y\ b
\\
=&(fb)\circ x+(fb)\circ y
\\
(fb)\circ(px)
=&f\circ(px)\ b
=p\ f\circ x\ b
\\
=&p\ (fb)\circ x
\end{align*}
}
{linear map times constant, 1, algebra}

\DefEq
{
\begin{align*}
(pf)\circ(x_1,...,x_i+y_i,...,x_n)
=&p\ f\circ(x_1,...,x_i+y_i,...,x_n)
\\
=&p\ (f\circ (x_1,...,x_i,...,x_n)+f\circ (x_1,...,y_i,...,x_n))
\\
=&p\ f\circ (x_1,...,x_i,...,x_n)+p\ f\circ (x_1,...,y_i,...,x_n)
\\
=&(pf)\circ (x_1,...,x_i,...,x_n)+(pf)\circ (x_1,...,y_i,...,x_n)
\\
(pf)\circ(x_1,...,qx_i,...,x_n)
=&p\ f\circ(x_1,...,qx_i,...,x_n)
=pq\ f\circ (x_1,...,x_i,...,x_n)
\\
=&qp\ f\circ (x_1,...,x_n)
=q\ (pf)\circ (x_1,...,x_n)
\\
(p(qf))\circ (x_1,...,x_n)=&p\ (qf)\circ (x_1,...,x_n)
=p\ (q\ f\circ (x_1,...,x_n))
\\
=&(pq)\ f\circ (x_1,...,x_n)=((pq)f)\circ (x_1,...,x_n)
\\
((p+q)f)\circ (x_1,...,x_n)=&(p+q)\ f\circ (x_1,...,x_n)
\\
=&p\ f\circ (x_1,...,x_n)+q\ f\circ (x_1,...,x_n)
\\
=&(pf)\circ (x_1,...,x_n)+(qf)\circ (x_1,...,x_n)
\end{align*}
}
{polylinear map times scalar, 1, algebra}

\DefEq
{
$\mathcal L(A_1,...,A_n;S)$
}
{module of polylinear mappings}

\DefEq
{
\[
(\Vector e_{\gi i}\Vector e_{\gi j})\Vector e_{\gi k}
=
\Vector e_{\gi i}(\Vector e_{\gi j}\Vector e_{\gi k})
\]
}
{associative product in algebra}

\DefEq
{
\[
(
f:R_1\rightarrow R_2
,
g:A_1\rightarrow A_2
)
\]
}
{linear map from A1 R1 to A2 R2}

\DefEq
{
\[
f:A_1\rightarrow A_2
\]
}
{linear map from A1 to A2}

\DefEq
{
$f\circ 0=0$.
}
{linear map, 0, D algebra}

\DefEq
{
\[
f(a+0)
=f(a)+f(0)
\]
}
{linear map, 0, D algebra, 1}

\DefEq
{
\begin{align}
g\circ(a+b)&=g\circ a+g\circ b
\EqLabel{linear map from A1 R1 to A2 R2, 1 1}
\\
g\circ(pa)&=(f\circ p)(g\circ a)
\EqLabel{linear map from A1 R1 to A2 R2, 1 2}
\\
f\circ(pq)&=(f\circ p)(f\circ q)
\EqLabel{linear map from A1 R1 to A2 R2, 1 3}
\end{align}
\[
\begin{matrix}
a,b\in A_1
&
p,q\in R_1
\end{matrix}
\]
}
{linear map from A1 R1 to A2 R2, 1}

\DefEq
{
\begin{align*}
g(a+b)&=g(a)+g(b)
\\
g(pa)&=f(p)g(a)
\\
f(pq)&=f(p)f(q)
\end{align*}
\[
\begin{matrix}
a,b\in A_1
&
p,q\in R_1
\end{matrix}
\]
}
{linear map from A1 R1 to A2 R2, 1, old}

\DefEq
{
\begin{align*}
(af)\circ x&=a\ f\circ x
\\
(fb)\circ x&=f\circ x\ b
\end{align*}
}
{linear map times constant, 0, algebra}

\DefEq
{
\begin{align*}
((a\otimes b)\circ(f_1+f_2))\circ x
&=(a(f_1+f_2)b)\circ x
=a((f_1+f_2)\circ x)b
\\
&=a(f_1\circ x+f_2\circ x)b
=a(f_1\circ x)b+a(f_2\circ x)b
\\
&=(af_1b)\circ x+(af_2b)\circ x
\\
&=(a\otimes b)\circ f_1\circ x+(a\otimes b)\circ f_2\circ x
\\
&=((a\otimes b)\circ f_1+(a\otimes b)\circ f_2)\circ x
\\
((a\otimes b)\circ(pf))\circ x
&=(a(pf)b)\circ x
=a((pf)\circ x)b
\\
&=a(p\ f\circ x)b
=pa(f\circ x)b
\\
&=p\ (afb)\circ x
=p\ ((a\otimes b)\circ f)\circ x
\\
&=(p((a\otimes b)\circ f))\circ x
\end{align*}
}
{representation A2 in LA, 2}

\DefEq
{
\begin{align*}
&((a_0\otimes...\otimes a_n,\sigma)\circ(f_1,...,f_i+g_i,...,f_n))\circ (x_1,...,x_n)
\\
=&(a_0\sigma(f_1)a_1...\sigma(f_i+g_i)...a_{n-1}\sigma(f_n)a_n)\circ (x_1,...,x_n)
\\
=&a_0\sigma(f_1\circ x_1)a_1...\sigma((f_i+g_i)\circ x_i)...a_{n-1}\sigma(f_n\circ x_n)a_n
\\
=&a_0\sigma(f_1\circ x_1)a_1...\sigma(f_i\circ x_i+g_i\circ x_i)...a_{n-1}\sigma(f_n\circ x_n)a_n
\\
=&a_0\sigma(f_1\circ x_1)a_1...\sigma(f_i\circ x_i)...a_{n-1}\sigma(f_n\circ x_n)a_n
\\
+&a_0\sigma(f_1\circ x_1)a_1...\sigma(g_i\circ x_i)...a_{n-1}\sigma(f_n\circ x_n)a_n
\\
=&(a_0\sigma(f_1)a_1...\sigma(f_i)...a_{n-1}\sigma(f_n)a_n)\circ (x_1,...,x_n)
\\
+&(a_0\sigma(f_1)a_1...\sigma(g_i)...a_{n-1}\sigma(f_n)a_n)\circ (x_1,...,x_n)
\\
=&((a_0\otimes ...\otimes a_n,\sigma)\circ(f_1,...,f_i,...,f_n))\circ (x_1,...,x_n)
\\
+&((a_0\otimes ...\otimes a_n,\sigma)\circ(f_1,...,g_i,...,f_n))\circ (x_1,...,x_n)
\\
=&((a_0\otimes ...\otimes a_n,\sigma)\circ(f_1,...,f_i,...,f_n)
\\
+&(a_0\otimes ...\otimes a_n,\sigma)\circ(f_1,...,g_i,...,f_n))\circ (x_1,...,x_n)
\\[4pt]
&((a_0\otimes...\otimes a_n,\sigma)\circ(f_1,...,pf_i,...,f_n))\circ (x_1,...,x_n)
\\
=&(a_0\sigma(f_1)a_1,...\sigma(pf_i)...a_{n-1}\sigma(f_n) a_n)\circ (x_1,...,x_n)
\\
=&a_0\sigma(f_1\circ x_1)a_1,...\sigma((pf_i)\circ x_i)...a_{n-1}\sigma(f_n\circ x_n) a_n
\\
=&a_0\sigma(f_1\circ x_1)a_1,...\sigma(p(f_i\circ x_i))...a_{n-1}\sigma(f_n\circ x_n) a_n
\\
=&p(a_0\sigma(f_1\circ x_1)a_1,...\sigma(f_i\circ x_i)...a_{n-1}\sigma(f_n\circ x_n) a_n)
\\
=&p(((a_0\otimes ...\otimes a_n,\sigma)\circ(f_1,..., f_i, ..., f_n))\circ (x_1,...,x_n))
\\
=&(p((a_0\otimes ...\otimes a_n,\sigma)\circ(f_1,..., f_i, ..., f_n)))\circ (x_1,...,x_n)
\end{align*}
}
{representation An in LAAnA, 3}

\DefEquation
{
(c\otimes d)\circ(a\otimes b)=(ca)\otimes(bd)
}
{product in algebra AA}

\DefEquation
{
(b_0\otimes b_1\otimes...\otimes b_n)\circ(a_0\otimes a_1\otimes...\otimes a_n)
=(b_0a_0)\otimes(a_1b_1)\otimes...\otimes(a_nb_n)
}
{product in algebra An}

\DefEq
{
\[
(pf)\circ x=p\ f\circ x
\]
}
{linear map times scalar, 0, algebra}

\DefEq
{
\[
(f+g)\circ a=f\circ a+g\circ a
\]
}
{sum of transformations of Abelian group}

\DefEq
{
\[
f(ab)=f(a)\circ f(b)
\]
}
{product of transformations of representation}

\DefEquation
{
f(a)\circ m=f(b)\circ m
}
{representation of ring, 1}

\DefEq
{
\[
g:A_2\times A_2\rightarrow\mathcal L(A_1;A_2)
\]
\[
g(a,b)\circ f=afb
\]
}
{bilinear mapping A LA}

\DefEq
{
\[
h:\ATwo\rightarrow \mathcal L(A_1;A_2)
\]
}
{linear map AA LAA}

\DefEq
{
\[
h:\AOn\times S_n\rightarrow \LAAnA
\]
}
{linear map An LAAnA}

\DefEq
{
\[
h:\mathcal L(A_1;A)\otimes...\otimes\mathcal L(A_n;A)\rightarrow \LAAnA
\]
}
{linear mapping LA LA}

\DefEq
{
\[
h:A_1\rightarrow A_2
\]
}
{linear map h A1 A2}

\DefEq
{
\[
f:A_1\rightarrow A_2
\]
}
{linear map f A1 A2}

\DefEquation
{
h=
(
a\pC{s}{0}
\otimes
a\pC{s}{1}
)
\circ f
=a\pC{s}{0}
f
a\pC{s}{1}
}
{h generated by f, associative algebra}

\DefEq
{
\[
\begin{matrix}
f\pC{tr}{q+1}^k=h\pC{t}{q}^{k-1}
&q=1, ..., k-1
\\
f\pC{tr}{q}^k=g\pC{tr}{q}
&q=0, 1
\end{matrix}
\]
}
{polylinear map, induction on n, 5, associative algebra}

\DefEquation
{
f_{k_1...k_n\cdot s\cdot p}^n=\Vector e_{\gi i}f_{k_1...k_n\cdot s\cdot p}^{n\gi i}
}
{polylinear map, associative algebra, components extention}

\DefEquation
{
\Vector e'_{\gi i}=\Vector e_{\gi j}h_{\gi i}^{\gi j}
}
{polylinear map over ring, associative algebra, change basis}

\DefEq
{
\[
f\circ (a_1,...,a_n)=
f\circ (\Vector e_{\gi{i_1}}a^{\gi{i_1}}_1,...,\Vector e_{\gi{i_n}}a^{\gi{i_n}}_n)
=a^{\gi{i_1}}_1...a^{\gi{i_n}}_nf\circ(\Vector e_{\gi{i_1}},...,\Vector e_{\gi{i_n}})
\]
}
{polylinear map over ring, 1, associative algebra}

\DefEquation
{
\begin{array}{r@{\ }l}
a_j&=\Vector e_{\gi i}a^{\gi i}_j
\\
f_{\gi{i_1...i_n}}&=f\circ(\Vector e_{\gi{i_1}},...,\Vector e_{\gi{i_n}})
\end{array}
}
{polylinear map over ring, coordinates, associative algebra}

\DefEquation
{
\begin{array}{r@{\ }l}
f'_{\gi{i_1...i_n}}&=f\circ(\Vector e'_{\gi{i_1}},...,\Vector e'_{\gi{i_n}})
\\
&=f\circ(\Vector e_{\gi{j_1}}h_{\gi{i_1}}^{\gi{j_1}} ,...,
\Vector e'_{\gi{j_n}}h_{\gi{i_n}}^{\gi{j_n}})
\\
&=h_{\gi{i_1}}^{\gi{j_1}} ...
h_{\gi{i_n}}^{\gi{j_n}} f\circ(\Vector e_{\gi{j_1}},...,\Vector e_{\gi{j_n}})
\\
&=h_{\gi{i_1}}^{\gi{j_1}} ... h_{\gi{i_n}}^{\gi{j_n}} f_{\gi{j_1...j_n}}
\end{array}
}
{polylinear map over ring, associative algebra, change coordinates}

\DefEquation
{
f_{\gi{i_1,...,i_n}}=f_{\sigma(\gi{i_1}),...,\sigma(\gi{i_n})}
}
{polylinear map symmetric, associative algebra}

\DefEq
{
\begin{align*}
a^{\gi{i_1}}_1 ... a^{\gi{i_n}}_n f_{\gi{i_1...i_n}}
=&f\circ(a_1,...,a_n)
\\
=&f\circ(\sigma(a_1),...,\sigma(a_n))
\\
=&a_1^{\gi{i_1}} ... a_n^{\gi{i_n}} f_{\sigma(\gi{i_1})...\sigma(\gi{i_n})}
\end{align*}
}
{polylinear map symmetric, 1, associative algebra}

\DefEq
{
\begin{align}
f_{\gi{l_1...l_n}}
=&f_{t\cdot k_1...k_n}^{\gi{i_0...i_n}}
I_{k_1\cdot}{}^{\gi{j_1}}_{\gi{l_1}}...I_{k_n\cdot}{}^{\gi{j_n}}_{\gi{l_n}}
C_{\gi{i_0}\sigma_t(\gi{j_1})}^{\gi{k_1}}
C_{\gi{k_1i_1}}^{\gi{l_1}}
 ... B{}_{\gi{l_{n-1}}\sigma_t(\gi{j_n})}^{\gi{k_n}}
C_{\gi{k_ni_n}}^{\gi{l_n}}\Vector e_{\gi{l_n}}
\EqLabel{coordinates of polylinear map, associative algebra over ring}
\\
f_{\gi{l_1...l_n}}^{\gi p}
=&f_{t\cdot k_1...k_n}^{\gi{i_0...i_n}}
I_{k_1\cdot}{}^{\gi{j_1}}_{\gi{l_1}}...I_{k_n\cdot}{}^{\gi{j_n}}_{\gi{l_n}}
C_{\gi{i_0}\sigma_t(\gi{j_1})}^{\gi{k_1}}
C_{\gi{k_1i_1}}^{\gi{l_1}}
... C_{\gi{l_{n-1}}\sigma_t(\gi{j_n})}^{\gi{k_n}}
C_{\gi{k_ni_n}}^{\gi p}
\EqLabel{coordinates of polylinear map, associative algebra over ring, 1}
\end{align}
}
{coordinates of polylinear map, associative algebra over ring}

\DefEquation
{
\begin{array}{r@{\ }l}
f\circ(a_1,...,a_n)
=&f_{t\cdot k_1...k_n}^{\gi{i_0...i_n}}
\ \Vector e_{\gi{i_0}}
\sigma_t(a_1^{\gi{l_1}}I_{k_1\cdot}{}^{\gi{j_1}}_{\gi{l_1}}\Vector e_{\gi{j_1}})
\Vector e_{\gi{i_1}} ... 
\sigma_t(a_n^{\gi{l_n}}I_{k_n\cdot}{}^{\gi{j_n}}_{\gi{l_n}}\Vector e_{\gi{j_n}}) \Vector e_{\gi{i_n}}
\\[4pt]
=&a_1^{\gi{l_1}} ...a_n^{\gi{l_n}}
f_{t\cdot k_1...k_n}^{\gi{i_0...i_n}}
I_{k_1\cdot}{}^{\gi{j_1}}_{\gi{l_1}}...I_{k_n\cdot}{}^{\gi{j_n}}_{\gi{l_n}}
\Vector e_{\gi{i_0}} \sigma_t(\Vector e_{\gi{j_1}})
\Vector e_{\gi{i_1}} ... \sigma_t(\Vector e_{\gi{j_n}}) \Vector e_{\gi{i_n}}
\\[4pt]
=&a_1^{\gi{l_1}} ...a_n^{\gi{l_n}}
f_{t\cdot k_1...k_n}^{\gi{i_0...i_n}}
I_{k_1\cdot}{}^{\gi{j_1}}_{\gi{l_1}}...I_{k_n\cdot}{}^{\gi{j_n}}_{\gi{l_n}}
C_{\gi{i_0}\sigma_t(\gi{j_1})}^{\gi{k_1}}
C_{\gi{k_1i_1}}^{\gi{l_1}}
\\[4pt]
&...C_{\gi{l_{n-1}}\sigma_t(\gi{j_n})}^{\gi{k_n}}
C_{\gi{k_ni_n}}^{\gi{l_n}} \Vector e_{\gi{l_n}}
\end{array}
}
{polylinear map, associative algebra, standard representation, 1a}

\DefEquation
{
f\circ(a_1,...,a_n)
=\Vector e_{\gi p}
f_{\gi{i_1...i_n}}^{\gi p}a^{\gi{i_1}}_1...a^{\gi{i_n}}_n
}
{polylinear map over field, associative algebra, 1a}

\DefEq
{
\[
\begin{matrix}
f_i\in\mathcal L(A_i;A)
&i=1,...,n
\end{matrix}
\]
}
{set of mappings, L(A;A)}

\DefEq
{
\[
h:A^{n+1}\rightarrow \LAAnA
\]
}
{n+1 linear mapping A LA}

\DefEq
{
\[
h:\prod_{i=1}^n\mathcal L(A_i;A)\rightarrow \LAAnA
\]
}
{n linear mapping LA LA}

\DefEq
{
$a\in\AOn$
}
{tensor over A}

\DefEq
{
\[
(a_0,...,a_n,\sigma)\circ(f_1,...,f_n)
=a_0\sigma(f_1)a_1...a_{n-1}\sigma(f_n)a_n
\]
}
{n+1 linear mapping A LA 1}

\DefEq
{
\[
(a_0\otimes...\otimes a_n,\sigma)\circ(f_1,...,f_n)
=a_0\sigma(f_1)a_1...a_{n-1}\sigma(f_n)a_n
\]
}
{n linear mapping LA LA 1}

\DefEquation
{
\begin{array}{r@{\ }l}
&((a_0,...,a_n)\circ\sigma(f_1,...,f_n))\circ(x_1,...,x_n)
\\=&
(a_0\sigma(f_1)a_1...a_{n-1}\sigma(f_n)a_n)\circ(x_1,...,x_n)
\\=&
a_0\sigma(f_1\circ x_1)a_1...a_{n-1}\sigma(f_n\circ x_n)a_n
\end{array}
}
{n linear mapping A LA}

\DefEq
{
\[
\begin{matrix}
a_j\in A&j=0,...,n
\end{matrix}
\]
}
{set of elements A}

\DefEq
{
\[
I_{k_i}\circ a_i=\Vector e_{\gi{j_i}}I_{k_i\cdot}{}^{\gi{j_i}}_{\gi{l_i}}a_i^{\gi{l_i}}
\]
}
{coordinates of polylinear map, 1, associative algebra over ring}

\DefEq
{
\begin{align*}
a^{\gi{i_1}}_1 ... a^{\gi{i_n}}_n f_{\gi{i_1...i_n}}
=&f\circ(a_1,...,a_n)
\\
=&|\sigma|f\circ(\sigma(a_1),...,\sigma(a_n))
\\
=&a_1^{\gi{i_1}} ... a_n^{\gi{i_n}} |\sigma| f_{\sigma(\gi{i_1})...\sigma(\gi{i_n})}
\end{align*}
}
{polylinear map skew symmetric, 1, associative algebra}

\DefEquation
{
f_{\gi{i_1,...,i_n}}=|\sigma|f_{\sigma(\gi{i_1}),...,\sigma(\gi{i_n})}
}
{polylinear map skew symmetric, associative algebra}

\DefEq
{
\[
f\circ(a_1,...,a_n)=|\sigma|f\circ(\sigma(a_1),...,\sigma(a_n))
\]
}
{polylinear map skew symmetric, associative algebra, definition}

\DefEq
{
\[
f\circ(a_1,...,a_n)=f\circ(\sigma(a_1),...,\sigma(a_n))
\]
}
{polylinear map symmetric, associative algebra, definition}

\DefEquation
{
f(a_1,...,a_n)=a^{\gi{i_1}}_1...a^{\gi{i_n}}_nf_{\gi{i_1...i_n}}
}
{polylinear map over ring, associative algebra}

\DefEquation
{
f_{t\cdot k_1...k_n}^{\gi{j_0...j_n}}=
f_{k_1...k_n\cdot s\cdot 0}^{n\gi{j_1}}\ ...f_{k_1...k_n\cdot s\cdot n}^{n\gi{j_n}}
}
{polylinear map, associative algebra, standard representation, 2}

\DefEquation
{
a=f_{k_1...k_n\cdot s\cdot 0}^{n\gi{j_1}}\ \Vector e_{\gi{j_1}}\ \sigma_s(I_{k_1}\circ a_1)
\ f_{k_1...k_n\cdot s\cdot 1}^{n\gi{j_2}}\ \Vector e_{\gi{j_2}}\ ...
\ \sigma_s(I_{k_n}\circ a_n)\ f_{k_1...k_n\cdot s\cdot n}^{n\gi{j_n}}\ \Vector e_{\gi{j_n}}
}
{polylinear map, associative algebra, standard representation, 1}

\DefEq
{
\symb{f_{t\cdot k_1...k_n}^{\gi{i_0...i_n}}}0
{standard component of polylinear map, associative algebra}
\begin{equation}
f\circ (a_1,...,a_n)=
\ShowSymbol{standard component of polylinear map, associative algebra}
\ \Vector e_{\gi{i_0}}\ \sigma_t(I_{k_1}\circ a_1)\ 
\Vector e_{\gi{i_1}}
\ ...\ \sigma_t(I_{k_n}\circ a_n)\ \Vector e_{\gi{i_n}}
\EqLabel{polylinear map, associative algebra, standard representation}
\end{equation}
}
{polylinear map, associative algebra, standard representation}

\DefEq
{
\symb{f\pC{s}{p}^n}1
{component of polylinear map, associative algebra}
}
{component of polylinear map, associative algebra}

\DefEq
{
\[
\sigma_s=\sigma_{tr}=
\left(
\begin{array}{cccc}
a_k&a_1&...&a_{k-1}
\\
a_k&\sigma_t(a_1)&...&\sigma_t(a_{k-1})
\end{array}
\right)
\]
}
{polylinear map, induction on n, 4, associative algebra}

\DefEq
{
\[
h\pC{t}{0}^{k-1}
=g\pC{tr}{0}
\ I_{k\cdot r}\circ a_k\ g\pC{tr}{1}
\]
}
{polylinear map, induction on n, 3, associative algebra}

\DefEq
{
\[
a
=a\pC{s}{0}
\otimes
a\pC{s}{1}
\]
}
{h generated by f 1, associative algebra}

\DefEquation
{
l(a)\circ x=ax
}
{left shift algebra}

\DefEquation
{
l(a)\circ l(b)=l(ab)-(a,b)_1
}
{left shift algebra, 1}

\DefEq
{
\[
(a,b)_1\circ x=(a,b,x)
\]
}
{left shift algebra, 2}

\DefEquation
{
\begin{array}{r@{\ }l}
(l(a)\circ l(b))\circ x
&=l(a)\circ(l(b)\circ x)
\\
&=a(bx)=(ab)x-(a,b,x)
\\
&=l(ab)\circ x-(a,b)_1\circ x
\end{array}
}
{left shift algebra, 3}

\DefEq
{
\begin{align*}
l(a)\circ(l(b)\circ l(c))\circ x
&=l(a)\circ(l(bc)\circ x-(b,c,x))
\\
&=l(a)\circ l(bc)\circ x-l(a)\circ(b,c,x)
\\
&=l(a(bc))\circ x-(a,bc,x)-l(a)\circ(b,c,x)
\end{align*}
}
{left shift algebra, 2a}

\DefEquation
{
r(a)\circ x=xa
}
{right shift algebra}

\DefEquation
{
\begin{array}{r@{\ }l}
(r(a)\circ r(b))\circ x
&=r(a)\circ(r(b)\circ x)
\\
&=(xb)a=x(ba)+(x,b,a)
\\
&=r(ba)\circ x+(x,b,a)
\end{array}
}
{right shift algebra, 3}

\DefEquation
{
r(a)\circ r(b)=r(ba)+(b,a)_2
}
{right shift algebra, 1}

\DefEq
{
\[
(b,a)_2\circ x=(x,b,a)
\]
}
{right shift algebra, 2}

\def\AcF{a\circ f}
\def\GAcF{g=\AcF}

\DefEquation
{
\GAcF
}
{standard representation of mapping A1 A2, 1, associative algebra}

\DefEquation
{
\GAcF
}
{standard representation of mapping A1 A2, 1, nonassociative algebra}

\DefEq
{
\[
\GAcF
\]
}
{coordinates of mapping A1 A2, 1}

\DefEq
{
\begin{align}
\EqLabel{coordinates of mapping f, associative algebra}
f\circ x&=f^{\gi i}_{\gi j}x^{\gi j}\Vector e_{2\cdot\gi i}
\\
\EqLabel{coordinates of mapping g, associative algebra}
g\circ x&=g^{\gi i}_{\gi j}x^{\gi j}\Vector e_{2\cdot\gi i}
\end{align}
}
{coordinates of mappings f g, associative algebra}

\DefEq
{
\begin{align}
\EqLabel{coordinates of mapping f, nonassociative algebra}
f\circ x&=f^{\gi i}_{\gi j}x^{\gi j}\Vector e_{2\cdot\gi i}
\\
\EqLabel{coordinates of mapping g, nonassociative algebra}
g\circ x&=g^{\gi i}_{\gi j}x^{\gi j}\Vector e_{2\cdot\gi i}
\end{align}
}
{coordinates of mappings f g, nonassociative algebra}

\DefEq
{
\begin{equation}
\mathcal B=
\left(
\mathcal C^{\cdot}{}_{\gi m}^{\gi k}{}_{\cdot\gi{ij}}
\right)
=
\left(
C_{2\cdot}{}^{\gi p}_{\gi{im}}C_{2\cdot}{}^{\gi k}_{\gi{pj}}
\right)
\EqLabel{linear map over ring, matrix}
\end{equation}
}
{linear map over ring, matrix}

\DefEq
{
${}^{\cdot}{}_{\gi m}^{\gi k}$
}
{linear map over ring, row of matrix}

\DefEq
{
${}_{\cdot\gi{ij}}$
}
{linear map over ring, column of matrix}

\DefEquation
{
\rank
\begin{pmatrix}
\mathcal C^{\cdot}{}_{\gi m}^{\gi k}{}_{\cdot\gi{ij}}
&g_{\gi m}^{\gi k}
\end{pmatrix}
=\rank\mathcal C
}
{linear map over ring, determinant=0, 1}

\DefEquation
{
\rank
\begin{pmatrix}
\mathcal C^{\cdot}{}_{\gi m}^{\gi k}{}_{\cdot\gi{ij}}
&g_{\gi m}^{\gi k}
&f_{\gi m}^{\gi k}
\end{pmatrix}
=\rank\mathcal C
}
{linear map over ring, determinant=0, 2}

\DefEquation
{
\rank
\begin{pmatrix}
\mathcal C^{\cdot}{}_{\gi m}^{\gi k}{}_{\cdot\gi{ij}}
&f_{\gi m}^{\gi k}
\end{pmatrix}
=\rank\mathcal C
}
{linear map over ring, determinant=0, 3}

\DefEq
{
$f=\delta$
}
{linear map over ring, f=1}

\DefEq
{
$g_{\gi k}^{\gi l}$
}
{linear map over ring, left side}

\DefEquation
{
\begin{array}{r@{}l}
g^{\gi k}_{\gi l}x^{\gi l}\Vector e_{2\cdot\gi k}
&=
a^{\gi{ij}}(\Vector e_{2\cdot\gi i}
(f^{\gi m}_{\gi l}x^{\gi l}\Vector e_{2\cdot\gi m}))
\Vector e_{2\cdot\gi j}
\\
&\VirtVar=
a^{\gi{ij}}
f^{\gi m}_{\gi l}x^{\gi l}
C_{2\cdot}{}^{\gi p}_{\gi{im}}
C_{2\cdot}{}^{\gi k}_{\gi{pj}}
\Vector e_{2\cdot\gi k}
\end{array}
}
{coordinates of mapping A1 A2, 3, nonassociative algebra}

\DefEquation
{
\begin{array}{r@{}l}
g^{\gi k}_{\gi l}x^{\gi l}\Vector e_{2\cdot\gi k}
&=
a^{\gi{ij}}\Vector e_{2\cdot\gi i}
f^{\gi m}_{\gi l}x^{\gi l}\Vector e_{2\cdot\gi m}
\Vector e_{2\cdot\gi j}
\\
&\VirtVar=
a^{\gi{ij}}
f^{\gi m}_{\gi l}x^{\gi l}
C_{2\cdot}{}^{\gi p}_{\gi{im}}
C_{2\cdot}{}^{\gi k}_{\gi{pj}}
\Vector e_{2\cdot\gi k}
\end{array}
}
{coordinates of mapping A1 A2, 3, associative algebra}

\DefEquation
{
g^{\gi k}_{\gi l}=
f^{\gi m}_{\gi l}g^{\gi{ij}}
C_{2\cdot}{}^{\gi p}_{\gi{im}}C_{2\cdot}{}^{\gi k}_{\gi{pj}}
}
{coordinates of mapping A1 A2, 2, associative algebra}

\DefEquation
{
g^{\gi k}_{\gi l}
=f^{\gi m}_{\gi l}g^{\gi{ij}}
C_{2\cdot}{}^{\gi p}_{\gi{im}}C_{2\cdot}{}^{\gi k}_{\gi{pj}}
}
{coordinates of mapping A1 A2, 2, nonassociative algebra}

\DefEquation
{
g^{\gi k}_{\gi l}
=f^{\gi m}_{\gi l}g^{\gi{ij}}
C_{2\cdot}{}^{\gi k}_{\gi{ip}}C_{2\cdot}{}^{\gi p}_{\gi{mj}}
}
{coordinates of mapping A1 A2, 21, nonassociative algebra}

\DefEquation
{
g^{\gi k}_{\gi l}=
f^{\gi m}_{\gi l}g^{\gi{ij}}
C^{\gi p}_{\gi{im}}C^{\gi k}_{\gi{pj}}
}
{coordinates of mapping A}

\DefEquation
{
f^{\gi k}_{\gi l}=
f^{\gi{ij}}
C^{\gi p}_{\gi{il}}C^{\gi k}_{\gi{pj}}
}
{coordinates of mapping A1 A2, 4f}

\DefEquation
{
f=f^{\gi{ij}}\Vector e_{\gi i}\otimes\Vector e_{\gi j}
}
{coordinates of mapping A1 A2, 5f}

\DefEquation
{
g^{\gi k}_{\gi l}=
g^{\gi{ij}}
C^{\gi p}_{\gi{il}}C^{\gi k}_{\gi{pj}}
}
{coordinates of mapping A1 A2, 4g}

\DefEquation
{
g=g^{\gi{ij}}\Vector e_{\gi i}\otimes\Vector e_{\gi j}
}
{coordinates of mapping A1 A2, 5g}

\DefEquation
{
g=
a^{\gi{ij}}(\Vector e_{\gi i}\otimes\Vector e_{\gi j})\circ f
=
a^{\gi{ij}}\Vector e_{\gi i}f\Vector e_{\gi j}
}
{standard representation of mapping A1 A2, 2, associative algebra}

\DefEquation
{
g
=a^{\gi{ij}}(\Vector e_{\gi i}\otimes\Vector e_{\gi j})\circ f
=
a^{\gi{ij}}(\Vector e_{\gi i}f)\Vector e_{\gi j}
}
{standard representation of mapping A1 A2, 2, nonassociative algebra}

\DefEquation
{
a=a^{\gi{ij}}\Vector e_{\gi i}\otimes\Vector e_{\gi j}
}
{standard representation of mapping A1 A2, 3, associative algebra}

\DefEquation
{
h:\ATwo\rightarrow {}^*\mathcal L(A_1;A_2)
}
{representation A2 in LA}

\DefEquation
{
h:A\otimes A\rightarrow {}^*\mathcal L(A;A)
}
{representation AA in LA}

\DefEq
{
\[
h:\AOn\times S_n\rightarrow {}^*\LAnA
\]
}
{representation An in LAnA}

\DefEq
{
$\GAcF$
}
{orbit of linear mapping, 01}

\DefEq
{
$a\in A_2\otimes A_2$
}
{orbit of linear mapping, 02}

\DefEq
{
$h\in(\ATwo)\circ g$,
}
{orbit of linear mapping, 1}

\DefEq
{
$h\in(\ATwo)\circ f$,
}
{orbit of linear mapping, 5}

\DefEquation
{
(\ATwo)\circ g\subset(\ATwo)\circ f
}
{orbit of linear mapping, g in f}

\DefEquation
{
(\ATwo)\circ f\subset(\ATwo)\circ g
}
{orbit of linear mapping, f in g}

\DefEquation
{
(\ATwo)\circ f=(\ATwo)\circ g
}
{orbit of linear mapping, f = g}

\DefEquation
{
f=a^{-1}\circ g
}
{orbit of linear mapping, 6}

\DefEq
{
$b\in\ATwo$
}
{orbit of linear mapping, 2}

\DefEq
{
$h=b\circ g$.
}
{orbit of linear mapping, 3}

\DefEquation
{
h=b\circ f
}
{orbit of linear mapping, 7}

\DefEquation
{
h=b\circ(\AcF)=(b\circ a)\circ f
}
{orbit of linear mapping, 4}

\DefEq
{
\[
h=b\circ(a^{-1}\circ g)=(b\circ a^{-1})\circ g
\]
}
{orbit of linear mapping, 8}

\DefEq
{
\symb{(\ATwo)\circ f}
0{orbit of linear mapping}
\[
\ShowSymbol{orbit of linear mapping}
=\{g=d\circ f:d\in \ATwo\}
\]
}
{orbit of linear mapping}

\DefEquation
{
\begin{matrix}
(a_0\otimes...\otimes a_n,\sigma)\circ(f_1\otimes...\otimes f_n)
=a_0\sigma(f_1)a_1...a_{n-1}\sigma(f_n)a_n
\\
\begin{matrix}
a_0,...,a_n\in A&\sigma\in S_n&f_1,...,f_n\in\mathcal L(A_n;A)
\end{matrix}
\end{matrix}
}
{representation An in LAnA, 1}

\DefEquation
{
b_1\sigma(I_{1\cdot s}\circ x_1)c_1b_2...c_{n-1}b_n\sigma(I_{n\cdot s}\circ x_n)c_n
}
{n linear mapping A LA, 1}

\DefEq
{
\[
\begin{matrix}
b_1=f\pC{s}{0}^n&b_2=...=b_n=e&c_1=f\pC{s}{1}^n&...&c_n=f\pC{s}{n}^n
\end{matrix}
\]
}
{n linear mapping A LA, 2}

\DefEq
{
\[
\begin{matrix}
f_i=\sigma^{-1}(b_i)I_{i\cdot s}\sigma^{-1}(c_i)&i=1,...,n
\end{matrix}
\]
}
{n linear mapping A LA, 3}

\DefEquation
{
\begin{matrix}
(a\otimes b)\circ f=afb
&a,b\in A_2&f\in\mathcal L(A_1;A_2)
\end{matrix}
}
{representation A2 in LA, 1}

\DefEquation
{
\begin{matrix}
(a\otimes b)\circ f=afb
&a,b\in A&f\in\mathcal L(A;A)
\end{matrix}
}
{representation AA in LA, 1}

\DefEquation
{
\xymatrix{
&B\ar[rd]^g
\\
A\ar[ru]^f\ar[rr]^{g\circ f}&&C
}
}
{product of linear mapping, algebra}

\DefEq
{
\begin{align}
(f+g)\circ a&=f\circ a+g\circ a
\\
(pf)\circ a&=p\,f\circ a
\end{align}
}
{operations over linear map}

\DefEquation
{
a^{\gi j}=A_{\gi i}^{\gi j}a'^{\gi i}
}
{free module over ring, change basis, 3}

\DefEquation
{
\Vector e'_{\gi i}=\Vector e_{\gi j}A_{\gi i}^{\gi j}
}
{free module over ring, change basis}

\DefEquation
{
a=\Vector e'_{\gi i}a'^{\gi i}=\Vector e_{\gi j}a^{\gi j}
}
{free module over ring, change basis, 1}

\DefEquation
{
\Vector e_{\gi j}a^{\gi j}=\Vector e'_{\gi i}a'^{\gi i}
=\Vector e_{\gi j}A_{\gi i}^{\gi j}a'^{\gi i}
}
{free module over ring, change basis, 2}

\DefEq
{
\[
\xymatrix{
&A\ar[rd]^g
\\
A\ar[ru]^f\ar[rr]^{g\circ f}&&A
}
\]
}
{product of linear mapping, algebra 1}

\DefEquation
{
f^*:g\in\mathcal L(B;C)\rightarrow g\circ f\in\mathcal L(A;C)
}
{product of linear mapping, algebra, f}

\DefEquation
{
\xymatrix{
&B\ar[rd]^g_{}="G"
\\
A\ar[ru]^f\ar[rr]_{g\circ f}^{}="GF"&&C
\ar @{=>} "G";"GF"
}
}
{product of linear mapping, algebra, f diagram}

\DefEquation
{
\xymatrix{
&B\ar[rd]^g
\\
A\ar[ru]^f_{}="F"\ar[rr]_{g\circ f}^{}="GF"&&C
\ar @{=>} "F";"GF"
}
}
{product of linear mapping, algebra, g diagram}

\DefEquation
{
g_*:f\in\mathcal L(A;B)\rightarrow g\circ f\in\mathcal L(A;C)
}
{product of linear mapping, algebra, g}

\DefEquation
{
\circ:(g,f)\in\mathcal L(B;C)\times\mathcal L(A;B)
\rightarrow g\circ f\in\mathcal L(A;C)
}
{product of linear mapping, algebra, gf}

\DefEquation
{
\circ:(g,f)\in\mathcal L(A;A)\times\mathcal L(A;A)
\rightarrow g\circ f\in\mathcal L(A;A)
}
{module L(A;A) is algebra}

\DefEq
{
\begin{align*}
(g\circ f)\circ(a+b)&=g\circ(f\circ(a+b))
=g\circ (f\circ a+f\circ b)
\\
&=g\circ (f\circ a)+g\circ (f\circ b)
=(g\circ f)\circ a+(g\circ f)\circ b
\\
(g\circ f)\circ(pa)&=g\circ(f\circ(pa))
=g\circ (p\ f\circ a)
=p\ g\circ (f\circ a)
\\
&=p\ (g\circ f)\circ a
\end{align*}
}
{product of linear mapping, algebra, 1}

\DefEq
{
\begin{align*}
((g_1+g_2)\circ f)\circ a&=(g_1+g_2)\circ(f\circ a)
=g_1\circ(f\circ a)+g_2\circ(f\circ a)
\\
&=(g_1\circ f)\circ a+(g_2\circ f)\circ a
\\
&=(g_1\circ f+g_2\circ f)\circ a
\\
((pg)\circ f)\circ a&=(pg)\circ(f\circ a)
=p\ g\circ(f\circ a)
=p\ (g\circ f)\circ a
\\
&=(p(g\circ f))\circ a
\end{align*}
}
{product of linear mapping, algebra, f1}

\DefEq
{
\begin{align*}
(g\circ (f_1+f_2))\circ a&=g\circ((f_1+f_2)\circ a)
=g\circ(f_1\circ a+f_2\circ a)
\\
&=g\circ (f_1\circ a)+g\circ (f_2\circ a)
=(g\circ f_1)\circ a+(g\circ f_2)\circ a
\\
&=(g\circ f_1+g\circ f_2)\circ a
\\
(g\circ (pf))\circ a&=g\circ((pf)\circ a)
=g\circ(p\ (f\circ a))
=p\ g\circ(f\circ a)
\\
&=p\ (g\circ f)\circ a
=(p(g\circ f))\circ a
\end{align*}
}
{product of linear mapping, algebra, g1}

\DefEquation
{
(a\otimes b)\circ f=afb
}
{linear map AA LAA, 1}

\DefEquation
{
\begin{array}{r@{\ }l}
(a_0\otimes...\otimes a_n,{\sigma})\circ (f_1,...,f_n)
&=(a_0,..., a_n,\sigma)\circ (f_1,...,f_n)
\\
&=a_0\sigma (f_1)a_1...a_{n-1}\sigma (f_n) a_n
\end{array}
}
{linear map An LAAnA, 1}

\DefEquation
{
(a_0\otimes...\otimes a_n,{\sigma})\circ (f_1\otimes...\otimes f_n)=
(a_0\otimes...\otimes a_n,{\sigma})\circ (f_1,...,f_n)
}
{linear mapping LA LA, 1}

\DefEquation
{
x\rightarrow (d\circ f)\circ x
}
{tensor and mapping in A, algebra}

\DefEquation
{
((a\otimes b)\circ\delta)\circ x
=(a\delta b)\circ x=a\ (\delta\circ x)\ b=axb
}
{tensor and mapping in A, 1, algebra}

\DefEquation
{
\begin{array}{r@{\ }l}
((a_0\otimes...\otimes a_n,\sigma)\circ(\delta,...,\delta))\circ (x_1,...,x_n)
&=a_0\ (\delta\circ x_1)\ ...\ (\delta\circ x_n)\ a_n
\\
&=a_0\ x_1\ ...\ x_n\ a_n
\end{array}
}
{n tensor and mapping in A, 1, algebra}

\DefEquation
{
((a\otimes b)\circ\delta)\circ x
=(a\otimes b)\circ(\delta\circ x)
=(a\otimes b)\circ x
}
{tensor and mapping in A, 2, algebra}

\DefEquation
{
\begin{array}{r@{\ }l}
&((a_0\otimes...\otimes a_n,\sigma)\circ(\delta,...,\delta))\circ (x_1,...,x_n)
\\
=&(a_0\otimes...\otimes a_n,\sigma)\circ(\delta\circ x_1,...,\delta\circ x_n)
\\
=&(a_0\otimes...\otimes a_n,\sigma)\circ (x_1,...,x_n)
\end{array}
}
{n tensor and mapping in A, 2, algebra}

\DefEq
{
$d=a_0\otimes...\otimes a_n$%
}
{n tensor and mapping in A, 3, algebra}%

\DefEq
{
\begin{align*}
((a_1+a_2)fb)\circ x
&=(a_1+a_2)\ f\circ x\ b
=a_1\ f\circ x\ b+a_2\ f\circ x\ b
\\
&=(a_1fb)\circ x+(a_2fb)\circ x
=(a_1fb+a_2fb)\circ x
\\
((pa)fb)\circ x
&=(pa)\ f\circ x\ b
=p(a\ f\circ x\ b)
=p((afb)\circ x)
=(p(afb))\circ x
\\
(af(b_1+b_2))\circ x
&=a\ f\circ x\ (b_1+b_2)
=a\ f\circ x\ b_1+a\ f\circ x\ b_2
\\
&=(afb_1)\circ x+(afb_2)\circ x
=(afb_1+afb_2)\circ x
\\
(af(pb))\circ x
&=a\ f\circ x\ (pb)
=p(a\ f\circ x\ b)
=p((afb)\circ x)
=(p(afb))\circ x
\end{align*}
}
{bilinear mapping A LA, 1}

\DefEq
{
\begin{align*}
&((a_0,...,a_i+b_i,...a_n,\sigma)\circ(f_1,...,f_n))\circ(x_1,...,x_n)
\\
=&a_0\sigma(f_1\circ x_1)a_1...(a_i+b_i)...a_{n-1}\sigma(f_n\circ x_n)a_n
\\
=&a_0\sigma(f_1\circ x_1)a_1...a_i...a_{n-1}\sigma(f_n\circ x_n)a_n
+a_0\sigma(f_1\circ x_1)a_1...b_i...a_{n-1}\sigma(f_n\circ x_n)a_n
\\
=&((a_0,...,a_i,...,a_n,\sigma)\circ(f_1,...,f_n))\circ(x_1,...,x_n)
\\
+&((a_0,...,b_i,...,a_n,\sigma)\circ(f_1,...,f_n))\circ(x_1,...,x_n)
\\
=&((a_0,...,a_i,...,a_n,\sigma)\circ(f_1,...,f_n)
+(a_0,...,b_i,...,a_n,\sigma)\circ(f_1,...,f_n))\circ (x_1,...,x_n)
\\[4pt]
&((a_0,...,pa_i,...a_n,\sigma)\circ(f_1,...,f_n))\circ(x_1,...,x_n)
\\
=&a_0\sigma(f_1\circ x_1)a_1...pa_i...a_{n-1}\sigma(f_n\circ x_n)a_n
\\
=&p(a_0\sigma(f_1\circ x_1)a_1...a_i...a_{n-1}\sigma(f_n\circ x_n)a_n)
\\
=&p(((a_0,...,a_i,...,a_n,\sigma)\circ(f_1,...,f_n))\circ(x_1,...,x_n))
\\
=&(p((a_0,...,a_i,...,a_n,\sigma)\circ(f_1,...,f_n)))\circ (x_1,...,x_n)
\end{align*}
}
{n+1 linear mapping A LA, 1}

\DefEq
{
\begin{align*}
&((a_0,...,a_n,\sigma)\circ(f_1,...,f_n))\circ(x_1,...,x_i+y_i,...,x_n)
\\
=&a_0\sigma(f_1\circ x_1)a_1...\sigma(f_i\circ(x_i+y_i))...a_{n-1}\sigma(f_n\circ x_n)a_n
\\
=&a_0\sigma(f_1\circ x_1)a_1...\sigma(f_i\circ x_i+f_i\circ y_i)...a_{n-1}\sigma(f_n\circ x_n)a_n
\\
=&a_0\sigma(f_1\circ x_1)a_1...\sigma(f_i\circ x_i)...a_{n-1}\sigma(f_n\circ x_n)a_n
\\
+&a_0\sigma(f_1\circ x_1)a_1...\sigma(f_i\circ y_i)...a_{n-1}\sigma(f_n\circ x_n)a_n
\\
=&((a_0,...,a_n,\sigma)\circ(f_1,...,f_n))\circ(x_1,...,x_i,...,x_n)
\\
+&((a_0,...,a_n,\sigma)\circ(f_1,...,f_n))\circ(x_1,...,y_i,...,x_n)
\\[4pt]
&((a_0,...,a_n,\sigma)\circ(f_1,...,f_n))\circ(x_1,...,px_i,...,x_n)
\\
=&a_0\sigma(f_1\circ x_1)a_1...\sigma(f_i\circ(px_i))...a_{n-1}\sigma(f_n\circ x_n)a_n
\\
=&a_0\sigma(f_1\circ x_1)a_1...\sigma(p(f_i\circ x_i))...a_{n-1}\sigma(f_n\circ x_n)a_n
\\
=&p(a_0\sigma(f_1\circ x_1)a_1...\sigma(f_i\circ x_i)...a_{n-1}\sigma(f_n\circ x_n)a_n)
\\
=&p(((a_0,...,a_n,\sigma)\circ(f_1,...,f_n))\circ(x_1,...,x_i,...,x_n))
\end{align*}
}
{polylinear mapping A 1n A, 1}

\DefEq
{
\[
f(a-b)\circ m=0
\]
}
{representation of ring, 2}

\DefEq
{
\[
f(a+b)\circ x=f(a)\circ x+f(b)\circ x
\]
}
{sum of transformations of Abelian group, 1}

\DefEquation
{
f:D\rightarrow {}^\star A
}
{representation of the ring}

\DefEq
{
\begin{align*}
p(qf)&=(pq)f
\\
(p+q)f&=pf+qf
\end{align*}
}
{linear map times scalar, pf, algebra}

\DefEq
{
\symb{\mathcal L(A_1;A_2)}
1{set linear mappings, algebra}
}
{set linear mappings, algebra}

\DefEq
{
$f$, $g$, $h\in\mathcal L(A_1,...,A_2;S)$.
}
{module of polylinear mappings, 1}

\DefEq
{
$a=(a_1,...,a_n)$, $a_1\in A_1$, ..., $a_n\in A_n$,
\begin{align*}
(f+g)\circ a=&f\circ a+g\circ a=g\circ a+f\circ a
\\
=&(g+f)\circ a
\\
((f+g)+h)\circ a=&
(f+g)\circ a+h\circ a=(f\circ a+g\circ a)+h\circ a
\\
=&f\circ a+(g\circ a+h\circ a)=f\circ a+(g+h)\circ a
\\
=&(f+(g+h))\circ a
\end{align*}
}
{module of polylinear mappings, 2}

\DefEq
{
\[
z\circ a=0
\]
}
{module of polylinear mappings, 3}

\DefEq
{
\[
(z+f)\circ a=z\circ a +f\circ a=0+f\circ a=f\circ a
\]
}
{module of polylinear mappings, 4}

\DefEq
{
\[
g\circ a=-f\circ a
\]
}
{module of polylinear mappings, 5}

\DefEq
{
\[
f+g=z
\]
}
{module of polylinear mappings, 6}

\DefEq
{
\[
(f+g)\circ a=f\circ a+g\circ a=f\circ a-f\circ a=0
\]
}
{module of polylinear mappings, 7}

\DefEq
{
\begin{align*}
f&:A_1\rightarrow A_2
\\
g&:A_1\rightarrow A_2
\end{align*}
}
{sum of linear maps, algebra}

\DefEq
{
\begin{align*}
f&:A_1\times...\times A_n\rightarrow S
\\
g&:A_1\times...\times A_n\rightarrow S
\end{align*}
}
{sum of polylinear maps, algebra}

\DefEq
{
\begin{align}
f\circ(a+b)&=f\circ a+f\circ b
\EqLabel{linear map from A1 to A2, 1 1}
\\
f\circ(pa)&=p(f\circ a)
\EqLabel{linear map from A1 to A2, 1 2}
\end{align}
\[
\begin{matrix}
a,b\in A_1
&
p\in D
\end{matrix}
\]
}
{linear map from A1 to A2, 1}

\DefEq
{
\begin{align*}
f(a+b)&=f(a)+f(b)
\\
f(pa)&=pf(a)
\end{align*}
\[
\begin{matrix}
a,b\in A_1
&
p\in D
\end{matrix}
\]
}
{linear map from A1 to A2, 1 old}

\DefEquation
{
\left\{
\begin{array}{r@{\ }l}
f\circ(a+b)&=f\circ a+f\circ b
\\
f\circ(pa)&=pf\circ a
\\
\multicolumn{2}{c}
{
\begin{matrix}
a,b\in A_1
&
p\in D
\end{matrix}
}
\end{array}
\right.
}
{linear map from A1 to A2, algebra}

\DefEq
{
\[
f:A\times A\rightarrow A
\]
}
{product in algebra, definition 1}

\DefEq
{
\[
\begin{matrix}
f:A\rightarrow A
&
f=(ax)b
\end{matrix}
\]
}
{f over A}

\DefEquation
{
f\circ x=f^{\gi{ij}}\ (\Vector e_{\gi i}x)\Vector e_{\gi j}
}
{linear map, standard representation, nonassociative algebra}

\DefEq
{
\[
f\circ x=f^{\gi{ij}}\Vector e_{\gi i}(x\Vector e_{\gi j})
\]
}
{linear map, standard representation, nonassociative algebra, 1}

\DefEq
{
\[
\begin{matrix}
g:A\rightarrow A
&
g=(cf)d
\end{matrix}
\]
}
{g over A}

\DefEquation
{
ab=f\circ(a,b)
}
{product in algebra, definition 2}

\DefEq
{
\symb{A^*}1
{opposite algebra}
}
{opposite algebra}

\DefEq
{
\[
ba=f\circ(a,b)
\]
}
{opposite algebra, product}

\DefEq
{
\symb{[a,b]}
0{commutator of algebra}
\[
\ShowSymbol{commutator of algebra}=ab-ba
\]
}
{commutator of algebra}

\DefEq
{
\symb{(a,b,c)}
0{associator of algebra}
\begin{equation}
\ShowSymbol{associator of algebra}=(ab)c-a(bc)
\EqLabel{associator of algebra}
\end{equation}
}
{associator of algebra}

\DefEq
{
\[
[a,b]=0
\]
}
{commutative D algebra}

\DefEq
{
\[
(a,b,c)=0
\]
}
{associative D algebra}

\DefEquation
{
a(b,c,d)+(a,b,c)d=(ab,c,d)-(a,bc,d)+(a,b,cd)
}
{associator of algebra, 1}

\DefEq
{
\begin{align*}
a(b,c,d)+(a,b,c)d
&=a((bc)d-b(cd))+((ab)c-a(bc))d
\\
&=a((bc)d)-a(b(cd))+((ab)c)d-(a(bc))d
\\
&=((ab)c)d-(ab)(cd)+(ab)(cd)
\\
&+a((bc)d)-a(b(cd))-(a(bc))d
\\
&=(ab,c,d)-(a(bc))d+a((bc)d)+(ab)(cd)-a(b(cd))
\\
&=(ab,c,d)-(a,(bc),d)+(a,b,cd)
\end{align*}
}
{associator of algebra, 2}

\DefEq
{
\symb{C^{\gi k}_{\gi{ij}}}1
{structural constants of algebra}
}
{structural constants of algebra}

\DefEquation
{
ab=a^{\gi i}b^{\gi j}\Eij
}
{product in algebra, 1}

\DefEquation
{
ab=a^{\gi i}b^{\gi j}
C^{\gi k}_{\gi{ij}}\Vector e_{\gi k}
}
{product in algebra, 2}

\DefEq
{
\[
\begin{matrix}
a=a^{\gi i}e_{\gi i}&b=b^{\gi i}e_{\gi i}&a, b\in A
\end{matrix}
\]
}
{a b in basis of algebra}

\DefEq
{
\begin{align*}
(a+b)c&=f\circ(a+b,c)=f\circ(a,c)+f\circ(b,c)=ac+bc
\\
a(b+c)&=f\circ(a,b+c)=f\circ(a,b)+f\circ(a,c)=ab+ac
\end{align*}
}
{product distributive in algebra}

%% file: Tensor.Algebra.English.tex
%auto-ignore
\input{\FilePrefix Tensor.Algebra.Eq}

%\chapter{Tensor Product}

\ePrints{1011.3102}
\ifx\Semafor\ValueOff
\Section{Tensor Product of Algebras}
\else
\section{Tensor Product of Algebras}
\fi
\label{Section: Tensor Product of Algebras}

\ePrints{2010.05.23}
\ifx\Semafor\ValueOff
\begin{definition}
\label{definition: tensor product of algebras}
Let $A_1$, ..., $A_n$ be
free algebras over commutative ring $D$.\footnote{I give definition
of tensor product of $D$\Hyph algebras
following to definition in \citeBib{Serge Lang}, p. 601 - 603.
\ePrints{Calculus.Paper}
\Items{1011.3102}
\ifx\Semafor\ValueOn
This subsection
is written on the base of the section
\xRef{8433-5163}{Section: Tensor Product of Algebras}.
\fi
}
Let us consider category $\mathcal A$ whose objects are
polylinear over commutative ring $D$ mappings
\ShowEq{polylinear mappings category}
where $S_1$, $S_2$ are modules over ring $D$,
We define morphism $f\rightarrow g$
to be linear over commutative ring $D$ mapping $h:S_1\rightarrow S_2$
making diagram
\ShowEq{polylinear mappings category, diagram}
commutative.
Universal object
\ShowEq{tensor product of algebras}
of category $\mathcal A$
is called \AddIndex{tensor product of algebras $A_1$, ..., $A_n$}
{tensor product of algebras}.
\qed
\end{definition}

\begin{definition}
\label{definition: tensor power of algebra}
Tensor product
\ShowEq{tensor power of algebra}
is called
\AddIndex{tensor power of algebra}{tensor power of algebra} $A$.
\qed
\end{definition}

\begin{theorem}
\label{theorem: tensor product of algebras is module}
There exists tensor product of algebras.
\end{theorem}
\begin{proof}
Let $M$ be module over ring $D$ generated
by product $\Times$
of $D$\Hyph algebras $A_1$, ..., $A_n$.
Injection
\ShowEq{map i, algebra, tensor product, 1}
is defined according to rule
\ShowEq{map i, algebra, tensor product}
Let $N\subset M$ be submodule generated by elements of the following type
\ShowEq{kernel, algebra, tensor product}
where $d_i\in A_i$, $c_i\in A_i$, $a\in D$.
Let
\[
j:M\rightarrow M/N
\]
be canonical mapping on factor module.
Consider commutative diagram
\ShowEq{diagram top, algebra, tensor product}
Since elements \EqRef{kernel 1, algebra, tensor product}
and \EqRef{kernel 2, algebra, tensor product}
belong to kernel of linear mapping $j$,
then, from equation
\EqRef{map i, algebra, tensor product},
it follows
\ShowEq{f, algebra, tensor product}
From equations \EqRef{f 1, algebra, tensor product}
and \EqRef{f 2, algebra, tensor product}
it follows that mapping $f$ is polylinear over ring $D$.
Since $M$ is module with basis
$\Times$, then, according to theorem \citeBib{Serge Lang}-4.1
on p. 135,
for any module $V$ and any
polylinear over $D$ map
\ShowEq{diagram bottom, 1, algebra, tensor product}
there exists a unique homomorphism $k:M\rightarrow V$,
for which following diagram is commutative
\ShowEq{diagram bottom, algebra, tensor product}
Since $g$ is polylinear over $D$, then
\ShowEq{ker k subseteq N}
According to statement 
on p. \citeBib{Serge Lang}-119,
map $j$ is universal in the category of homomorphisms of
vector space $M$
whose kernel contains $N$.
Therefore, we have homomorphism
\[
h:M/N\rightarrow V
\]
which makes the following diagram commutative
\begin{equation}
\EqLabel{diagram right, algebra, tensor product}
\xymatrix{
&M/N\ar[dd]^h
\\
M\ar[dr]^k\ar[ur]_j
\\
&V
}
\end{equation}
We join diagrams
\EqRef{diagram top, algebra, tensor product},
\EqRef{diagram bottom, algebra, tensor product},
\EqRef{diagram right, algebra, tensor product},
and get commutative diagram
\ShowEq{diagram, algebra, tensor product}
Since $\mathrm{Im}f$ generates $M/N$,
then map $h$ is uniquely determined.
\end{proof}

According to proof of theorem
\ref{theorem: tensor product of algebras is module}
\[A_1\otimes...\otimes A_n=M/N\]
If $d_i\in A_i$, we write
\ShowEq{map j, algebra, tensor product}

\begin{theorem}
\label{theorem: tensor product and polylinear mapping}
Let $A_1$, ..., $A_n$ be
algebras over commutative ring $D$.
Let
\ShowEq{map f, 1, algebra, tensor product}
be polylinear mapping defined by equation
\ShowEq{map f, algebra, tensor product}
Let
\ShowEq{map g, algebra, tensor product}
be polylinear mapping into $D$\Hyph module $V$.
There exists an $D$\Hyph linear mapping
\ShowEq{map h, algebra, tensor product}
such that the diagram
\ShowEq{map gh, algebra, tensor product}
is commutative.
\end{theorem}
\begin{proof}
Equation
\EqRef{map f, algebra, tensor product}
follows
from equations
\EqRef{map i, algebra, tensor product} and
\EqRef{map j, algebra, tensor product}.
An existence of the mapping $h$ follows from the definition
\ref{definition: tensor product of algebras}
and constructions made in the proof of the theorem
\ref{theorem: tensor product of algebras is module}.
\end{proof}

We can write equations \EqRef{f 1, algebra, tensor product}
and \EqRef{f 2, algebra, tensor product}
as
\ShowEq{tensors, algebra, tensor product}

\begin{theorem}
\label{theorem: tensor product and product in algebra}
Let $A$ be
algebra over commutative ring $D$.
There exists a linear mapping
\ShowEq{tensor product and product in algebra}
\end{theorem}
\begin{proof}
The theorem is corollary of the theorem
\ref{theorem: tensor product and polylinear mapping}
and the definition
\ref{definition: algebra over ring}.
\end{proof}
\else
\newpage
For $a_i\in A_i$, $i=1$, ..., $n$,
$a_1\otimes...\otimes a_n\in A_1\otimes...\otimes A_n$.
According to definition
\ShowEq{tensors, algebra, tensor product}

\ePrints{2010.05.23}
\ifx\Semafor\ValueOn
\newpage
\fi
\begin{remark}
\label{remark: representation of tensor, ring}
Representation of tensor as
\ShowEq{tensor representation, algebra}
is ambiguous.
Since $r\in D$, then
\[
(d_1 r)
\otimes
d_2
=d_1
\otimes
(rd_2 )
\]
We can increase or decrease number of summands
using algebraic operations.
Following transformation
\begin{align*}
&d_1
\otimes
d_2
+c_1
\otimes
c_2
\\
=&d_1
\otimes
(d_2-c_2+c_2)
+c_1
\otimes
c_2
\\
=&d_1
\otimes
(d_2-c_2)
+d_1
\otimes
c_2
+c_1
\otimes
c_2
\\
=&d_1
\otimes
(d_2-c_2)
+(d_1
+c_1)
\otimes
c_2
\end{align*}
is between possible transformations.
\qed
\end{remark}
\fi

\ePrints{2010.05.23}
\ifx\Semafor\ValueOn
\newpage
\fi
\begin{theorem}
\label{theorem: standard component of tensor, algebra}
Tensor product $\Tensor A$ of free
finite dimensional algebras
$A_1$, ..., $A_n$ over the commutative ring $D$ is free
finite dimensional algebra.

Let
\ShowEq{tensor product of algebras, basis i}
be the basis of algebra $A_i$ over ring $D$.
We can represent any tensor $a\in\Tensor A$ in the following form
\ShowEq{tensor canonical representation, algebra}
%Expression \EqRef{tensor canonical representation, algebra}
%is defined unique up to selected basis.
Expression
$\ShowSymbol{standard component of tensor, algebra}$
is called
\AddIndex{standard component of tensor}
{standard component of tensor, algebra}.
\ePrints{2010.05.23}
\ifx\Semafor\ValueOn
\qed
\end{theorem}
\newpage
\else
\end{theorem}
\begin{proof}
\ePrints{Calculus.Paper}
\Items{1011.3102}
\ifx\Semafor\ValueOn
See the proof of the theorem
\xRef{8433-5163}{theorem: standard component of tensor, algebra}.
\else
Algebras
$A_1$, ..., $A_n$
are modules over the ring $D$.
According to theorem
\ref{theorem: tensor product of algebras is module},
$\Tensor A$ is module. 

Vector $a_i\in A_i$ has expansion
\ShowEq{expansion i relative basis}
relative to basis $\Basis e_i$.
From equations
\EqRef{tensors 1, algebra, tensor product},
\EqRef{tensors 2, algebra, tensor product},
it follows
\ShowEq{expansion i relative basis, 1}
Since set of tensors $\Tensor a$
is the generating set of module $\Tensor A$,
then we can write tensor
$a\in\Tensor A$
in form
\ShowEq{tensor canonical representation 1, algebra}
where
\ShowEq{tensor canonical representation 2, algebra}
Let
\ShowEq{tensor canonical representation 3, algebra}
Then equation \EqRef{tensor canonical representation 1, algebra}
has form
\EqRef{tensor canonical representation, algebra}.

Therefore, set of tensors $\TensorBasis i$
is the generating set of module $\Tensor A$.
Since the dimension of module $A_i$, $i=1$, ..., $n$, is finite,
then the set of tensors $\TensorBasis i$ is finite.
Therefore, the set of tensors $\TensorBasis i$
contains a basis of module $\Tensor A$, and the module $\Tensor A$
is free module over the ring $D$.

We define the product of tensors like $\Tensor a$ componentwise
\ShowEq{product of tensors, algebra, tensor product}
In particular, if for any $i$, $i=1$, ..., $n$,
$a_i\in A_i$ has inverse,
then tensor
\ShowEq{inverse tensor, algebra, 2}
is inverse tensor to tensor
\ShowEq{inverse tensor, algebra, 1}

The definition of the product
\EqRef{product of tensors, algebra, tensor product}
agreed with equation
\EqRef{tensors 2, algebra, tensor product}
because
\ShowEq{product of tensors, algebra, tensor product, 1}
The distributive property of multiplication over
addition
\ShowEq{product of tensors, algebra, tensor product, 2}
follows from the equation
\EqRef{tensors 1, algebra, tensor product}.
The equation
\EqRef{product of tensors, algebra, tensor product, 2}
allows us to define the product for any tensors $a$, $b$.
\fi
\end{proof}

\ePrints{Calculus.Paper}
\Items{1011.3102}
\ifx\Semafor\ValueOff
\begin{remark}
\label{remark: tensor product of algebras}
According to the remark
\ref{remark: opposite algebra},
we can define different structures of algebra
in the tensor product of algebras.
For instance, algebras $A_1\otimes A_2$,
$A_1\otimes A_2^*$, $A_1^*\otimes A_2$
are defined in the same module.
\qed
\end{remark}
\fi
\fi

\ePrints{Calculus.Paper}
\Items{1011.3102}
\ifx\Semafor\ValueOff
\begin{theorem}
Let $\Basis e_i$ be the basis of the algebra $A_i$ over the ring $D$.
Let $B_{i\cdot}{}^{\gi j}_{\gi{kl}}$ be
structural constants of the algebra $A_i$ relative the basis $\Basis e_i$.
Structural constants of the tensor product $\Tensor A$
relative to the basis $\TensorBasis i$ have form
\ShowEq{structural constants tensor algebra}
\ePrints{2010.05.23}
\ifx\Semafor\ValueOn
\qed
\end{theorem}
\newpage
\else
\end{theorem}
\begin{proof}
Direct multiplication of tensors
$\TensorBasis i$ has form
\ShowEq{structural constants tensor algebra, proof, 1}
According to the definition of structural constants
\ShowEq{structural constants tensor algebra, proof, 2}
The equation
\EqRef{structural constants tensor algebra}
follows from comparison
\EqRef{structural constants tensor algebra, proof, 1},
\EqRef{structural constants tensor algebra, proof, 2}.

From the chain of equations
\ShowEq{structural constants tensor algebra, proof, 3}
it follows that definition of product
\EqRef{structural constants tensor algebra, proof, 2}
with structural constants
\EqRef{structural constants tensor algebra}
agreed with the definition of product
\EqRef{product of tensors, algebra, tensor product}.
\end{proof}

\begin{theorem}
For tensors $a$, $b\in\Tensor A$,
standard components of product satisfy to equation
\ShowEq{standard components of product, tensor product}
\end{theorem}
\begin{proof}
According to the definition
\ShowEq{standard components of product, tensor product, 1}
At the same time
\ShowEq{standard components of product, tensor product, 2}
The equation
\EqRef{standard components of product, tensor product}
follows from equations
\EqRef{standard components of product, tensor product, 1},
\EqRef{standard components of product, tensor product, 2}.
\end{proof}

\begin{theorem}
If the algebra $A_i$, $i=1$, ..., $n$, is associative,
then the tensor product $\Tensor A$ is associative algebra.
\end{theorem}
\begin{proof}
Since
\ShowEq{tensor product, associative algebra, 1}
then
\ShowEq{tensor product, associative algebra, 2}
\end{proof}
\fi
\fi

%% file: Tensor.Algebra.Eq.tex
%auto-ignore

\def\Times{A_1\times...\times A_n}

\newcommand{\Tensor}[1]{#1_1\otimes...\otimes #1_n}

\DefEquation
{
a^s d_{s\cdot 1}\otimes...\otimes d_{s\cdot n}
}
{tensor representation, algebra}

\DefEq
{
\begin{equation}
\EqLabel{map j, algebra, tensor product}
j\circ(d_1,..., d_n)=\Tensor d
\end{equation}
}
{map j, algebra, tensor product}

\DefEq
{
\[
h:a\otimes b\in A\otimes A\rightarrow ab\in A
\]
}
{tensor product and product in algebra}

\DefEq
{
\symb{A^{\otimes n}}
0{tensor power of algebra}
\[
\begin{matrix}
\ShowSymbol{tensor power of algebra}
=A_1\otimes...\otimes A_n
&
A_1=...=A_n=A
\end{matrix}
\]
}
{tensor power of algebra}

\DefEquation
{
f\circ(d_1,...,d_n)=\Tensor d
}
{map f, algebra, tensor product}

\DefEq
{
\[
f:A_1\times...\times A_n\rightarrow\Tensor A
\]
}
{map f, 1, algebra, tensor product}

\DefEq
{
\[
g:A_1\times...\times A_n\rightarrow V
\]
}
{map g, algebra, tensor product}

\DefEq
{
\[
h:\Tensor A\rightarrow V
\]
}
{map h, algebra, tensor product}

\DefEq
{
\begin{equation}
\begin{array}{r@{}l}
&a_1\otimes...\otimes(a_i+b_i)\otimes...\otimes a_n
\\
\VirtVar
=&a_1\otimes...\otimes a_i\otimes...\otimes a_n+
a_1\otimes...\otimes b_i\otimes...\otimes a_n
\end{array}
\EqLabel{tensors 1, algebra, tensor product}
\end{equation}
\begin{equation}
a_1\otimes...\otimes(ca_i)\otimes...\otimes a_n=
c(a_1\otimes...\otimes a_i\otimes...\otimes a_n)
\EqLabel{tensors 2, algebra, tensor product}
\end{equation}
\[
\begin{matrix}
a_i\in A_i&b_i\in A_i&c\in D
\end{matrix}
\]
}
{tensors, algebra, tensor product}

\DefEq
{
\begin{equation}
\EqLabel{product of tensors, algebra, tensor product, 2}
\begin{array}{r@{}l}
&(a_1\otimes...\otimes a_i\otimes...\otimes a_n)
\\
*&
((b_1\otimes...\otimes b_i\otimes...\otimes b_n)+
(b_1\otimes...\otimes c_i\otimes...\otimes b_n))
\\
=&
(a_1\otimes...\otimes a_i\otimes...\otimes a_n)
(b_1\otimes...\otimes (b_i+c_i)\otimes...\otimes b_n)
\\
=&
(a_1b_1)\otimes...\otimes (a_i(b_i+c_i))\otimes...\otimes (a_nb_n)
\\
=&
(a_1b_1)\otimes...\otimes (a_ib_i+a_ic_i)\otimes...\otimes (a_nb_n)
\\
=&
(a_1b_1)\otimes...\otimes (a_ib_i)\otimes...\otimes (a_nb_n)
\\
+&
(a_1b_1)\otimes...\otimes (a_ic_i)\otimes...\otimes (a_nb_n)
\\
=&
(a_1\otimes...\otimes a_i\otimes...\otimes a_n)
(b_1\otimes...\otimes b_i\otimes...\otimes b_n)
\\
+&
(a_1\otimes...\otimes a_i\otimes...\otimes a_n)
(b_1\otimes...\otimes c_i\otimes...\otimes b_n)
\end{array}
\end{equation}
}
{product of tensors, algebra, tensor product, 2}

\DefEq
{
\begin{align*}
&(a_1\otimes...\otimes (ca_i)\otimes...\otimes a_n)
(b_1\otimes...\otimes b_i\otimes...\otimes b_n)
\\
=&
(a_1b_1)\otimes...\otimes (ca_i)b_i\otimes...\otimes (a_nb_n)
\\
=&
c((a_1b_1)\otimes...\otimes (a_ib_i)\otimes...\otimes (a_nb_n))
\\
=&
c((\Tensor a)(\Tensor b))
\end{align*}
}
{product of tensors, algebra, tensor product, 1}

\DefEq
{
\[
\xymatrix{
f:\Times\ar[r]&S_1
&
g:\Times\ar[r]&S_2
}
\]
}
{polylinear mappings category}

\DefEq
{
\[
\xymatrix{
&S_1\ar[dd]^h
\\
\Times
\ar[ru]^f\ar[rd]_g
\\
&S_2
}
\]
}
{polylinear mappings category, diagram}

\DefEquation
{
\xymatrix{
&\Tensor A\ar[dd]^h
\\
\Times
\ar[ru]^f\ar[rd]_g
\\
&V
}
}
{map gh, algebra, tensor product}

\DefEq
{
\symb{\Tensor A}1{tensor product of algebras}
}
{tensor product of algebras}

\DefEq
{
\[
\xymatrix{
i:\Times
\ar[r]&M
}
\]
}
{map i, algebra, tensor product, 1}

\DefEq
{
$\mathrm{ker}\ k\subseteq N$.
}
{ker k subseteq N}

\DefEq
{
\begin{align}
&(d_1,...,d_i+c_i,...,d_n)-
(d_1,..., d_i,..., d_n)-
(d_1,..., c_i,..., d_n)
\EqLabel{kernel 1, algebra, tensor product}
\\
&(d_1,...,ad_i,..., d_n)-
a(d_1,..., d_i,..., d_n)
\EqLabel{kernel 2, algebra, tensor product}
\end{align}
}
{kernel, algebra, tensor product}

\DefEq
{
\[
\xymatrix{
g:\Times
\ar[r]&V
}
\]
}
{diagram bottom, 1, algebra, tensor product}

\DefEq
{
\begin{equation}
\EqLabel{diagram, algebra, tensor product}
\xymatrix{
&&&M/N\ar[dd]^h
\\
\Times\ar[rrru]^f\ar[rrrd]_g\ar[rr]_(.7)i
&&M\ar[ur]_j\ar[dr]^k
\\
&&&V
}
\end{equation}
}
{diagram, algebra, tensor product}

\DefEq
{
\[
a_i=a_i^{\gi k}\Vector e_{i\cdot\gi k}
\]
}
{expansion i relative basis}

\newcommand\TensorBasis[1]
{\Vector e_{1\cdot\gi{#1_1}}
\otimes...\otimes\Vector e_{n\cdot\gi{#1_n}}}

\DefEq
{
\begin{equation}
\EqLabel{structural constants tensor algebra}
C^{\cdot\gi{j_1...j_n}}_{\cdot\gi{k_1...k_n}\cdot\gi{l_1...l_n}}=
C_{1\cdot}{}^{\gi{j_1}}_{\gi{k_1l_1}}...C_{n\cdot}{}^{\gi{j_n}}_{\gi{k_nl_n}}
\end{equation}
}
{structural constants tensor algebra}

\DefEq
{
\begin{align}
\nonumber
&
(\TensorBasis k)
(\TensorBasis l)
\\
\nonumber
=&
(\Vector e_{1\cdot\gi{k_1}}\Vector e_{1\cdot\gi{l_1}})
\otimes...\otimes
(\Vector e_{n\cdot\gi{k_n}}\Vector e_{n\cdot\gi{l_n}})
\\
\EqLabel{structural constants tensor algebra, proof, 1}
=&
(\Vector e_{1\cdot\gi{k_1}}\Vector e_{1\cdot\gi{l_1}})
\otimes...\otimes
(\Vector e_{n\cdot\gi{k_n}}\Vector e_{n\cdot\gi{l_n}})
\\
\nonumber
=&
(C_{1\cdot}{}^{\gi{j_1}}_{\gi{k_1l_1}}\Vector e_{1\cdot\gi{j_1}})
\otimes...\otimes
(C_{n\cdot}{}^{\gi{j_n}}_{\gi{k_nl_n}}\Vector e_{n\cdot\gi{j_n}})
\\
\nonumber
=&
C_{1\cdot}{}^{\gi{j_1}}_{\gi{k_1l_1}}...C_{n\cdot}{}^{\gi{j_n}}_{\gi{k_nl_n}}
\Vector e_{1\cdot\gi{j_1}}\otimes...\otimes\Vector e_{n\cdot\gi{j_n}}
\end{align}
}
{structural constants tensor algebra, proof, 1}

\DefEq
{
\begin{equation}
\EqLabel{structural constants tensor algebra, proof, 2}
(\TensorBasis k)
(\TensorBasis l)
=
C^{\cdot\gi{j_1...j_n}}_{\cdot\gi{k_1...k_n}\cdot\gi{l_1...l_n}}
(\TensorBasis j)
\end{equation}
}
{structural constants tensor algebra, proof, 2}

\DefEq
{
\begin{align*}
&(\Tensor a)(\Tensor b)
\\
=&
(a_1^{\gi{k_1}}\Vector e_{1\cdot\gi{k_1}}
\otimes...\otimes
a_n^{\gi{k_n}}\Vector e_{n\cdot\gi{k_n}})
(b_1^{\gi{l_1}}\Vector e_{1\cdot\gi{l_1}}
\otimes...\otimes
b_n^{\gi{l_n}}\Vector e_{n\cdot\gi{l_n}})
\\
=&
a_1^{\gi{k_1}}...a_n^{\gi{k_n}}
b_1^{\gi{l_1}}...b_n^{\gi{l_n}}
(\Vector e_{1\cdot\gi{k_1}}
\otimes...\otimes
\Vector e_{n\cdot\gi{k_n}})
(\Vector e_{1\cdot\gi{l_1}}
\otimes...\otimes
\Vector e_{n\cdot\gi{l_n}})
\\
=&
a_1^{\gi{k_1}}...a_n^{\gi{k_n}}
b_1^{\gi{l_1}}...b_n^{\gi{l_n}}
C^{\cdot\gi{j_1...j_n}}_{\cdot\gi{k_1...k_n}\cdot\gi{l_1...l_n}}
(\TensorBasis j)
\\
=&
a_1^{\gi{k_1}}...a_n^{\gi{k_n}}
b_1^{\gi{l_1}}...b_n^{\gi{l_n}}
C_{1\cdot}{}^{\gi{j_1}}_{\gi{k_1l_1}}\ ...C_{n\cdot}{}^{\gi{j_n}}_{\gi{k_nl_n}}
(\TensorBasis j)
\\
=&
(a_1^{\gi{k_1}}b_1^{\gi{l_1}}\ C_{1\cdot}{}^{\gi{j_1}}_{\gi{k_1l_1}}
\Vector e_{1\cdot\gi{j_1}})
\otimes...\otimes
(a_n^{\gi{k_n}}b_n^{\gi{l_n}}C_{n\cdot}{}^{\gi{j_n}}_{\gi{k_nl_n}}
\Vector e_{n\cdot\gi{j_n}})
\\
=&
(a_1b_1)
\otimes...\otimes
(a_nb_n)
\end{align*}
}
{structural constants tensor algebra, proof, 3}

\DefEq
{
\begin{equation}
\EqLabel{standard components of product, tensor product}
(ab)^{\gi{j_1...j_n}}=
C^{\cdot\gi{j_1...j_n}}_{\cdot\gi{k_1...k_n}\cdot\gi{l_1...l_n}}
a^{\gi{k_1...k_n}}b^{\gi{l_1...l_n}}
\end{equation}
}
{standard components of product, tensor product}

\DefEq
{
\[
\Tensor a=a_1^{\gi{i_1}}...a_n^{\gi{i_n}}
\TensorBasis i
\]
}
{expansion i relative basis, 1}

\DefEq
{
\begin{equation}
\EqLabel{standard components of product, tensor product, 1}
ab=(ab)^{\gi{j_1...j_n}}\TensorBasis j
\end{equation}
}
{standard components of product, tensor product, 1}

\DefEq
{
\begin{equation}
\EqLabel{standard components of product, tensor product, 2}
\begin{array}{r@{}l}
ab
&=
a^{\gi{k_1...k_n}}\TensorBasis k
b^{\gi{k_1...k_n}}\TensorBasis l
\\
\VirtVar&=
a^{\gi{k_1...k_n}}b^{\gi{k_1...k_n}}
C^{\cdot\gi{j_1...j_n}}_{\cdot\gi{k_1...k_n}\cdot\gi{l_1...l_n}}
\TensorBasis j
\end{array}
\end{equation}
}
{standard components of product, tensor product, 2}

\DefEq
{
\begin{align*}
&((\TensorBasis i)(\TensorBasis j))(\TensorBasis k)
\\
=&
((\Vector e_{1\cdot\gi{i_1}}\Vector e_{1\cdot\gi{j_1}})
\otimes...\otimes
(\Vector e_{n\cdot\gi{i_n}}\Vector e_{1\cdot\gi{j_n}}))
(\TensorBasis k)
\\
=&
((\Vector e_{1\cdot\gi{i_1}}\Vector e_{1\cdot\gi{j_1}})
\Vector e_{1\cdot\gi{k_1}})
\otimes...\otimes
((\Vector e_{n\cdot\gi{i_n}}\Vector e_{1\cdot\gi{j_n}})
\Vector e_{1\cdot\gi{k_n}})
\\
=&
(\Vector e_{1\cdot\gi{i_1}}(\Vector e_{1\cdot\gi{j_1}}
\Vector e_{1\cdot\gi{k_1}}))
\otimes...\otimes
(\Vector e_{n\cdot\gi{i_n}}(\Vector e_{1\cdot\gi{j_n}}
\Vector e_{1\cdot\gi{k_n}}))
\\
=&
(\TensorBasis i)
((\Vector e_{1\cdot\gi{j_1}}\Vector e_{1\cdot\gi{k_1}})
\otimes...\otimes
(\Vector e_{1\cdot\gi{j_n}}\Vector e_{1\cdot\gi{k_n}}))
\\
=&
(\TensorBasis i)
((\TensorBasis j)(\TensorBasis k))
\end{align*}
}
{tensor product, associative algebra, 1}

\DefEq
{
\begin{align*}
(ab)c
=&
a^{\gi{i_1...i_n}}b^{\gi{j_1...j_n}}c^{\gi{k_1...k_n}}
\\
&((\TensorBasis i)(\TensorBasis j))(\TensorBasis k)
\\
=&
a^{\gi{i_1...i_n}}b^{\gi{j_1...j_n}}c^{\gi{k_1...k_n}}
\\
&(\TensorBasis i)((\TensorBasis j)(\TensorBasis k))
\\
=&
a(bc)
\end{align*}
}
{tensor product, associative algebra, 2}

\DefEq
{
\symb{a^{\gi{i_1...i_n}}}0
{standard component of tensor, algebra}
\begin{equation}
a=
\ShowSymbol{standard component of tensor, algebra}
\TensorBasis i
\EqLabel{tensor canonical representation, algebra}
\end{equation}
}
{tensor canonical representation, algebra}

\DefEq
{
$\Basis e_i$
}
{tensor product of algebras, basis i}

\DefEq
{
\begin{equation}
a=a^sa_{s\cdot}{}_1^{\gi{i_1}}...a_{s\cdot}{}_n^{\gi{i_n}}
\TensorBasis i
\EqLabel{tensor canonical representation 1, algebra}
\end{equation}
}
{tensor canonical representation 1, algebra}

\DefEq
{
$a^s$, $a_{s\cdot}{}_1^{\gi{i_1}}$, $a_{s\cdot}{}_n^{\gi{i_n}}\in F$.
}
{tensor canonical representation 2, algebra}

\DefEq
{
\begin{equation}
a^sa_{s\cdot}{}_1^{\gi{i_1}}...a_{s\cdot}{}_n^{\gi{i_n}}
=a^{\gi{i_1...i_n}}
\EqLabel{tensor canonical representation 3, algebra}
\end{equation}
}
{tensor canonical representation 3, algebra}

\DefEq
{
\[
\Tensor a\in\Tensor A
\]
}
{inverse tensor, algebra, 1}

\DefEq
{
\[
(\Tensor a)^{-1}=(a_1)^{-1}\otimes...\otimes (a_n)^{-1}
\]
}
{inverse tensor, algebra, 2}

\DefEq
{
\begin{equation}
\EqLabel{diagram bottom, algebra, tensor product}
\xymatrix{
\Times\ar[rrrd]_g\ar[rr]_(.7)i
&&M\ar[dr]^k
\\
&&&V
}
\end{equation}
}
{diagram bottom, algebra, tensor product}

\DefEquation
{
\xymatrix{
&&&M/N
\\
\Times\ar[rrru]^f\ar[rr]_(.7)i
&&M\ar[ur]_j
}
}
{diagram top, algebra, tensor product}

\DefEq
{
\begin{equation}
(\Tensor d)(\Tensor c)=(d_1c_1)\otimes...\otimes(d_nc_n)
\EqLabel{product of tensors, algebra, tensor product}
\end{equation}
}
{product of tensors, algebra, tensor product}

\DefEquation
{
i\circ(d_1,...,d_n)=(d_1,...,d_n)
}
{map i, algebra, tensor product}

\DefEq
{
\begin{align}
f\circ(d_1,...,d_i+c_i,..., d_n)
=&f\circ(d_1,..., d_i,..., d_n)+
f\circ(d_1,..., c_i,..., d_n)
\EqLabel{f 1, algebra, tensor product}
\\
f\circ(d_1,...,ad_i,..., d_n)=&
a\ f\circ(d_1,..., d_i,..., d_n)
\EqLabel{f 2, algebra, tensor product}
\end{align}
}
{f, algebra, tensor product}

%% file: Derivative.English.tex
%auto-ignore

\input{\FilePrefix Derivative.Eq}

\Chapter{Differentiable Mappings}
\label{chapter: Differentiable maps}

\ePrints{Calculus.Paper}
\ifx\Semafor\ValueOn
In this section we explore derivative and differential of the mapping
into $D$\hyph algebra.
Complex field is the algebra over real field.
In the calculus of functions of complex variable, we consider
linear mappings generated by the mapping $$I_0\circ z=z$$
In this section, we also consider linear mappings generated by the
mapping $$I_0\circ a=a$$
\fi%Calculus.Paper

\ePrints{2010.05.23}
\ifx\Semafor\ValueOff
\Section{Topological Ring}
\label{section: Topological Ring}

\begin{definition}
Ring $D$ is called 
\AddIndex{topological ring}{topological ring}\footnote{I
made definition according to definition
from \citeBib{Pontryagin: Topological Group},
chapter 4}
if $D$ is topological space and the algebraic operations
defined in $D$ are continuous in the topological space $D$.
\qed
\end{definition}

According to definition, for arbitrary elements $a$, $b\in D$
and for arbitrary neighborhoods $W_{a-b}$ of the element $a-b$,
$W_{ab}$ of the element $ab$ there exists neighborhoods $W_a$
of the element $a$ and $W_b$ of the element $b$ such
that $W_a-W_b\subset W_{a-b}$,
$W_aW_b\subset W_{ab}$.

\begin{definition}
\label{definition: norm on ring}
\AddIndex{Norm on ring}
{norm on ring} $D$\footnote{I
made definition according to definition
from \citeBib{Bourbaki: General Topology: Chapter 5 - 10},
IX, \S 3.2 and definition \citeBib{Arnautov Glavatsky Mikhalev}-1.1.12,
p. 23.} is a mapping
\[d\in D\rightarrow |d|\in R\]
which satisfies the following axioms
\begin{itemize}
\item $|a|\ge 0$
\item $|a|=0$ if, and only if, $a=0$
\item $|ab|=|a|\ |b|$
\item $|a+b|\le|a|+|b|$
\end{itemize}

Ring $D$, endowed with the structure defined by a given absolute value on
$D$, is called \AddIndex{normed ring}{normed ring}.
\qed
\end{definition}

Invariant distance on additive group of ring $D$
\[d(a,b)=|a-b|\]
defines topology of metric space,
compatible with ring structure of $D$.

\begin{definition}
\label{definition: limit of sequence, normed ring}
Let $D$ be normed ring.
Element $a\in D$ is called 
\AddIndex{limit of a sequence}{limit of sequence, normed ring}
\ShowEq{limit of sequence, normed ring}
if for every $\epsilon\in R$, $\epsilon>0$
there exists positive integer $n_0$ depending on $\epsilon$ and such,
that $|a_n-a|<\epsilon$ for every $n>n_0$.
\qed
\end{definition}

\begin{theorem}
\label{theorem: limit of sequence, normed ring, product on scalar}
Let $D$ be normed ring of characteristic $0$ and let $d\in D$.
Let $a\in D$ be limit of a sequence $\{a_n\}$.
Then
\[
\lim_{n\rightarrow\infty}(a_nd)=ad
\]
\[
\lim_{n\rightarrow\infty}(da_n)=da
\]
\end{theorem}
\begin{proof}
Statement of the theorem is trivial, however I give this proof
for completeness sake. 
Since $a\in D$ is limit of the sequence $\{a_n\}$,
then according to definition
\ref{definition: limit of sequence, normed ring}
for given $\epsilon\in R$, $\epsilon>0$,
there exists positive integer $n_0$ such,
that $|a_n-a|<\epsilon/|d|$ for every $n>n_0$.
According to definition \ref{definition: norm on ring}
the statement of theorem follows from inequalities
\[
|a_nd-ad|=|(a_n-a)d|=|a_n-a||d|<\epsilon/|d||d|=\epsilon
\]
\[
|da_n-da|=|d(a_n-a)|=|d||a_n-a|<|d|\epsilon/|d|=\epsilon
\]
for any $n>n_0$.
\end{proof}

\begin{definition}
Let $D$ be normed ring.
The sequence $\{a_n\}$, $a_n\in D$ is called 
\AddIndex{fundamental}{fundamental sequence, normed ring}
or \AddIndex{Cauchy sequence}{Cauchy sequence, normed ring},
if for every $\epsilon\in R$, $\epsilon>0$
there exists positive integer $n_0$ depending on $\epsilon$ and such,
that $|a_p-a_q|<\epsilon$ for every $p$, $q>n_0$.
\qed
\end{definition}

\begin{definition}
Normed ring $D$ is called
\AddIndex{complete}{complete ring}
if any fundamental sequence of elements
of ring $D$ converges, i.e.
has limit in ring $D$.
\qed
\end{definition}

Later on, speaking about normed ring of characteristic $0$,
we will assume that homeomorphism of field of rational numbers $Q$
into ring $D$ is defined.

\begin{theorem}
\label{theorem: complete ring contains real number}
Complete ring $D$ of characteristic $0$
contains as subfield an isomorphic image of the field $R$ of
real numbers.
It is customary to identify it with $R$.
\end{theorem}
\begin{proof}
Let us consider fundamental sequence of rational numbers $\{p_n\}$.
Let $p'$ be limit of this sequence in division ring $D$.
Let $p$ be limit of this sequence in field $R$.
Since immersion of field $Q$ into division ring $D$ is homeomorphism,
then we may identify $p'\in D$ and $p\in R$.
\end{proof}

\begin{theorem}
\label{theorem: complete ring and real number}
Let $D$ be complete ring of characteristic $0$ and let $d\in D$.
Then any real number $p\in R$ commute with $d$.
\end{theorem}
\begin{proof}
Let us represent real number $p\in R$ as
fundamental sequence of rational numbers $\{p_n\}$.
Statement of theorem follows from chain of equations
\[
pd=\lim_{n\rightarrow\infty}(p_nd)=\lim_{n\rightarrow\infty}(dp_n)=dp
\]
based on statement of theorem
\ref{theorem: limit of sequence, normed ring, product on scalar}.
\end{proof}

\Section{Topological \texorpdfstring{$D$}{D}-Algebra}

\begin{definition}
Given a topological commutative ring $D$ and
$D$\Hyph algebra $A$ such that $A$ has
a topology compatible with the structure of the additive
group of $A$ and mappings
\ShowEq{topological D algebra}
are continuous,
then $\Vector V$ is called a
\AddIndex{topological $D$\Hyph algebra}
{topological D algebra}\footnote{I
made definition according to definition
from \citeBib{Bourbaki: Topological Vector Space},
p. TVS I.1}.
\qed
\end{definition}

\begin{definition}
\label{definition: norm on d algebra}
\AddIndex{Norm on $D$\Hyph algebra}
{norm on D algebra} $A$
over normed commutative ring $D$\footnote{I
made definition according to definition
from \citeBib{Bourbaki: General Topology: Chapter 5 - 10},
IX, \S 3.3} is a map
\ShowEq{norm on d algebra}
which satisfies the following axioms
\begin{itemize}
\ShowEq{norm on d algebra 1}
\ShowEq{norm on d algebra 2, 1}
if, and only if,
\ShowEq{norm on d algebra 2, 2}
\ShowEq{norm on d algebra 3}
\end{itemize}

If $D$ is a normed commutative ring,
$D$\Hyph algebra $A$,
endowed with the structure defined by a given norm on
$A$, is called
\AddIndex{normed $D$\Hyph algebra}{normed D algebra}.
\qed
\end{definition}

\begin{definition}
\label{definition: limit of sequence, normed algebra}
Let $A$ be normed $D$\Hyph algebra.
Element $a\in A$ is called 
\AddIndex{limit of a sequence}{limit of sequence, normed algebra}
\ShowEq{limit of sequence, normed algebra}
if for every $\epsilon\in R$, $\epsilon>0$
there exists positive integer $n_0$ depending on $\epsilon$ and such,
that $|a_n-a|<\epsilon$ for every $n>n_0$.
\qed
\end{definition}

\begin{definition}
Let $A$ be normed $D$\Hyph algebra.
The sequence $\{a_n\}$, $a_n\in A$ is called 
\AddIndex{fundamental}{fundamental sequence, normed algebra}
or \AddIndex{Cauchy sequence}{Cauchy sequence, normed algebra},
if for every $\epsilon\in R$, $\epsilon>0$
there exists positive integer $n_0$ depending on $\epsilon$ and such,
that $|a_p-a_q|<\epsilon$ for every $p$, $q>n_0$.
\qed
\end{definition}

\begin{definition}
Normed $D$\Hyph algebra $A$ is called
\AddIndex{Banach $D$\Hyph algebra}{Banach algebra}
if any fundamental sequence of elements
of algebra $A$ converges, i.e.
has limit in algebra $A$.
\qed
\end{definition}

\begin{definition}
Let $A$ be Banach $D$\Hyph algebra.
Set of elements
\ShowEq{unit sphere in algebra}
is called
\AddIndex{unit sphere in algebra}{unit sphere in algebra} $A$.
\qed
\end{definition}

\begin{definition}
\label{definition: continuous function, algebra}
Mapping
\ShowEq{f A1 A2}
of Banach $D_1$\Hyph algebra $A_1$ with norm $|x|_1$
into Banach $D_2$\Hyph algebra $A_2$ with norm $|y|_2$
is called \AddIndex{continuous}{continuous function, algebra}, if
for every as small as we please $\epsilon>0$
there exist such $\delta>0$, that
\[
|x'-x|_1<\delta
\]
implies
\[
|f(x')-f(x)|_2<\epsilon
\]
\qed
\end{definition}

\begin{definition}
\label{definition: norm of map, algebra}
Let
\ShowEq{f A1 A2}
mapping of Banach $D_1$\Hyph algebra $A_1$ with norm $|x|_1$
into Banach $D_2$\Hyph algebra $A_2$ with norm $|y|_2$.
Value
\ShowEq{norm of map, algebra}
is called
\AddIndex{norm of mapping $f$}
{norm of map, algebra}.
\qed
\end{definition}

\begin{theorem}
Let
\ShowEq{f A1 A2}
linear mapping of Banach $D_1$\Hyph algebra $A_1$ with norm $|x|_1$
into Banach $D_2$\Hyph algebra $A_2$ with norm $|y|_2$.
Then
\ShowEq{norm of linear map, algebra}
\end{theorem}
\begin{proof}
From the definition
\ref{definition: linear map from A1 to A2, algebra}
and the theorem
\ref{theorem: complete ring contains real number},
it follows that
\ShowEq{linear map and real number}
From the equation \EqRef{linear map and real number}
and the definition \ref{definition: norm on d algebra}
it follows that
\ShowEq{linear map and real number, 1}
Assuming $\displaystyle r=\frac 1{|x|_1}$, we get
\ShowEq{norm of linear map, algebra, 1}
Equation \EqRef{norm of linear map, algebra}
follows from equations \EqRef{norm of linear map, algebra, 1}
and \EqRef{norm of map, algebra}.
\end{proof}

\begin{theorem}
Let
\ShowEq{f A1 A2}
linear mapping of Banach $D_1$\Hyph algebra $A_1$ with norm $|x|_1$
into Banach $D_2$\Hyph algebra $A_2$ with norm $|y|_2$.
Since $\|f\|<\infty$, then map $f$ is continuous.
\end{theorem}
\begin{proof}
Since map $f$ is linear, then
according to definition \ref{definition: norm of map, algebra}
\[
|f(x)-f(y)|_2=|f(x-y)|_2\le \|f\|\ |x-y|_1
\]
Let us assume arbitrary $\epsilon>0$. Assume $\displaystyle\delta=\frac\epsilon{\|f\|}$.
Then
\[
|f(x)-f(y)|_2\le \|f\|\ \delta=\epsilon
\]
follows from inequality
\[
|x-y|_1<\delta
\]
According to definition \ref{definition: continuous function, algebra}
map $f$ is continuous.
\end{proof}

\ifx\texFuture\Defined
\begin{theorem}
\label{theorem: continuity of map, norm, division ring}
Let
\[f:D_1\rightarrow D_2\]
mapping of complete division ring $D_1$
в complete division ring $D_2$.
map $f$ is continuous, if $\|f\|<\infty$.
\end{theorem}
\begin{proof}
Since map $f$ is additive, then
according to definition \ref{definition: norm of map, algebra}
\[
|f(x)-f(y)|_2=|f(x-y)|_2\le \|f\|\ |x-y|_1
\]
Let $\epsilon>0$ be arbitrary. Assume $\displaystyle\delta=\frac\epsilon{\|f\|}$.
Then from inequality
\[
|x-y|_1<\delta
\]
follows
\[
|f(x)-f(y)|_2\le \|f\|\ \delta=\epsilon
\]
According to definition \ref{definition: continuous function, algebra}
mapping $f$ is continuous.
\end{proof}
\fi
\else
\newpage
\fi

\Section{The Derivative of Mapping in Algebra}
\label{section: The Derivative of Mapping in Algebra}

\begin{definition}
\label{definition: function differentiable in Gateaux sense, algebra}
Let $A$ be Banach $D$\Hyph algebra.
The function
\ShowEq{f:A->A}
is called \AddIndex{differentiable in the G\^ateaux sense}
{function differentiable in Gateaux sense, algebra}
on the set $U\subset A$,
if at every point $x\in U$
the increment of the function $f$ can be represented as
\ShowEq{Gateaux derivative of map, algebra}
where
\AddIndex{the G\^ateaux derivative
$\ShowSymbol{Gateaux derivative of map, algebra}$
of map $f$}
{Gateaux derivative of map, algebra}
is linear map of increment $a$ and
\ShowEq{Gateaux derivative of map, algebra, 1}
is such continuous map that
\ShowEq{Gateaux derivative of map, algebra, 2}
\qed
\end{definition}
\ePrints{2010.05.23}
\ifx\Semafor\ValueOn
\newpage
\fi

\begin{remark}
\label{remark: differential L(A,A)}
According to definition
\ref{definition: function differentiable in Gateaux sense, algebra}
for given $x$, the G\^ateaux derivative
\ShowEq{differential L(A,A), 1}
Therefore, the G\^ateaux derivative of map $f$ is map
\ShowEq{differential L(A,A), 2}
Expressions $\partial f(x)$ and $\displaystyle\frac{\partial f(x)}{\partial x}$
are different notations for the same function.
We will use notation $\displaystyle\frac{\partial f(x)}{\partial x}$
to underline that this is  the G\^ateaux derivative with respect to variable $x$.
\qed
\end{remark}
\ePrints{2010.05.23}
\ifx\Semafor\ValueOn
\newpage
\fi

\begin{theorem}
\label{theorem: function differentiable in Gateaux sense, algebra}
It is possible to represent \AddIndex{the G\^ateaux differential
\ShowEq{Gateaux differential of map, algebra}
of mapping $f$}
{Gateaux differential of map, algebra}
as\footnote{Formally, we have to write the differential of the mapping in the form
\ShowEq{Gateaux differential, representation, algebra, 1}
However, for instance, in the theory of functions of complex variable we consider
only linear mappings generated by mapping $I_0\circ z=z$.
Therefore, exploring derivatives, we also restrict ourselves
to linear mappings generated by the mapping $I_0$.
To write expressions in the general case is not difficult.
}
\ShowEq{Gateaux differential, representation, algebra}
Expression
\ShowEq{component of Gateaux derivative of map, algebra}
is called
\AddIndex{component of the G\^ateaux derivative of map $f(x)$}
{component of Gateaux derivative of map, algebra}.
\ePrints{2010.05.23}
\ifx\Semafor\ValueOn
\qed
\end{theorem}
\newpage
\else
\end{theorem}
\begin{proof}
Corollary of definitions
\ref{definition: function differentiable in Gateaux sense, algebra}
and theorem
\ref{theorem: linear mapping in L(A,A), associative algebra}.
\end{proof}

From definitions
\ref{definition: linear map from A1 to A2, algebra},
\ref{definition: function differentiable in Gateaux sense, algebra}
and the theorem
\ref{theorem: complete ring contains real number}
it follows
\ShowEq{Gateaux differential is multiplicative over field R, algebra}
Combining equation
\EqRef{Gateaux differential is multiplicative over field R, algebra}
and definition
\ref{definition: function differentiable in Gateaux sense, algebra},
we get known definition of the G\^ateaux
differential
\ShowEq{Gateaux differential, representation in algebra}
Definitions of the G\^ateaux derivative
\EqRef{Gateaux derivative of map, algebra}
and \EqRef{Gateaux differential, representation in algebra}
are equivalent. Using this equivalence we tell
that map $f$
is called differentiable in the G\^ateaux sense on the set $U\subset D$,
if at every point $x\in U$
the increment of the function $f$ can be represented as
\ShowEq{Gateaux differential of map, t, algebra}
where
$o:R\rightarrow A$ is such continuous map that
\[
\lim_{t\rightarrow 0}\frac{|o(t)|}{|t|}=0
\]
\fi

\begin{theorem}
\label{theorem: Gateaux differential, standard form, algebra}
Let $A$ be Banach $D$\Hyph algebra.
Let $\Basis e$ be basis of algebra $A$ over ring $D$.
\AddIndex{Standard representation of the G\^ateaux derivative of mapping}
{Gateaux derivative, standard representation, algebra}
\ShowEq{f:A->A}
has form
\ShowEq{Gateaux derivative, algebra, standard representation}
Expression
\ShowEq{standard component of Gateaux derivative, algebra}
in equation
\EqRef{Gateaux derivative, algebra, standard representation}
is called
\AddIndex{standard component of the G\^ateaux derivative of mapping $f$}
{standard component of Gateaux derivative, algebra}.
\ePrints{2010.05.23}
\ifx\Semafor\ValueOn
\qed
\end{theorem}
\newpage
\else
\end{theorem}
\begin{proof}
Statement of theorem is corollary of
statement \eqref{f in L(A,A), 2, associative algebra} of the theorem
\ref{theorem: linear mapping in L(A,A), associative algebra}.
\end{proof}

\begin{theorem}
\label{theorem: Gateaux differential and jacobian, algebra}
Let $A$ be Banach $D$\Hyph algebra.
Let $\Basis e$ be basis of algebra $A$ over ring $D$.
Then it is possible to represent the G\^ateaux differential of mapping
\[
f:D\rightarrow D
\]
as
\ShowEq{Gateaux differential and jacobian, algebra}
where $dx\in A$ has expansion
\ShowEq{Gateaux differential and jacobian, 1, algebra}
relative to basis $\Basis e$ and Jacobian matrix of map $f$ has form
\ShowEq{standard components and Jacobian, algebra}
\end{theorem}
\begin{proof}
Statement of theorem is corollary of theorem
\ref{theorem: linear mapping in L(A,A), associative algebra}.
\end{proof}

\begin{theorem}
\label{theorem: derivative of the sum}
Let $A$ be Banach $D$\Hyph algebra.
Let $f$, $g$ be differentiable mappings
\ShowEq{bilinear map and derivative, 1}
The mapping
\ShowEq{derivative of the sum, 2}
is differentiable and
the G\^ateaux derivative satisfies to relationship
\ShowEq{derivative of the sum, 3}
\end{theorem}
\begin{proof}
According to the definition
\EqRef{Gateaux differential, representation in algebra},
\ShowEq{derivative of the sum, 4}
The equation \EqRef{derivative of the sum, 3}
follows from the equation \EqRef{derivative of the sum, 4}.
\end{proof}

\begin{theorem}
\label{theorem: bilinear map and derivative}
Let $A$ be Banach $D$\Hyph algebra.
Let
\ShowEq{bilinear map and derivative}
be continuous bilinear mapping.
Let $f$, $g$ be differentiable mappings
\ShowEq{bilinear map and derivative, 1}
The mapping
\ShowEq{bilinear map and derivative, 2}
is differentiable and
the G\^ateaux differential satisfies to relationship
\ShowEq{bilinear map and derivative, 3}
\end{theorem}
\begin{proof}
Equation \EqRef{bilinear map and derivative, 3}
follows from chain of equations
\ShowEq{bilinear map and derivative, 4}
based on definition
\EqRef{Gateaux differential, representation in algebra}.
\end{proof}
\fi

\begin{theorem}
Let $A$ be Banach $D$\Hyph algebra.
Let $f$, $g$ be differentiable mappings
\ShowEq{bilinear map and derivative, 1}
The G\^ateaux differential satisfies to relationship
\ShowEq{Gateaux differential of product, algebra}
\ePrints{2010.05.23}
\ifx\Semafor\ValueOn
\qed
\end{theorem}
\newpage
{\bf Table of Derivatives of Mappings of Associative Algebra}
\else
\end{theorem}
\begin{proof}
The theorem is corollary of theorem
\ref{theorem: bilinear map and derivative}
and definition \ref{definition: algebra over ring}.
\end{proof}

\begin{theorem}
Let $A$ be Banach $D$\Hyph algebra.
Suppose the G\^ateaux derivative of mapping
\ShowEq{f:A->A}
has expansion
\ShowEq{Gateaux derivative of f, algebra}
Suppose the G\^ateaux derivative of mapping
\ShowEq{g:A->A}
has expansion
\ShowEq{Gateaux derivative of g, algebra}
The G\^ateaux derivative of mapping $f(x)g(x)$
have form
\ShowEq{Gateaux derivative, fg, algebra}
\end{theorem}
\begin{proof}
Let us substitute \EqRef{Gateaux derivative of f, algebra}
and \EqRef{Gateaux derivative of g, algebra}
into equation \EqRef{Gateaux differential of product, algebra}
\begin{align}
\EqLabel{Gateaux differential of fg, division ring}
\partial (f(x) g(x))(a)
&=\partial f(x)(a)\ g(x)+f(x)\ \partial g(x)(a)
\\ \nonumber
&=\frac{\pC{s}{0}\partial f(x)}{\partial x}a
\frac{\pC{s}{1}\partial f(x)}{\partial x}
g(x)
+f(x)\frac{\pC{t}{0}\partial g(x)}{\partial x}a
\frac{\pC{t}{1}\partial g(x)}{\partial x}
\end{align}
Based \EqRef{Gateaux differential of fg, division ring},
we define equations
\EqRef{0 component of Gateaux derivative, fg, algebra},
\EqRef{1 component of Gateaux derivative, fg, algebra}.
\end{proof}

\begin{theorem}
Let $A$ be Banach $D$\Hyph algebra.
If the G\^ateaux derivative
$\partial f(x)$
exists in point $x$ and has finite norm,
then function $f$ is continuous at point $x$.
\end{theorem}
\begin{proof}
From definition \ref{definition: norm of map, algebra}
it follows
\ShowEq{derivative and continuos function, algebra, 1}
From \EqRef{Gateaux derivative of map, algebra},
\EqRef{derivative and continuos function, algebra, 1}
it follows
\begin{equation}
\EqLabel{derivative and continuos function, algebra, 2}
|f(x+a)-f(x)|<|a|\ \|\partial f(x)\|
\end{equation}
Let us assume arbitrary $\epsilon>0$. Assume
\[\delta=\frac\epsilon{\|\partial f(x)\|}\]
Then from inequality
\[
|a|<\delta
\]
it follows
\[
|f(x+a)-f(x)|\le \|\partial f(x)\|\ \delta=\epsilon
\]
According to definition \ref{definition: continuous function, algebra}
map $f$ is continuous at point $x$.
\end{proof}

\begin{theorem}
\label{theorem: Gateaux derivative, 0, D algebra}
Let $A$ be Banach $D$\Hyph algebra.
Let mapping
\ShowEq{f:A->A}
be differentiable in the G\^ateaux sense at point
$x$.
Then
\ShowEq{Gateaux derivative, 0, D algebra}
\end{theorem}
\begin{proof}
Corollary of definitions
\ref{definition: function differentiable in Gateaux sense, algebra}
and theorem
\ref{theorem: linear map, 0, D algebra}.
\end{proof}

\begin{theorem}
\label{theorem: composite map, derivative, D algebra}
Let $A$ be Banach $D$\Hyph algebra.
Let mapping
\ShowEq{f:A->A}
be differentiable in the G\^ateaux sense at point
$x$
and norm of the G\^ateaux derivative of mapping $f$
be finite
\ShowEq{composite map, norm f, D algebra}
Let mapping
\ShowEq{g:A->A}
be differentiable in the G\^ateaux sense at point
\ShowEq{composite map, y fx, D algebra}
and norm of the G\^ateaux derivative of mapping $g$
be finite
\ShowEq{composite map, norm g, D algebra}
Mapping
\ShowEq{composite map, gfx, D algebra}
is differentiable in the G\^ateaux sense at point
$x$
\ShowEq{composite map, derivative, D algebra}
\ShowEq{composite map, Gateaux derivative 01, D algebra}
\end{theorem}
\begin{proof}
According to definition
\ref{definition: function differentiable in Gateaux sense, algebra}
\ShowEq{Gateaux derivative of map g, D algebra}
where
\ShowEq{o1:A->A}
is such continuous map that
\ShowEq{o1:A->A, 1}
According to definition
\ref{definition: function differentiable in Gateaux sense, algebra}
\ShowEq{Gateaux derivative of map f, D algebra}
where
\ShowEq{o2:A->A}
is such continuous map that
\ShowEq{o2:A->A, 1}
According to
\EqRef{Gateaux derivative of map f, D algebra}
increment $a$ of value $x\in A$
leads to increment
\ShowEq{composite map, b fxa, D algebra}
of value $y$.
Using \EqRef{composite map, y fx, D algebra},
\EqRef{composite map, b fxa, D algebra}
in equation \EqRef{Gateaux derivative of map g, D algebra},
we get
\ShowEq{Gateaux derivative of map gf, 1, D algebra}
According to definitions
\ref{definition: function differentiable in Gateaux sense, algebra},
\ref{definition: linear map from A1 to A2, algebra}
from equation \EqRef{Gateaux derivative of map gf, 1, D algebra}
it follows
\ShowEq{Gateaux derivative of map gf, 2, D algebra}
According to definition \ref{definition: norm on d algebra}
\ShowEq{Gateaux derivative of map gf, 3, D algebra}
From \EqRef{composite map, norm g, D algebra}
it follows that
\ShowEq{Gateaux derivative of map gf, 4, D algebra}
From \EqRef{composite map, norm f, D algebra}
it follows that
\ShowEq{Gateaux derivative of map gf, 6, D algebra}
According to the theorem \ref{theorem: Gateaux derivative, 0, D algebra}
\ShowEq{Gateaux derivative of map gf, 7, D algebra}
Therefore,
\ShowEq{Gateaux derivative of map gf, 5, D algebra}
From equations
\EqRef{Gateaux derivative of map gf, 3, D algebra},
\EqRef{Gateaux derivative of map gf, 4, D algebra},
\EqRef{Gateaux derivative of map gf, 5, D algebra}
it follows
\ShowEq{Gateaux derivative of map gf, 8, D algebra}
According to definition
\ref{definition: function differentiable in Gateaux sense, algebra}
\ShowEq{Gateaux derivative of map gf, 9, D algebra}
where
\ShowEq{Gateaux derivative of map, algebra, 1}
is such continuous mapping that
\ShowEq{Gateaux derivative of map, algebra, 2}
Equation \EqRef{composite map, derivative, D algebra}
follows from
\EqRef{Gateaux derivative of map gf, 2, D algebra},
\EqRef{Gateaux derivative of map gf, 8, D algebra},
\EqRef{Gateaux derivative of map gf, 9, D algebra}.

From equation \EqRef{composite map, derivative, D algebra}
and theorem
\ref{theorem: function differentiable in Gateaux sense, algebra}
it follows that
\ShowEq{composite map, derivative 2, D algebra}
\EqRef{composite map, Gateaux derivative 01, D algebra}
follow from equation
\EqRef{composite map, derivative 2, D algebra}.
\end{proof}
\fi

\Section{Table of Derivatives of Map of Associative Algebra}

\begin{theorem}
\label{theorem: Gateaux derivative of const, algebra}
Let $D$ be the complete commutative ring of characteristic $0$.
Let $A$ be associative $D$\Hyph algebra.
Then for any $b\in A$
\ShowEq{Gateaux derivative of const, algebra}
\ePrints{2010.05.23}
\ifx\Semafor\ValueOn
\qed
\end{theorem}
\newpage
\else
\end{theorem}
\begin{proof}
Immediate corollary of definition
\ref{definition: function differentiable in Gateaux sense, algebra}.
\end{proof}
\fi

\begin{theorem}
\label{theorem: Gateaux derivative, product over constant, algebra}
Let $D$ be the complete commutative ring of characteristic $0$.
Let $A$ be associative $D$\Hyph algebra.
Then for any $b$, $c\in A$
\ShowEq{Gateaux derivative, product over constant, algebra}
\ePrints{2010.05.23}
\ifx\Semafor\ValueOn
\qed
\end{theorem}
\newpage
\else
\end{theorem}
\begin{proof}
Immediate corollary of equations
\EqRef{Gateaux differential of product, algebra},
\EqRef{0 component of Gateaux derivative, fg, algebra},
\EqRef{1 component of Gateaux derivative, fg, algebra}
because $\partial b=\partial c=0$.
\end{proof}
\fi

\begin{theorem}
\label{theorem: Gateaux derivative, fx=axb, algebra}
Let $D$ be the complete commutative ring of characteristic $0$.
Let $A$ be associative $D$\Hyph algebra.
Then for any $b$, $c\in A$
\ShowEq{Gateaux derivative, fx=axb, algebra}
\ePrints{2010.05.23}
\ifx\Semafor\ValueOn
\qed
\end{theorem}
\newpage
\else
\end{theorem}
\begin{proof}
Corollary of theorem
\ref{theorem: Gateaux derivative, product over constant, algebra},
when $f(x)=x$.
\end{proof}
\fi

\ePrints{2010.05.23}
\ifx\Semafor\ValueOff
\begin{theorem}
\label{theorem: Gateaux derivative of linear map}
Let $D$ be the complete commutative ring of characteristic $0$.
Let $A$ be associative $D$\Hyph algebra.
Let $f$ be linear mapping
\ShowEq{derivative linear map associative algebra}
Then
\ShowEq{derivative linear map associative algebra, 1}
\end{theorem}
\begin{proof}
Corollary of theorems
\ref{theorem: Gateaux derivative, fx=axb, algebra},
\ref{theorem: derivative of the sum}.
\end{proof}
\fi

\begin{corollary}
Let $D$ be the complete commutative ring of characteristic $0$.
Let $A$ be associative $D$\Hyph algebra.
Then for any $b\in A$
\ShowEq{Gateaux derivative, xb-bx, algebra}
\qed
\end{corollary}
\ePrints{2010.05.23}
\ifx\Semafor\ValueOn
\newpage
\fi

\begin{theorem}
Let $D$ be the complete commutative ring of characteristic $0$.
Let $A$ be associative $D$\Hyph algebra.
Then\footnote{The statement of the theorem is similar to example VIII,
\citeBib{Hamilton Elements of Quaternions 1}, p. 451.
If product is commutative, then the equation
\EqRef{derivative x2, algebra}
gets form
\ShowEq{derivative x2, field}}
\ShowEq{derivative x2, algebra}
\ePrints{2010.05.23}
\ifx\Semafor\ValueOn
\qed
\end{theorem}
\newpage
\else
\end{theorem}
\begin{proof}
Consider increment of map $f(x)=x^2$.
\ShowEq{derivative x2, algebra, 1}
\EqRef{derivative x2, algebra}
follows from equations \EqRef{Gateaux differential, representation, algebra},
\EqRef{derivative x2, algebra, 1}.
\end{proof}
\fi

\begin{theorem}
Let $D$ be the complete commutative ring of characteristic $0$.
Let $A$ be associative division $D$\Hyph algebra.
Then\footnote{The statement of the theorem is similar to example IX,
\citeBib{Hamilton Elements of Quaternions 1}, p. 451.
If product is commutative, then the equation
\EqRef{derivative x power -1, algebra}
gets form
\ShowEq{derivative x power -1, field}}
\ShowEq{derivative x power -1, algebra}
\ePrints{2010.05.23}
\ifx\Semafor\ValueOn
\qed
\end{theorem}
\newpage
\else
\end{theorem}
\begin{proof}
Let us substitute $f(x)=x^{-1}$ in definition
\EqRef{Gateaux differential, representation in algebra}.
\ShowEq{derivative x power -1, algebra, 1}
Equation \EqRef{derivative x power -1, algebra}
follows from chain of equations
\EqRef{derivative x power -1, algebra, 1}.
\end{proof}
\fi

\begin{theorem}
Let $D$ be the complete commutative ring of characteristic $0$.
Let $A$ be associative division $D$\Hyph algebra.
Then\footnote{If product is commutative, then
\[
y=xax^{-1}=a
\]
Accordingly, the derivative is $0$.}
\ShowEq{derivative xax power -1, algebra}
\ePrints{2010.05.23}
\ifx\Semafor\ValueOn
\qed
\end{theorem}
\newpage
\else
\end{theorem}
\begin{proof}
Equation \EqRef{derivative xax power -1, algebra}
is corollary of equations
\EqRef{Gateaux differential of product, algebra},
\EqRef{Gateaux derivative, fx=axb, algebra}.
\end{proof}
\fi

%% file: Derivative.Eq.tex
%auto-ignore

\def\LAA{\mathcal L(A;A)}

\DefEq
{
$\partial f(x)\in\LAA$.
}
{differential L(A,A), 1}

\DefEq
{
\[
\begin{array}{r@{\ }l}
\partial f&=f
\\
\partial f\circ dx&=f\circ dx
\end{array}
\]
}
{derivative linear map associative algebra, 1}

\DefEq
{
\[
f\circ x=(a_{s\cdot 0}\otimes a_{s\cdot 1})\circ x
=a_{s\cdot 0}\ x\  a_{s\cdot 1}
\]
}
{derivative linear map associative algebra}

\DefEq
{
$o:A\rightarrow A$
}
{Gateaux derivative of map, algebra, 1}

\DefEq
{
\[
f:A\rightarrow A
\]
}
{f:A->A}

\DefEquation
{
\|\partial f(x)\|=F\le\infty
}
{composite map, norm f, D algebra}

\DefEquation
{
y=f(x)
}
{composite map, y fx, D algebra}

\DefEquation
{
\|\partial g(y)\|=G\le\infty
}
{composite map, norm g, D algebra}

\DefEquation
{
\left\{
\begin{array}{r@{\,}l}
\partial(g\circ f)(x)
&=\partial g(y)\circ\partial f(x)
\\
\partial(g\circ f)(x)\circ a
&=\partial g(y)\circ\partial f(x)\circ a
\end{array}
\right.
}
{composite map, derivative, D algebra}

\DefEquation
{
\left\{
\begin{array}{r@{\,}l}
\displaystyle
\frac{\partial\pC{st}{0} (g\circ f)(x)}{\partial x}&=
\displaystyle
\frac{\partial\pC{s}{0} g(f(x))}{\partial f(x)}
\frac{\partial\pC{t}{0} f(x)}{\partial x}
\\
\displaystyle
\frac{\partial\pC{st}{1} (g\circ f)(x)}{\partial x}&=
\displaystyle
\frac{\partial\pC{t}{1} f(x)}{\partial x}
\ \frac{\partial\pC{s}{1} g(f(x))}{\partial f(x)}
\end{array}
\right.
}
{composite map, Gateaux derivative 01, D algebra}

\DefEq
{
$o_1:A\rightarrow A$
}
{o1:A->A}

\DefEq
{
$o_2:A\rightarrow A$
}
{o2:A->A}

\DefEq
{
\[
\lim_{b\rightarrow 0}\frac{|o_1(b)|}{|b|}=0
\]
}
{o1:A->A, 1}

\DefEq
{
\[
\lim_{a\rightarrow 0}\frac{|o_2(a)|}{|a|}=0
\]
}
{o2:A->A, 1}

\DefEquation
{
b=\partial f(x)\circ a+o_2(a)
}
{composite map, b fxa, D algebra}

\DefEquation
{
\begin{array}{r@{\,}l}
&g(f(x+a))-g(f(x))
\\
=&g(f(x)+\partial f(x)\circ a
+o_2(a))
-g(f(x))
\\
=&\partial g(f(x))\circ(\partial f(x)\circ a
+o_2(a))
-o_1(\partial f(x)\circ a
+o_2(a))
\end{array}
}
{Gateaux derivative of map gf, 1, D algebra}

\DefEquation
{
\begin{array}{r@{\,}l}
&g(f(x+a))-g(f(x))
\\
=&\partial g(f(x))\circ\partial f(x)\circ a
+\partial g(f(x))\circ o_2(a)
-o_1(\partial f(x)\circ a
+o_2(a))
\end{array}
}
{Gateaux derivative of map gf, 2, D algebra}

\DefEquation
{
\begin{array}{r@{\,}l}
&\displaystyle
\lim_{a\rightarrow 0}
\frac
{|\partial g(f(x))\circ o_2(a)
-o_1(\partial f(x)\circ a
+o_2(a))|}
{|a|}
\\
\le
&\displaystyle
\lim_{a\rightarrow 0}
\frac
{|\partial g(f(x))\circ o_2(a)|}
{|a|}
+\lim_{a\rightarrow 0}
\frac
{|o_1(\partial f(x)\circ a
+o_2(a))|}
{|a|}
\end{array}
}
{Gateaux derivative of map gf, 3, D algebra}

\DefEquation
{
\lim_{a\rightarrow 0}
\frac
{|\partial g(f(x))\circ o_2(a)|}
{|a|}
\le
G\lim_{a\rightarrow 0}
\frac
{|o_2(a)|}
{|a|}=0
}
{Gateaux derivative of map gf, 4, D algebra}

\DefEq
{
\begin{align*}
&\lim_{a\rightarrow 0}
\frac
{|o_1(\partial f(x)\circ a
+o_2(a))|}
{|a|}
\\
=&
\lim_{a\rightarrow 0}
\frac
{|o_1(\partial f(x)\circ a
+o_2(a))|}
{|\partial f(x)\circ a
+o_2(a)|}
\lim_{a\rightarrow 0}
\frac
{|\partial f(x)\circ a
+o_2(a)|}
{|a|}
\\
\le&
\lim_{a\rightarrow 0}
\frac
{|o_1(\partial f(x)\circ a
+o_2(a))|}
{|\partial f(x)\circ a
+o_2(a)|_2}
\lim_{a\rightarrow 0}
\frac
{\|\partial f(x)\||a|
+|o_2(a)|}
{|a|}
\\
=&
\lim_{a\rightarrow 0}
\frac
{|o_1(\partial f(x)\circ a
+o_2(a))|}
{|\partial f(x)\circ a
+o_2(a)|}
\|\partial f(x)\|
\end{align*}
}
{Gateaux derivative of map gf, 6, D algebra}

\DefEq
{
\[
\lim_{a\rightarrow 0}
(\partial f(x)\circ a)
+o_2(a))=0
\]
}
{Gateaux derivative of map gf, 7, D algebra}

\DefEquation
{
\lim_{a\rightarrow 0}
\frac
{|\partial g(f(x))\circ o_2(a)
-o_1(\partial f(x)\circ a
+o_2(a))|}
{|a|}=0
}
{Gateaux derivative of map gf, 8, D algebra}

\DefEquation
{
(g\circ f)(x+a)-(g\circ f)(x)
=\partial (g\circ f)(x)\circ a
+o(a)
}
{Gateaux derivative of map gf, 9, D algebra}

\DefEquation
{
\begin{array}{r@{\,}l}
&\displaystyle
\frac{\partial\pC{st}{0} (g\circ f)(x)}{\partial x}
\ a
\ \frac{\partial\pC{st}{1} (g\circ f)(x)}{\partial x}
\\[10pt]
=&\displaystyle
\frac{\partial\pC{s}{0} g(f(x))}{\partial f(x)}
\ (\partial f(x)\circ a)
\ \frac{\partial\pC{s}{1} g(f(x))}{\partial f(x)}
\\[10pt]
=&\displaystyle
\frac{\partial\pC{s}{0} g(f(x))}{\partial f(x)}
\frac{\partial\pC{t}{0} f(x)}{\partial x}
\ a
\ \frac{\partial\pC{t}{1} f(x)}{\partial x}
\ \frac{\partial\pC{s}{1} g(f(x))}{\partial f(x)}
\end{array}
}
{composite map, derivative 2, D algebra}

\DefEq
{
\[
\lim_{a\rightarrow 0}\frac{|o(a)|}{|a|}=0
\]
}
{Gateaux derivative of map, algebra, 2}

\DefEquation
{
\lim_{a\rightarrow 0}
\frac
{|o_1(\partial f(x)\circ a
+o_2(a))|}
{|a|}=0
}
{Gateaux derivative of map gf, 5, D algebra}

\DefEq
{
\[
\partial f(x)\circ 0=0
\]
}
{Gateaux derivative, 0, D algebra}

\DefEquation
{
f(x+a)-f(x)
=\partial f(x)\circ a
+o_2(a)
}
{Gateaux derivative of map f, D algebra}

\DefEquation
{
g(y+b)-g(y)
=\partial g(y)\circ b
+o_1(b)
}
{Gateaux derivative of map g, D algebra}

\DefEq
{
\[
(g\circ f)(x)=g(f(x))
\]
}
{composite map, gfx, D algebra}

\DefEq
{
\[
g:A\rightarrow A
\]
}
{g:A->A}

\DefEquation
{
\partial f(x)
=\frac{\partial\pC{s}{0} f(x)}{\partial x}
\otimes
\frac{\partial\pC{s}{1} f(x)}{\partial x}
}
{Gateaux derivative of f, algebra}

\DefEquation
{
\partial g(x)
=\frac{\partial\pC{t}{0} g(x)}{\partial x}
\otimes
\frac{\partial\pC{t}{1} g(x)}{\partial x}
}
{Gateaux derivative of g, algebra}

\DefEq
{
\begin{equation}
\EqLabel{Gateaux derivative, fg, algebra}
\partial (f(x)g(x))
=\frac{\partial\pC{s}{0} f(x)}{\partial x}
\otimes
\left(
\frac{\partial\pC{s}{1} f(x)}{\partial x}g(x)
\right)
+\left(
f(x)
\frac{\partial\pC{t}{0} g(x)}{\partial x}
\right)
\otimes
\frac{\partial\pC{t}{1} g(x)}{\partial x}
\end{equation}
\begin{align}
\EqLabel{0 component of Gateaux derivative, fg, algebra}
\frac{\partial\pC{s}{0} f(x)g(x)}{\partial x}
&
=\frac{\partial\pC{s}{0} f(x)}{\partial x}
&
\frac{\partial\pC{t}{0} f(x)g(x)}{\partial x}
&
=f(x)\frac{\partial\pC{t}{0} g(x)}{\partial x}
\\
\EqLabel{1 component of Gateaux derivative, fg, algebra}
\frac{\partial\pC{s}{1} f(x)g(x)}{\partial x}
&
=\frac{\partial\pC{s}{1} f(x)}{\partial x}
g(x)
&
\frac{\partial\pC{t}{1} f(x)g(x)}{\partial x}
&
=\frac{\partial\pC{t}{1} g(x)}{\partial x}
\end{align}
}
{Gateaux derivative, fg, algebra}

\DefEquation
{
\left\{
\begin{array}{r@{\ }l}
\partial (bf(x)c)
&=b\ \partial f(x)\ c
\\
\partial (bf(x)c)\circ dx
&=b(\partial f(x)\circ dx)c
\\
\displaystyle\VirtFrac
\frac{\partial\pC{s}{0} bf(x)c}{\partial x}
&\displaystyle
=b\frac{\partial\pC{s}{0} f(x)}{\partial x}
\\
\displaystyle\VirtFrac
\frac{\partial\pC{s}{1} bf(x)c}{\partial x}
&\displaystyle
=\frac{\partial\pC{s}{1} f(x)}{\partial x}c
\end{array}
\right.
}
{Gateaux derivative, product over constant, algebra}

\DefEquation
{
\left\{
\begin{array}{r@{\ }lr@{\ }l}
\partial (bxc)&=b\otimes c
&
\partial (bxc)\circ dx&=b\ dx\ c
\\
\displaystyle
\frac{\partial\pC{1}{0} bxc}{\partial x}&=b
&
\displaystyle
\frac{\partial\pC{1}{1} bxc}{\partial x}&=c
\end{array}
\right.
}
{Gateaux derivative, fx=axb, algebra}

\DefEquation
{
\left\{
\begin{array}{cc}
\multicolumn{2}{c}
{
\begin{array}{r@{\ }l}
\partial x^2&=x\otimes 1+1\otimes x
\\
\partial x^2\circ dx&=x\ dx+dx\ x
\end{array}
}
\\
\displaystyle\VirtFrac
\frac{\partial\pC{1}{0} x^2}{\partial x}=x
&
\displaystyle\frac{\partial\pC{1}{1} x^2}{\partial x}=1
\\
\displaystyle\VirtFrac
\frac{\partial\pC{2}{0} x^2}{\partial x}=1
&
\displaystyle\frac{\partial\pC{2}{1} x^2}{\partial x}=x
\end{array}
\right.
}
{derivative x2, algebra}

\DefEquation
{
\left\{
\begin{array}{cc}
\multicolumn{2}{c}
{
\begin{array}{r@{\ }l}
\partial x^{-1}&=-x^{-1}\otimes x^{-1}
\\
\partial x^{-1}\circ dx&=-x^{-1}\ dx\ x^{-1}
\end{array}
}
\\
\displaystyle\frac{\partial\pC{1}{0} x^{-1}}{\partial x}=-x^{-1}
&
\displaystyle\frac{\partial\pC{1}{1} x^{-1}}{\partial x}=x^{-1}
\end{array}
\right.
}
{derivative x power -1, algebra}

\DefEq
{
\[
\begin{matrix}
f:A\rightarrow A&g:A\rightarrow A
\end{matrix}
\]
}
{bilinear map and derivative, 1}

\DefEquation
{
\partial h(f(x),g(x))\circ dx
=h(\partial f(x)\circ dx,g(x))
+h(f(x),\partial g(x)\circ dx)
}
{bilinear map and derivative, 3}

\DefEquation
{
\partial (f+g)(x)
=\partial f(x)
+\partial g(x)
}
{derivative of the sum, 3}

\DefEq
{
\begin{align*}
\partial h(f(x),g(x))\circ a
&=\lim_{t\rightarrow 0}(t^{-1}(h(f(x+ta),g(x+ta))-h(f(x),g(x))))
\\
&=\lim_{t\rightarrow 0}(t^{-1}(h(f(x+ta),g(x+ta))-h(f(x),g(x+ta))))
\\
&+\lim_{t\rightarrow 0}(t^{-1}(h(f(x),g(x+ta))-h(f(x),g(x))))
\\
&=h(\lim_{t\rightarrow 0}t^{-1}(f(x+ta)-f(x)),g(x))
\\
&+h(f(x),\lim_{t\rightarrow 0}t^{-1}(g(x+ta)-g(x)))
\end{align*}
}
{bilinear map and derivative, 4}

\DefEq
{
\[
h(f,g):A\rightarrow A
\]
}
{bilinear map and derivative, 2}

\DefEq
{
\[
f+g:A\rightarrow A
\]
}
{derivative of the sum, 2}

\DefEq
{
\[
h:A\times A\rightarrow A
\]
}
{bilinear map and derivative}

\DefEquation
{
\left\{
\begin{array}{ll}
\multicolumn{2}{c}
{
\begin{array}{r@{\ }l}
\partial (xax^{-1})=1\otimes ax^{-1}-xax^{-1}\otimes x^{-1}
\\
\partial (xax^{-1})\circ dx=dx\ ax^{-1}-xax^{-1}\ dx\ x^{-1}
\end{array}
}
\\
\displaystyle\frac{\partial\pC{1}{0} x^{-1}}{\partial x}=1
&
\displaystyle\frac{\partial\pC{1}{1} x^{-1}}{\partial x}=ax^{-1}
\\
\displaystyle\frac{\partial\pC{2}{0} x^{-1}}{\partial x}=-xax^{-1}
&
\displaystyle\frac{\partial\pC{2}{1} x^{-1}}{\partial x}=x^{-1}
\end{array}
\right.
}
{derivative xax power -1, algebra}

\DefEq
{
\begin{align}
\partial f(x)\circ h
&=\lim_{t\rightarrow 0,\ t\in R}(t^{-1}((x+th)^{-1}-x^{-1}))
\nonumber
\\
&=\lim_{t\rightarrow 0,\ t\in R}(t^{-1}((x+th)^{-1}-x^{-1}(x+th)(x+th)^{-1}))
\nonumber
\\
&=\lim_{t\rightarrow 0,\ t\in R}(t^{-1}(1-x^{-1}(x+th))(x+th)^{-1})
\EqLabel{derivative x power -1, algebra, 1}
\\
&=\lim_{t\rightarrow 0,\ t\in R}(t^{-1}(1-1-x^{-1}th)(x+th)^{-1})
\nonumber
\\
&=\lim_{t\rightarrow 0,\ t\in R}(-x^{-1}h(x+th)^{-1})
\nonumber
\end{align}
}
{derivative x power -1, algebra, 1}

\DefEq
{
\[
\left\{
\begin{array}{r@{\ }lr@{\ }l}
\multicolumn{4}{c}
{
\begin{array}{r@{\ }l}
\partial (xb-bx)
&=1\otimes b-b\otimes 1
\\
\partial (xb-bx)\circ dx
&=dx\ b-b\ dx
\end{array}
}
\\
\displaystyle\VirtFrac
\frac{\partial\pC{1}{0} (xb-bx)}{\partial x}
&=1
&\displaystyle
\frac{\partial\pC{1}{1} (xb-bx)}{\partial x}
&=b
\\
\displaystyle\VirtFrac
\frac{\partial\pC{2}{0} (xb-bx)}{\partial x}
&=-b
&\displaystyle
\frac{\partial\pC{2}{1} (xb-bx)}{\partial x}
&=1
\end{array}
\right.
\]
}
{Gateaux derivative, xb-bx, algebra}

\DefEq
{
\[
\partial b=0
\]
}
{Gateaux derivative of const, algebra}

\DefEq
{
\begin{align*}
(a,v)\in D\times A&\rightarrow av\in A
\\
(v,w)\in A\times A&\rightarrow vw\in A
\end{align*}
}
{topological D algebra}

\DefEq
{
\item $|a+b|\le |a|+|b|$
\item $|ab|=|a|\ |b|$
\item $|da|=|d|\ |a|$, $d\in D$, $a\in A$
}
{norm on d algebra 3}

\DefEq
{
$a\in A$, $|a|=1$,
}
{unit sphere in algebra}

\DefEq
{
$a=0$
}
{norm on d algebra 2, 2}

\DefEq
{
\item $|a|=0$
}
{norm on d algebra 2, 1}

\DefEq
{
\item $|a|\ge 0$
}
{norm on d algebra 1}

\DefEq
{
\[a\in A\rightarrow |a|\in R\]
}
{norm on d algebra}

\DefEq
{
$\{a_n\}$
\symb{\lim_{n\rightarrow\infty}a_n}0{limit of sequence, normed ring}
\[
a=\ShowSymbol{limit of sequence, normed ring}
\]
}
{limit of sequence, normed ring}

\DefEq
{
$\{a_n\}$
\symb{\lim_{n\rightarrow\infty}a_n}0{limit of sequence, normed algebra}
\[
a=\ShowSymbol{limit of sequence, normed algebra}
\]
}
{limit of sequence, normed algebra}

\DefEq
{
\symb{\partial f(x)\circ dx}
0{Gateaux differential of map, algebra}
$\ShowSymbol{Gateaux differential of map, algebra}$
}
{Gateaux differential of map, algebra}

\DefEq
{
\begin{equation}
\EqLabel{Gateaux differential is multiplicative over field R, algebra}
\partial f(x)\circ(ra)
=
r\partial f(x)\circ a
\end{equation}
\[
\begin{matrix}
r\in R& r\ne 0&a\in A&a\ne 0
\end{matrix}
\]
}
{Gateaux differential is multiplicative over field R, algebra}

\DefEquation
{
f\circ(x+ta)-f\circ x
=t\partial f(x)\circ a
+o(t)
}
{Gateaux differential of map, t, algebra}

\DefEquation
{
\begin{array}{r@{\ }l}
\partial f(x)\circ a
&=\displaystyle
\lim_{t\rightarrow 0,\ t\in R}(t^{-1}((f+g)\circ(x+ta)-(f+g)\circ x))
\\
&=\displaystyle
\lim_{t\rightarrow 0,\ t\in R}(t^{-1}(f\circ(x+ta)+g\circ(x+ta)-f\circ x-g\circ x))
\\
&=\displaystyle
\lim_{t\rightarrow 0,\ t\in R}(t^{-1}(f\circ(x+ta)-f\circ x))
\\
&+\displaystyle
\lim_{t\rightarrow 0,\ t\in R}(t^{-1}(g\circ(x+ta)-g\circ x))
\\
&=\displaystyle
\partial f(x)\circ a
+\partial g(x)\circ a
\end{array}
}
{derivative of the sum, 4}

\DefEquation
{
\partial f(x)\circ a=
\lim_{t\rightarrow 0,\ t\in R}(t^{-1}(f\circ(x+ta)-f\circ x))
}
{Gateaux differential, representation in algebra}

\DefEquation
{
\ShowSymbol{Gateaux differential of map, algebra}
=
\left(
\frac{\partial_{s\cdot 0} f(x)}{\partial x}
\otimes
\frac{\partial_{s\cdot 1} f(x)}{\partial x}
\right)
\circ dx
=\frac{\partial_{s\cdot 0} f(x)}{\partial x}
dx
\frac{\partial_{s\cdot 1} f(x)}{\partial x}
}
{Gateaux differential, representation, algebra}

\DefEquation
{
\ShowSymbol{Gateaux differential of map, algebra}
=
\left(
\frac{\partial_{k\cdot s\cdot 0} f(x)}{\partial x}
\otimes
\frac{\partial_{k\cdot s\cdot 1} f(x)}{\partial x}
\right)
\circ I_k\circ dx
=\frac{\partial_{k\cdot s\cdot 0} f(x)}{\partial x}
(I_k\circ dx)
\frac{\partial_{k\cdot s\cdot 1} f(x)}{\partial x}
}
{Gateaux differential, representation, algebra, 1}

\DefEquation
{
\partial (f(x)g(x))\circ dx
=(\partial f(x)\circ dx)\ g(x)+f(x)\ (\partial g(x)\circ dx)
}
{Gateaux differential of product, algebra}

\DefEquation
{
\partial f(x)\circ dx=dx^{\gi i}\frac{\partial f^{\gi j}}{\partial x^{\gi i}}
\Vector e_{\gi j}
}
{Gateaux differential and jacobian, algebra}

\DefEquation
{
f(x+h)-f(x)
=(x+h)^2-x^2
=xh+hx+h^2
=xh+hx+o(h)
}
{derivative x2, algebra, 1}

\DefEquation
{
|\partial f(x)\circ a|
\le\|\partial f(x)\|\ |a|
}
{derivative and continuos function, algebra, 1}

\DefEquation
{
\frac{\partial f^{\gi j}}{\partial x^{\gi i}}=
\StandPartial{f(x)}{x}{kr} C_{\gi{ki}}^{\gi p}C_{\gi{pr}}^{\gi j}
}
{standard components and Jacobian, algebra}

\DefEq
{
\[
\begin{matrix}
dx=dx^{\gi i}\ \Vector e_{\gi i}&&dx^{\gi i}\in D
\end{matrix}
\]
}
{Gateaux differential and jacobian, 1, algebra}

\DefEq
{
$\displaystyle
\ShowSymbol{standard component of Gateaux derivative, algebra}$
}
{standard component of Gateaux derivative, algebra}

\DefEq
{
\symb{\StandPartial{f(x)}{x}{ij}}0
{standard component of Gateaux derivative, algebra}
\begin{equation}
\partial f(x)
=\ShowSymbol{standard component of Gateaux derivative, algebra}
\Vector e_{\gi i}\otimes \Vector e_{\gi j}
\EqLabel{Gateaux derivative, algebra, standard representation}
\end{equation}
}
{Gateaux derivative, algebra, standard representation}

\DefEq
{
\symb{\frac{\partial\pC{s}{p} f(x)}{\partial x}}0
{component of Gateaux derivative of map, algebra}
$\displaystyle\ShowSymbol{component of Gateaux derivative of map, algebra}$,
$p=0$, $1$,
}
{component of Gateaux derivative of map, algebra}

\DefEq
{
\[
\partial f:A\rightarrow \LAA
\]
}
{differential L(A,A), 2}

\DefEquation
{
\frac{|f(x)|_2}{|x|_1}=\left|f\left(\frac x{|x|_1}\right)\right|_2
}
{norm of linear map, algebra, 1}

\DefEquation
{
\begin{matrix}
f(rx)=rf(x)
&
r\in R
\end{matrix}
}
{linear map and real number}

\DefEq
{
\[
\frac{|f(rx)|_2}{|rx|_1}
=\frac{|r|\ |f(x)|_2}{|r|\ |x|_1}
=\frac{|f(x)|_2}{|x|_1}
\]
}
{linear map and real number, 1}

\DefEquation
{
\|f\|=\text{sup}\{|f(x)|_2:|x|_1=1\}
}
{norm of linear map, algebra}

\DefEq
{
\symb{\|f\|}0{norm of map, algebra}
\begin{equation}
\ShowSymbol{norm of map, algebra}=
\text{sup}\frac{|f(x)|_2}{|x|_1}
\EqLabel{norm of map, algebra}
\end{equation}
}
{norm of map, algebra}

\DefEq
{
\[f:A_1\rightarrow A_2\]
}
{f A1 A2}

\DefEq
{
\symb{\partial f(x)}
0{Gateaux derivative of map, algebra}
\symb{\frac{\partial f(x)}{\partial x}}
0{Gateaux derivative of map, fraction, algebra}
\begin{equation}
\EqLabel{Gateaux derivative of map, algebra}
f(x+a)-f(x)
=\ShowSymbol{Gateaux derivative of map, algebra}\circ a
+o(a)
=\ShowSymbol{Gateaux derivative of map, fraction, algebra}
\circ a
+o(a)
\end{equation}
}
{Gateaux derivative of map, algebra}

\DefEq
{
\begin{align*}
\partial x^{-1}\circ dx&=-x^{-2}dx
\\
\frac{dx^{-1}}{dx}&=-x^{-2}
\end{align*}
}
{derivative x power -1, field}

\DefEq
{
\begin{align*}
\partial x^2\circ dx&=2x\ dx
\\
\frac{dx^2}{dx}&=2x
\end{align*}
}
{derivative x2, field}

%% file: Second.Derivative.English.tex
%auto-ignore

\input{\FilePrefix Second.Derivative.Eq}

\Chapter{Derivative of Second Order of Map of Division Ring}
\label{chapter: Derivative of Second Order of Map of Division Ring}

\ePrints{2010.05.23}
\ifx\Semafor\ValueOn
\newpage
{\bf Derivative of Second Order of Map of Algebra}
\fi
\Section{Derivative of Second Order of Map of Algebra}

Let $D$ be the complete commutative ring of characteristic $0$.
Let $A$ be associative $D$\Hyph algebra.
Let
\ShowEq{f:A->A}
function differentiable in the G\^ateaux sense.
\ePrints{2010.05.23}
\ifx\Semafor\ValueOff
According to remark \ref{remark: differential L(A,A)}
the G\^ateaux derivative is map
\ShowEq{differential L(A,A), 2}
According to theorem
\ShowEq{theorem: module L(A;A) is algebra}
and definition \ref{definition: norm of map, algebra}
set $\LAA$ is
Banach $D$\Hyph algebra.
Therefore, we may consider the question,
if map $\partial f$ is differentiable in the G\^ateaux sense.

According to definition
\ref{definition: function differentiable in Gateaux sense, algebra}
\ShowEq{Gateaux differential of map df, algebra}
where
\ShowEq{Gateaux differential of map df, algebra, o2}
is such continuous map, that
\ShowEq{Gateaux differential of map df, algebra, o2 lim}
According to definition
\ref{definition: function differentiable in Gateaux sense, algebra}
the mapping
\ShowEq{increment of derivative, algebra}
is linear map of variable $a_2$. From equation
\EqRef{Gateaux differential of map df, algebra}
it follows that mapping
\ShowEq{increment of derivative, algebra}
is linear mapping of variable $a_1$.
Therefore, the mapping
\ShowEq{increment of derivative, algebra}
is bilinear mapping.
\fi

\begin{definition}
\label{definition: Gateaux derivative of Second Order, algebra}
Polylinear map
\ShowEq{Gateaux derivative of Second Order, algebra}
is called
\AddIndex{the G\^ateaux derivative of second order of map $f$}
{Gateaux derivative of Second Order, algebra}.
\qed
\end{definition}

\ePrints{2010.05.23}
\ifx\Semafor\ValueOn
\newpage
\fi
\begin{remark}
According to definition
\ref{definition: Gateaux derivative of Second Order, algebra}
for given $x$ the G\^ateaux derivative of second order
\ShowEq{differential L(A,A;A), 1}.
Therefore, the G\^ateaux derivative of second order of map $f$ is mapping
\ShowEq{differential L(A,A;A)}
According to the theorem
\ShowEq{theorem: tensor product and polylinear mapping}
we may consider also expression
\ShowEq{differential AA A,A}
Then
\ShowEq{differential L(AA;A)}
We use the same notation for mapping
because of the nature of the argument it is clear what kind of mapping we consider.
\qed
\end{remark}

\ePrints{2010.05.23}
\ifx\Semafor\ValueOn
\newpage
\fi
\begin{theorem}
It is possible to represent
\AddIndex{the G\^ateaux differential of second order of map $f$}
{Gateaux differential of Second Order, algebra}
as
\ShowEq{Differential of Second Order, algebra, representation}
Expression
\ePrints{2010.05.23}
\ifx\Semafor\ValueOff
\footnote{We suppose
\ShowEq{component of Gateaux derivative of Second Order, 1, algebra}
}
\fi
\ShowEq{component of Gateaux derivative of Second Order, algebra}
is called
\AddIndex{component of the G\^ateaux derivative of second order of map $f(x)$}
{component of Gateaux derivative of Second Order, algebra}.
\ePrints{2010.05.23,1006.2597}
\ifx\Semafor\ValueOn
\qed
\end{theorem}
\else
\end{theorem}
\begin{proof}
Corollary of definition
\ref{definition: Gateaux derivative of Second Order, algebra}
and theorem
\ref{theorem: polylinear map, algebra}.
\end{proof}
\fi

\ePrints{2010.05.23}
\ifx\Semafor\ValueOn
\newpage
\fi
By induction, assuming that we defined the G\^ateaux derivative
$\partial^{n-1} f(x)$ of order $n-1$, we define
\ShowEq{Gateaux derivative of Order n, algebra}
\AddIndex{the G\^ateaux derivative of order $n$ of map $f$}
{Gateaux derivative of Order n, algebra}.
We also assume $\partial^0 f(x)=f(x)$.

\ePrints{2010.05.23}
\ifx\Semafor\ValueOn
\newpage
{\bf Taylor Series}
\fi
\Section{Taylor Series}
\label{section: Taylor Series}

Let $D$ be the complete commutative ring of characteristic $0$.
Let $A$ be associative $D$\Hyph algebra.
Let $p_k(x)$ be the monomial of power $k$, $k>0$,
in one variable over $D$\Hyph algebra $A$.

It is evident that monomial of power $0$ has form $a_0$, $a_0\in A$.
For $k>0$,
\[
p_k(x)=p_{k-1}(x)xa_k
\]
where $a_k\in A$.
\ePrints{2010.05.23}
\ifx\Semafor\ValueOff
Actually, last factor of monomial $p_k(x)$ is either $a_k\in A$,
or has form $x^l$, $l\ge 1$.
In the later case we assume $a_k=1$.
Factor preceding $a_k$ has form $x^l$, $l\ge 1$.
We can represent this factor as $x^{l-1}x$.
Therefore, we proved the statement.

\fi
In particular, monomial of power $1$ has form $p_1(x)=a_0xa_1$.

\ePrints{2010.05.23}
\ifx\Semafor\ValueOff
Without loss of generality, we assume $k=n$.

\begin{theorem}
\label{theorem: Gateaux derivative of f(x)x, algebra}
For any $m>0$ the following equation is true
\ShowEq{Gateaux derivative of f(x)x, algebra}
where symbol $\widehat{h^i}$ means absense of variable $h^i$ in the list.
\end{theorem}
\begin{proof}
For $m=1$, this is corollary  of equation
\EqRef{Gateaux differential of product, algebra}
\ShowEq{Gateaux derivative of f(x)x, 1, algebra}
Assume, \EqRef{Gateaux derivative of f(x)x, algebra} is true
for $m-1$. Then
\ShowEq{Gateaux derivative of f(x)x, m-1, algebra}
Using equations \EqRef{Gateaux differential of product, algebra}
and \EqRef{Gateaux derivative, product over constant, algebra}
we get
\ShowEq{Gateaux derivative of f(x)x, m, algebra}
The difference between
equations \EqRef{Gateaux derivative of f(x)x, algebra} and
\EqRef{Gateaux derivative of f(x)x, m, algebra}
is only in form of presentation.
We proved the theorem.
\end{proof}
\fi

\begin{theorem}
\label{theorem: Gateaux derivative of pn is symmetric, m < n, algebra}
The G\^ateaux derivative
\ShowEq{Gateaux derivative of pn is symmetric}
is symmetric polynomial with respect to variables $h_1$, ..., $h_m$.
\ePrints{2010.05.23}
\ifx\Semafor\ValueOn
\qed
\end{theorem}
\newpage
\else
\end{theorem}
\begin{proof}
To prove the theorem we consider algebraic properties
of the G\^ateaux derivative and give equivalent definition.
We start from construction of monomial.
For any monomial $p_n(x)$ we build symmetric
polynomial $r_n(x)$ according to following rules
\begin{itemize}
\item If $p_1(x)=a_0xa_1$, then $r_1(x_1)=a_0x_1a_1$
\item If $p_n(x)=p_{n-1}(x)a_n$, then
\[
r_n(x_1,...,x_n)=r_{n-1}(x_{[1},...,x_{n-1})x_{n]}a_n
\]
where square brackets express symmetrization of expression
with respect to variables $x_1$, ..., $x_n$.
\end{itemize}
It is evident that
\[
p_n(x)=r_n(x_1,...,x_n)\ \ \ x_1=...=x_n=x
\]
We define the G\^ateaux derivative of power $k$ according to rule
\ShowEq{Gateaux derivative of monomial, algebraic definition, algebra}
According to construction, polynomial $r_n(h_1,...,h_k,x_{k+1},...,x_n)$
is symmetric with respect to variables $h_1$, ..., $h_k$, $x_{k+1}$, ..., $x_n$.
Therefore, polynomial
\EqRef{Gateaux derivative of monomial, algebraic definition, algebra}
is symmetric with respect to variables $h_1$, ..., $h_k$.

For $k=1$, we will prove that definition
\EqRef{Gateaux derivative of monomial, algebraic definition, algebra}
of the G\^ateaux derivative coincides with definition
\EqRef{Gateaux derivative of map, algebra}.

For $n=1$, $r_1(h_1)=a_0h_1a_1$. This expression coincides
with expression of the G\^ateaux derivative in theorem
\ref{theorem: Gateaux derivative, fx=axb, algebra}.

Let the statement be true for $n-1$.
The following equation is true
\ShowEq{Gateaux derivative of monomial, n, algebra}
Assume $x_2=...=x_n=x$. 
According to suggestion of induction, from equations
\EqRef{Gateaux derivative of monomial, algebraic definition, algebra},
\EqRef{Gateaux derivative of monomial, n, algebra}
it follows that
\ShowEq{Gateaux derivative of monomial, n1, algebra}

According to theorem
\ref{theorem: Gateaux derivative of f(x)x, algebra}
\ShowEq{Gateaux derivative of pn is symmetric, 1}
This proves the equation
\EqRef{Gateaux derivative of monomial, algebraic definition, algebra}
for $k=1$.

Let us prove now that definition
\EqRef{Gateaux derivative of monomial, algebraic definition, algebra}
of the G\^ateaux derivative coincides with definition
\EqRef{Gateaux derivative of Order n, algebra} for $k>1$.

Let equation
\EqRef{Gateaux derivative of monomial, algebraic definition, algebra}
be true for $k-1$.
Let us consider arbitrary monomial of polynomial $r_n(h_1,...,h_{k-1},x_k,...,x_n)$.
Identifying variables $h_1$, ..., $h_{k-1}$
with elements of division ring $D$, we consider polynomial
\begin{equation}
R_{n-k}(x_k,...,x_n)=r_n(h_1,...,h_{k-1},x_k,...,x_n)
\EqLabel{reduced polynom, algebra}
\end{equation}
Assume $P_{n-k}(x)=R_{n-k}(x_k,...,x_n)$, $x_k=...=x_n=x$.
Therefore
\ShowEq{Gateaux derivative of pn is symmetric, 2}
According to definition of the G\^ateaux derivative
\EqRef{Gateaux derivative of Order n, algebra}
\ShowEq{Gateaux derivative of Order n, 1, algebra}
According to definition
\EqRef{Gateaux derivative of monomial, algebraic definition, algebra}
of the G\^ateaux derivative
\ShowEq{Gateaux derivative of monomial, algebraic definition, 1, algebra}
According to definition \EqRef{reduced polynom, algebra},
from equation
\EqRef{Gateaux derivative of monomial, algebraic definition, 1, algebra}
it follows that
\begin{equation}
\partial P_{n-k}(x)(h_k)=r_n(h_1,...,h_k,x_{k+1},...,x_n)\ \ \ x_{k+1}=x_n=x
\EqLabel{Gateaux derivative of monomial, algebraic definition, 2, algebra}
\end{equation}
From comparison of equations \EqRef{Gateaux derivative of Order n, 1, algebra}
and \EqRef{Gateaux derivative of monomial, algebraic definition, 2, algebra}
it follows that
\[
\partial^k p_n(x)(h_1;...;h_k)=r_n(h_1,...,h_k,x_{k+1},...,x_n)\ \ \ x_{k+1}=x_n=x
\]
Therefore equation
\EqRef{Gateaux derivative of monomial, algebraic definition, algebra}
is true for any $k$ and $n$.

We proved the statement of theorem.
\end{proof}

\begin{theorem}
\label{theorem: Gateaux derivative of pn = 0, algebra}
For any $n\ge 0$ following equation is true
\ShowEq{Gateaux derivative n1 of pn, algebra}
\end{theorem}
\begin{proof}
Since $p_0(x)=a_0$, $a_0\in D$, then for $n=0$ theorem is corollary
of theorem \ref{theorem: Gateaux derivative of const, algebra}.
Let statement of theorem is true for $n-1$. According to theorem
\ref{theorem: Gateaux derivative of f(x)x, algebra}
when $f(x)=p_{n-1}(x)$ we get
\begin{align*}
\partial^{n+1} p_n(x)(h_1;...;h_{n+1})
=&\partial^{n+1}(p_{n-1}(x)xa_n)(h_1;...;h_{n+1})
\\
=&\partial^{n+1} p_{n-1}(x)(h_1;...;h_m)xa_n
\\
+&\partial^n p_{n-1}(x)(h_1;...;h_{m-1})h_ma_n
\\
+&\partial^n p_{n-1}(x)(\widehat{h_1};...;h_{m-1};h_m)h_1a_n
+...
\\
+&\partial^{m-1} p_{n-1}(x)(h_1;...;\widehat{h_{m-1}};h_m)h_{m-1}a_n
\end{align*}
According to suggestion of induction all monomials are equal $0$.
\end{proof}

\begin{theorem}
\label{theorem: Gateaux derivative of pn = 0, m < n, algebra}
If $m<n$, then following equation is true
\ShowEq{m Gateaux derivative of polinom pn, algebra}
\end{theorem}
\begin{proof}
For $n=1$ following equation is true
\[
\partial^0 p_1(0)=a_0xa_1=0
\]
Assume that statement is true
for $n-1$. Then according to theorem \ref{theorem: Gateaux derivative of f(x)x, algebra}
\begin{align*}
&\partial^m (p_{n-1}(x)xa_n)(h_1;...;h_m)
\\
=&\partial^m p_{n-1}(x)(h_1;...;h_m)xa_n
+\partial^{m-1} p_{n-1}(x)(h_1;...;h_{m-1})h_ma_n
\\
+&\partial^{m-1} p_{n-1}(x)(\widehat{h_1};...;h_{m-1};h_m)h_1a_n
+...
\\
+&\partial^{m-1} p_{n-1}(x)(h_1;...;\widehat{h_{m-1}};h_m)h_{m-1}a_n
\end{align*}
First term equal $0$ because $x=0$.
Because $m-1<n-1$, then
rest terms equal $0$ according to suggestion of induction.
We proved the statement of theorem.
\end{proof}
\fi

When $h_1=...=h_n=h$, we assume
\ShowEq{dfxh}
This notation does not create ambiguity, because we can determine function
according to number of arguments.

\begin{theorem}
\label{theorem: n Gateaux derivative of polinom pn, algebra}
For any $n>0$ following equation is true
\ShowEq{n Gateaux derivative of polinom pn, algebra}
\ePrints{2010.05.23}
\ifx\Semafor\ValueOn
\qed
\end{theorem}
\newpage
\else
\end{theorem}
\begin{proof}
For $n=1$ following equation is true
\ShowEq{n Gateaux derivative of polinom pn, 2, algebra}
Assume the statement is true
for $n-1$. Then according to theorem \ref{theorem: Gateaux derivative of f(x)x, algebra}
\ShowEq{n Gateaux derivative of polinom pn, 1, algebra}
First term equal $0$ according to theorem
\ref{theorem: Gateaux derivative of pn = 0, algebra}.
The rest $n$ terms equal, and according to suggestion of induction
from equation \EqRef{n Gateaux derivative of polinom pn, 1, algebra}
it follows
\ShowEq{n Gateaux derivative of polinom pn, 3, algebra}
Therefore, statement of theorem is true for any $n$.
\end{proof}

Let $p(x)$ be polynomial of power $n$.\footnote{I consider
Taylor polynomial for polynomials by analogy with
construction of Taylor polynomial in \citeBib{Fihtengolts: Calculus volume 1}, p. 246.}
\[
p(x)=p_0+p_{1i_1}(x)+...+p_{ni_n}(x)
\]
We assume sum by index $i_k$ which enumerates terms
of power $k$.
According to theorem \ref{theorem: Gateaux derivative of pn = 0, algebra},
\ref{theorem: Gateaux derivative of pn = 0, m < n, algebra},
\ref{theorem: n Gateaux derivative of polinom pn, algebra}
\ShowEq{Taylor polynomial 1, algebra}
Therefore, we can write
\ShowEq{Taylor polynomial 2, algebra}
This representation of polynomial is called
\AddIndex{Taylor polynomial}{Taylor polynomial, algebra}.
If we consider substitution of variable $x=y-y_0$, then considered above construction
remain true for polynomial
\[
p(y)=p_0+p_{1i_1}(y-y_0)+...+p_{ni_n}(y-y_0)
\]
Therefore
\ShowEq{Taylor polynomial 3, algebra}
\fi

Assume that function $f(x)$ is differentiable in the G\^ateaux sense at point $x_0$
up to any order.\footnote{I explore
construction of Taylor series by analogy with
construction of Taylor series in \citeBib{Fihtengolts: Calculus volume 1}, p. 248, 249.}

\begin{theorem}
\label{theorem: n Gateaux derivative equal 0, algebra}
If function $f(x)$ holds
\ShowEq{n Gateaux derivatives of function, algebra}
then for $t\rightarrow 0$ expression $f(x+th)$ is infinitesimal of order
higher than $n$ with respect to $t$
\[
f(x_0+th)=o(t^n)
\]
\ePrints{2010.05.23}
\ifx\Semafor\ValueOn
\qed
\end{theorem}
\newpage
\else
\end{theorem}
\begin{proof}
When $n=1$ this statement follows from equation
\EqRef{Gateaux differential of map, t, algebra}.

Let statement be true for $n-1$.
Map
\ShowEq{n Gateaux derivative equal 0, 1, algebra}
satisfies to condition
\ShowEq{n Gateaux derivative equal 0, 2, algebra}
According to suggestion of induction
\[
f_1(x_0+th)=o(t^{n-1})
\]
Then equation
\EqRef{Gateaux differential, representation in algebra}
gets form
\[
o(t^{n-1})=\lim_{t\rightarrow 0,\ t\in R}(t^{-1}f(x+th))
\]
Therefore,
\[
f(x+th)=o(t^n)
\]
\end{proof}
\fi

Let us form polynomial
\ShowEq{Taylor polynomial, f(x), algebra}
According to theorem \ref{theorem: n Gateaux derivative equal 0, algebra}
\[
f(x_0+t(x-x_0))-p(x_0+t(x-x_0))=o(t^n)
\]
Therefore, polynomial $p(x)$ is good approximation of map $f(x)$.

If the mapping $f(x)$ has the G\^ateaux derivative of any order,
then passing to the limit
$n\rightarrow\infty$, we get expansion into series
\ShowEq{Taylor series, f(x), algebra}
which is called \AddIndex{Taylor series}{Taylor series, algebra}.

\ePrints{2010.05.23}
\ifx\Semafor\ValueOn
\newpage
{\bf Integral}
\fi
\Section{Integral}

\ePrints{2010.05.23}
\ifx\Semafor\ValueOff
Concept of integral has different aspect. In this section we consider
integration as operation inverse to differentiation.
As a matter of fact, we consider procedure of solution of ordinary differential
equation
\ShowEq{differential equation, algebra}

\begin{example}
I start from example of differential equation over real field.
\begin{equation}
\EqLabel{differential equation y=xx, 1, real number}
y'=3x^2
\end{equation}
\begin{equation}
\EqLabel{differential equation y=xx, initial, real number}
\begin{matrix}
x_0=0&y_0=0
\end{matrix}
\end{equation}
Differentiating one after another equation
\EqRef{differential equation y=xx, 1, real number},
we get the chain of equations
\begin{align}
\EqLabel{differential equation y=xx, 2, real number}
y''&=6x
\\
\EqLabel{differential equation y=xx, 3, real number}
y'''&=6
\\
\EqLabel{differential equation y=xx, 4, real number}
y^{(n)}&=0&n&>3
\end{align}
From equations
\EqRef{differential equation y=xx, 1, real number},
\EqRef{differential equation y=xx, initial, real number},
\EqRef{differential equation y=xx, 2, algebra},
\EqRef{differential equation y=xx, 3, algebra},
\EqRef{differential equation y=xx, 4, algebra}
it follows expansion into Taylor series
\[
y=x^3
\]
\qed
\end{example}

\begin{example}
\fi
Let us consider similar differential equation over algebra
\ShowEq{differential equation y=xx, algebra}
\ePrints{2010.05.23}
\ifx\Semafor\ValueOn
\newpage
\fi
Differentiating one after another equation
\EqRef{differential equation y=xx, 1, algebra},
we get the chain of equations
\ShowEq{differential equation y=xx, 234, algebra}
From equations
\EqRef{differential equation y=xx, 1, algebra},
\EqRef{differential equation y=xx, initial, algebra},
\EqRef{differential equation y=xx, 2, algebra},
\EqRef{differential equation y=xx, 3, algebra},
\EqRef{differential equation y=xx, 4, algebra}
expansion into Taylor series follows
\[
y=x^3
\]
\ePrints{2010.05.23}
\ifx\Semafor\ValueOn
\newpage
\else
\qed
\end{example}
\fi

\begin{remark}
I will write following equations
to show how derivative works. 
\ShowEq{differential equation y=xx, algebra, 1}
\qed
\end{remark}

\ePrints{2010.05.23}
\ifx\Semafor\ValueOn
\newpage
\fi
\begin{remark}
Differential equation
\ShowEq{differential equation y=xx, a, algebra}
also leads to answer $y=x^3$. it is evident that this
map does not satisfies differential equation.
However, contrary to theorem
\ref{theorem: Gateaux derivative of pn is symmetric, m < n, algebra}
second derivative is
not symmetric polynomial. This means that
equation \EqRef{differential equation y=xx, 1, a, algebra}
does not possess a solution.
\qed
\end{remark}

\ePrints{2010.05.23}
\ifx\Semafor\ValueOff
\begin{example}
It is evident that, if function satisfies to differential equation
\ShowEq{differential equation, additive function, algebra}
then The G\^ateaux derivative of second order
\ShowEq{differential equation, additive function, 1, algebra}
In that case, if initial condition is $y(0)=0$,
then differential equation
\EqRef{differential equation, additive function, algebra} has solution
\ShowEq{differential equation, additive function, solution, algebra}
\qed
\end{example}

\Section{Exponent}

In this section we consider one of possible models
of exponent.

In a field we can define exponent as solution of differential equation
\ShowEq{exponent over field}
It is evident that we cannot write such equation for division ring.
However we can use equation
\ShowEq{derivative over field}
From equations
\EqRef{exponent over field}, \EqRef{derivative over field}
it follows
\ShowEq{exponent derivative over field}
This equation is closer to our goal,
however there is the question: in which order
we should multiply $y$ and $h$?
To answer this question we change equation
\else
\newpage
Consider equation
\fi
\ShowEq{exponent derivative over division ring}
Hence, our goal is to solve differential equation
\EqRef{exponent derivative over division ring}
with initial condition $y(0)=1$.

\ePrints{2010.05.23}
\ifx\Semafor\ValueOff
For the statement and proof of the theorem
\ref{theorem: exponent derivative n over division ring}
I introduce following notation.
Let
\ShowEq{exponent derivative, transposition n}
be transposition of the tuple of variables
\ShowEq{exponent transposition tuple 1}
Let $p_{\sigma}(h_i)$ be position that
variable $h_i$ gets in the tuple
\ShowEq{exponent transposition tuple 2}
For instance, if transposition $\sigma$ has form
\ShowEq{exponent transposition tuple 3}
then following tuples equal
\ShowEq{exponent transposition tuple 4}

\begin{theorem}
\label{theorem: exponent derivative n over division ring}
If function $y$ is solution of differential equation
\EqRef{exponent derivative over division ring}
then the G\^ateaux derivative of order $n$ of function $y$
has form
\ShowEq{exponent derivative n over division ring}
where sum is over transpositions
\ShowEq{exponent derivative, transposition n}
of the set of variables $y$, $h_1$, ..., $h_n$.
Transposition $\sigma$ has following properties
\begin{enumerate}
\item If there exist $i$, $j$, $i\ne j$, such that
$p_{\sigma}(h_i)$ is situated in product
\EqRef{exponent derivative n over division ring}
on the left side of $p_{\sigma}(h_j)$ and $p_{\sigma}(h_j)$ is situated
on the left side of $p_{\sigma}(y)$, then $i<j$.
\label{enumerate: exponent derivative, transposition n, left}
\item If there exist $i$, $j$, $i\ne j$, such that
$p_{\sigma}(h_i)$ is situated in product
\EqRef{exponent derivative n over division ring}
on the right side of $p_{\sigma}(h_j)$ and $p_{\sigma}(h_j)$ is situated
on the right side of $p_{\sigma}(y)$, then $i>j$.
\label{enumerate: exponent derivative, transposition n, right}
\end{enumerate}
\end{theorem}
\begin{proof}
We prove this statement by induction.
For $n=1$ the statement is true because this is differential equation
\EqRef{exponent derivative over division ring}.
Let the statement be true for $n=k-1$.
Hence
\ShowEq{exponent derivative n=k-1 over division ring}
where the sum is over transposition
\ShowEq{exponent derivative, transposition n=k-1}
of the set of variables $y$, $h_1$, ..., $h_{k-1}$.
Transposition $\sigma$ satisfies to conditions
\eqref{enumerate: exponent derivative, transposition n, left},
\eqref{enumerate: exponent derivative, transposition n, right}
in theorem.
According to definition
\EqRef{Gateaux derivative of Order n, algebra}
the G\^ateaux derivative of order $k$ has form
\ShowEq{exponent derivative n=k over division ring, 1}
From equations
\EqRef{exponent derivative over division ring},
\EqRef{exponent derivative n=k over division ring, 1}
it follows that
\ShowEq{exponent derivative n=k over division ring, 2}
It is easy to see that arbitrary transposition $\sigma$
from sum
\EqRef{exponent derivative n=k over division ring, 2}
forms two transpositions
\ShowEq{exponent derivative, transposition n=k}
From \EqRef{exponent derivative n=k over division ring, 2}
and \EqRef{exponent derivative, transposition n=k}
it follows that
\ShowEq{exponent derivative n=k over division ring, 3}
In expression
\EqRef{exponent derivative n=k over division ring, 3}
$p_{\tau_1}(h_k)$ is written immediately before $p_{\tau_1}(y)$.
Since $k$ is smallest value of index
then transposition $\tau_1$ satisfies to conditions
\eqref{enumerate: exponent derivative, transposition n, left},
\eqref{enumerate: exponent derivative, transposition n, right}
in the theorem.
In expression
\EqRef{exponent derivative n=k over division ring, 3}
$p_{\tau_2}(h_k)$ is written immediately after $p_{\tau_2}(y)$.
Since $k$ is largest value of index
than transposition $\tau_2$ satisfies to conditions
\eqref{enumerate: exponent derivative, transposition n, left},
\eqref{enumerate: exponent derivative, transposition n, right}
in the theorem.

It remains to show that in the expression
\EqRef{exponent derivative n=k over division ring, 3}
we get all transpositions $\tau$ that satisfy to conditions
\eqref{enumerate: exponent derivative, transposition n, left},
\eqref{enumerate: exponent derivative, transposition n, right}
in the theorem.
Since $k$ is largest index then according to conditions
\eqref{enumerate: exponent derivative, transposition n, left},
\eqref{enumerate: exponent derivative, transposition n, right}
in the theorem $\tau(h_k)$
is written either immediately before or immediately after $\tau(y)$.
Therefore, any transposition $\tau$ has
either form $\tau_1$ or form $\tau_2$.
Using equation
\EqRef{exponent derivative, transposition n=k},
we can find corresponding
transposition $\sigma$ for given transposition $\tau$.
Therefore, the statement of theorem is true for $n=k$.
We proved the theorem.
\end{proof}
\fi

\begin{theorem}
The solution of differential equation
\EqRef{exponent derivative over division ring}
with initial condition $y(0)=1$
is exponent
\ShowEq{exponent over division ring}
that has following Taylor series expansion
\ShowEq{exponent Taylor series over division ring}
\ePrints{2010.05.23}
\ifx\Semafor\ValueOn
\qed
\end{theorem}
\newpage
\else
\end{theorem}
\begin{proof}
The G\^ateaux derivative of order $n$ has $2^n$ items.
In fact, the G\^ateaux derivative of order $1$ has $2$ items,
and each differentiation increase number of items twice.
From initial condition $y(0)=1$ and theorem
\ref{theorem: exponent derivative n over division ring}
it follows that the G\^ateaux derivative of order $n$
of required solution has form
\ShowEq{exponent derivative n over division ring, x=0}
Taylor series expansion
\EqRef{exponent Taylor series over division ring}
follows
from equation
\EqRef{exponent derivative n over division ring, x=0}.
\end{proof}
\fi

\begin{theorem}
\label{theorem: exponent of sum}
The equation
\ShowEq{exponent of sum}
is true iff
\ShowEq{exponent of sum, 1}
\ePrints{2010.05.23}
\ifx\Semafor\ValueOn
\qed
\end{theorem}
\newpage
\else
\end{theorem}
\begin{proof}
To prove the theorem it is enough to consider Taylor series
\ShowEq{exponent a b a+b}
Let us multiply expressions
\EqRef{exponent a} and \EqRef{exponent b}.
The sum of monomials of order $3$ has form
\ShowEq{exponent ab 3}
and in general does not equal expression
\ShowEq{exponent a+b 3}
The proof of statement that \EqRef{exponent of sum} follows from
\EqRef{exponent of sum, 1}
is trivial.
\end{proof}

The meaning of the theorem \ref{theorem: exponent of sum}
becomes more clear if we recall that there exist two
models of design of exponent. First model is
the solution of differential equation
\EqRef{exponent derivative over division ring}.
Second model is exploring of one parameter group of transformations.
For field both models lead to the same function.
I cannot state this now for general case. This is the subject of separate
research. However if we recall that quaternion is analogue
of transformation of three dimensional space then the statement of the theorem
becomes evident.
\fi

%% file: Second.Derivative.Eq.tex
%auto-ignore

\def\DfTwo{\partial(\partial f(x)\circ a_1)\circ a_2}

\DefEq
{
$o_2:A\rightarrow \LAA$
}
{Gateaux differential of map df, algebra, o2}

\DefEq
{
\ePrints{1006.2597}
\Items{Calculus.Paper}
\ifx\Semafor\ValueOn
\ref{theorem: module L(A;A) is algebra}
\else
\xRef{1003.1544}{theorem: module L(A;A) is algebra}
\fi
}
{theorem: module L(A;A) is algebra}

\DefEq
{
\ePrints{Calculus.Paper}
\ifx\Semafor\ValueOn
\ref{theorem: tensor product and polylinear mapping}.
\else
\xRef{1003.1544}{theorem: tensor product and polylinear mapping}.
\fi
}
{theorem: tensor product and polylinear mapping}

\DefEq
{
\[
\lim_{a_2\rightarrow 0}\frac{\|o_2(a_2)\|}{|a_2|}=0
\]
}
{Gateaux differential of map df, algebra, o2 lim}

\DefEq
{
\[
\partial^2 f(x)\circ(a_1\otimes a_2)=
\partial^2 f(x)\circ(a_1;a_2)
\]
}
{differential AA A,A}

\DefEq
{
$\partial^2 f(x)\in\mathcal L(A,A;A)$%
}
{differential L(A,A;A), 1}

\DefEq
{
\[
\partial^2 f:A\rightarrow \mathcal L(A,A;A)
\]
}
{differential L(A,A;A)}

\DefEq
{
\begin{align*}
\partial^2 f(x)\in&\ \mathcal L(A\otimes A;A)
\\
\partial^2 f:A\rightarrow&\ \mathcal L(A\otimes A;A)
\end{align*}
}
{differential L(AA;A)}

\DefEq
{
\symb{\partial^2 f(x)\circ(a_1; a_2)}
0{Gateaux differential of Second Order, algebra}
\begin{equation}
\begin{array}{r@{}l}
\ShowSymbol{Gateaux differential of Second Order, algebra}
&\displaystyle
=
\left(
\frac{\partial\pC{s}{0}^2 f(x)}{\partial x^2}
\otimes
\frac{\partial\pC{s}{1}^2 f(x)}{\partial x^2}
\otimes
\frac{\partial\pC{s}{2}^2 f(x)}{\partial x^2},\sigma_s
\right)
\circ(a_1; a_2)
\\
&\displaystyle\VirtFrac
=\frac{\partial\pC{s}{0}^2 f(x)}{\partial x^2}
\sigma_s(a_1)
\frac{\partial\pC{s}{1}^2 f(x)}{\partial x^2}
\sigma_s(a_2)
\frac{\partial\pC{s}{2}^2 f(x)}{\partial x^2}
\end{array}
\EqLabel{Differential of Second Order, algebra, representation}
\end{equation}
}
{Differential of Second Order, algebra, representation}

\DefEquation
{
\begin{array}{r@{\ }l}
\partial P_{n-k}(x)\circ h_k
=&\partial(\partial^{k-1} p_n(x)\circ(h_1;...;h_{k-1}))\circ h_k
\\
=&\partial^k p_n(x)\circ(h_1;...;h_{k-1}; h_k)
\end{array}
}
{Gateaux derivative of Order n, 1, algebra}

\DefEq
{
\[
\partial (f(x)x)\circ h_1=(\partial f(x)\circ h_1)x+
f(x)h_1
\]
}
{Gateaux derivative of f(x)x, 1, algebra}

\DefEq
{
\[
f_1(x_0)=\partial f_1(x_0)\circ h=...=\partial^{n-1} f_1(x_0)\circ h=0
\]
}
{n Gateaux derivative equal 0, 2, algebra}

\DefEq
{
\[
f_1(x)=\partial f(x)\circ h
\]
}
{n Gateaux derivative equal 0, 1, algebra}

\DefEq
{
\[
\partial^n p_n(x)\circ h
=n(\partial^{n-1} p_{n-1}(x)\circ h)ha_n=n(n-1)!p_{n-1}(h)ha_n=n!p_n(h)
\]
}
{n Gateaux derivative of polinom pn, 3, algebra}

\DefEquation
{
\begin{array}{r@{\ }l}
\partial^n p_n(x)\circ h
=&(\partial^n p_{n-1}(x)\circ h)xa_n
+(\partial^{n-1} p_{n-1}(x)\circ h)ha_n
\\
+&...
+(\partial^{n-1} p_{n-1}(x)\circ h)ha_n
\end{array}
}
{n Gateaux derivative of polinom pn, 1, algebra}

\DefEq
{
\[
\partial p_1(x)\circ h=\partial (a_0xa_1)\circ h=a_0ha_1=1!p_1(h)
\]
}
{n Gateaux derivative of polinom pn, 2, algebra}

\DefEq
{
\[
P_{n-k}(x)=\partial^{k-1} p_n(x)\circ(h_1;...;h_{k-1})
\]
}
{Gateaux derivative of pn is symmetric, 2}

\DefEquation
{
r_n(h_1,x_2,...,x_n)=r_{n-1}(h_1,x_{[2},...,x_{n-1})x_{n]}a_n+r_{n-1}(x_2,...,x_n)h_1a_n
}
{Gateaux derivative of monomial, n, algebra}

\DefEq
{
\[
r_n(h_1,x_2,...,x_n)=\partial p_n(x)\circ h_1
\]
}
{Gateaux derivative of pn is symmetric, 1}

\DefEq
{
\[
\begin{array}{r@{\ }l}
&\partial^{m-1} (f(x)x)\circ(h_1;...;h_{m-1})
\\
=&\partial^{m-1} f(x)\circ(h_1;...;h_{m-1})x
+\partial^{m-2} f(x)\circ(h_1;...;h_{m-2})h_{m-1}
\\
+&\partial^{m-2} f(x)\circ(\widehat{h_1};...;h_{m-2};h_{m-1})h_1
+...
\\
+&\partial^{m-2} f(x)\circ(h_1;...;\widehat{h_{m-2}};h_{m-1})h_{m-2}
\end{array}
\]
}
{Gateaux derivative of f(x)x, m-1, algebra}

\DefEquation
{
\begin{array}{r@{\ }l}
&\partial^m (f(x)x)\circ(h_1;...;h_{m-1};h_m)
\\
=&\partial^m f(x)\circ(h_1;...;h_{m-1};h_m)x
\\
+&\partial^{m-1} f(x)\circ(h_1;...;h_{m-2};h_{m-1})h_m
\\
+&\partial^{m-1} f(x)\circ(h_1;...;h_{m-2};\widehat{h_{m-1}};h_m)h_{m-1}
\\
+&\partial^{m-2} f(x)\circ(\widehat{h_1};...;h_{m-2};h_{m-1};h_m)h_1
+...
\\
+&\partial^{m-2} f(x)\circ(h_1;...;\widehat{h_{m-2}};h_{m-1};h_m)h_{m-2}
\end{array}
}
{Gateaux derivative of f(x)x, m, algebra}

\DefEquation
{
\begin{array}{r@{\ }l}
&\partial^m (f(x)x)\circ(h_1;...;h_m)
\\
=&\partial^m f(x)\circ(h_1;...;h_m)x
+\partial^{m-1} f(x)\circ(h_1;...;h_{m-1})h_m
\\
+&\partial^{m-1} f(x)\circ(\widehat{h_1};...;h_{m-1};h_m)h_1
+...
\\
+&\partial^{m-1} f(x)\circ(h_1;...;\widehat{h_{m-1}};h_m)h_{m-1}
\end{array}
}
{Gateaux derivative of f(x)x, algebra}

\DefEq
{
$\partial^m p_n(x)\circ(h_1;...;h_m)$
}
{Gateaux derivative of pn is symmetric}

\DefEq
{
\[
\frac{\partial\pC{s}{p}^2 f(x)}{\partial x^2}
=\frac{\partial\pC{s}{p}^2 f(x)}{\partial x\partial x}
\]}
{component of Gateaux derivative of Second Order, 1, algebra}

\DefEq
{
\[
\partial^k p(0)\circ(h_1;...;h_k)=k!p_{ki_k}(x)
\]
}
{Taylor polynomial 1, algebra}

\DefEq
{
\[
p(x)=p_0
+(1!)^{-1}\partial p(0)\circ x
+(2!)^{-1}\partial^2 p(0)\circ x+...+(n!)^{-1}\partial^n p(0)\circ x
\]
}
{Taylor polynomial 2, algebra}

\DefEq
{
\[
f(x_0)=\partial f(x_0)\circ h=...=\partial^n f(x_0)\circ h=0
\]
}
{n Gateaux derivatives of function, algebra}

\DefEq
{
\[
\partial f(x)\circ h=F(x;h)
\]
}
{differential equation, algebra}

\DefEq
{
\[
f(x)=\sum_{n=0}^{\infty}(n!)^{-1}\partial^n f(x_0)\circ(x-x_0)
\]
}
{Taylor series, f(x), algebra}

\DefEq
{
\[
p(x)=f(x_0)+(1!)^{-1}\partial f(x_0)\circ(x-x_0)
+...+(n!)^{-1}\partial^n f(x_0)\circ(x-x_0)
\]
}
{Taylor polynomial, f(x), algebra}

\DefEq
{
\[
p(y)=p_0+(1!)^{-1}\partial p(y_0)\circ(y-y_0)
+(2!)^{-1}\partial^2 p(y_0)\circ(y-y_0)+...
+(n!)^{-1}\partial^n p(y_0)\circ(y-y_0)
\]
}
{Taylor polynomial 3, algebra}

\DefEq
{
\[
\partial^n p_n(x)\circ h=n!p_n(h)
\]
}
{n Gateaux derivative of polinom pn, algebra}

\DefEq
{
\[
\partial^nf(x)\circ h=\partial^nf(x)\circ(h_1;...;h_n)
\]
}
{dfxh}

\DefEq
{
\begin{align}
\EqLabel{differential equation y=xx, 1, algebra}
\partial y&=1\otimes x^2+x\otimes x+x^2\otimes 1
\\
\EqLabel{differential equation y=xx, initial, algebra}
&\begin{matrix}
x_0=0&y_0=0
\end{matrix}
\end{align}
}
{differential equation y=xx, algebra}

\DefEq
{
\begin{align*}
\partial y\circ h&=hx^2+xhx+x^2h
\\
\partial^2 y\circ(h_1;h_2)
&=h_1h_2x+h_1xh_2+h_2h_1x
\\
&+xh_1h_2+h_2xh_1+xh_2h_1
\\
\partial^3 y\circ(h_1;h_2;h_3)
&=h_1h_2h_3+h_1h_3h_2+h_2h_1h_3
\\
&+h_3h_1h_2+h_2h_3h_1+h_3h_2h_1
\end{align*}
}
{differential equation y=xx, algebra, 1}

\DefEq
{
\[
\partial^m p_n(0)=0
\]
}
{m Gateaux derivative of polinom pn, algebra}

\DefEq
{
\[
\partial^{n+1} p_n(x)=0
\]
}
{Gateaux derivative n1 of pn, algebra}

\DefEq
{
\symb{\partial^n f(x)}
0{Gateaux derivative of Order n, algebra}
\symb{\frac{\partial^n f(x)}{\partial x^n}}
0{Gateaux derivative of Order n, fraction, algebra}
\begin{equation}
\EqLabel{Gateaux derivative of Order n, algebra}
\begin{array}{r@{}l}
\ShowSymbol{Gateaux derivative of Order n, algebra}
\circ(a_1;...;a_n)=
\\
\displaystyle\VirtFrac
\ShowSymbol{Gateaux derivative of Order n, fraction, algebra}
\circ(a_1;...;a_n)=
&
\ \partial(\partial^{n-1} f(x)
\circ(a_1;...;a_{n-1}))\circ a_n
\end{array}
\end{equation}
}
{Gateaux derivative of Order n, algebra}

\DefEquation
{
\begin{matrix}
\partial^k p_n(x)\circ(h_1;...;h_k)=r_n(h_1,...,h_k,x_{k+1},...,x_n)
&
x_{k+1}=...=x_n=x
\end{matrix}
}
{Gateaux derivative of monomial, algebraic definition, algebra}

\DefEq
{
\symb{\frac{\partial\pC{s}{p}^2 f(x)}{\partial x^2}}0
{component of Gateaux derivative of Second Order, algebra}
\[
\begin{matrix}
\displaystyle
\ShowSymbol{component of Gateaux derivative of Second Order, algebra}
&p=0, 1, 2
\end{matrix}
\]
}
{component of Gateaux derivative of Second Order, algebra}

\DefEq
{
\begin{equation}
\EqLabel{differential equation y=xx, 2, algebra}
\begin{array}{r@{}l}
\partial^2 y
&
=1\otimes_1 1\otimes_2x+1\otimes_1x\otimes_21
+1\otimes_21\otimes_1x
\\
&+x\otimes_11\otimes_21+1\otimes_2x\otimes_11
+x\otimes_21\otimes_11
\end{array}
\end{equation}
\begin{equation}
\EqLabel{differential equation y=xx, 3, algebra}
\begin{array}{r@{}l}
\partial^3 y
&=1\otimes_11\otimes_21\otimes_31+1\otimes_11\otimes_31\otimes_21
+1\otimes_21\otimes_11\otimes_31
\\
&+1\otimes_31\otimes_11\otimes_21+1\otimes_21\otimes_31\otimes_11
+1\otimes_31\otimes_21\otimes_11
\end{array}
\end{equation}
\begin{equation}
\EqLabel{differential equation y=xx, 4, algebra}
\partial^n y=0\ \ \ n>3
\end{equation}
}
{differential equation y=xx, 234, algebra}

\DefEquation
{
\partial P_{n-k}(x)\circ h_k=R_{n-k}(h_k,x_{k+1},...,x_n)\ \ \ x_{k+1}=...=x_n=x
}
{Gateaux derivative of monomial, algebraic definition, 1, algebra}

\DefEq
{
\[
r_n(h_1,x_2,...,x_n)=(\partial p_{n-1}(x)\circ h_1)xa_n+p_{n-1}(x)h_1a_n
\]
}
{Gateaux derivative of monomial, n1, algebra}

\DefEq
{
\begin{equation}
\EqLabel{differential equation y=xx, 1, a, algebra}
\partial y=3\otimes x^2
\end{equation}
\begin{equation}
\EqLabel{differential equation y=xx, initial, a, algebra}
\begin{matrix}
x_0=0&y_0=0
\end{matrix}
\end{equation}
}
{differential equation y=xx, a, algebra}

\DefEq
{
\begin{equation}
\partial y=f\pC{s}{0}\otimes f\pC{s}{1}
\EqLabel{differential equation, additive function, algebra}
\end{equation}
\[
\begin{matrix}
f\pC{s}{0}\in A&f\pC{s}{1}\in A
\end{matrix}
\]
}
{differential equation, additive function, algebra}

\DefEq
{
\[
\partial^2 f(x)=0
\]
}
{differential equation, additive function, 1, algebra}

\DefEq
{
\[
y=f\pC{s}{0}\ x\ f\pC{s}{1}
\]
}
{differential equation, additive function, solution, algebra}

\DefEq
{
\begin{equation}
y'=y
\EqLabel{exponent over field}
\end{equation}
}
{exponent over field}

\DefEq
{
\begin{equation}
\partial(y)\circ h=y'h
\EqLabel{derivative over field}
\end{equation}
}
{derivative over field}

\DefEq
{
\begin{equation}
\partial(y)\circ h=yh
\EqLabel{exponent derivative over field}
\end{equation}
}
{exponent derivative over field}

\DefEq
{
\begin{equation}
\partial(y)\circ h=\frac 12(yh+hy)
\EqLabel{exponent derivative over division ring}
\end{equation}
}
{exponent derivative over division ring}

\DefEq
{
\[
\sigma=
\begin{pmatrix}
y&h_1&...&h_n
\\
\sigma(y)&\sigma(h_1)&...&\sigma(h_n)
\end{pmatrix}
\]
}
{exponent derivative, transposition n}

\DefEq
{
\[
\begin{pmatrix}
y&h_1&...&h_n
\end{pmatrix}
\]
}
{exponent transposition tuple 1}

\DefEq
{
\[
\begin{pmatrix}
\sigma(y)&\sigma(h_1)&...&\sigma(h_n)
\end{pmatrix}
\]
}
{exponent transposition tuple 2}

\DefEq
{
\[
\begin{pmatrix}
y&h_1&h_2&h_3
\\
h_2&y&h_3&h_1
\end{pmatrix}
\]
}
{exponent transposition tuple 3}

\DefEq
{
\begin{align*}
\begin{pmatrix}
\sigma(y)&\sigma(h_1)&\sigma(h_2)&\sigma(h_3)
\end{pmatrix}
=&
\begin{pmatrix}
h_2&y&h_3&h_1
\end{pmatrix}
\\
=&
\begin{pmatrix}
p_{\sigma}(h_2)&p_{\sigma}(y)&p_{\sigma}(h_3)&p_{\sigma}(h_1)
\end{pmatrix}
\end{align*}
}
{exponent transposition tuple 4}

\DefEquation
{
\partial^n(y)\circ(h_1,...,h_n)=\frac 1{2^n}
\sum_{\sigma} \sigma(y)\sigma(h_1) ... \sigma(h_n)
}
{exponent derivative n over division ring}

\DefEquation
{
\partial^{k-1}(y)\circ (h_1,...,h_{k-1})=\frac 1{2^{k-1}}
\sum_{\sigma} \sigma(y)\sigma(h_1) ... \sigma(h_{k-1})
}
{exponent derivative n=k-1 over division ring}

\DefEq
{
\[
\sigma=
\begin{pmatrix}
y&h_1&...&h_{k-1}
\\
\sigma(y)&\sigma(h_1)&...&\sigma(h_{k-1})
\end{pmatrix}
\]
}
{exponent derivative, transposition n=k-1}

\DefEquation
{
\begin{array}{r@{\ }l}
\partial^k(y)\circ(h_1,...,h_k)
=&\partial(\partial^{k-1}(y)\circ(h_1,...,h_{k-1}))\circ h_k
\\
=&\displaystyle\frac 1{2^{k-1}}\partial
\left(
\sum_{\sigma} \sigma(y)\sigma(h_1) ... \sigma(h_{k-1})
\right)
\circ h_k
\end{array}
}
{exponent derivative n=k over division ring, 1}

\DefEquation
{
\begin{array}{r@{\ }l}
&\partial^k(y)\circ(h_1,...,h_k)
\\
=&\displaystyle\frac 1{2^{k-1}}\frac 12
\left(
\sum_{\sigma} \sigma(yh_k)\sigma(h_1) ... \sigma(h_{k-1})
+
\sum_{\sigma} \sigma(h_ky)\sigma(h_1) ... \sigma(h_{k-1})
\right)
\end{array}
}
{exponent derivative n=k over division ring, 2}

\DefEq
{
\begin{equation}
\begin{array}{rl}
\tau_1
&=
\begin{pmatrix}
y&h_1&...&h_{k-1}&h_k
\\
\tau_1(y)&\tau_1(h_1)&...&\tau_1(h_{k-1})&\tau_1(h_k)
\end{pmatrix}
\\
&=
\begin{pmatrix}
h_ky&h_1&...&h_{k-1}&
\\
\sigma(h_ky)&\sigma(h_1)&...&\sigma(h_{k-1})
\end{pmatrix}
\\
\tau_2
&=
\begin{pmatrix}
y&h_1&...&h_{k-1}&h_k
\\
\tau_2(y)&\tau_2(h_1)&...&\tau_2(h_{k-1})&\tau_2(h_k)
\end{pmatrix}
\\
&=
\begin{pmatrix}
yh_k&h_1&...&h_{k-1}&
\\
\sigma(yh_k)&\sigma(h_1)&...&\sigma(h_{k-1})
\end{pmatrix}
\end{array}
\EqLabel{exponent derivative, transposition n=k}
\end{equation}
}
{exponent derivative, transposition n=k}

\DefEq
{
\begin{align}
&\partial^k(y)\circ(h_1,...,h_k)
\nonumber
\\
=&\frac 1{2^k}
\Bigg(
\sum_{\tau_1} \tau_1(y)\tau_1(h_1) ... \tau_1(h_{k-1})\tau_1(h_k)
\EqLabel{exponent derivative n=k over division ring, 3}
\\
+&
\sum_{\tau_2} \tau_2(y)\tau_2(h_1) ... \tau_2(h_{k-1})\tau_2(h_k)
\Bigg)
\nonumber
\end{align}
}
{exponent derivative n=k over division ring, 3}

\DefEq
{
$y=e^x$
}
{exponent over division ring}

\DefEq
{
\begin{equation}
e^x=\sum_{n=0}^{\infty}
\frac 1{n!}x^n
\EqLabel{exponent Taylor series over division ring}
\end{equation}
}
{exponent Taylor series over division ring}

\DefEquation
{
\partial^n(0)\circ(h,...,h)=1
}
{exponent derivative n over division ring, x=0}

\DefEq
{
\begin{equation}
e^{a+b}=e^ae^b
\EqLabel{exponent of sum}
\end{equation}
}
{exponent of sum}

\DefEq
{
\begin{equation}
ab=ba
\EqLabel{exponent of sum, 1}
\end{equation}
}
{exponent of sum, 1}

\DefEq
{
\begin{align}
e^a&=\sum_{n=0}^{\infty}\frac 1{n!}a^n
\EqLabel{exponent a}
\\
e^b&=\sum_{n=0}^{\infty}\frac 1{n!}b^n
\EqLabel{exponent b}
\\
e^{a+b}&=\sum_{n=0}^{\infty}\frac 1{n!}(a+b)^n
\EqLabel{exponent a+b}
\end{align}
}
{exponent a b a+b}

\DefEq
{
\begin{equation}
\frac 16a^3+\frac 12a^2b+\frac 12ab^2+\frac 16b^3
\EqLabel{exponent ab 3}
\end{equation}
}
{exponent ab 3}

\DefEq
{
\begin{equation}
\frac 16(a+b)^3=\frac 16a^3+\frac 16a^2b+\frac 16aba+\frac 16ba^2
+\frac 16ab^2+\frac 16bab+\frac 16b^2a+\frac 16b^3
\EqLabel{exponent a+b 3}
\end{equation}
}
{exponent a+b 3}

\DefEq
{
$\DfTwo$
}
{increment of derivative, algebra}

\DefEq
{
\symb{\partial^2 f(x)}
0{Gateaux derivative of Second Order, algebra}
\symb{\frac{\partial^2 f(x)}{\partial x^2}}
0{Gateaux derivative of Second Order, fraction, algebra}
\begin{equation}
\EqLabel{Gateaux derivative of Second Order, algebra}
\ShowSymbol{Gateaux derivative of Second Order, algebra}\circ(a_1;a_2)
=\ShowSymbol{Gateaux derivative of Second Order, fraction, algebra}
\circ(a_1;a_2)
=\DfTwo
\end{equation}
}
{Gateaux derivative of Second Order, algebra}

\DefEq
{
\begin{equation}
\EqLabel{Gateaux differential of map df, algebra}
(\partial f\circ(x+a_2))\circ a_1-(\partial f\circ x)\circ a_1
=\DfTwo
+o_2(a_2)
\end{equation}
}
{Gateaux differential of map df, algebra}

%% file: Biblio.English.tex
%auto-ignore
\OpenBiblio

%22. A.G. Kurosh, Lectures on General Algebra, Chelsea, New York, 1963

\BiblioItem{Einstein: Electrodynamics of Moving Bodies}
{
Albert Einstein,
On the Electrodynamics of Moving Bodies, 1905,\\
Hendrik Antoon Lorentz, The Principle of Relativity, 37 - 65,
Translated by W Perrett, G B Jeffery
Courier Dover Publications, 1952
}%

\BiblioItem{Einstein: Foundations of general relativity}
{
Albert Einstein,
Die Grundlage der allgemeinen Relativit\"atstheorie,
Ann. Phys., 1916, {\bf 49}, 769 - 822,\\
Einstein's Annalen Papers: The Complete Collection 1901-1922,
edited by J\"urgen Renn, 517 - 571,\\
Wiley-VCH Verlag GmbH \& Co. KGaA, 2005
}%

\BiblioItem{Einstein: Geometry and Experience}
{
Albert Einstein, Geometry and Experience, (1921)\\
Albert Einstein, Sidelights on Relativity, 25 - 56,\\
Courier Dover Publications, 1983
}%

\BiblioItem{Einstein: Main problems of general relativity}
{
Albert Einstein,
Grundgedanken und Probleme der Relativit\"atstheorie, (1923),\\
Nobelstiftelsen, Les Prix Nobel en 1921 - 1922,
Imprimerie Royale, Stockholm, 1923
}%

\BiblioItem{Einstein: Noneuclidean Geometry and Physics}
{
Albert Einstein,
Nichtenklidische Geometrie in der Physik Neue Rundschan, (1925)
Berlin, S. 16 - 20
}%

\BiblioItem{Einstein: Isaak Newton}
{
Albert Einstein,
Isaak Newton, 1927,
Out of My Later Years, 
Citadel Press, 1995, 219 - 223
}%

\BiblioItem{Einstein: On Science}
{
Albert Einstein,
On Science, 
Cosmic Religion, with Other Opinions and Aphorisms,142 - 146,
New York, 1931, 97 - 103
}%

\BiblioItem{Einstein: Autobiographical Notes}
{
Albert Einstein,
Autobiographical Notes, 1949,\\
Paul A. Schilpp, editor, Albert Einstein: Philosopher-Scientist,
Evanston, 
Illinois, The Library of Living Philosophers, 1949, 1 - 95
}%

\BiblioItem{Cite: 104}
{
Cite 104, Source unknown
}%

\BiblioItem{Ghez}
{
Ghez et al.,
The First Measurement of Spectral Lines in a Short-Period Star Bound to the Galaxy's Central Black Hole: A Paradox of Youth,
\href{http://www.journals.uchicago.edu/ApJ/journal/issues/ApJL/v586n2/16990/brief/16990.abstract.html}{ApJL, 586, L127} (2003),
eprint \href{http://arxiv.org/abs/astro-ph/0302299}{arXiv:astro-ph/0302299} (2003)
}%

\BiblioItem{Schodel}
{
R. Sch\"odel et al.,
A star in a 15.2-year orbit around the supermassive black hole at the centre of the Milky Way,
\href{http://www.nature.com/cgi-taf/DynaPage.taf?file=/nature/journal/v419/n6908/abs/nature01121_fs.html}{Nature 419, 694} (2002)
}%

\BiblioItem{Mielke}
{
Eckehard W. Mielke, Affine generalization of the Komar complex of general relativity,
\href{http://prola.aps.org/searchabstract/PRD/v63/i4/e044018}{Phys. Rev. D 63, 044018} (2001)
}%

\BiblioItem{Obukhov}
{
Yu. N. Obukhov and J. G. Pereira, Metric\hyph affine approach to teleparallel gravity,
\href{http://scitation.aip.org/getabs/servlet/GetabsServlet?prog=normal&id=PRVDAQ000067000004044016000001&idtype=cvips&gifs=Yes}
{Phys. Rev. D 67, 044016} (2003),
eprint \href{http://arxiv.org/abs/gr-qc/0212080}{arXiv:gr-qc/0212080} (2002)
}%

\BiblioItem{Sardanashvily}
{
Giovanni Giachetta, Gennadi Sardanashvily, Dirac Equation in Gauge and Affine-Metric Gravitation Theories,
eprint \href{http://arxiv.org/abs/gr-qc/9511035}{arXiv:gr-qc/9511035} (1995)
}%

\BiblioItem{Gauge}
{
Frank Gronwald and Friedrich W. Hehl, On the Gauge Aspects of Gravity, eprint
\href{http://arxiv.org/abs/gr-qc/9602013}{arXiv:gr-qc/9602013} (1996)
}%

\BiblioItem{Neeman}
{
Yuval Neeman, Friedrich W. Hehl, Test Matter in a Spacetime with Nonmetricity, eprint
\href{http://arxiv.org/abs/gr-qc/9604047}{arXiv:gr-qc/9604047} (1996)
}%

\BiblioItem{torsion}
{
F. W. Hehl, P. von der Heyde, G. D. Kerlick, and J. M. Nester,
General relativity with spin and torsion: Foundations and prospects,\\
\href{http://prola.aps.org/abstract/RMP/v48/i3/p393_1}{Rev. Mod. Phys. 48, 393} (1976)
}%

\BiblioItem{Megged}
{
O. Megged, Post-Riemannian Merger of Yang-Mills Interactions with Gravity,
eprint \href{http://arxiv.org/abs/hep-th/0008135}{arXiv:hep-th/0008135} (2001)
}%

%\BiblioItem{Hehl}
%{
%Friedrich W. Hehl, Uwe Muench,
%eprint \href{http://arxiv.org/abs/gr-qc/9708007}{arXiv:gr-qc/9708007} (1997)
%}%

\BiblioItem{gr-qc-9604027}
{
Yu.N. Obukhov, E.J. Vlachynsky, W. Esser, R. Tresguerres and F.W. Hehl,
An exact solution of the metric\hyph affine gauge theory with dilation, shear, and spin charges,
eprint \href{http://arxiv.org/abs/gr-qc/9604027}{arXiv:gr-qc/9604027} (1996)
}%

\BiblioItem{4419-7514}
{
Mari\'an Fabian, Petr Habala, Petr H\'ajek, Vicente Montesinos, V\'aclav Zizler.
Banach Space Theory: The Basis for Linear and Nonlinear Analysis.
\\
Springer; New York, 2010; ISBN-13: 978-1441975140
}%

\BiblioItem{Weinberg I}
{
Steven Weinberg.
The Quantum Theory of Fields. Volume I. Foundations.
Cambridge university press, 1995
}%

\BiblioItem{Weinberg II}
{
Steven Weinberg.
The Quantum Theory of Fields. Volume II. Modern applications.
Cambridge university press, 1996
}%

\BiblioItem{Reinhardt}
{
Walter Greiner, Joachim Reinhardt. Field Quantization. Springer.
}%

\BiblioItem{978-3540875604}
{
Walter Greiner, Joachim Reinhardt. Quantum Electrodynamics. Springer, 2009.
}%

\BiblioItem{978-1898563020}
{
H. Robert Mills. Practical Astronomy. Woodhead Publishing, 1994. ISBN-13: 978-1898563020.
}%

\BiblioItem{Landau}
{
L. D. Landau, E. M. Lifshich, The classical theory of fields.
Oxford, New York, Pergamon Press
}%

\BiblioItem{Wheeler}
{
Ignazio Ciufolini, John Wheeler. Gravitation and Inertia.
Princeton university press.
}%

\BiblioItem{Anderson02}
{
J. D. Anderson, P. A. Laing, E. L. Lau, A. S. Liu, M. M. Nieto, and S. G. Turyshev,
Study of the anomalous acceleration of Pioneer 10 and 11,
\href{http://prola.aps.org/searchabstract/PRD/v65/i8/e082004}{Phys. Rev. D 65, 082004, 50 pp.}, (2002),
eprint \href{http://arxiv.org/abs/gr-qc/0104064}{arXiv:gr-qc/0104064} (2001)
}%

\BiblioItem{Anderson98}
{
J. D. Anderson, P. A. Laing, E. L. Lau, A. S. Liu, M. M. Nieto, and S. G. Turyshev,
Indication, from Pioneer 10/11, Galileo, and Ulysses Data, of an Apparent Anomalous, Weak, Long-Range Acceleration,
\href{http://prola.aps.org/abstract/PRL/v81/i14/p2858_1}{Phys. Rev. Lett. 81, 2858}, (1998),
eprint \href{http://arxiv.org/abs/gr-qc/9808081}{arXiv:gr-qc/9808081} (1998)
}%

%\BiblioItem{Havas} Peter Havas, The Classical Equations of Motion of Point Particles, I,
%{
%\href{http://prola.aps.org/abstract/PR/v87/i2/p309_1}{Phys. Rev. 87, 309} (1952)
%}%

\BiblioItem{H. Aslaksen}
{
H. Aslaksen.  Quaternionic determinants \textit{Math.
Intelligencer} {\bf 18}(3), pp.57-65, (1996).
}%

\BiblioItem{L. Chen: Definition of determinant}
{
L. Chen, Definition of determinant and Cramer solutions over
quaternion field, \textit{Acta Math. Sinica (N.S.)} {\bf 7},
pp.171-180, (1991).
}%

\BiblioItem{L. Chen: Inverse matrix}
{
L. Chen,
Inverse matrix and properties of double determinant over quaternion
field, \textit{Sci. China, Ser. A} {\bf 34}, pp.528-540, (1991).
}%

\BiblioItem{N. Cohen S. De Leo}
{
N. Cohen, S. De Leo, The quaternionic determinant, \textit{The Electronic Journal Linear
Algebra} {\bf 7}, pp.100-111, (2000).
}%

\BiblioItem{Dyson: Quaternion determinants}
{
F. J. Dyson, Quaternion determinants, \textit{Helvetica Phys.
Acta} {\bf 45}, pp. 289-302, (1972).
}%

\BiblioItem{Melvin Hausner}
{
Melvin Hausner,
A Vector Space Approach to Geometry,
Dover Publications, 1998
}%

\BiblioItem{Serge Lang}
{
Serge Lang,
Algebra, Springer, 2002
}%

\BiblioItem{Burris Sankappanavar}
{
S. Burris, H.P. Sankappanavar,
A Course in Universal Algebra, Springer-Verlag (March, 1982),
\\eprint
\href{http://www.math.uwaterloo.ca/~snburris/htdocs/ualg.html}
{http://www.math.uwaterloo.ca/~snburris/htdocs/ualg.html}
\\(The Millennium Edition)
}%

\BiblioItem{Shilov}
{
G. E. Shilov,
Calculus, Multivariable Functions,
Moscow, Nauka, 1972
}%

\BiblioItem{Kolmogorov Fomin}
{
A. N. Kolmogorov and S. V. Fomin,
Elements of the Theory of Functions and Functional Analysis,
Courier Dover Publication, 1999
}%

\BiblioItem{Lebedev Vorovich}
{
L. P. Lebedev, I. I. Vorovich,
Functional Analysis in Mechanics,
Springer, 2002
}%

\BiblioItem
{Rashevsky}
{
P. K. Rashevsky, Riemann Geometry and Tensor Calculus,\\
Moscow, Nauka, 1967
}%

\BiblioItem
{Kurosh: High Algebra}
{
A. G. Kurosh, High Algebra,
Moscow, Nauka, 1968
}%

\BiblioItem
{Kurosh: General Algebra}
{
A. G. Kurosh, Lectures on General Algebra,
Chelsea Pub Co, 1965 
}%

\BiblioItem
{Sabinin: Smooth Quasigroups}
{
Lev V. Sabinin, Smooth Quasigroups and Loops,
Kluwer Academic Publisher, 1999 
}%

\BiblioItem{Dubrovin Fomenko Novikov part 1}
{
B. A. Dubrovin, A. T. Fomenko, S. P. Novikov,
Modern Geometry - Methods and Applications,\\
Part 1, The Geometry of Surfaces, Transformation Groups, and Fields,\\
Translated by Robert G. Burns,\\
Springer - New York, 1992
}%

\BiblioItem{Korn}
{
Granino A. Korn, Theresa M. Korn,
Mathematical Handbook for Scientists and Engineer,
McGraw-Hill Book Company, New York, San Francisco,
Toronto, London, Sydney, 1968
}%

\BiblioItem{Hocking Young Topology}
{
John G. Hocking, Gail S. Young,
Topology,\\
Courier Dover Publications, 1988
}%

\BiblioItem{Olver: Lie groups to differential equations}
{
Peter J. Olver,
Applications of Lie groups to differential equations,\\
Springer, 2000
}%

\BiblioItem{Tartaglia}
{
Angelo Tartaglia and Matteo Luca Ruggiero,
Angular Momentum Effects in Michelson\Hyph Morley Type Experiments,
Gen.Rel.Grav. 34, 1371-1382 (2002),\\
eprint \href{http://arxiv.org/abs/gr-qc/0110015}{arXiv:gr-qc/0110015} (2001)
}%

\BiblioItem{Tomozawa}
{
Yukio Tomozawa, Speed of Light in Gravitational Fields, eprint
\href{http://arxiv.org/abs/astro-ph/0303047}{arXiv:astro-ph/0303047} (2004)
}%

\BiblioItem{Magueijo}
{
Joao Magueijo,
Covariant and locally Lorentz-invariant varying speed of light theories,
\href{http://prola.aps.org/abstract/PRD/v62/i10/e103521}{Phys. Rev. D 62, 103521} (2000),
eprint \href{http://arxiv.org/abs/gr-qc/0007036}{arXiv:gr-qc/0007036} (2000)
}%

\BiblioItem{Bassett}
{
Bruce A. Bassett, Stefano Liberati, Carmen Molina-Paris, and Matt Visser,
Geometrodynamics of variable-speed-of-light cosmologies,
\href{http://prola.aps.org/abstract/PRD/v62/i10/e103518}{Phys. Rev. D 62}, 103518 (2000),
eprint \href{http://arxiv.org/abs/astro-ph/0001441}{arXiv:astro-ph/0001441} (2000)
}%

\BiblioItem{C.A. Deavours The Quaternion Calculus}
{
C.A. Deavours, The Quaternion Calculus, 
American Mathematical Monthly, {\bf 80} (1973), pp. 995 - 1008
}%

\BiblioItem{Straumann}
{
Lochlainn O'Raifeartaigh and Norbert Straumann,
Gauge theory: Historical origins and some modern developments,
\href{http://prola.aps.org/abstract/RMP/v72/i1/p1_1}{Rev. Mod. Phys. 72, 1} (2000)
}%

\BiblioItem{Lammerzahl}
{
Claus L\"ammerzahl, Mark P. Haugan,
On the interpretation of Michelson\Hyph Morley experiments,
%\href{http://www.sciencedirect.com/science?\_ob=ArticleURL&\_udi=B6TVM-42WP7CR-1&\_user=10&\_handle=W-WA-A-A-AZ-MsSAYZW-UUW-AUDDYZYZAU-WZCBYCEDW-AZ-U&\_fmt=summary&\_coverDate=04%2F23%2F2001&\_rdoc=1&\_orig=browse&\_srch=%23toc%235538%232001%23997179995%23246657!&\_cdi=5538&view=c&\_acct=C000050221&\_version=1&\_urlVersion=0&\_userid=10&md5=385478cda8c5568dea1aeaf0c43669da}
{Phys. Lett. A282 223-229} (2001),\\
eprint \href{http://arxiv.org/abs/gr-qc/0103052}{arXiv:gr-qc/0103052} (2001)
}%

\BiblioItem{0305117}
{
Holger Mueller, Sven Herrmann, Claus Braxmaier, Stephan Schiller, Achim Peters.
Modern Michelson-Morley Experiment using Cryogenic Optical Resonators.
eprint \href{http://arxiv.org/abs/physics/0305117}{arXiv:physics/0305117} (2003)
\\
Phys. Rev. Lett. 91:020401, 2003
}%

\BiblioItem{0706.2031}
{
Holger Mueller, Paul Louis Stanwix, Michael Edmund Tobar,
Eugene Ivanov, Peter Wolf, Sven Herrmann, Alexander Senger,
Evgeny Kovalchuk, Achim Peters.
Relativity tests by complementary rotating Michelson-Morley experiments.
eprint \href{http://arxiv.org/abs/0706.2031}{arXiv:0706.2031 [physics.class-ph]} (2006)
\\
Phys. Rev. Lett. 99:050401, 2007
}%

\BiblioItem{1008.1205}
{
M. Nagel, K. M\"ohle, K. D\"oringshoff, S. Herrmann, A. Senger, E.V. Kovalchuk, A. Peters.
Testing Lorentz Invariance by Comparing Light Propagation in Vacuum and Matter.
eprint \href{http://arxiv.org/abs/1008.1205}{arXiv:1008.1205 [physics.ins-det]} (2010)
}%

\BiblioItem{1109.4897}
{
The OPERA Collaboration.
Measurement of the neutrino velocity with the OPERA detector in the CNGS beam.
eprint \href{http://arxiv.org/abs/1109.4897}{arXiv:1109.4897 [hep-ex]} (2011)
}%

\BiblioItem{Ranada}
{
Antonio F. Ranada,
Pioneer acceleration and variation of light speed: experimental situation,
eprint \href{http://arxiv.org/abs/gr-qc/0402120}{arXiv:gr-qc/0402120} (2004)
}%

\BiblioItem{Gelfand Minlos: rotation and Lorentz groups}
{
Izrail Moiseevich Gelfand, Robert Adolfovich Minlos,
Representations of the rotation and Lorentz groups and their applications;\\
Engl. transl. ed. H. K. Farahat; Transl. by G. Cummins and T. Boddongton;\\
Pergamon Press, 1963
}%

\BiblioItem{math.QA-0208146}
{
I. Gelfand, S. Gelfand, V. Retakh, R. Wilson,
Quasideterminants,\\
eprint \href{http://arxiv.org/abs/math.QA/0208146}{arXiv:math.QA/0208146} (2002)
}%

\BiblioItem{q-alg-9705026}
{
I.Gelfand, V.Retakh,
Quasideterminants, I,\\
eprint \href{http://arxiv.org/abs/q-alg/9705026}{arXiv:q-alg/9705026} (1997)
}%

\BiblioItem{Gelfand Retakh 1991}
{
I. Gelfand and V. Retakh, Determinants of Matrices over Noncommutative Rings, Funct.
Anal. Appl. 25 (1991), no. 2, 91-102
}%

\BiblioItem{Gelfand Retakh 1992}
{
I. Gelfand and V. Retakh, A Theory of Noncommutative Determinants and Characteristic
Functions of Graphs, Funct. Anal. Appl. 26 (1992), no. 4, 1-20
}%

\BiblioItem{hep-th-9407124}
{
I. M. Gelfand, D. Krob, A. Lascoux, B. Leclerc, V.S. Retakh and J.-Y. Thibon,
Noncommutative symmetric functions,\\
eprint \href{http://arxiv.org/abs/hep-th/9407124}{arXiv:hep-th/9407124} (1994)
}%

\BiblioItem{Naimark Shtern: Theory of group representations}
{
Mark Aronovich Naimark, Aleksandr Isaakovich Shtern,
Theory of group representations;\\
Heidelberg, 1982
}%

\BiblioItem{Barut Raczka: Theory of group representations}
{
Asim Orhan Barut; Ryszard R\c{a}czka;
Theory of group representations and applications;\\
World Scientific Publishing Co. Pre. Ltd., 1986
}%

\BiblioItem{Mihalev Pilz: concise handbook of algebra}
{
Aleksandr Vasilevich Mikhalev; G\"{u}nter Pilz;
The concise handbook of algebra;\\
Kluwer Academic Publishers, 2002
}%

\BiblioItem{Shafarevich: Basic notions of algebra}
{
I. R. Shafarevich,
Basic notions of algebra,\\
Translated from the Russian by M. Reid,\\
Springer, 2005
}%

\BiblioItem{Elsgolts: Differential Equations}
{
Lev Elsgolts,
Differential Equations and the Calculus of Variations,\\
University Press of the Pacific, 2003 
}%

\BiblioItem{Baez Huerta: algebra of grand unified theories}
{
John Baez; John Huerta;
The algebra of grand unified theories;\\
Bull. Amer. Math. Soc. {\bf 47} (2010), 483-552
}%

\BiblioItem{J. Fan: Determinants}
{
J. Fan, Determinants and multiplicative functionals
on quaternion matrices, \textit{Linear Algebra and Its
Applications} {\bf 369}, pp. 193-201, (2003).
}%

\BiblioItem{Carl Faith 1}
{
Carl Faith, Algebra: Rings, Modules and Categories I,
Springer - Verlag, Berlin - Heidelberg - New York, 1973
}%

\BiblioItem{Gilson Nimmo Ohta}
{
 C.R.Gilson, J.J.C.Nimmo, Y.Ohta, Quasideterminant solutions of a non-Abelian Hirota-Miwa
 equation, \textit{Journal of Physics A: Mathematical and Theoretical} {\bf 40}(42), pp.
 12607-12617,(2007).
}%

\BiblioItem{Haider Hassan}
{
B. Haider, M. Hassan, Quasideterminant solutions of an integrable chiral model in two
 dimensions, \textit{Journal of Physics A: Mathematical and Theoretical} {\bf 42} (35), art. no.
 355211, (2009).
}%

%\BiblioItem{Pareigis}
%{
%Bodo Pareigis, Categories and Functors,
%Academic Press - New York - London, 1970
%}%

%\BiblioItem{Beachy}
%{
%John A. Beachy, Introductory Lectures on Rings i Modules,
%Cambridge University Press, 1999
%}%

\BiblioItem{0702447}
{
I.I. Kyrchei, Cramer's rule for quaternion systems of linear equations,
\textit{Journal of Mathematical Sciences} {\bf 155}(6), 839-858, (2008).
 Translated from  \textit{Fundamental and Appl. Math.}
 {\bf 13}(4), pp.67-94, (2007). (in Russian)\\
eprint
\href{http://arxiv.org/abs/math/0702447}{arXiv:math.RA/0702447}
(2007)
}%

\BiblioItem{1004.4380}
{
I.I. Kyrchei, Cramer's rule for some quaternion matrix
    equations,  \textit{Applied Mathematics and Computation} {\bf 217}(5), pp.2024-2030, (2010).\\eprint
\href{http://arxiv.org/abs/1004.4380
}{arXiv:math.RA/arXiv:1004.4380 } (2010)
}%

\BiblioItem{1005.0736}
{
I.I. Kyrchei,Determinantal representations of the Moore-Penrose inverse
 over the quaternion skew field and corresponding Cramer's rules,
 \\
eprint
\href{http://arxiv.org/abs/1005.0736}{arXiv:math.RA/1005.0736}
(2010)
}%

\BiblioItem{0412.391}
{
Aleks Kleyn,
Basis Manifold,
eprint \href{http://arxiv.org/abs/math.DG/0412391}{arXiv:math.DG/0412391} (2007)
}%

\BiblioItem{0405.027}
{
Aleks Kleyn,
Reference Frame in General Relativity,\\
eprint \href{http://arxiv.org/abs/gr-qc/0405027}{arXiv:gr-qc/0405027} (2008)
}%

\BiblioItem{0405.028}
{
Aleks Kleyn, Metric\hyph Affine Manifold,\\
eprint \href{http://arxiv.org/abs/gr-qc/0405028}{arXiv:gr-qc/0405028} (2008)
}%

\BiblioItem{0612.111}
{
Aleks Kleyn,
Biring of Matrices,\\
eprint \href{http://arxiv.org/abs/math.OA/0612111}{arXiv:math.OA/0612111} (2007)
}%

\BiblioItem{0701.238}
{
Aleks Kleyn,
Lectures on Linear Algebra over Division Ring,\\
eprint \href{http://arxiv.org/abs/math.GM/0701238}{arXiv:math.GM/0701238} (2010)
}%

\BiblioItem{0702.561}
{
Aleks Kleyn,
Fibered $\mathfrak{F}$\Hyph Algebra,\\
eprint \href{http://arxiv.org/abs/math.DG/0702561}{arXiv:math.DG/0702561} (2007)
}%

\BiblioItem{math.RA-0501237}
{
Aleks Kleyn,
Vector Space Over Division Ring,\\
eprint \href{http://arxiv.org/abs/math.RA/0412391}{arXiv:math.RA/0501237} (2007)
}%

\BiblioItem{math.RA-0501237v1}
{
Aleks Kleyn,
Module Over Division Ring, version 1,\\
eprint \href{http://arxiv.org/abs/math/0501237v1}{arXiv:math.RA/0501237v1} (2005)
}%

\BiblioItem{0707.2246}
{
Aleks Kleyn,
Fibered Correspondence,\\
eprint \href{http://arxiv.org/abs/0707.2246}{arXiv:0707.2246} (2007)
}%

\BiblioItem{0803.2620}
{
Aleks Kleyn,
Morphism of \Ts Representations,\\
eprint \href{http://arxiv.org/abs/0803.2620}{arXiv:0803.2620} (2008)
}%

\BiblioItem{0803.3276}
{
Aleks Kleyn,
Lorentz Transformation and General Covariance Principle,\\
eprint \href{http://arxiv.org/abs/0803.3276}{arXiv:0803.3276} (2009)
}%

\BiblioItem{0812.4763}
{
Aleks Kleyn,
Introduction into Calculus over Division Ring,\\
eprint \href{http://arxiv.org/abs/0812.4763}{arXiv:0812.4763} (2010)
}%

\BiblioItem{0906.0135}
{
Aleks Kleyn,
Introduction into Geometry over Division Ring,\\
eprint \href{http://arxiv.org/abs/0906.0135}{arXiv:0906.0135} (2010)
}%

\BiblioItem{0909.0855}
{
Aleks Kleyn,
Quaternion Rhapsody,\\
eprint \href{http://arxiv.org/abs/0909.0855}{arXiv:0909.0855} (2010)
}%

\BiblioItem{0912.3315}
{
Aleks Kleyn,
Representation of Universal Algebra,\\
eprint \href{http://arxiv.org/abs/0912.3315}{arXiv:0912.3315} (2009)
}%

\BiblioItem{0912.4061}
{
Aleks Kleyn,
Linear Equation in Finite Dimensional Algebra,\\
eprint \href{http://arxiv.org/abs/0912.4061}{arXiv:0912.4061} (2010)
}%

\BiblioItem{1001.4852}
{
Aleks Kleyn,
The Matrix of Linear Mappings,\\
eprint \href{http://arxiv.org/abs/1001.4852}{arXiv:1001.4852} (2010)
}%

\BiblioItem{1003.1544}
{
Aleks Kleyn,
Linear Mappings of Free Algebra,\\
eprint \href{http://arxiv.org/abs/1003.1544}{arXiv:1003.1544} (2010)
}%

\BiblioItem{1011.3102}
{
Aleks Kleyn,
Polylinear Mapping of Free Algebra,\\
eprint \href{http://arxiv.org/abs/1011.3102}{arXiv:1011.3102} (2010)
}%

\BiblioItem{1104.5197}
{
Aleks Kleyn,
$C^*$-Rhapsody,\\
eprint \href{http://arxiv.org/abs/1104.5197}{arXiv:1104.5197} (2011)
}%

\BiblioItem{1105.4307}
{
Aleks Kleyn,
Algebra with Conjugation,\\
eprint \href{http://arxiv.org/abs/1105.4307}{arXiv:1105.4307} (2011)
}%

\BiblioItem{1107.1139}
{
Aleks Kleyn,
Linear Mappings of Quaternion Algebra,\\
eprint \href{http://arxiv.org/abs/1107.1139}{arXiv:1107.1139} (2011)
}%

\BiblioItem{1107.5037}
{
Aleks Kleyn,
Orthogonal Basis and Motion in Finsler Geometry,\\
eprint \href{http://arxiv.org/abs/1107.5037}{arXiv:1107.5037} (2011)
}%

\BiblioItem{8433-5163}
{
Aleks Kleyn,
Linear Mappings of Free Algebra: First Steps in Noncommutative Linear Algebra,\\
Lambert Academic Publishing, 2010
}%

\BiblioItem{8443-0072}
{
Aleks Kleyn,
Representation Theory: Representation of Universal Algebra,\\
Lambert Academic Publishing, 2011
}%

\BiblioItem{Lauve: Quantum coordinates}
{
A. Lauve, Quantum- and quasi-Plucker coordinates,
\textit{Journal of Algebra} {\bf 296}(2), pp.440-461,
(2006).
}%

\BiblioItem{Lewis D. W. Quaternion algebras}
{
Lewis D. W. Quaternion algebras and the algebraic legacy
of Hamilton's quaternions, \textit{Irish Math. Soc. Bulletin} {\bf
57}, pp. 41-64, (2006).
}%

\BiblioItem{0812.2865}
{
Jos\'e Miguel Figueroa-O'Farrill,
Three lectures on 3-algebras,
eprint \href{http://arxiv.org/abs/0812.2865}{arXiv:0812.2865} (2008)
}%

\BiblioItem{Li Nimmo: Darboux transformations}
{
C.X.Li, J.J.C. Nimmo, Darboux transformations for a twisted
derivation and quasideterminant solutions to the super KdV
equation, \textit{Proceedings of the Royal Society A:
Mathematical, Physical and Engineering Sciences} {\bf 466} (2120),
pp. 2471-2493, (2010)
}%

\BiblioItem{Schiebold: Cauchy-type determinants}
{
C. Schiebold, Cauchy-type determinants and integrable
systems, \textit{Linear Algebra and Its Applications} {\bf 433}
(2), pp. 447-475, (2010)
}%

\BiblioItem{Suzuki: Noncommutative spectral decomposition}
{
T. Suzuki, Noncommutative
spectral decomposition with qua\-si\-de\-ter\-mi\-nant, \textit{Advances in
Mathematics} {\bf 217}(5), pp. 2141-2158, (2008)
}%

\BiblioItem{1105.3456}
{
C. W. F. Everitt, D. B. DeBra, B. W. Parkinson, J. P. Turneaure, J. W. Conklin,
M. I. Heifetz, G. M. Keiser, A. S. Silbergleit, T. Holmes, J. Kolodziejczak,
M. Al-Meshari, J. C. Mester, B. Muhlfelder, V. Solomonik, K. Stahl, P. Worden,
W. Bencze, S. Buchman, B. Clarke, A. Al-Jadaan, H. Al-Jibreen, J. Li, J. A. Lipa,
J. M. Lockhart, B. Al-Suwaidan, M. Taber, S. Wang,\\
Gravity Probe B: Final Results of a Space Experiment to Test General Relativity,\\
eprint \href{http://arxiv.org/abs/1105.3456}{arXiv:1105.3456[gr-qc]} (2011)
}%

\BiblioItem{0009305}
{
G. S. Asanov.
Can Neutrinos and High-Energy Particles Test Finsler Metric of Space-Time?\\
eprint \href{http://arxiv.org/abs/hep-ph/0009305}{arXiv:hep-ph/0009305} (2000)
}%

\BiblioItem{Asanov 2004}
{
G. S. Asanov.
Finsleroid - space supplemented by angle and scalar product.\\
Hypercomplex Numbers in Geometry and Physics, {\bf 1}, 2004, p. 40 - 62
}%

\BiblioItem{1004.3007}
{
Sergiu I. Vacaru,
Principles of Einstein-Finsler Gravity and Perspectives in Modern Cosmology,\\
eprint \href{http://arxiv.org/abs/1004.3007}{arXiv:1004.3007[math-ph]} (2010)
}%

\BiblioItem{1012.4148}
{
Sergiu I. Vacaru.
Principles of Einstein-Finsler Gravity and Cosmology.\\
eprint \href{http://arxiv.org/abs/1012.4148}{arXiv:1012.4148[physics.gen-ph]} (2010)
}%

\BiblioItem{1112.5641}
{
Christian Pfeifer, Mattias N.R. Wohlfarth.
Finsler geometric extension of Einstein gravity.\\
eprint \href{http://arxiv.org/abs/1112.5641}{arXiv:1112.5641[gr-qc]} (2011)
}%

\BiblioItem{0711.0056}
{
Zhe Chang, Xin Li.
Lorentz Invariance Violation and Symmetry in Randers\Hyph Finsler Spaces.\\
eprint \href{http://arxiv.org/abs/0711.0056}{arXiv:0711.0056[hep-th]} (2011)
}%

\BiblioItem{Rund Finsler geometry}
{
Hanno Rund,
The differential geometry of Finsler spaces.
\\
Springer - Verlag, Berlin - G\"ottingen - Heidelberg, 1959
}%

\BiblioItem{Beem Dostoglou Ehrlich}
{
John K. Beem, Stamatis A. Dostoglou, Paul E. Ehrlich,
Advances in differential geometry and general relativity.
\\
American Mathematical Society, 2004
}%

\BiblioItem{978-0719033414}
{
Malcolm Pemberton, Nicholas Rau,
Mathematics for economists: an introductory textbook.
\\
Manchester University Press, November 2001; ISBN-13: 978-0719033414
}%

\BiblioItem{0 521 59180 5}
{
Cyrus D. Cantrell,
Modern mathematical methods for physicists and engineers.
\\
Cambridge University Press, 2000
}%

\BiblioItem{Arveson spectral theory}
{
William Arveson,
A short course on spectral theory.
\\
Springer - Verlag, New York, 2002
}%

\BiblioItem{Robert Hermann}
{
Robert Hermann,
Topics in the mathematics of quantum mechanics.
\\
Math Sci Press, 1973
}%

\BiblioItem{9705.009}%q-alg-9705009
{
John C. Baez,
An Introduction to n-Categories,\\
eprint \href{http://arxiv.org/abs/q-alg/9705009}{arXiv:q-alg/9705009} (1997)
}%

\BiblioItem{0105.155}
{
John C. Baez,
The Octonions,\\
eprint \href{http://arxiv.org/abs/math.RA/0105155}{arXiv:math.RA/0105155} (2002)
}%

\BiblioItem{John Baez: Math Blogs}
{
John C. Baez,
What do mathematicians need to know about blogging?,\\
Notices of the American Mathematical Society,
(2010), 3, {\bf 57}, 333,\\
\url{http://www.ams.org/notices/201003/rtx100300333p.pdf}
}%

\BiblioItem{Tolstoi about Anna Karenina}
{
Tolstoi about Anna Karenina,
in book A Karenina Companion, by C. J. G. Turner,
published by Wilfrid Laurier University Press (August 1993)
}%

\BiblioItem
{Cohn: Universal Algebra}
{
Paul M. Cohn,
Universal Algebra,
Springer, 1981
}%

\BiblioItem
{Maunder: Algebraic Topology}
{
C. R. F. Maunder,
Algebraic Topology,
Dover Publications, Inc, Mineola, New York, 1996
}%

\BiblioItem{Pommaret: Partial Differential Equations}
{
J.-F. Pommaret,
Partial Differential Equations and Group Theory,
Springer, 1994
}%

\BiblioItem{Bourbaki: Set Theory}
{
N. Bourbaki,
Theory of sets,
Springer, 2004
}%

\BiblioItem{Bourbaki: Algebra 1}
{
N. Bourbaki,
Algebra 1,
Springer, 2004
}%

\BiblioItem
{Bourbaki: General Topology 1}
{
N. Bourbaki,
General Topology, Chapters 1 - 4,
Springer, 1989
}

\BiblioItem{Bourbaki: General Topology: Chapter 5 - 10}
{
N. Bourbaki,
General Topology, Chapters 5 - 10,
Springer, 1989
}

\BiblioItem{Bourbaki: Topological Vector Space}
{
N. Bourbaki,
Topological Vector Spaces, Chapters 1 - 5,
Transl. by H. G. Eggleston $\&$ S. Madan,
Springer, 2003
}

\BiblioItem{Bourbaki: Real Group Lie}
{
N. Bourbaki,
Lie Groups and Lie Algebras, Chapters 7 - 9,
Translator Andrew Pressley,
Springer, 2005
}

\BiblioItem{Shabat: Complex Analysis}
{
Shabat B. V.,
Introduction to Complex Analysis,
\\ \url{http://www.math.uchicago.edu/~ryzhik/shabat-all.pdf},
\\Translated from Russian by L.Ryzhik, 2003
(Moscow, Nauka, 1969)
}

\BiblioItem{Pontryagin: Topological Group}
{
L. S. Pontryagin,
Selected Works, Volume Two, Topological Groups,
Gordon and Breach Science Publishers, 1986
}

\BiblioItem
{Eisenhart: Riemannian Geometry}
{
Eisenhart,
Riemannian Geometry,
Princeton University Press, Princeton, 1949
}

\BiblioItem
{Eisenhart: Continuous Groups of Transformations}
{
Eisenhart,
Continuous Groups of Transformations,
Dover Publications, New York, 1961
}

\BiblioItem
{Condon Odabasi}
{
Edward Uhler Condon, Halis Odabasi,
Atomic Structure,
CUP Archive, 1980
}

\BiblioItem{Postnikov: Differential Geometry}
{
Postnikov M. M.,
Geometry IV: Differential geometry,
Moscow, Nauka, 1983
}

\BiblioItem{Fihtengolts: Calculus volume 1}
{
Fihtengolts G. M.,
Differential and Integral Calculus Course, volume 1,
Moscow, Nauka, 1969
}

\BiblioItem{Hatcher: Algebraic Topology}
{
Allen Hatcher,
Algebraic Topology,
Cambridge University Press, 2002
}

\BiblioItem{geometry of differential equations}
{
Vinogradov, A. M., Krasil'shchik, I. S., and Lychagin, V. V.,
Introduction to geometry of nonlinear differential equations,
Nauka, Moscow, 1986
}

\BiblioItem{cohomological analysis}
{
A. M. Vinogradov,
Cohomological Analysis of Partial Differential Equations
and Secondary Calculus,
American Mathematical Society, 2001
}

\BiblioItem{0801.1734}
{
Brandon S. DiNunno, Richard A. Matzner,
The Volume Inside a Black Hole,\\
eprint \href{http://arxiv.org/abs/0801.1734v1}{arXiv:0801.1734v1} (2008)
}

\BiblioItem{0702.447}
{
Ivan Kyrchei,
Cramer's rule for some quaternion matrix equations,\\
eprint \href{http://arxiv.org/abs/math/0702447}{arXiv:math.RA/0702447} (2007)
}

\BiblioItem{Izrail M. Gelfand: Quaternion Groups}
{
I. M. Gelfand, M. I. Graev,
Representation of Quaternion Groups over Localy Compact and
Functional Fields,\\
Funct. Anal. Appl. {\bf 2} (1968) 19 - 33;\\
Izrail Moiseevich Gelfand, Semen Grigorevich Gindikin,\\
Izrail M. Gelfand: Collected Papers, volume II, 435 - 449,\\
Springer, 1989
}

\BiblioItem{Richard D. Schafer}
{
Richard D. Schafer,
An Introduction to Nonassociative Algebras,
Dover Publications, Inc., New York, 1995
}

\BiblioItem{Bamberg Sternberg}
{
Paul Bamberg, Shlomo Sternberg,
A course in mathematics for students of physics,
Cambridge University Press, 1991
}

\BiblioItem{Conway Smith}
{
John Horton Conway, Derek Alan Smith,
On quaternions and octonions: their geometry, arithmetic, and symmetry,
A K Peters, Natick, Massachussets, 2003
}

\BiblioItem{Fueter}
{
Fueter, R.
Die Funktionentheorie der Differentialgleichungen $\Delta u = 0$ und
$\Delta \Delta u = 0$ mit vier reellen Variablen.
Comment. Math. Helv. {\bf 7} (1935), 307-330
}

\BiblioItem{Sudbery Quaternionic Analysis}
{
A. Sudbery,
Quaternionic Analysis,
Math. Proc. Camb. Phil. Soc. (1979), {\bf 85}, 199 - 225
}

\BiblioItem{0902.4771}
{
Fabrizio Colombo, Graziano Gentili, Irene Sabadini,
A Cauchy kernel for slice regular functions,\\
eprint \href{http://arxiv.org/abs/0902.4771v1}{arXiv:0902.4771v1} (2009)
}

\BiblioItem{Vadim Komkov}
{
Vadim Komkov,
Variational Principles of Continuum Mechanics with Engineering Applications: Critical Points Theory,\\
Springer, 1986
}

\BiblioItem{Alain Connes 1994}
{
Alain Connes,
Noncommutative Geometry,\\
Academic Press, 1994
}

\BiblioItem{Hamilton papers 3}
{
Sir William Rowan Hamilton,
The Mathematical Papers, Vol. III, Algebra,\\
Cambridge at the University Press, 1967
}

\BiblioItem{Hamilton Elements of Quaternions 1}
{
Sir William Rowan Hamilton,
Elements of Quaternions, Volume I,\\
Longmans, Green, and Co., London, New York, and Bombay, 1899
}

\BiblioItem{Cartan geometry in reper}
{
Elie Cartan, Vladislav V. Goldberg, Serge\u{i} Pavlovich Finikov,\\
Riemannian geometry in an orthogonal frame:
from lectures delivered by Elie Cartan at the Sorbonne in 1926-1927,\\
translated by Vladislav V. Goldberg,\\
World Scientific, 2001
}

\BiblioItem{Arnautov Glavatsky Mikhalev}
{
V. I. Arnautov, S. T. Glavatsky, A. V. Mikhalev,\\
Introduction to the theory of topological rings and modules,
Volume 1995,\\
Marcel Dekker, Inc, 1996
}

\BiblioItem{Moore Yaqub}
{
Hal G. Moore, Adil Yaqub,
A first course in linear algebra with applications,
Edition 3, Academic Press, 1998 
}

\BiblioItem{math.CV-0405471}
{
S. V. Ludkovsky,
Differentiable functions of Cayley-Dickson numbers,\\
eprint \href{http://arxiv.org/abs/math.CV/0405471}{arXiv:math.CV/0405471} (2004)
}%

\BiblioItem{W.Bertram H.Glockner K.Neeb}
{
W.Bertram, H.Glockner, K.Neeb,
Differential Calculus over General Base Fields and Rings,
Expositiones Mathematicae (2004), Volume 22, Issue 3, Pages 213-282
}

\CloseBiblio

%% file: Index.English.tex
%auto-ignore
\OpenIndex
\SetIndexSpace%
\Index%1%1%1-rcd form
   {$1$-\rcd form}%
   {1-rcd form, vector spaces}%
\SetIndexSpace%
\Index%2%2%2- ary fibered relation
   {$2$\Hyph ary fibered relation}%
   {2 ary fibered relation}%
\SetIndexSpace%
\Index%3%A%a b-CR quasideterminant
   {$(^a_b)$\hyph \CR quasideterminant}%
   {a b cr-quasideterminant}%
\Index%668%A%A CRcirc- basis for module
   {$A\CRcirc$\Hyph basis for module}%
   {A CRcirc basis, module over algebra}%
\Index%667%A%A CRcirc- linearly dependent set of vectors
   {$A\CRcirc$\Hyph linearly dependent set of vectors}%
   {CRcirc linearly dependent, Astar module over algebra}%
\Index%666%A%A CRcirc- linearly independent set of vectors
   {$A\CRcirc$\Hyph linearly independent set of vectors}%
   {CRcirc linearly independent, Astar module over algebra}%
\Index%693%A%A linear mapping of modules
   {$A$\Hyph linear mapping of modules}%
   {A linear map of modules}%
\Index%603%A%A(A) mapping
   {$\mathcal A(A)$\Hyph mapping}%
   {A(A) mapping}%
\Index%685%A%A- module
   {$A$\Hyph module}%
   {module over algebra}%
\Index%5%A%A- valued function
   {$A$\Hyph valued function}%
   {A valued function}%
\Index%645%A%Abelian Omega group
   {Abelian $\Omega$\Hyph group}%
   {Abelian Omega group}%
\Index%284%A%absolute value on division ring
   {absolute value on division ring}%
   {absolute value on division ring}%
\Index%793%A%absolute value on ring
   {absolute value on ring}%
   {absolute value on ring}%
\Index%680%A%Acr linear mapping of modules
   {\Acr linear mapping of modules}%
   {Acr linear map of modules}%
\Index%764%A%ACRcirc linear combination
   {$A\CRcirc$\Hyph linear combination}%
   {ACRcirc linear combination}%
\Index%100%A%active representation
   {active representation}%
   {active representation}%
\Index%569%A%active representation of group in basis manifold of representation
   {active representation of group $G(f)$ in basis manifold of representation}%
   {active representation in basis manifold}%
\Index%716%A%active representation of group in basis manifold of tower of representations
   {active representation of group $G(\Vector f)$ in basis manifold of tower of representations}%
   {active representation in basis manifold, tower of representations}%
\Index%99%A%active starT- representation
   {active \sT representation}%
   {active representation, vector space}%
\Index%326%A%active transformation of basis manifold of representation
   {active transformation of basis manifold of representation}%
   {active transformation of basis, representation}%
\Index%715%A%active transformation of basis manifold of tower of representations
   {active transformation of basis manifold of tower of representations}%
   {active transformation of basis, tower of representations}%
\Index%101%A%active transformation on basis manifold
   {active transformation on basis manifold}%
   {active transformation}%
\Index%102%A%active transformation on the set of rcd bases
   {active transformation on the set of \rcd bases}%
   {active transformation, vector space}%
\Index%110%A%affine basis
   {affine basis}%
   {Affine Basis}%
\Index%733%A%affine structure on set
   {affine structure on set}%
   {affine structure on set}%
\Index%141%A%affine transformation group
   {affine transformation group}%
   {drc affine transformation group}%
\Index%142%A%affine transformation group
   {affine transformation group}%
   {affine transformation group}%
\Index%109%A%affine transformation on basis manifold
   {affine transformation on basis manifold}%
   {affine transformation}%
\Index%607%A%algebra over ring
   {algebra over ring}%
   {algebra over ring}%
\Index%780%A%algebra with conjugation
   {algebra with conjugation}%
   {algebra with conjugation}%
\Index%105%A%alternative representation of matrix
   {alternative representation of matrix}%
   {Alternative representation}%
\Index%275%A%anholonomic coordinate
   {anholonomic coordinate}%
   {anholonomic coordinate}%
\Index%278%A%anholonomic coordinates of connection
   {anholonomic coordinates of connection}%
   {anholonomic coordinates of connection}%
\Index%276%A%anholonomic coordinates of vector
   {anholonomic coordinates of vector}%
   {vector anholonomic coordinates}%
\Index%277%A%anholonomic coordinates on manifold
   {anholonomic coordinates on manifold}%
   {anholonomic coordinates on manifold}%
\Index%293%A%anholonomity object
   {anholonomity object}%
   {anholonomity object}%
\Index%106%A%antihomomorphism of fibered groups
   {antihomomorphism of fibered groups}%
   {antihomomorphism of fibered groups}%
\Index%107%A%antisymmetric 2- ary fibered relation
   {antisymmetric $2$\Hyph ary fibered relation}%
   {antisymmetric 2 ary fibered relation}%
\Index%582%A%Arc basis for vector space
   {$A\RCstar$\Hyph basis for vector space}%
   {Arc basis, vector space}%
\Index%581%A%Arc linearly dependent vectors
   {$A\RCstar$\Hyph linearly dependent vectors}%
   {linearly dependent, A vector space}%
\Index%580%A%Arc linearly independent vectors
   {$A\RCstar$\Hyph linearly independent vectors}%
   {linearly independent, A vector space}%
\Index%108%A%arity of operation
   {arity of operation}%
   {arity of operation}%
\Index%617%A%associative D algebra
   {associative $D$\Hyph algebra}%
   {associative D algebra}%
\Index%585%A%associative law for Astar linear mappings of vector spaces
   {associative law for $A\star$\Hyph linear mappings of vector spaces}%
   {associative law for Astar linear mappings of vector spaces}%
\Index%662%A%associative law for Astar- module
   {associative law for $A\star$\Hyph module}%
   {associative law, Astar module over algebra}%
\Index%576%A%associative law for Astar- vector space
   {associative law for $A\star$\Hyph vector space}%
   {associative law, Astar vector space}%
\Index%165%A%associative law for covariant starT- representation
   {associative law for covariant \sT representation}%
   {associative law for covariant starT representation}%
\Index%166%A%associative law for covariant Tstar- representation
   {associative law for covariant \Ts representation}%
   {associative law for covariant Tstar representation}%
\Index%162%A%associative law for Drc linear maps of vector bundles
   {associative law for \Drc linear maps of vector bundles}%
   {associative law for drc linear maps of vector bundles}%
\Index%654%A%associative law for Dstar- module
   {associative law for $D\star$\Hyph module}%
   {associative law, Dstar module}%
\Index%164%A%associative law for Dstar- vector fields
   {associative law for $\mathcal D\star$\Hyph vector fields}%
   {associative law, Dstar vector fields}%
\Index%163%A%associative law for Dstar- vector space
   {associative law for $D\star$\Hyph vector space}%
   {associative law, Dstar vector space}%
\Index%161%A%associative law for rcd linear maps of vector spaces
   {associative law for \rcd linear maps of vector spaces}%
   {associative law for rcd linear maps of vector spaces}%
\Index%671%A%associative law for starA- module
   {associative law for $\star A$\Hyph module}%
   {associative law, starA module over algebra}%
\Index%167%A%associative law for twin representations
   {associative law for twin representations}%
   {associative law for twin representations}%
\Index%168%A%associative law of composition of fibered correspondences
   {associative law of composition of fibered correspondences}%
   {associative law, composition of fibered correspondences}%
\Index%616%A%associator of R algebra
   {associator of $R$\Hyph algebra}%
   {associator of algebra}%
\Index%789%A%Astar antilinear mapping of algebra with conjugation
   {$A\star$\Hyph antilinear mapping of algebra with conjugation}%
   {antilinear mapping of algebra with conjugation}%
\Index%583%A%Astar linear map of vector spaces
   {$A\star$\Hyph linear map of vector spaces}%
   {Astar linear map of vector spaces}%
\Index%575%A%Astar vector space
   {$A\star$\Hyph vector space}%
   {Astar vector space}%
\Index%660%A%Astar- module
   {$A\star$\Hyph module}%
   {Astar-module}%
\Index%665%A%Astar- product of vector over scalar
   {$A\star$\Hyph product of vector over scalar}%
   {Astar product of vector over scalar, Astar module}%
\Index%579%A%Astar-product of vector over scalar
   {$A\star$\hyph product of vector over scalar}%
   {Astar product of vector over scalar, vector space}%
\Index%95%A%auto parallel line
   {auto parallel line}%
   {auto parallel line}%
\Index%94%A%automorphism of representation of Omega algebra
   {automorphism of representation of $\Omega$\Hyph algebra}%
   {automorphism of representation}%
\Index%560%A%automorphism of tower of representations
   {automorphism of tower of representations}%
   {automorphism of tower of representations}%
\Index%282%A%norm of quaternion
   {norm of quaternion}%
   {norm of quaternion}%
\SetIndexSpace%
\Index%503%B%Banach D- algebra
   {Banach $D$\Hyph algebra}%
   {Banach algebra}%
\Index%114%B%base of fibered correspondence
   {base of fibered correspondence}%
   {base of fibered correspondence}%
\Index%113%B%base of mapping
   {base of mapping}%
   {base of map}%
\Index%256%B%basis manifold of affine space
   {basis manifold of affine space}%
   {Basis Manifold, Affine Space}%
\Index%261%B%basis manifold of central affine space
   {basis manifold of central affine space}%
   {Basis Manifold, Central Affine Space, division ring}%
\Index%262%B%basis manifold of central affine space
   {basis manifold of central affine space}%
   {Basis Manifold, Central Affine Space}%
\Index%259%B%basis manifold of Euclid space
   {basis manifold of Euclid space}%
   {Basis Manifold, Euclid Space}%
\Index%260%B%basis manifold of Euclid space
   {basis manifold of Euclid space}%
   {Basis Manifold, Euclid Space, division ring}%
\Index%257%B%basis manifold of rcd affine space
   {basis manifold of \rcd affine space}%
   {Basis Manifold, rcd Affine Space, division ring}%
\Index%255%B%basis manifold of rcd vector space
   {basis manifold of \rcd vector space}%
   {basis manifold of rcd vector space}%
\Index%536%B%basis manifold of representation
   {basis manifold of representation}%
   {basis manifold representation F algebra}%
\Index%714%B%basis manifold of tower of representations
   {basis manifold of tower of representations}%
   {basis manifold tower of representations}%
\Index%258%B%basis manifold of vector space
   {basis manifold of vector space}%
   {basis manifold of vector space}%
\Index%689%B%basis of A- module
   {basis of $A$\Hyph module}%
   {basis of A module}%
\Index%766%B%basis of ACRcirc module
   {basis of $A\CRcirc$\Hyph module}%
   {basis of ACRcirc module}%
\Index%653%B%basis of algebra L(A,A)
   {basis of algebra $\mathcal L(A;A)$}%
   {basis of algebra L(A,A)}%
\Index%767%B%basis of CRcircA module
   {basis of $\RCcirc A$\Hyph module}%
   {basis of CRcircA module}%
\Index%329%B%basis of representation
   {basis of representation}%
   {basis of representation}%
\Index%567%B%basis of tower of representations
   {basis of tower of representations}%
   {basis of tower of representations}%
\Index%115%B%basis of vector space
   {basis of vector space}%
   {Basis}%
\Index%117%B%basis vector of representation of Lie group
   {basis vector of representation of Lie group}%
   {basis vector of representation of Lie group}%
\Index%116%B%basis vector of representation of Lie group over algebra A
   {basis vector of representation of Lie group over algebra $A$}%
   {basis vector of representation of Lie group over algebra A}%
\Index%122%B%biring
   {biring}%
   {biring}%
\Index%419%B%bundle of level 2
   {bundle of level $2$}%
   {bundle of level 2}%
\Index%418%B%bundle of level n
   {bundle of level $n$}%
   {bundle of level n}%
\SetIndexSpace%
\Index%6%C%c row of matrix
   {\subs row of matrix}%
   {c row}%
\Index%29%C%c rows rcd vector space
   {\sups rows \rcd vector space}%
   {sups rows rcd vector space}%
\Index%12%C%c-row of matrix
   {$c$\hyph row of matrix}%
   {c-row}%
\Index%443%C%Cartan connection
   {Cartan connection}%
   {Cartan connection}%
\Index%227%C%Cartan curvature
   {Cartan curvature}%
   {Cartan curvature}%
\Index%401%C%Cartan derivative
   {Cartan derivative}%
   {Cartan derivative}%
\Index%444%C%Cartan symbol
   {Cartan symbol}%
   {Cartan symbol}%
\Index%338%C%Cartan transport
   {Cartan transport}%
   {Cartan transport}%
\Index%149%C%Cartesian power A of bundle B
   {Cartesian power $\Bundle A$ of bundle $\Bundle B$}%
   {Cartesian power A of bundle B}%
\Index%148%C%Cartesian power A of set B
   {Cartesian power $A$ of set $B$}%
   {Cartesian power of set}%
\Index%150%C%Cartesian power n of bundle E
   {Cartesian power $n$ of bundle $\Bundle E$}%
   {Cartesian power n of bundle E}%
\Index%594%C%Cartesian power n of H algebra
   {Cartesian power $n$ of $\mathfrak{H}$\Hyph algebra}%
   {Cartesian power of algebra}%
\Index%175%C%category of drc vector spaces
   {category of \drc vector spaces}%
   {category of drc vector spaces}%
\Index%179%C%category of fibered correspondences over diagonal
   {category of fibered correspondences over diagonal}%
   {category of fibered correspondences over diagonal}%
\Index%178%C%category of reduced fibered correspondences
   {category of reduced fibered correspondences}%
   {category of reduced fibered correspondences}%
\Index%176%C%category of Tstar- representations of Omega1 algebra A
   {category of \Ts representations of $\Omega_1$\Hyph algebra $A$}%
   {category of Tstar representations of Omega1 algebra}%
\Index%177%C%category of Tstar- representations of Omega1 algebra from category A
   {category of \Ts representations of $\Omega_1$\Hyph algebra from category $\mathcal A$}%
   {category of Tstar representations of Omega1 algebra from category}%
\Index%271%C%Cauchy sequence in normed algebra
   {Cauchy sequence in normed algebra}%
   {Cauchy sequence, normed algebra}%
\Index%241%C%Cauchy sequence in normed ring
   {Cauchy sequence in normed ring}%
   {Cauchy sequence, normed ring}%
\Index%351%C%Cauchy sequence in valued division ring
   {Cauchy sequence in valued division ring}%
   {Cauchy sequence, valued division ring}%
\Index%797%C%Cauchy sequence in valued ring
   {Cauchy sequence in valued ring}%
   {Cauchy sequence, valued ring}%
\Index%620%C%center of algebra
   {center of an $R$\Hyph algebra $A$}%
   {center of algebra}%
\Index%507%C%center of ring D
   {center of ring $D$}%
   {center of ring}%
\Index%508%C%central affine basis
   {central affine basis}%
   {Central Affine Basis, division ring}%
\Index%509%C%central affine basis
   {central affine basis}%
   {Central Affine Basis}%
\Index%750%C%column determinant
   {column determinant}%
   {column determinant}%
\Index%124%C%column vector
   {column vector}%
   {column vector}%
\Index%615%C%commutative D algebra
   {commutative $D$\Hyph algebra}%
   {commutative D algebra}%
\Index%190%C%commutative diagram of correspondences
   {commutative diagram of correspondences}%
   {commutative diagram of correspondences}%
\Index%614%C%commutator of R algebra
   {commutator of $R$\Hyph algebra}%
   {commutator of algebra}%
\Index%191%C%compact-open topology
   {compact\hyph open topology}%
   {compact open topology}%
\Index%349%C%complete division ring
   {complete division ring}%
   {complete division ring}%
\Index%270%C%complete ring
   {complete ring}%
   {complete ring}%
\Index%348%C%complete system of linear partial differential equations
   {complete system of linear partial differential equations}%
   {Complete System of Linear Partial Differential Equations}%
\Index%128%C%completely integrable system
   {completely integrable system}%
   {completely integrable system}%
\Index%199%C%component of Gateaux derivative of map \Vector f(\Vector x)
   {component of the G\^ateaux derivative of map $\Vector f(\Vector x)$}%
   {component of Gateaux derivative of map, D vector space}%
\Index%198%C%component of Gateaux derivative of map f(x)
   {component of the G\^ateaux derivative of map $f(x)$}%
   {component of Gateaux derivative of map, division ring}%
\Index%197%C%component of Gateaux derivative of second order of map  of division ring
   {component of the G\^ateaux derivative of second order of map  of division ring}%
   {component of Gateaux derivative of Second Order, division ring}%
\Index%196%C%component of Gateaux derivative of second order of map \Vector f(\Vector x)
   {component of the G\^ateaux derivative of second order of map $\Vector f(\Vector x)$}%
   {component of Gateaux derivative of Second Order, D vector space}%
\Index%193%C%component of linear map f of division ring
   {component of linear map $f$ of division ring}%
   {component of linear map, division ring}%
\Index%192%C%component of linear map of D- vector space
   {component of linear map of $D$\Hyph vector space}%
   {component of linear map, D vector space}%
\Index%647%C%component of polylinear map into associative algebra
   {component of polylinear map into associative algebra}%
   {component of polylinear map, associative algebra}%
\Index%195%C%component of polylinear map of division ring
   {component of polylinear map of division ring}%
   {component of polylinear map, division ring}%
\Index%194%C%component of polylinear mapping \Vector A
   {component of polylinear mapping $\Vector A$}%
   {component of polyadditive map, D vector space}%
\Index%629%C%component of the Gateaux derivative of map f(x) of algebra
   {component of the G\^ateaux derivative of map $f(x)$ of algebra}%
   {component of Gateaux derivative of map, algebra}%
\Index%633%C%component of the Gateaux derivative of second order of map f(x) of algebra
   {component of the G\^ateaux derivative of second order of map $f(x)$ of algebra}%
   {component of Gateaux derivative of Second Order, algebra}%
\Index%394%C%composition of fibered correspondences
   {composition of fibered correspondences}%
   {composition of fibered correspondences}%
\Index%393%C%composition of reduced fibered correspondences
   {composition of reduced fibered correspondences}%
   {composition of reduced fibered correspondences}%
\Index%495%C%condition of reducibility of products
   {condition of reducibility of products}%
   {condition of reducibility of products}%
\Index%746%C%conjugate of quaternion
   {conjugate of quaternion $x$}%
   {conjugate of quaternion}%
\Index%788%C%conjugation in algebra
   {conjugation in algebra}%
   {conjugation in algebra}%
\Index%777%C%conjugation in ring
   {conjugation in ring}%
   {conjugation in ring}%
\Index%538%C%connection coefficients in D affine space
   {connection coefficients in $D$\Hyph affine space}%
   {connection coefficients, D affine space}%
\Index%281%C%continuous correspondence
   {continuous correspondence}%
   {continuous correspondence}%
\Index%280%C%continuous function of division ring
   {continuous function of division ring}%
   {continuous function, division ring}%
\Index%313%C%continuous function over algebra
   {continuous function over $D$\Hyph algebra}%
   {continuous function, algebra}%
\Index%200%C%contravariant starT- representation of fibered group
   {contravariant \sT representation of fibered group}%
   {contravariant starT representation of fibered group}%
\Index%202%C%contravariant starT- representation of group
   {contravariant \sT representation of group}%
   {contravariant starT representation of group}%
\Index%203%C%contravariant Tstar- representation of fibered group
   {contravariant \Ts representation of fibered group}%
   {contravariant Tstar representation of fibered group}%
\Index%201%C%contravariant Tstar- representation of group
   {contravariant \Ts representation of group}%
   {contravariant Tstar representation of group}%
\Index%210%C%coordinate Drc vector bundle
   {coordinate \Drc vector bundle}%
   {coordinate drc vector bundle}%
\Index%215%C%coordinate isomorphism
   {coordinate isomorphism}%
   {coordinate isomorphism}%
\Index%206%C%coordinate matrix of set of vectors in dcr rows vector space
   {coordinate matrix of set of vectors in \dcr rows vector space}%
   {coordinate matrix of set of vectors, dcr vector space}%
\Index%207%C%coordinate matrix of set of vectors in rcd rows vector space
   {coordinate matrix of set of vectors in \rcd rows vector space}%
   {coordinate matrix of set of vectors, rcd vector space}%
\Index%205%C%coordinate matrix of vector field in rcD basis
   {coordinate matrix of vector field in \rcD basis}%
   {coordinate matrix of vector field in drc basis}%
\Index%691%C%coordinate matrix of vector in basis
   {coordinate matrix of vector in basis}%
   {coordinate matrix of vector in basis, module}%
\Index%204%C%coordinate matrix of vector in drc basis
   {coordinate matrix of vector in \drc basis}%
   {coordinate matrix of vector in rcd basis}%
\Index%214%C%coordinate rcd isomorphism
   {coordinate \rcd isomorphism}%
   {coordinate rcd isomorphism}%
\Index%209%C%coordinate rcd vector space
   {coordinate \rcd vector space}%
   {coordinate rcd vector space}%
\Index%208%C%coordinate reference frame
   {coordinate reference frame}%
   {coordinate reference frame}%
\Index%697%C%coordinate representation in Omega_2- algebra
   {coordinate representation in $\Omega_2$\Hyph algebra}%
   {coordinate representation, Omega_2 algebra}%
\Index%212%C%coordinate representation in rcd vector space
   {coordinate representation in \rcd vector space}%
   {coordinate representation, rcd vector space}%
\Index%720%C%coordinate representation in tuple of Omega- algebras
   {coordinate representation in tuple of $\VX\Omega$\Hyph algebras}%
   {coordinate tower of representations, Omega algebra}%
\Index%213%C%coordinate representation of group in vector space
   {coordinate representation of group in vector space}%
   {coordinate representation, vector space}%
\Index%211%C%coordinate vector space
   {coordinate vector space}%
   {coordinate vector space}%
\Index%701%C%coordinates of a geometric object in Omega_2- algebra
   {coordinates of a geometric object in $\Omega_2$\Hyph algebra $M$}%
   {coordinates of geometric object, representation g}%
\Index%723%C%coordinates of a geometric object in tuple of Omega- algebras
   {coordinates of a geometric object in tuple of $\VX\Omega$\Hyph algebras}%
   {coordinates of geometric object, tower of representations g}%
\Index%325%C%coordinates of basis of representation
   {coordinates of basis of representation}%
   {coordinates of basis relative to basis, representation}%
\Index%558%C%coordinates of element of representation relative to set
   {coordinates of element $m$ of representation $f$ relative to set $X$}%
   {coordinates of element relative to set, representation}%
\Index%696%C%coordinates of endomorphism of representation
   {coordinates of endomorphism of representation}%
   {coordinates of endomorphism, representation}%
\Index%711%C%coordinates of endomorphism of tower of representations
   {coordinates of endomorphism of tower of representations}%
   {coordinates of endomorphism, tower of representations}%
\Index%218%C%coordinates of geometric object
   {coordinates of geometric object}%
   {coordinates of geometric object, vector space}%
\Index%220%C%coordinates of geometric object in coordinate rcd vector space
   {coordinates of geometric object in coordinate \rcd vector space}%
   {coordinates of geometric object, coordinate rcd vector space}%
\Index%221%C%coordinates of geometric object in coordinate representation
   {coordinates of geometric object in coordinate representation}%
   {coordinates of geometric object, coordinate vector space}%
\Index%699%C%coordinates of geometric object in coordinate space of representation
   {coordinates of geometric object in coordinate space of representation}%
   {coordinates of geometric object, coordinate representation}%
\Index%721%C%coordinates of geometric object in coordinate space of tower of representations
   {coordinates of geometric object in coordinate space of tower of representations}%
   {coordinates of geometric object, coordinate tower of representations}%
\Index%219%C%coordinates of geometric object in rcd vector space
   {coordinates of geometric object in \rcd vector space}%
   {coordinates of geometric object, rcd vector space}%
\Index%730%C%coordinates of point of affine space relative to basis)$
   {coordinates of point $A$ of affine space $\overset{\circ}{A}$ relative to basis $(O,\Basis e)$}%
   {coordinates in affine space}%
\Index%737%C%coordinates of point of rcd affine space relative to basis
   {coordinates of point of \rcd affine space relative to basis}%
   {coordinates in rcd affine space}%
\Index%224%C%coordinates of representation
   {coordinates of representation}%
   {coordinates of representation, drc vector space}%
\Index%225%C%coordinates of representation
   {coordinates of representation}%
   {coordinates of representation}%
\Index%222%C%coordinates of set of vectors in dcr vector space
   {coordinates of set of vectors in \dcr vector space}%
   {coordinates of set of vectors, dcr vector space}%
\Index%223%C%coordinates of set of vectors in rcd vector space
   {coordinates of set of vectors in \rcd vector space}%
   {coordinates of set of vectors, rcd vector space}%
\Index%217%C%coordinates of vector field in Drc basis
   {coordinates of vector field in \Drc basis}%
   {coordinates of vector field in drc basis}%
\Index%692%C%coordinates of vector in basis
   {coordinates of vector in basis}%
   {coordinates of vector in basis, module}%
\Index%216%C%coordinates of vector in rcd basis
   {coordinates of vector in \rcd basis}%
   {coordinates of vector in rcd basis}%
\Index%792%C%coordinates of vector relative to Schauder basis
   {coordinates of vector relative to Schauder basis}%
   {coordinates of vector, Schauder basis}%
\Index%455%C%correspondence continuous on the set
   {correspondence continuous on the set}%
   {correspondence continuous on the set}%
\Index%454%C%correspondence of homomorphism
   {correspondence of homomorphism}%
   {correspondence of homomorphism}%
\Index%185%C%covariant starT- representation of fibered group
   {covariant \sT representation of fibered group}%
   {covariant starT representation of fibered group}%
\Index%184%C%covariant starT- representation of group
   {covariant \sT representation of group}%
   {covariant starT representation of group}%
\Index%187%C%covariant Tstar- representation of fibered group
   {covariant \Ts representation of fibered group}%
   {covariant Tstar representation of fibered group}%
\Index%186%C%covariant Tstar- representation of group
   {covariant \Ts representation of group}%
   {covariant Tstar representation of group}%
\Index%8%C%CR inverse element of biring
   {\CR inverse element of biring}%
   {cr-inverse element}%
\Index%7%C%CR matrix group
   {\CR matrix group}%
   {cr-matrix group}%
\Index%10%C%CR power
   {\CR power}%
   {cr power}%
\Index%9%C%CR product of matrices
   {\CR product of matrices}%
   {cr-product of matrices}%
\Index%589%C%cr product of matrices of mappings
   {$\CRcirc$\Hyph product of matrices of mappings}%
   {cr product of matrices of mappings}%
\Index%11%C%crd vector space
   {\crd vector space}%
   {crd vector space}%
\Index%768%C%Cstar-algebra
   {$C^*$\Hyph algebra}%
   {Cstar-algebra}%
\Index%537%C%curvilinear coordinates of point in affine space
   {curvilinear coordinates of point in affine space}%
   {curvilinear coordinates of point in affine space}%
\Index%19%C%r rows drc vector space
   {\subs rows \drc vector space}%
   {subs rows drc vector space}%
\SetIndexSpace%
\Index%23%D%d affine space
   {$D$\Hyph affine space}%
   {d affine space}%
\Index%549%D%D- affine connection on manifold with affine connections
   {$D$\Hyph affine connection on manifold with affine connections}%
   {D affine connection, affine manifold}%
\Index%48%D%D- valued variable
   {$D$\Hyph valued variable}%
   {D valued variable}%
\Index%14%D%D- vector function
   {$D$\Hyph vector function}%
   {d vector function}%
\Index%550%D%D-affine connection coefficients on manifold
   {$D$\Hyph affine connection coefficients on manifold}%
   {D affine connection coefficients, manifold}%
\Index%15%D%D-vector space
   {$D$\hyph vector space}%
   {D vector space}%
\Index%18%D%dcr vector space
   {\dcr vector space}%
   {dcr vector space}%
\Index%303%D%determinant of matrix
   {determinant of matrix}%
   {determinant}%
\Index%309%D%deviation of trajectories
   {deviation of trajectories}%
   {deviation of trajectories}%
\Index%153%D%diagonal in bundle
   {diagonal in bundle}%
   {diagonal in bundle}%
\Index%154%D%diagram of correspondences
   {diagram of correspondences}%
   {diagram of correspondences}%
\Index%725%D%diagram of representations
   {diagram of representations}%
   {diagram of representations}%
\Index%417%D%dimension of rcd vector space
   {dimension of \rcd vector space}%
   {dimension of vector space}%
\Index%151%D%direct product of bundles
   {direct product of bundles}%
   {Cartesian product of bundles}%
\Index%410%D%direct product of D- vector spaces
   {direct product of $D$\Hyph vector spaces}%
   {direct product of D vector spaces}%
\Index%415%D%direct product of division rings
   {direct product of division rings}%
   {direct product of division rings}%
\Index%411%D%direct product of rcd vector spaces
   {direct product of \rcd vector spaces}%
   {direct product, rcd vector space}%
\Index%414%D%direct product of representations of fibered group
   {direct product of representations of fibered group}%
   {direct product of representations of fibered group}%
\Index%413%D%direct product of representations of group
   {direct product of representations of group}%
   {direct product of representations of group}%
\Index%152%D%direct product of total spaces
   {direct product of total spaces}%
   {Cartesian product of total spaces}%
\Index%412%D%direct product of Tstar- representations of group
   {direct product of \Ts representations of group}%
   {direct product of Tstar representations of group}%
\Index%409%D%direct sum of representations
   {direct sum of representations}%
   {direct sum of representations}%
\Index%663%D%distributive law for Astar- module
   {distributive law for $A\star$\Hyph module}%
   {distributive law, Astar module over algebra}%
\Index%577%D%distributive law for Astar- vector space
   {distributive law for $A\star$\Hyph vector space}%
   {distributive law, Astar vector space}%
\Index%655%D%distributive law for Dstar- module
   {distributive law for $D\star$\Hyph module}%
   {distributive law, Dstar module}%
\Index%170%D%distributive law for Dstar- vector fields
   {distributive law for $\mathcal D\star$\Hyph vector fields}%
   {distributive law, Dstar vector fields}%
\Index%169%D%distributive law for Dstar- vector space
   {distributive law for $D\star$\Hyph vector space}%
   {distributive law, Dstar vector space}%
\Index%672%D%distributive law for starA- module
   {distributive law for $\star A$\Hyph module}%
   {distributive law, starA module over algebra}%
\Index%755%D%double determinant
   {double determinant}%
   {double determinant}%
\Index%26%D%Drc basis for vector  bundle
   {\Drc basis for vector  bundle}%
   {drc basis, vector bundle}%
\Index%16%D%drc basis of r rows vector space
   {\drc basis of \subs rows vector space}%
   {drc basis, r rows vector space}%
\Index%36%D%Drc linear map of vector bundles
   {\Drc linear map of vector bundles}%
   {drc linear map of vector bundles}%
\Index%20%D%Drc linearly dependent vector fields
   {\Drc linearly dependent vector fields}%
   {linearly dependent vector fields}%
\Index%33%D%Drc linearly independent vector fields
   {\Drc linearly independent vector fields}%
   {linearly independent vector fields}%
\Index%70%D%drc representation of group
   {\drc representation of group}%
   {drc linear representation of group}%
\Index%17%D%drc vector
   {\drc vector}%
   {drc vector}%
\Index%30%D%drc vector function
   {\drc vector function}%
   {drc vector function}%
\Index%28%D%drc vector space
   {\drc vector space}%
   {drc vector space}%
\Index%734%D%Dstar affine space
   {$D\star$\Hyph affine space}%
   {Dstar affine space}%
\Index%771%D%Dstar antilinear homomorphism
   {$D\star$\Hyph antilinear homomorphism}%
   {Dstar antilinear homomorphism}%
\Index%779%D%Dstar antilinear mapping of ring with conjugation
   {$D\star$\Hyph antilinear mapping of ring with conjugation}%
   {antilinear mapping of ring with conjugation}%
\Index%42%D%Dstar component of coordinates of vector \Vector r
   {\Ds component of coordinates of vector $\Vector r$}%
   {Dstar component of coordinates of vector, D vector space}%
\Index%769%D%Dstar linear homomorphism
   {$D\star$\Hyph linear homomorphism}%
   {Dstar linear homomorphism}%
\Index%658%D%Dstar linearly dependent vectors of Dstar module
   {$D\star$\Hyph linearly dependent vectors of $D\star$\Hyph module}%
   {Dstar linearly dependent, module}%
\Index%657%D%Dstar linearly independent vectors of Dstar module
   {$D\star$\Hyph linearly independent vectors of $D\star$\Hyph module}%
   {Dstar linearly independent, module}%
\Index%659%D%Dstar- basis for module
   {$D\star$\Hyph basis for module}%
   {D basis, module}%
\Index%71%D%Dstar- module
   {$D\star$\Hyph module}%
   {Dstar-module}%
\Index%41%D%Dstar- vector bundle
   {$\mathcal D\star$\Hyph vector bundle}%
   {Dstar vector bundle}%
\Index%39%D%Dstar- vector field
   {$\mathcal D\star$\Hyph vector field}%
   {Dstar vector field}%
\Index%40%D%Dstar- vector space
   {$D\star$\hyph  vector space}%
   {Dstar vector space}%
\Index%43%D%Dstar-linear composition of vector fields
   {$\mathcal D\star$\hyph linear composition of vector fields}%
   {linear composition of vector fields}%
\Index%45%D%Dstar-product of vector field over scalar
   {$\mathcal D\star$\hyph product of vector field over scalar}%
   {Dstar product of vector field over scalar, vector space}%
\Index%44%D%Dstar-product of vector over scalar
   {$D\star$\hyph product of vector over scalar}%
   {Dstar product of vector over scalar, vector space}%
\Index%159%D%dual space of rcd vector space
   {dual space of \rcd vector space}%
   {dual space of rcd vector space}%
\Index%382%D%duality principle for biring
   {duality principle for biring}%
   {duality principle for biring}%
\Index%383%D%duality principle for biring of matrices
   {duality principle for biring of matrices}%
   {duality principle for biring of matrices}%
\Index%25%D%rcd basis for c rows vector space
   {\rcd basis for \sups rows vector space}%
   {rcd basis, c rows vector space}%
\Index%32%D%rcd linear span in vector space
   {\rcd linear span in vector space}%
   {linear span, vector space}%
\SetIndexSpace%
\Index%521%E%effective representation of division ring
   {effective representation of division ring}%
   {effective representation of division ring}%
\Index%519%E%effective representation of fibered Omega- algebra
   {effective representation of fibered $\Omega$\Hyph algebra}%
   {effective representation of fibered Omega-algebra}%
\Index%518%E%effective representation of group
   {effective representation of group}%
   {effective representation of group}%
\Index%517%E%effective representation of Omega- algebra A
   {effective representation of $\Omega$\Hyph algebra $A$}%
   {effective representation of algebra}%
\Index%622%E%effective representation of ring
   {effective representation of ring}%
   {effective representation of ring}%
\Index%514%E%effective Tstar- representation of fibered division ring
   {effective \Ts representation of fibered division ring}%
   {effective representation of fibered division ring}%
\Index%520%E%effective Tstar- representation of fibered group
   {effective \Ts representation of fibered group}%
   {effective representation of fibered group}%
\Index%516%E%effective Tstar- representation of group
   {effective \Ts representation of group}%
   {effective Tstar representation of group}%
\Index%557%E%endomorphism of representation of Omega algebra
   {endomorphism of representation of $\Omega$\Hyph algebra}%
   {endomorphism of representation}%
\Index%327%E%endomorphism of representation regular on generating set X
   {endomorphism of representation regular on generating set $X$}%
   {endomorphism of representation, regular on set}%
\Index%328%E%endomorphism of representation singular on generating set X
   {endomorphism of representation singular on generating set $X$}%
   {endomorphism of representation, singular on set}%
\Index%559%E%endomorphism of tower of representations
   {endomorphism of tower of representations}%
   {endomorphism of tower of representations}%
\Index%564%E%endomorphism of tower of representations regular on tuple of generating sets
   {endomorphism of tower of representations regular on tuple of generating sets}%
   {endomorphism of representation, regular on tuple}%
\Index%565%E%endomorphism of tower of representations singular on tuple of generating sets
   {endomorphism of tower of representations singular on tuple of generating sets}%
   {endomorphism of representation, singular on tuple}%
\Index%50%E%enhanced Lie group
   {enhanced Lie group}%
   {enhanced Lie group}%
\Index%706%E%equivalence generated by representation
   {equivalence generated by representation $f$}%
   {equivalence of representation}%
\Index%51%E%essential parameters in a set of functions
   {essential parameters in a set of functions}%
   {essential parameters}%
\Index%539%E%Euclidean metric on division ring
   {Euclidean metric on division ring}%
   {Euclidean metric on division ring}%
\Index%547%E%Euclidean scalar product in D- vector space
   {Euclidean scalar product in $D$\Hyph vector space}%
   {Euclidean scalar product, vector space}%
\Index%542%E%Euclidean scalar product on division ring
   {Euclidean scalar product on division ring}%
   {Euclidean scalar product on division ring}%
\Index%437%E%extended matrix of drc linear equations
   {extended matrix of \drc linear equations}%
   {extended matrix, system of drc linear equations}%
\Index%438%E%extended matrix of rcd linear equations
   {extended matrix of \rcd linear equations}%
   {extended matrix, system of rcd linear equations}%
\Index%386%E%extension of correspondence
   {extension of correspondence}%
   {extension of correspondence}%
\Index%515%E%extreme line
   {extreme line}%
   {extreme line}%
\SetIndexSpace%
\Index%434%F%fibered coordinate Drc isomorphism
   {fibered coordinate \Drc isomorphism}%
   {fibered coordinate drc isomorphism}%
\Index%432%F%fibered correspondence from A to B
   {fibered correspondence from $\Bundle A$ to $\Bundle B$}%
   {fibered correspondence from A to B}%
\Index%430%F%fibered correspondence in A
   {fibered correspondence in $\Bundle{A}$}%
   {fibered correspondence in A}%
\Index%431%F%fibered correspondence of homomorphism
   {fibered correspondence of homomorphism}%
   {fibered correspondence of homomorphism}%
\Index%427%F%fibered equivalence
   {fibered equivalence}%
   {fibered equivalence}%
\Index%422%F%fibered group
   {fibered group}%
   {fibered group}%
\Index%436%F%fibered identification morphism
   {fibered identification morphism}%
   {fibered identification morphism}%
\Index%424%F%fibered little group
   {fibered little group}%
   {fibered little group}%
\Index%435%F%fibered morphism from bundle A into B
   {fibered morphism from bundle $\Bundle A$ into $\Bundle B$}%
   {fibered morphism from A into B}%
\Index%433%F%fibered natural morphism
   {fibered natural morphism}%
   {fibered natural morphism}%
\Index%420%F%fibered Omega- algebra
   {fibered $\Omega$\Hyph algebra}%
   {fibered Omega-algebra}%
\Index%421%F%fibered Omega- subalgebra
   {fibered $\Omega$\Hyph subalgebra}%
   {fibered Omega-subalgebra}%
\Index%426%F%fibered ordering
   {fibered ordering}%
   {fibered ordering}%
\Index%425%F%fibered preordering
   {fibered preordering}%
   {fibered preordering}%
\Index%428%F%fibered ring
   {fibered ring}%
   {fibered ring}%
\Index%423%F%fibered stability group
   {fibered stability group}%
   {fibered stability group}%
\Index%429%F%fibered subset
   {fibered subset}%
   {fibered subset}%
\Index%479%F%field-strength tensor
   {field-strength tensor}%
   {field-strength tensor}%
\Index%497%F%filter F converges to A
   {filter $\mathfrak{F}$ converges to $A$}%
   {filter converges}%
\Index%808%F%Finsler metric
   {Finsler metric}%
   {Finsler metric}%
\Index%798%F%Finsler space
   {Finsler space}%
   {Finsler space}%
\Index%807%F%Finsler structure
   {Finsler structure}%
   {Finsler structure}%
\Index%337%F%first Newton law
   {first Newton law}%
   {First Newton law}%
\Index%690%F%free A- module
   {free $A$\Hyph module}%
   {free A module}%
\Index%606%F%free algebra over ring
   {free algebra over ring}%
   {free algebra over ring}%
\Index%612%F%free module over ring
   {free module over ring}%
   {free module over ring}%
\Index%442%F%free Tstar- representation of fibered group
   {free \Ts representation of fibered group}%
   {free representation of fibered group}%
\Index%441%F%free Tstar- representation of group
   {free \Ts representation of group}%
   {free representation of group}%
\Index%339%F%Frenet transport
   {Frenet transport}%
   {Frenet transport}%
\Index%506%F%function continuous with respect to set of arguments
   {function continuous with respect to set of arguments}%
   {function continuous with respect to set of arguments}%
\Index%805%F%function homogeneous of degree
   {function homogeneous of degree $k$}%
   {function homogeneous}%
\Index%502%F%function of \gi n D- valued variables
   {function of $\gi n$ $D$\Hyph valued variables}%
   {function of n D valued variables}%
\Index%625%F%function of algebra differentiable in the Gateaux sense
   {function of algebra differentiable in the G\^ateaux sense}%
   {function differentiable in Gateaux sense, algebra}%
\Index%501%F%function of D- vector space V to D- vector space W differentiable in Gateaux sense
   {function of $D$\Hyph vector space $\Vector{V}$ to $D$|Hyph vector space $\Vector W$ differentiable in the G\^ateaux sense}%
   {function differentiable in Gateaux sense, D vector space}%
\Index%505%F%function of division ring differentiable in Gateaux sense
   {function of division ring differentiable in the G\^ateaux sense}%
   {function differentiable in Gateaux sense, division ring}%
\Index%504%F%function of division ring Dstar differentiable in the Fr\'echet sense
   {function of division ring \Ds differentiable in the Fr\'echet sense}%
   {function Dstar differentiable in Frechet sense, division ring}%
\Index%243%F%fundamental sequence in normed algebra
   {fundamental sequence in normed algebra}%
   {fundamental sequence, normed algebra}%
\Index%642%F%fundamental sequence in normed ring
   {fundamental sequence in normed ring}%
   {fundamental sequence, normed ring}%
\Index%500%F%fundamental sequence in valued division ring
   {fundamental sequence in valued division ring}%
   {fundamental sequence, valued division ring}%
\Index%796%F%fundamental sequence in valued ring
   {fundamental sequence in valued ring}%
   {fundamental sequence, valued ring}%
\SetIndexSpace%
\Index%52%G%G- reference frame
   {$G$\Hyph reference frame}%
   {G reference frame}%
\Index%53%G%G-basis of vector space
   {$G$\Hyph basis of vector space}%
   {G-basis}%
\Index%54%G%G-coordinates of basis
   {$G$\Hyph coordinates of basis}%
   {G-coordinates}%
\Index%55%G%G-space
   {$G$\Hyph space}%
   {GSpace}%
\Index%13%G%Gateaux crd derivative of map \Vector f of D-vector space \Vector V to D-vector space \Vector W
   {the G\^ateaux \crd derivative of map $\Vector f$ of $D$\hyph vector space $\Vector V$ to $D$\hyph vector space $\Vector W$}%
   {Gateaux crd derivative of map, D vector space}%
\Index%397%G%Gateaux derivative of map f
   {the G\^ateaux derivative of map $f$}%
   {Gateaux derivative of map, division ring}%
\Index%398%G%Gateaux derivative of map Vector f of normed D- vector space Vector V to normed D- vector space Vector W
   {the G\^ateaux derivative of map $\Vector f$ of normed $D$\Hyph vector space $\Vector{V}$ to normed $D$\Hyph vector space $\Vector{W}$}%
   {Gateaux derivative of map, D vector space}%
\Index%627%G%Gateaux derivative of mapping of algebra
   {the G\^ateaux derivative of mapping of algebra}%
   {Gateaux derivative of map, algebra}%
\Index%399%G%Gateaux derivative of order n of map \Vector f
   {the G\^ateaux derivative of order $n$ of map $\Vector f$}%
   {Gateaux derivative of Order n, D vector space}%
\Index%400%G%Gateaux derivative of order n of map f of division ring
   {the G\^ateaux derivative of order $n$ of map $f$ of division ring}%
   {Gateaux derivative of Order n, division ring}%
\Index%634%G%Gateaux derivative of order n of mapping f of algebra
   {the G\^ateaux derivative of order $n$ of map $f$ of algebra}%
   {Gateaux derivative of Order n, algebra}%
\Index%626%G%Gateaux derivative of second order of map of algebra
   {the G\^ateaux derivative of second order of mapping of algebra}%
   {Gateaux derivative of Second Order, algebra}%
\Index%396%G%Gateaux derivative of second order of map of division ring
   {the G\^ateaux derivative of second order of map of division ring}%
   {Gateaux derivative of Second Order, division ring}%
\Index%157%G%Gateaux differential of map \Vector f of normed D- vector space V to normed D- vector space W
   {the G\^ateaux differential of map $\Vector f$ of normed $D$\Hyph vector space $\Vector{V}$ to normed $D$\Hyph vector space $\Vector{W}$}%
   {Gateaux differential of map, D vector space}%
\Index%158%G%Gateaux differential of map f
   {the G\^ateaux differential of map $f$}%
   {Gateaux differential of map, division ring}%
\Index%628%G%Gateaux differential of mapping of algebra
   {the G\^ateaux differential of mapping $f$ of algebra}%
   {Gateaux differential of map, algebra}%
\Index%155%G%Gateaux differential of second order of map of division ring
   {the G\^ateaux differential of second order of map of division ring}%
   {Gateaux differential of Second Order, division ring}%
\Index%156%G%Gateaux differential of second order of mapping \Vector f
   {the G\^ateaux differential of second order of mapping $\Vector f$}%
   {Gateaux differential of Second Order, D vector space}%
\Index%632%G%Gateaux differential of second order of mapping f of algebra
   {the G\^ateaux differential of second order of mapping $f$ of algebra}%
   {Gateaux differential of Second Order, algebra}%
\Index%38%G%Gateaux drc derivative of map \Vector f of D- vector space \Vector V to D- vector space \Vector W
   {the G\^ateaux \drc derivative of map $\Vector f$ of $D$\Hyph vector space $\Vector V$ to $D$\Hyph vector space $\Vector W$}%
   {Gateaux drc derivative of map, D vector space}%
\Index%46%G%Gateaux Dstar derivative of map f of division ring D
   {the G\^ateaux \Ds derivative of map $f$ of division ring $D$}%
   {Gateaux Dstar derivative of map, division ring}%
\Index%453%G%Gateaux mixed partial derivative of map f^j with respect to variables v^i, v^j
   {the G\^ateaux mixed partial derivative of map $f^j$ with respect to variables $v^i$, $v^j$}%
   {Gateaux partial derivative of Second Order, D vector space}%
\Index%512%G%Gateaux partial derivative of map f^j with respect to variable v^i
   {the G\^ateaux partial derivative of map $f^j$ with respect to variable $v^i$}%
   {Gateaux partial derivative, D vector space}%
\Index%511%G%Gateaux partial drc derivative of map f^b with respect to variable x^a
   {the G\^ateaux partial \drc derivative of map $f^b$ with respect to variable $x^a$}%
   {Gateaux partial drc derivative of map with respect to variable, D vector space}%
\Index%510%G%Gateaux partial drc derivative of map fb with respect to variable xa
   {the G\^ateaux partial \drc derivative of map $f^b$ with respect to variable $x^a$}%
   {Gateaux partial crd derivative of map with respect to variable, D vector space}%
\Index%76%G%Gateaux starD derivative of map f of division ring D
   {the G\^ateaux \sD derivative of map $f$ of division ring $D$}%
   {Gateaux starD derivative of map, division ring}%
\Index%332%G%generating set of representation
   {generating set of representation}%
   {generating set of representation}%
\Index%331%G%generating set of subrepresentation
   {generating set of subrepresentation}%
   {generating set of subrepresentation}%
\Index%289%G%generator of linear mapping
   {generator of linear mapping}%
   {generator of linear map, division ring}%
\Index%700%G%geometric object defined in Omega_2- algebra
   {geometric object defined in $\Omega_2$\Hyph algebra $M$}%
   {geometric object, representation g}%
\Index%136%G%geometric object defined in rcd vector space
   {geometric object defined in \rcd vector space}%
   {geometric object, rcd vector space}%
\Index%722%G%geometric object defined in tuple of Omega- algebras
   {geometric object defined in tuple of $\VX\Omega$\Hyph algebras $\VX A$}%
   {geometric object, tower of representations g}%
\Index%133%G%geometric object in coordinate representation
   {geometric object in coordinate representation}%
   {geometric object, coordinate vector space}%
\Index%698%G%geometric object in coordinate representation defined in Omega_2- algebra
   {geometric object in coordinate representation defined in $\Omega_2$\Hyph algebra $M$}%
   {geometric object, coordinate representation g}%
\Index%134%G%geometric object in coordinate representation defined in rcd vector space
   {geometric object in coordinate representation defined in \rcd vector space}%
   {geometric object, coordinate rcd vector space}%
\Index%719%G%geometric object in coordinate representation defined in tuple of Omega- algebras
   {geometric object in coordinate representation defined in tuple of $\VX\Omega$\Hyph algebras $\VX A$}%
   {geometric object, coordinate tower of representations g}%
\Index%132%G%geometric object in vector space
   {geometric object in vector space}%
   {geometric object, vector space}%
\Index%704%G%geometric object of type
   {geometric object of type $H$}%
   {geometric object of type H, representation g}%
\Index%135%G%geometric object of type A
   {geometric object of type $A$ in vector space}%
   {geometric object of type A, vector space}%
\Index%146%G%group algebra
   {group algebra}%
   {group algebra}%
\Index%568%G%group of automorphisms of representation
   {group of automorphisms of representation}%
   {group of automorphisms of representation}%
\SetIndexSpace%
\Index%292%H%Hadamard inverse of matrix
   {Hadamard inverse of matrix}%
   {Hadamard inverse of matrix}%
\Index%783%H%Hamel basis
   {Hamel basis}%
   {Hamel basis}%
\Index%546%H%hermitian conjugated vector
   {hermitian conjugated vector}%
   {hermitian conjugated vector}%
\Index%541%H%hermitian conjugation in division ring
   {hermitian conjugation in division ring}%
   {hermitian conjugation, division ring}%
\Index%747%H%hermitian matrix
   {hermitian matrix}%
   {hermitian matrix}%
\Index%544%H%hermitian metric on division ring
   {hermitian metric on division ring}%
   {hermitian metric on division ring}%
\Index%350%H%hermitian scalar product in D- vector space
   {hermitian scalar product in $D$\Hyph vector space}%
   {hermitian scalar product, vector space}%
\Index%545%H%hermitian scalar product on division ring
   {hermitian scalar product on division ring}%
   {hermitian scalar product on division ring}%
\Index%138%H%holonomic coordinates of connection
   {holonomic coordinates of connection}%
   {holonomic coordinates of connection}%
\Index%137%H%holonomic coordinates of vector
   {holonomic coordinates of vector}%
   {vector holonomic coordinates}%
\Index%296%H%homogeneous bundle of fibered group
   {homogeneous bundle of fibered group}%
   {homogeneous bundle of fibered group}%
\Index%804%H%homogeneous linear geometric object
   {homogeneous linear geometric object}%
   {homogeneous linear geometric object}%
\Index%802%H%homogeneous linear representation of Lie group
   {homogeneous linear representation of Lie group}%
   {homogeneous Linear Representation of Lie Group}%
\Index%294%H%homogeneous map of degree k over field F
   {homogeneous map of degree $k$ over field $F$}%
   {homogeneous map of degree over field, D vector space}%
\Index%295%H%homogeneous space of group
   {homogeneous space of group}%
   {homogeneous space of group}%
\Index%140%H%homomorphism of fibered groups
   {homomorphism of fibered groups}%
   {homomorphism of fibered groups}%
\Index%139%H%homomorphism of fibered universal algebras
   {homomorphism of fibered universal algebras}%
   {homomorphism of fibered universal algebras}%
\SetIndexSpace%
\Index%121%I%infinitesimal generator
   {infinitesimal generator}%
   {infinitesimal generator}%
\Index%174%I%infinitesimal generators of group Lie
   {infinitesimal generators of group Lie}%
   {infinitesimal generators of group Lie}%
\Index%705%I%invariance principle in drc vector space
   {invariance principle in \drc vector space}%
   {invariance principle}%
\Index%384%I%invariance principle in representation of universal algebra
   {invariance principle in representation of universal algebra}%
   {invariance principle, representation g}%
\Index%710%I%invariance principle in tower of representations of universal algebras
   {invariance principle in tower of representations of universal algebras}%
   {invariance principle, tower of representations g}%
\Index%385%I%invariance principle in vector space
   {invariance principle in vector space}%
   {invariance principle, vector space}%
\Index%291%I%inverse fibered correspondence
   {inverse fibered correspondence}%
   {inverse fibered correspondence}%
\Index%290%I%inverse reduced fibered correspondence
   {inverse reduced fibered correspondence}%
   {inverse reduced fibered correspondence}%
\Index%745%I%involution in quaternion algebra
   {involution in quaternion algebra}%
   {involution, quaternion algebra}%
\Index%173%I%isomorphism of fibered Omega- algebras
   {isomorphism of fibered $\Omega$\Hyph algebras}%
   {isomorphism of fibered Omega-algebras}%
\Index%266%I%isomorphism of repesentations of Omega- algebra
   {isomorphism of repesentations of $\Omega$\Hyph algebra}%
   {isomorphism of repesentations of Omega algebra}%
\Index%810%I%isotropic vector
   {isotropic vector}%
   {isotropic vector}%
\SetIndexSpace%
\Index%596%J%j i $RCcirc quasideterminant
   {$(^j_i)$\hyph $\RCcirc$\Hyph quasideterminant}%
   {j i RCcirc-quasideterminant}%
\Index%592%J%Jacobian complete system of differential equations
   {Jacobian complete system of differential equations}%
   {Jacobian complete system of differential equations}%
\Index%593%J%Jacobian complete system of drc differential equations
   {Jacobian complete system of \drv differential equations}%
   {Jacobian complete system of drc differential equations}%
\Index%678%J%ji- quasideterminant
   {$(ji)$\hyph quasideterminant}%
   {j i quasideterminant}%
\Index%527%J%the Jacobi Gateaux matrix of map of D- vector space
   {the Jacobi\Hyph G\^ateaux matrix of map of $D$\Hyph vector space}%
   {Jacobi Gateaux matrix of map, D vector space}%
\SetIndexSpace%
\Index%526%K%kernel of inefficiency of representation of fibered group
   {kernel of inefficiency of representation of fibered group}%
   {kernel of inefficiency of representation of fibered group}%
\Index%525%K%kernel of inefficiency of representation of group
   {kernel of inefficiency of representation of group}%
   {kernel of inefficiency of representation of group}%
\Index%524%K%kernel of inefficiency of Tstar- representation of group G
   {kernel of inefficiency of \Ts representation of group $G$}%
   {kernel of inefficiency of Tstar representation of group}%
\Index%522%K%kernel of linear mapping of D- vector space
   {kernel of linear mapping of $D$\Hyph vector space}%
   {kernel of linear map, D vector space}%
\Index%523%K%kernel of linear mapping of division ring
   {kernel of linear mapping of division ring}%
   {kernel of linear map, division ring}%
\Index%493%K%Killing equation
   {Killing equation}%
   {Killing equation}%
\Index%494%K%Killing equation of second type
   {Killing equation of second type}%
   {Killing equation second type}%
\Index%123%K%Killing vector of second type
   {Killing vector of second type}%
   {Killing vector second type}%
\Index%445%K%Kronecker symbol
   {Kronecker symbol}%
   {Kronecker symbol}%
\SetIndexSpace%
\Index%754%L%left cofactor of entry of matrix
   {left cofactor of entry of matrix}%
   {left cofactor, matrix}%
\Index%301%L%left defined Lie algebra of Lie group
   {left defined Lie algebra of Lie group}%
   {left defined Lie algebra}%
\Index%756%L%left double cofactor of entry of matrix
   {left double cofactor of entry of matrix}%
   {left double cofactor}%
\Index%229%L%left invariant vector field
   {left invariant vector field}%
   {left invariant vector}%
\Index%661%L%left module over D- algebra A
   {left module over $D$\Hyph algebra $A$}%
   {left module over algebra}%
\Index%237%L%left module over ring
   {left module over a ring $D$}%
   {left module over ring}%
\Index%731%L%left shift of R- module
   {left shift of $R$\Hyph module}%
   {left shift of module}%
\Index%240%L%left shift on fibered group
   {left shift on fibered group}%
   {Tstar shift, fibered group}%
\Index%238%L%left shift on group
   {left shift on group}%
   {left shift}%
\Index%239%L%left shift on group
   {left shift on group}%
   {left shift, group}%
\Index%236%L%left structural constant of Lie algebra
   {left structural constant of Lie algebra}%
   {left structural constant of Lie algebra}%
\Index%228%L%left vector space
   {left vector space}%
   {left vector space}%
\Index%751%L%left-ordered cycle notation of permutation
   {left-ordered cycle notation of permutation}%
   {left-ordered cycle notation of permutation}%
\Index%231%L%left-side contravariant representation of group
   {left-side contravariant representation of group}%
   {left-side contravariant representation of group}%
\Index%230%L%left-side covariant representation of group
   {left-side covariant representation of group}%
   {left-side covariant representation of group}%
\Index%233%L%left-side representation of fibered Omega- algebra
   {left-side representation of fibered $\Omega$\Hyph algebra}%
   {left-side representation of fibered Omega-algebra}%
\Index%232%L%left-side representation of Omega1 algebra A in Omega2 algebra M
   {left-side representation of $\Omega_1$\Hyph algebra $A$ in $\Omega_2$\Hyph algebra $M$}%
   {left-side representation of algebra}%
\Index%234%L%left-side transformation
   {left-side transformation}%
   {left-side transformation}%
\Index%235%L%left-side transformation on bundle
   {left-side transformation on bundle}%
   {left-side transformation of bundle}%
\Index%104%L%Lie algebra of Lie group
   {Lie algebra of Lie group}%
   {algebra Lie group Lie}%
\Index%402%L%Lie derivative
   {Lie derivative}%
   {Lie derivative}%
\Index%404%L%Lie derivative of connection
   {Lie derivative of connection}%
   {Lie derivative of connection}%
\Index%403%L%Lie derivative of metric
   {Lie derivative of metric}%
   {Lie derivative of metric}%
\Index%118%L%Lie group basic operators
   {Lie group basic operators}%
   {Lie group basic operators}%
\Index%246%L%lift of correspondence
   {lift of correspondence}%
   {lift of correspondence}%
\Index%245%L%lift of mapping
   {lift of mapping}%
   {lift of map}%
\Index%365%L%limit of correspondence with respect to the filter
   {limit of correspondence with respect to the filter}%
   {limit of correspondence with respect to the filter}%
\Index%366%L%limit of filter
   {limit of filter}%
   {limit of filter}%
\Index%460%L%limit of sequence in normed algebra
   {limit of sequence in normed algebra}%
   {limit of sequence, normed algebra}%
\Index%467%L%limit of sequence in normed ring
   {limit of sequence in normed ring}%
   {limit of sequence, normed ring}%
\Index%795%L%limit of sequence in ring
   {limit of sequence in ring}%
   {limit of sequence, valued ring}%
\Index%364%L%limit of sequence in valued division ring
   {limit of sequence in valued division ring}%
   {limit of sequence, valued division ring}%
\Index%367%L%limit set of filter
   {limit set of filter}%
   {limit set of filter}%
\Index%686%L%linear combination of vectors of A- module
   {linear combination of vectors of $A$\Hyph module}%
   {linear combination in A module}%
\Index%688%L%linear dependent vectors of A- module
   {linear dependent vectors of $A$\Hyph module}%
   {linear dependent vectors, module}%
\Index%803%L%linear geometric object
   {linear geometric object}%
   {linear geometric object}%
\Index%687%L%linear independent vectors of A- module
   {linear independent vectors of $A$\Hyph module}%
   {linear independent vectors, module}%
\Index%242%L%linear map of D- vector spaces over field F
   {linear map of $D$\Hyph vector spaces over field $F$}%
   {linear map over field, vector space}%
\Index%610%L%linear mapping of $R$\Hyph algebra $A_1$ into R algebra
   {linear mapping of $R$\Hyph algebra $A_1$ into $R$\Hyph algebra $A_2$}%
   {linear mapping of R algebras}%
\Index%613%L%linear mapping of $R$\Hyph module $A_1$ into $R$\Hyph module $A_2$
   {linear mapping of $R$\Hyph module $A_1$ into $R$\Hyph module $A_2$}%
   {linear mapping of R modules}%
\Index%97%L%linear mapping of division ring
   {linear mapping of division ring}%
   {linear mapping of division ring}%
\Index%98%L%linear mapping of division ring generated by mapping G
   {linear mapping of division ring generated by mapping $G$}%
   {linear map generated by map, division ring}%
\Index%608%L%linear mapping of R_1 module into R_2 module
   {linear mapping of $R_1$\Hyph module $A_1$ into $R_2$\Hyph module $A_2$}%
   {linear mapping of R_1 module into R_2 module}%
\Index%244%L%linear representation of group
   {linear representation of group}%
   {linear representation of group}%
\Index%801%L%linear representation of Lie group
   {linear representation of Lie group}%
   {Linear Representation of Lie Group}%
\Index%740%L%linear transformation of rcd affine space
   {linear transformation of \rcd affine space}%
   {linear transformation, rcd affine space}%
\Index%250%L%little group
   {little group}%
   {little group}%
\Index%247%L%local reference frame
   {local reference frame}%
   {local reference frame}%
\Index%248%L%locally compact at point p space
   {locally compact at point $p$ space}%
   {locally compact at point space}%
\Index%249%L%locally compact space
   {locally compact space}%
   {locally compact space}%
\Index%373%L%Lorentz transformation
   {Lorentz transformation}%
   {Lorentz transformation}%
\SetIndexSpace%
\Index%56%M%m- dimensional parallelepiped
   {$m$\Hyph dimensional parallelepiped}%
   {m dimensional parallelepiped}%
\Index%57%M%m- vector
   {$m$\Hyph vector}%
   {m-vector}%
\Index%556%M%manifold with D- affine connections
   {manifold with $D$\Hyph affine connections}%
   {manifold with D- affine connections}%
\Index%316%M%map of type G on manifold
   {map of type $G$ on manifold}%
   {map of type G on manifold}%
\Index%572%M%map polylinear over finite dimensional algebras
   {map polylinear over finite dimensional algebras}%
   {map polylinear over finite dimensional algebras}%
\Index%312%M%mapping of rings polylinear over commutative ring
   {mapping of rings polylinear over commutative ring}%
   {map polylinear over commutative ring, ring}%
\Index%407%M%mapping space
   {mapping space}%
   {mapping space}%
\Index%681%M%matrix of Acr linear mapping of module
   {matrix of \Acr linear mapping of module}%
   {matrix of Acr linear map}%
\Index%772%M%matrix of antilinear homomorphism
   {matrix of antilinear homomorphism}%
   {matrix of antilinear homomorphism}%
\Index%584%M%matrix of Astar linear map
   {matrix of $A\star$\Hyph linear map}%
   {matrix of Astar linear map}%
\Index%252%M%matrix of bilinear function
   {matrix of bilinear function}%
   {matrix of bilinear function}%
\Index%595%M%matrix of endomorphisms of Omega algebra
   {matrix of endomorphisms of $\Omega$\Hyph algebra}%
   {matrix of endomorphisms of Omega algebra}%
\Index%253%M%matrix of fibered Drc linear map
   {matrix of fibered \Drc linear map}%
   {matrix of fibered drc linear map}%
\Index%770%M%matrix of linear homomorphism
   {matrix of linear homomorphism}%
   {matrix of linear homomorphism}%
\Index%591%M%matrix of linear mappings
   {matrix of linear mappings}%
   {matrix of linear mappings}%
\Index%588%M%matrix of mappings
   {matrix of mappings}%
   {matrix of mappings}%
\Index%552%M%matrix of quadratic map
   {matrix of quadratic map}%
   {matrix of quadratic map, division ring}%
\Index%251%M%matrix of rcd linear mapping
   {matrix of \rcd linear mapping}%
   {matrix of rcd linear map}%
\Index%96%M%metric tensor in Minkowski space
   {metric tensor in Minkowski space}%
   {metric tensor, Minkowski space}%
\Index%254%M%metric-affine manifold
   {metric-affine manifold}%
   {metric-affine manifold}%
\Index%635%M%Minkowski space
   {Minkowski space}%
   {Minkowski space, Finsler}%
\Index%611%M%module over ring
   {module over ring}%
   {module over ring}%
\Index%265%M%morphism from tower of Tstar- representations into tower of Tstar- representations
   {morphism from tower of \Ts representations into tower of \Ts representations}%
   {morphism from tower of representations into tower of representations}%
\Index%269%M%morphism of fibered Tstar- representations from F into G
   {morphism of fibered \Ts representations from $\Bundle F$ into $\Bundle G$}%
   {morphism of fibered representations from f into g}%
\Index%732%M%morphism of representation f
   {morphism of representation $f$}%
   {morphism of representation f}%
\Index%268%M%morphism of representations from f into g
   {morphism of representations from $f$ into $g$}%
   {morphism of representations from f into g}%
\Index%267%M%morphism of representations of Omega1 algebra in Omega2 algebra
   {morphism of representations of $\Omega_1$\Hyph algebra in $\Omega_2$\Hyph algebra}%
   {morphism of representations of Omega1 algebra in Omega2 algebra}%
\Index%264%M%morphism of Tstar- representations of fibered Omega- algebra
   {morphism of \Ts representations of fibered $\Omega$\Hyph algebra}%
   {morphism of representations of fibered Omega algebra}%
\Index%799%M%motion of Minkowski space
   {motion of Minkowski space}%
   {motion, Minkowski space}%
\Index%147%M%movement on basis manifold
   {movement on basis manifold}%
   {movement transformation}%
\SetIndexSpace%
\Index%726%N%n algebra over ring
   {$n$\Hyph algebra over the ring}%
   {n algebra over ring}%
\Index%58%N%n- ary fibered relation
   {$n$\Hyph ary fibered relation}%
   {fibered relation}%
\Index%279%N%nonmetricity
   {nonmetricity}%
   {nonmetricity}%
\Index%272%N%nonsingular bilinear function
   {nonsingular bilinear function}%
   {nonsingular bilinear function}%
\Index%273%N%nonsingular system of drc linear equations
   {nonsingular system of \rcd linear equations}%
   {nonsingular system of linear equations}%
\Index%643%N%nonsingular tensor
   {nonsingular tensor}%
   {nonsingular tensor, algebra}%
\Index%274%N%nonsingular transformation
   {nonsingular transformation}%
   {nonsingular transformation}%
\Index%742%N%norm in quaternion algebra
   {norm in quaternion algebra}%
   {norm, quaternion algebra}%
\Index%286%N%norm of map \Vector A of normed D-vector space
   {norm of map $\Vector A$ of normed $D$\hyph vector space}%
   {norm of map, D vector space}%
\Index%314%N%norm of mapping into D- algebra
   {norm of mapping into $D$\Hyph algebra}%
   {norm of map, algebra}%
\Index%285%N%norm of mapping of division ring
   {norm of mapping of division ring}%
   {norm of map, division ring}%
\Index%346%N%norm on D- algebra
   {norm on $D$\Hyph algebra}%
   {norm on D algebra}%
\Index%283%N%norm on D- vector space
   {norm on $D$\Hyph vector space}%
   {norm on D vector space}%
\Index%641%N%norm on ring
   {norm on ring}%
   {norm on ring}%
\Index%343%N%normed D- algebra
   {normed $D$\Hyph algebra}%
   {normed D algebra}%
\Index%287%N%normed D- vector space
   {normed $D$\Hyph vector space}%
   {normed D vector space}%
\Index%640%N%normed ring
   {normed ring}%
   {normed ring}%
\Index%813%N%not complete group
   {not complete group}%
   {not complete group}%
\Index%806%N%not complete Omega algebra
   {not complete $\Omega$\Hyph algebra}%
   {not complete Omega algebra}%
\Index%619%N%nucleus of algebra
   {nucleus of $R$\Hyph algebra $A$}%
   {nucleus of algebra}%
\SetIndexSpace%
\Index%639%O%octonion algebra
   {octonion algebra}%
   {octonion algebra}%
\Index%644%O%Omega group
   {$\Omega$\Hyph group}%
   {Omega group}%
\Index%646%O%Omega linear mapping
   {$\Omega$\Hyph linear mapping}%
   {Omega linear map}%
\Index%707%O%Omega_2- word of element of representation relative to generating set
   {$\Omega_2$\Hyph word of element of representation relative to generating set}%
   {word of element relative to generating set, representation}%
\Index%300%O%operation on bundle
   {operation on bundle}%
   {operation on bundle}%
\Index%623%O%opposite algebra to algebra
   {opposite algebra to algebra $P$}%
   {opposite algebra}%
\Index%408%O%opposite fibered preordering
   {opposite fibered preordering}%
   {opposite fibered preordering}%
\Index%624%O%orbit of linear mapping
   {orbit of linear mapping}%
   {orbit of linear mapping}%
\Index%306%O%orbit of representation of fibered group
   {orbit of representation of fibered group}%
   {orbit of representation of fibered group}%
\Index%305%O%orbit of representation of group
   {orbit of representation of group}%
   {orbit of representation of group}%
\Index%304%O%orbit of Tstar- representation of group
   {orbit of \Ts representation of group}%
   {orbit of Tstar  representation of group}%
\Index%729%O%origin of coordinate system of affine space
   {origin of coordinate system of affine space}%
   {origin of coordinate system of affine space}%
\Index%738%O%origin of coordinate system of star D affine space
   {origin of coordinate system of $\star D$\Hyph affine space}%
   {origin of coordinate system of starD affine space}%
\Index%637%O%orthogonal basis in Minkowski space
   {orthogonal basis in Minkowski space}%
   {orthogonal basis, Minkowski space}%
\Index%636%O%orthogonality in Minkowski space
   {orthogonality in Minkowski space}%
   {Minkowski orthogonality}%
\Index%638%O%orthonormal basis in Minkowski space
   {orthonormal basis in Minkowski space}%
   {orthonormal basis, Minkowski space}%
\Index%307%O%orthonornal basis
   {orthonornal basis}%
   {Orthonornal Basis}%
\Index%308%O%orthonornal basis
   {orthonornal basis}%
   {Orthonornal Basis, division ring}%
\SetIndexSpace%
\Index%695%P%active representation of group in basis manifold of representation
   {passive representation of group $G(f)$ in basis manifold of representation}%
   {passive representation in basis manifold}%
\Index%739%P%parallel shift of rcd affine space
   {parallel shift of \rcd affine space}%
   {parallel shift, rcd affine space}%
\Index%317%P%parallelogram
   {parallelogram}%
   {parallelogram}%
\Index%513%P%partial linear mapping of variable v^i
   {partial linear mapping of variable $v^i$}%
   {partial linear map of variable}%
\Index%333%P%passive representation
   {passive representation}%
   {passive representation}%
\Index%718%P%passive representation of group in basis manifold of tower of representations
   {passive representation of group $G(\Vector f)$ in basis manifold of tower of representations}%
   {passive representation in basis manifold, tower of representations}%
\Index%321%P%passive starT- representation
   {passive \sT representation}%
   {passive starT representation}%
\Index%570%P%passive transformation of the basis manifold of representation
   {passive transformation of the basis manifold of representation}%
   {passive transformation of basis, representation}%
\Index%717%P%passive transformation of the basis manifold of tower of representations
   {passive transformation of the basis manifold of tower of representations}%
   {passive transformation of basis, tower of representations}%
\Index%335%P%passive transformation on basis manifold
   {passive transformation on basis manifold}%
   {passive transformation}%
\Index%336%P%passive transformation on the set of rcd bases
   {passive transformation on the set of \rcd bases}%
   {passive transformation, vector space}%
\Index%744%P%permutability property of trace
   {permutability property of trace}%
   {permutability property of trace}%
\Index%416%P%pfaffian derivative
   {pfaffian derivative}%
   {pfaffian derivative}%
\Index%342%P%polylinear mapping of (n)-D-vector spaces
   {polylinear mapping of $(n)$\hyph $D$\hyph vector spaces}%
   {polylinear map of D vector spaces}%
\Index%609%P%polylinear mapping of algebras
   {polylinear mapping of algebras}%
   {polylinear map of algebras}%
\Index%621%P%polylinear mapping of modules
   {polylinear mapping of modules}%
   {polylinear map of modules}%
\Index%345%P%polylinear skew symmetric map
   {polylinear skew symmetric map}%
   {polylinear map skew symmetric, division ring}%
\Index%347%P%polylinear symmetric map
   {polylinear symmetric map}%
   {polylinear map symmetric, division ring}%
\Index%760%P%polymorphism of representations
   {polymorphism of representations}%
   {polymorphism of representations}%
\Index%344%P%polyvector
   {polyvector}%
   {polyvector}%
\Index%352%P%potential energy
   {potential energy}%
   {potential energy}%
\Index%388%P%product of geometric object and constant
   {product of geometric object and constant}%
   {product of geometric object and constant}%
\Index%389%P%product of geometric object and constant in vector space
   {product of geometric object and constant in vector space}%
   {product of geometric object and constant, vector space}%
\Index%390%P%product of groups
   {product of groups}%
   {product of groups}%
\Index%391%P%product of morphisms of representations of universal algebra
   {product of morphisms of representations of universal algebra}%
   {product of morphisms of representations of universal algebra}%
\Index%561%P%product of morphisms of tower of representations
   {product of morphisms of tower of representations}%
   {product of morphisms of tower of representations}%
\Index%392%P%product of morphisms of Tstar- representations of fibered Omega- algebra
   {product of morphisms of \Ts representations of fibered $\Omega$\Hyph algebra}%
   {product of morphisms of representations of fibered Omega algebra}%
\Index%395%P%product of objects in category
   {product of objects in category}%
   {product of objects in category}%
\Index%387%P%projection of bundle E along fiber E
   {projection of bundle $\Bundle E$ along fiber $E$}%
   {projection of bundle along fiber}%
\Index%540%P%pseudo-Euclidean metric on division ring
   {pseudo\Hyph Euclidean metric on division ring}%
   {pseudo-Euclidean metric on division ring}%
\Index%548%P%pseudo-Euclidean scalar product in D- vector space
   {pseudo\Hyph Euclidean scalar product in $D$\Hyph vector space}%
   {pseudo-Euclidean scalar product, vector space}%
\Index%543%P%pseudo-Euclidean scalar product on division ring
   {pseudo-Euclidean scalar product on division ring}%
   {pseudo-Euclidean scalar product on division ring}%
\SetIndexSpace%
\Index%555%Q%quadratic form in division ring
   {quadratic form in division ring}%
   {quadratic form, division ring}%
\Index%551%Q%quadratic map of division ring
   {quadratic map of division ring}%
   {Quadratic Map of Division Ring}%
\Index%180%Q%quasi affine transformation on basis manifold
   {quasi affine transformation on basis manifold}%
   {quasi affine transformation}%
\Index%181%Q%quasi affine transformation on basis manifold
   {quasi affine transformation on basis manifold}%
   {quasi affine drc transformation}%
\Index%182%Q%quasi movement on basis manifold
   {quasi movement on basis manifold}%
   {quasi movement, division ring}%
\Index%183%Q%quasi movement on basis manifold
   {quasi movement on basis manifold}%
   {quasi movement}%
\Index%604%Q%quasiclosed ring of mappings
   {quasiclosed ring of mappings}%
   {quasiclosed ring of mappings}%
\Index%679%Q%quasideterminant
   {quasideterminant}%
   {quasideterminant definition}%
\Index%800%Q%quasimotion of Minkowski space
   {quasimotion of Minkowski space}%
   {Quasimotion, Minkowski space}%
\Index%741%Q%quaternion algebra
   {quaternion algebra}%
   {quaternion algebra}%
\Index%103%Q%quaternion algebra E over the field F
   {quaternion algebra $E$ over the field $F$}%
   {quaternion algebra over the field}%
\Index%496%Q%quotient bundle
   {quotient bundle}%
   {quotient bundle}%
\SetIndexSpace%
\Index%49%R%drc basis dual to rcd basis of vector space
   {\drc basis dual to \rcd basis of vector space}%
   {basis dual to basis, rcd vector space}%
\Index%4%R%j i-RC quasideterminant
   {$(^j_i)$\hyph \RC quasideterminant}%
   {j i RC-quasideterminant}%
\Index%59%R%r row of matrix
   {\sups row of matrix}%
   {r row}%
\Index%728%R%R- module
   {$R$\Hyph module}%
   {R- module}%
\Index%72%R%r-row of matrix
   {$r$\hyph row of matrix}%
   {r-row}%
\Index%749%R%rank of Hermitian matrix by principal minors
   {rank of Hermitian matrix by principal minors}%
   {rank of Hermitian matrix by principal minors}%
\Index%554%R%rank of quadratic map of division ring
   {rank of quadratic map of division ring}%
   {rank of quadratic map, division ring}%
\Index%65%R%RC inverse element of biring
   {\RC inverse element of biring}%
   {rc-inverse element}%
\Index%61%R%RC major minor
   {\RC major minor}%
   {RC-major minor}%
\Index%63%R%RC matrix group
   {\RC matrix group}%
   {rc-matrix group}%
\Index%64%R%RC nonsingular matrix
   {\RC nonsingular matrix}%
   {RC nonsingular matrix}%
\Index%68%R%RC power
   {\RC power}%
   {rc power}%
\Index%66%R%RC product of matrices
   {\RC product of matrices}%
   {rc-product of matrices}%
\Index%590%R%rc product of matrices of mappings
   {$\RCcirc$\Hyph product of matrices of mappings}%
   {rc product of matrices of mappings}%
\Index%62%R%RC quasideterminant
   {\RC quasideterminant}%
   {RC-quasideterminant}%
\Index%67%R%RC rank of matrix
   {\RC rank of matrix}%
   {rc-rank of matrix}%
\Index%60%R%RC singular matrix
   {\RC singular matrix}%
   {RC singular matrix}%
\Index%677%R%RCcirc A- basis for module
   {$\RCcirc A$\Hyph basis for module}%
   {RCcirc A basis, module over algebra}%
\Index%676%R%RCcirc A- linearly dependent set of vectors
   {$\RCcirc A$\Hyph linearly dependent set of vectors}%
   {RCcirc linearly dependent, starA module over algebra}%
\Index%675%R%RCcirc A- linearly independent set of vectors
   {$\RCcirc A$\Hyph linearly independent set of vectors}%
   {RCcirc linearly independent, starA module over algebra}%
\Index%601%R%RCcirc nonsingular matrix of A(A) mappings
   {$\RCcirc$\Hyph nonsingular matrix of $\mathcal A(A)$\Hyph mappings}%
   {RCcirc nonsingular matrix of A(A) mappings}%
\Index%598%R%RCcirc nonsingular matrix of endomorphisms
   {$\RCcirc$\Hyph nonsingular matrix of endomorphisms}%
   {RCcirc nonsingular matrix of endomorphisms}%
\Index%602%R%RCcirc nonsingular system of additive equations
   {$\RCcirc$\Hyph nonsingular system of additive equations}%
   {RCcirc nonsingular system of additive equations}%
\Index%597%R%RCcirc quasideterminant
   {$\RCcirc$\Hyph quasideterminant}%
   {RCcirc-quasideterminant definition}%
\Index%600%R%RCcirc singular matrix of A(A) mappings
   {$\RCcirc$\Hyph singular matrix of $\mathcal A(A)$\Hyph mappings}%
   {RCcirc singular matrix of A(A) mappings}%
\Index%599%R%RCcirc singular matrix of endomorphisms
   {$\RCcirc$\Hyph singular matrix of endomorphisms}%
   {RCcirc singular matrix of endomorphisms}%
\Index%765%R%RCcircA linear combination
   {$\RCcirc A$\Hyph linear combination}%
   {RCcircA linear combination}%
\Index%111%R%rcd affine basis
   {\rcd affine basis}%
   {rcd affine basis, division ring}%
\Index%759%R%rcd affine plane
   {\rcd affine plane}%
   {rcd affine plane}%
\Index%736%R%rcd affine space
   {\rcd affine space}%
   {rcd affine space}%
\Index%758%R%rcd affine transformation
   {\rcd affine transformation}%
   {rcd affine transformation}%
\Index%22%R%rcd automorphism of vector space
   {\rcd automorphism of vector space}%
   {automorphism of vector space}%
\Index%24%R%rcd basis for vector space
   {\rcd basis for vector space}%
   {rcd basis, vector space}%
\Index%31%R%rcd isomorphism of vector spaces
   {\rcd isomorphism of vector spaces}%
   {isomorphism of vector spaces}%
\Index%35%R%rcd linear map of vector spaces
   {\rcd linear map of vector spaces}%
   {rcd linear map of vector spaces}%
\Index%37%R%rcd linear Tstar- representation of group
   {\rcd linear \Ts representation of group}%
   {rcd linear Tstar representation of group}%
\Index%21%R%rcd linearly dependent vectors
   {\rcd linearly dependent vectors}%
   {linearly dependent, vector space}%
\Index%34%R%rcd linearly independent vectors
   {\rcd linearly independent vectors}%
   {linearly independent, vector space}%
\Index%27%R%rcd vector
   {\rcd vector}%
   {rcd vector}%
\Index%69%R%rcd vector space
   {\rcd vector space}%
   {rcd vector space}%
\Index%377%R%reduced Cartesian product of bundles
   {reduced Cartesian product of bundles}%
   {reduced Cartesian product of bundles}%
\Index%378%R%reduced Cartesian product of total spaces
   {reduced Cartesian product of total spaces}%
   {reduced Cartesian product of total spaces}%
\Index%379%R%reduced fibered correspondence from A to B
   {reduced fibered correspondence from $\Bundle{A}$ to $\Bundle B$}%
   {reduced fibered correspondence from A to B}%
\Index%380%R%reduced fibered correspondence in A
   {reduced fibered correspondence in $\Bundle{A}$}%
   {reduced fibered correspondence in A}%
\Index%761%R%reduced polymorphism of representations
   {reduced polymorphism of representations}%
   {reduced polymorphism of representations}%
\Index%381%R%reducible biring
   {reducible biring}%
   {reducible biring}%
\Index%450%R%reference frame in event space
   {reference frame in event space}%
   {reference frame in event space}%
\Index%263%R%reference frame manifold
   {reference frame manifold}%
   {reference frame manifold}%
\Index%439%R%reflexive 2- ary fibered relation
   {reflexive $2$\Hyph ary fibered relation}%
   {reflexive 2 ary fibered relation}%
\Index%330%R%regular endomorphism of representation
   {regular endomorphism of representation}%
   {regular endomorphism of representation}%
\Index%566%R%regular endomorphism of tower of representations
   {regular endomorphism of tower of representations}%
   {regular endomorphism of tower of representations}%
\Index%553%R%regular quadratic map in division ring
   {regular quadratic map in division ring}%
   {regular quadratic map, division ring}%
\Index%372%R%representation of group
   {representation of group}%
   {representation of group}%
\Index%528%R%representation of Omega algebra in representation
   {representation of $\Omega$\Hyph algebra in representation}%
   {representation of Omega algebra in representation}%
\Index%529%R%representation of Omega algebra in tower of representations
   {representation of $\Omega$\Hyph algebra in tower of representations}%
   {representation of Omega algebra in tower of representations}%
\Index%370%R%representation of Omega- algebra A in category B
   {representation of $\Omega$\Hyph algebra $A$ in category $\mathcal B$}%
   {representation of Omega algebra in category}%
\Index%371%R%representation of Omega1 algebra A in Omega2 algebra M
   {representation of $\Omega_1$\Hyph algebra $A$ in $\Omega_2$\Hyph algebra $M$}%
   {representation of algebra}%
\Index%702%R%representative of geometric object in drc vector space
   {representative of geometric object in \drc vector space}%
   {representative of geometric object, drc vector space}%
\Index%703%R%representative of geometric object in Omega_2- algebra
   {representative of geometric object in $\Omega_2$\Hyph algebra}%
   {representative of geometric object, representation g}%
\Index%368%R%representative of geometric object in rcd vector space
   {representative of geometric object in \rcd vector space}%
   {representative of geometric object, rcd vector space}%
\Index%724%R%representative of geometric object in tuple of Omega- algebras
   {representative of geometric object in tuple of $\VX\Omega$\Hyph algebras}%
   {representative of geometric object, tower of representations g}%
\Index%369%R%representative of geometric object in vector space
   {representative of geometric object in vector space}%
   {representative of geometric object, vector space}%
\Index%474%R%restriction of correspondence \Phi to set C
   {restriction of correspondence $\Phi$ to set $C$}%
   {restriction of correspondence}%
\Index%753%R%right cofactor of entry of matrix
   {right cofactor of entry of matrix}%
   {right cofactor, matrix}%
\Index%302%R%right defined Lie algebra of Lie group
   {right defined Lie algebra of Lie group}%
   {right defined Lie algebra}%
\Index%757%R%right double cofactor of entry of matrix
   {right double cofactor of entry of matrix}%
   {right double cofactor}%
\Index%354%R%right invariant vector field
   {right invariant vector field}%
   {right invariant vector}%
\Index%670%R%right module over D- algebra A
   {right module over $D$\Hyph algebra $A$}%
   {right module over algebra}%
\Index%361%R%right module over ring
   {right module over a ring $D$}%
   {right module over ring}%
\Index%362%R%right shift on group
   {right shift on group}%
   {right shift}%
\Index%363%R%right shift on group
   {right shift on group}%
   {right shift, group}%
\Index%360%R%right structural constant of Lie algebra
   {right structural constant of Lie algebra}%
   {right structural constant of Lie algebra}%
\Index%353%R%right vector space
   {right vector space}%
   {right vector space}%
\Index%752%R%right-ordered cycle notation of permutation
   {right-ordered cycle notation of permutation}%
   {right-ordered cycle notation of permutation}%
\Index%356%R%right-side contravariant representation of group
   {right-side contravariant representation of group}%
   {right-side contravariant representation of group}%
\Index%355%R%right-side covariant representation of group
   {right-side covariant representation of group}%
   {right-side covariant representation of group}%
\Index%358%R%right-side representation of fibered Omega- algebra
   {right-side representation of fibered $\Omega$\Hyph algebra}%
   {right-side representation of fibered Omega-algebra}%
\Index%357%R%right-side representation of Omega1 algebra A in Omega2 algebra M
   {right-side representation of $\Omega_1$\Hyph algebra $A$ in $\Omega_2$\Hyph algebra $M$}%
   {right-side representation of algebra}%
\Index%359%R%right-side transformation
   {right-side transformation}%
   {right-side transformation}%
\Index%188%R%ring has characteristic 0
   {ring has characteristic $0$}%
   {ring has characteristic 0}%
\Index%189%R%ring has characteristic p
   {ring has characteristic $p$}%
   {ring has characteristic p}%
\Index%778%R%ring with conjugation
   {ring with conjugation}%
   {ring with conjugation}%
\Index%748%R%row determinant
   {row determinant}%
   {row determinant}%
\Index%127%R%row vector
   {row vector}%
   {row vector}%
\SetIndexSpace%
\Index%669%S%$starA- module
   {$\star A$\Hyph module}%
   {starA-module}%
\Index%77%S%Dstar- product of rcd linear map A over scalar
   {$D\star$\Hyph product of \rcd linear map $A$ over scalar}%
   {Dstar product of rcd linear map over scalar}%
\Index%73%S%RCstarS,RCstarT- linear map of vector spaces
   {$(\RCstar S,\RCstar T)$\Hyph linear map of vector spaces}%
   {rcs rct linear map of vector spaces}%
\Index%786%S%scalar algebra of algebra
   {scalar algebra of algebra}%
   {scalar algebra of algebra}%
\Index%775%S%scalar algebra of ring
   {scalar algebra of ring}%
   {scalar algebra of ring}%
\Index%784%S%scalar of element of algebra
   {scalar of element of algebra}%
   {scalar of algebra}%
\Index%773%S%scalar of element of ring
   {scalar of element of ring}%
   {scalar of ring}%
\Index%781%S%scalar of mapping
   {scalar of mapping}%
   {scalar of mapping}%
\Index%451%S%scalar potential
   {scalar potential}%
   {scalar potential}%
\Index%791%S%Schauder basis
   {Schauder basis}%
   {Schauder basis}%
\Index%129%S%second Newton law
   {second Newton law}%
   {Second Newton law}%
\Index%694%S%section of bundle
   {section of bundle}%
   {section of bundle}%
\Index%562%S%set of coordinates of representation
   {set of coordinates of representation}%
   {coordinate set of representation}%
\Index%571%S%set of Omega_2- words of representation
   {set of $\Omega_2$\Hyph words of representation}%
   {word set of representation}%
\Index%563%S%set of tuples of coordinates of tower of representations
   {set of tuples of coordinates of tower of representations}%
   {coordinate set of tower of representations}%
\Index%709%S%set of tuples of Omega- words of tower of representations
   {set of tuples of $\Vector\Omega$\Hyph words of tower of representations}%
   {word set of tower of representations}%
\Index%405%S%simple polyvector
   {simple polyvector}%
   {simple polyvector}%
\Index%299%S%single transitive representation of fibered Omega- algebra
   {single transitive representation of fibered $\Omega$\Hyph algebra}%
   {single transitive representation of fibered Omega-algebra}%
\Index%298%S%single transitive representation of group
   {single transitive representation of group}%
   {single transitive representation of group}%
\Index%297%S%single transitive representation of Omega- algebra A
   {single transitive representation of $\Omega$\Hyph algebra $A$}%
   {single transitive representation of algebra}%
\Index%130%S%singular linear mapping of D- vector space
   {singular linear mapping of $D$\Hyph vector space}%
   {singular linear map, D vector space}%
\Index%131%S%singular linear mapping of division ring
   {singular linear mapping of division ring}%
   {singular linear map, division ring}%
\Index%226%S%skew product of vectors
   {skew product of vectors}%
   {skew product of vectors}%
\Index%651%S%skew symmetric polylinear mapping into associative algebra
   {skew symmetric polylinear mapping into associative algebra}%
   {polylinear map skew symmetric, associative algebra}%
\Index%406%S%space of orbits of Tstar- representation
   {space of orbits of \Ts representation}%
   {space of orbits of Ts representation}%
\Index%812%S%spacelike vector
   {spacelike vector}%
   {spacelike vector}%
\Index%452%S%speed of deviation
   {speed of deviation}%
   {speed of deviation}%
\Index%74%S%SRCstar,TRCstar- linear map of vector bundles
   {$(\mathcal S\RCstar,\mathcal T\RCstar)$\Hyph linear map of vector bundles}%
   {src trc linear map of vector bundles}%
\Index%112%S%Sstar, star T-bimodule
   {($S\star$, $\star T$)\hyph bimodule}%
   {(Sstar,starT)-bimodule}%
\Index%145%S%stability group
   {stability group}%
   {stability group}%
\Index%457%S%stable set of representation
   {stable set of representation}%
   {stable set of representation}%
\Index%458%S%standard component of Gateaux differential of map f
   {standard component of the G\^ateaux differential of map $f$}%
   {standard component of Gateaux differential, division ring}%
\Index%462%S%standard component of polylinear map f of division ring
   {standard component of polylinear map $f$ of division ring}%
   {standard component of polylinear map, division ring}%
\Index%649%S%standard component of polylinear mapping into associative algebra
   {standard component of polylinear mapping into associative algebra}%
   {standard component of polylinear map, associative algebra}%
\Index%459%S%standard component of quadratic map f over field F
   {standard component of quadratic map $f$ over field $F$}%
   {standard component of quadratic map, division ring}%
\Index%463%S%standard component of tensor
   {standard component of tensor}%
   {standard component of tensor, division ring}%
\Index%574%S%standard component of tensor in tensor product of algebras
   {standard component of tensor in tensor product of algebras}%
   {standard component of tensor, algebra}%
\Index%630%S%standard component of the Gateaux derivative of mapping f
   {standard component of the G\^ateaux derivative of mapping $f$}%
   {standard component of Gateaux derivative, algebra}%
\Index%461%S%standard component over field F of bilitnear map f
   {standard component over field $F$ of bilitnear map $f$}%
   {standard component of bilinear map, division ring}%
\Index%472%S%standard coordinates of basis
   {standard coordinates of basis}%
   {standard coordinates of basis}%
\Index%471%S%standard coordinates of rcd basis
   {standard coordinates of \rcd basis}%
   {standard coordinates of rcd basis}%
\Index%456%S%standard F- component of linear mapping f
   {standard $F$\Hyph component of linear mapping $f$}%
   {standard component of linear map, division ring}%
\Index%464%S%standard F- representation of linear mapping of division ring
   {standard $F$\Hyph representation of linear mapping of division ring}%
   {linear map, standard representation, division ring}%
\Index%465%S%standard representation of Gateaux differential of map of division ring over field F
   {standard representation of the G\^ateaux differential of map of division ring over field $F$}%
   {Gateaux differential, standard representation, division ring}%
\Index%468%S%standard representation of matrix
   {standard representation of matrix}%
   {Standard representation}%
\Index%470%S%standard representation of polylinear map of division ring
   {standard representation of polylinear map of division ring}%
   {polylinear map, standard representation, division ring}%
\Index%648%S%standard representation of polylinear mapping into associative algebra
   {standard representation of polylinear mapping into associative algebra}%
   {polylinear map, standard representation, associative algebra}%
\Index%466%S%standard representation of quadratic map of division ring over field F
   {standard representation of quadratic map of division ring over field $F$}%
   {quadratic map, standard representation, division ring}%
\Index%469%S%standard representation over field F of bilinear map of division ring
   {standard representation over field $F$ of bilinear map of division ring}%
   {bilinear map, standard representation, division ring}%
\Index%683%S%star A- product of Acr linear map over scalar
   {$\star A$\Hyph product of \Acr linear map $\Vector f$ over scalar}%
   {starA product of Acr linear map over scalar}%
\Index%586%S%star A- product of Astar linear map over scalar
   {$\star A$\Hyph product of $A\star$\Hyph linear mapping over scalar}%
   {starA product of Astar linear map over scalar}%
\Index%78%S%star D-vector space
   {$\star D$\hyph vector space}%
   {starD-vector space}%
\Index%79%S%star R-module
   {$\star R$\hyph module}%
   {starR-module}%
\Index%674%S%starA- product of vector over scalar
   {$\star A$\Hyph product of vector over scalar}%
   {starA product of vector over scalar, starA module}%
\Index%735%S%starD affine space
   {$\star D$\Hyph affine space}%
   {starD affine space}%
\Index%75%S%starD component of coordinates of vector \Vector r
   {\sD component of coordinates of vector $\Vector r$}%
   {starD component of coordinates of vector, D vector space}%
\Index%84%S%starT- representation of fibered group
   {\sT representation of fibered group}%
   {starT representation of fibered group}%
\Index%83%S%starT- representation of fibered Omega- algebra
   {\sT representation of fibered $\Omega$\Hyph algebra}%
   {starT representation of fibered Omega-algebra}%
\Index%82%S%starT- representation of Omega1 algebra A in Omega2 algebra M
   {\sT representation of $\Omega_1$\Hyph algebra $A$ in $\Omega_2$\Hyph algebra $M$}%
   {starT representation of algebra}%
\Index%80%S%starT- shift
   {\sT shift}%
   {starT shift}%
\Index%81%S%starT- shift on fibered group
   {\sT shift on fibered group}%
   {starT shift, fibered group}%
\Index%85%S%starT- transformation
   {\sT transformation}%
   {starT transformation}%
\Index%86%S%starT- transformation on bundle
   {\sT transformation on bundle}%
   {starT transformation of bundle}%
\Index%618%S%structural constants of algebra $A$ over ring $R$
   {structural constants of algebra $P$ over ring $D$}%
   {structural constants of algebra}%
\Index%473%S%structural constants of division ring D over field F
   {structural constants of division ring $D$ over field $F$}%
   {structural constants of division ring over field}%
\Index%340%S%subbundle
   {subbundle}%
   {subbundle}%
\Index%341%S%subbundle of Dstar-vector space
   {subbundle of $\mathcal D\star$\hyph vector space}%
   {subbundle of Dstar vector bundle}%
\Index%322%S%subrepresentation generated by set X
   {subrepresentation generated by set $X$}%
   {subrepresentation generated by set}%
\Index%334%S%subrepresentation of representation
   {subrepresentation of representation}%
   {subrepresentation of representation}%
\Index%682%S%sum of Acr linear mappings of module
   {sum of \Acr linear mappings of module}%
   {sum of Acr linear maps, module}%
\Index%477%S%sum of geometric objects in vector space
   {sum of geometric objects in vector space}%
   {sum of geometric objects, vector space}%
\Index%476%S%sum of geometrical objects
   {sum of geometric objects}%
   {sum of geometric objects}%
\Index%475%S%sum of rcd linear maps
   {sum of \rcd linear maps}%
   {sum of rcd linear maps, rcd vector spaces}%
\Index%323%S%superposition of coordinates of the representation and the element
   {superposition of coordinates of the representation $f$ and the element $m$}%
   {superposition of coordinates, representation}%
\Index%712%S%superposition of coordinates of the tower of representations and the element
   {superposition of coordinates of the tower of representations $\Vector f$ and the element $\VX a$}%
   {superposition of coordinates, tower of representations}%
\Index%446%S%symmetric 2- ary fibered relation
   {symmetric $2$\Hyph ary fibered relation}%
   {symmetric 2 ary fibered relation}%
\Index%447%S%symmetric bilinear map of D- vector space to division ring
   {symmetric bilinear map of $D$\Hyph vector space to division ring}%
   {symmetric bilinear map, vector space to division ring}%
\Index%650%S%symmetric polylinear mapping into associative algebra
   {symmetric polylinear mapping into associative algebra}%
   {polylinear map symmetric, associative algebra}%
\Index%143%S%symmetry group
   {symmetry group}%
   {symmetry group}%
\Index%144%S%symmetry group
   {symmetry group}%
   {SymmetryGroup}%
\Index%448%S%synchronization of reference frame
   {synchronization of reference frame}%
   {synchronization of reference frame}%
\Index%605%S%system of additive equations
   {system of additive equations}%
   {system of additive equations}%
\Index%449%S%system of drc linear equations
   {system of \drc linear equations}%
   {system of drc linear equations}%
\Index%311%S%system of rcd linear equations
   {system of \rcd linear equations}%
   {system of rcd linear equations}%
\Index%631%S%tandard representation of Gateaux derivative of mapping over algebra
   {tandard representation of the G\^ateaux derivative of mapping over algebra}%
   {Gateaux derivative, standard representation, algebra}%
\SetIndexSpace%
\Index%499%T%Taylor polynomial
   {Taylor polynomial}%
   {Taylor polynomial, division ring}%
\Index%440%T%Taylor series
   {Taylor series}%
   {Taylor series, division ring}%
\Index%652%T%tensor power of algebra
   {tensor power of algebra}%
   {tensor power of algebra}%
\Index%763%T%tensor power of representation
   {tensor power of representation}%
   {tensor power of representation}%
\Index%573%T%tensor product of algebras
   {tensor product of algebras}%
   {tensor product of algebras}%
\Index%480%T%tensor product of D- vector spaces
   {tensor product of $D$\Hyph vector spaces}%
   {tensor product of D vector spaces}%
\Index%484%T%tensor product of division rings
   {tensor product of division rings}%
   {tensor product of division rings}%
\Index%481%T%tensor product of Dstar vector spaces
   {tensor product of \Ds vector spaces}%
   {tensor product of Dstar vector spaces}%
\Index%483%T%tensor product of representations
   {tensor product of representations}%
   {tensor product of representations}%
\Index%762%T%tensor product of representations
   {tensor product of representations}%
   {tensor product of representations}%
\Index%482%T%tensor product of rings over commutative ring
   {tensor product of rings over commutative ring}%
   {tensor product of rings}%
\Index%47%T%the Fr\'echet Dstar derivative of map f of division ring D at point x
   {the Fr\'echet \Ds derivative of map $f$ of division ring $D$ at point $x$}%
   {Frechet Dstar derivative of map, division ring}%
\Index%811%T%timelike vector
   {timelike vector}%
   {timelike vector}%
\Index%486%T%topological D- vector space
   {topological $D$\Hyph vector space}%
   {topological D vector space}%
\Index%310%T%topological D-algebra
   {topological $D$\Hyph algebra}%
   {topological D algebra}%
\Index%488%T%topological division ring
   {topological division ring}%
   {topological division ring}%
\Index%487%T%topological drc vector space
   {topological \drc vector space}%
   {topological drc vector space}%
\Index%790%T%topological ring
   {topological ring}%
   {topological ring}%
\Index%498%T%torsion form
   {torsion form}%
   {torsion form}%
\Index%478%T%torsion tensor
   {torsion tensor}%
   {torsion tensor}%
\Index%120%T%tower of bundles
   {tower of bundles}%
   {tower of bundles}%
\Index%324%T%tower of effective representations
   {tower of effective representations}%
   {tower of effective representations}%
\Index%119%T%tower of representations of Omega algebras
   {tower of representations of $\Vector{\Omega}$\Hyph algebras}%
   {tower of representations of algebras}%
\Index%531%T%tower of subrepresentations
   {tower of subrepresentations}%
   {tower of subrepresentations}%
\Index%532%T%tower of subrepresentations of tower of representations generated by tuple of sets
   {tower of subrepresentations of tower of representations $\Vector f$ generated by tuple of sets $\VX X$}%
   {subrepresentation generated by tuple of sets}%
\Index%743%T%trace of quaternion
   {trace of quaternion}%
   {trace, quaternion algebra}%
\Index%376%T%transformation coordinated with equivalence
   {transformation coordinated with equivalence}%
   {transformation coordinated with equivalence}%
\Index%374%T%transformation of universal algebra
   {transformation of universal algebra}%
   {transformation of universal algebra}%
\Index%375%T%transformation on bundle
   {transformation on bundle}%
   {transformation of bundle}%
\Index%489%T%transitive 2- ary fibered relation
   {transitive $2$\Hyph ary fibered relation}%
   {transitive 2 ary fibered relation}%
\Index%492%T%transitive representation of fibered Omega- algebra
   {transitive representation of fibered $\Omega$\Hyph algebra}%
   {transitive representation of fibered Omega-algebra}%
\Index%491%T%transitive representation of group
   {transitive representation of group}%
   {transitive representation of group}%
\Index%490%T%transitive representation of Omega- algebra A
   {transitive representation of $\Omega$\Hyph algebra $A$}%
   {transitive representation of algebra}%
\Index%88%T%Tstar- linear composition of  vectors
   {\Ts linear composition of  vectors}%
   {linear composition of  vectors}%
\Index%87%T%Tstar- matrices vector space
   {\Ts matrices vector space}%
   {matrices vector space}%
\Index%91%T%Tstar- representation of fibered Omega- algebra
   {\Ts representation of fibered $\Omega$\Hyph algebra}%
   {Tstar representation of fibered Omega-algebra}%
\Index%90%T%Tstar- representation of Omega1 algebra A in Omega2 algebra M
   {\Ts representation of $\Omega_1$\Hyph algebra $A$ in $\Omega_2$\Hyph algebra $M$}%
   {Tstar representation of algebra}%
\Index%89%T%Tstar- shift
   {\Ts shift}%
   {Tstar shift}%
\Index%92%T%Tstar- transformation
   {\Ts transformation}%
   {Tstar transformation}%
\Index%93%T%Tstar- transformation on bundle
   {\Ts transformation on bundle}%
   {Tstar transformation of bundle}%
\Index%535%T%tuple of coordinates of element relative to tuple of sets
   {tuple of coordinates of element $\Vector a$ relative to tuple of sets $\VX X$}%
   {coordinates of element, tower of representations}%
\Index%713%T%tuple of equivalence generated by tower of representations
   {tuple of equivalence generated by tower of representations $\Vector f$}%
   {tuple of equivalence of tower of representations}%
\Index%534%T%tuple of generating sets of tower of representations
   {tuple of generating sets of tower of representations}%
   {tuple of generating sets of tower of representations}%
\Index%533%T%tuple of generating sets of tower subrepresentations
   {tuple of generating sets of tower subrepresentations}%
   {tuple of generating sets of tower subrepresentations}%
\Index%708%T%tuple of Omega- words of elements of tower of representations relative to tuple of generating sets
   {tuple of $\Vector{\Omega}$\Hyph words of element of tower of representations relative to tuple of generating sets}%
   {tuple of words relative to tuple of generating sets, tower of representations}%
\Index%530%T%tuple of stable sets of tower of representation
   {tuple of stable sets of tower of representation}%
   {tuple of stable sets of tower of representations}%
\Index%587%T%twin representations of associative algebra
   {twin representations of associative algebra}%
   {twin representations of associative algebra}%
\Index%684%T%twin representations of D- algebra
   {twin representations of $D$\Hyph algebra}%
   {twin representations of D algebra}%
\Index%320%T%twin representations of division ring
   {twin representations of division ring}%
   {twin representations of division ring}%
\Index%319%T%twin representations of fibered group
   {twin representations of fibered group}%
   {twin representations of fibered group}%
\Index%318%T%twin representations of group
   {twin representations of group}%
   {twin representations of group}%
\SetIndexSpace%
\Index%315%U%unit sphere in algebra
   {unit sphere in $D$\Hyph algebra}%
   {unit sphere in algebra}%
\Index%160%U%unit sphere in division ring
   {unit sphere in division ring}%
   {unit sphere in division ring}%
\Index%809%U%unit vector
   {unit vector}%
   {unit vector}%
\Index%172%U%unitarity law for  Dstar- vector fields
   {unitarity law for  $\mathcal D\star$\Hyph vector fields}%
   {unitarity law, Dstar vector fields}%
\Index%664%U%unitarity law for Astar- module
   {unitarity law for $A\star$\Hyph module}%
   {unitarity law, Astar module over algebra}%
\Index%578%U%unitarity law for Astar- vector space
   {unitarity law for $A\star$\Hyph vector space}%
   {unitarity law, Astar vector space}%
\Index%656%U%unitarity law for Dstar- module
   {unitarity law for $D\star$\Hyph module}%
   {unitarity law, Dstar module}%
\Index%171%U%unitarity law for Dstar- vector space
   {unitarity law for $D\star$\Hyph vector space}%
   {unitarity law, Dstar vector space}%
\Index%673%U%unitarity law for starA- module
   {unitarity law for $\star A$\Hyph module}%
   {unitarity law, starA module over algebra}%
\SetIndexSpace%
\Index%288%V%valued division ring
   {valued division ring}%
   {valued division ring}%
\Index%794%V%valued ring
   {valued ring}%
   {valued ring}%
\Index%125%V%vector bundle
   {vector bundle}%
   {vector bundle}%
\Index%787%V%vector module of algebra
   {vector module of algebra}%
   {vector module of algebra}%
\Index%776%V%vector module of ring
   {vector module of ring}%
   {vector module of ring}%
\Index%785%V%vector of element of algebra
   {vector of element of algebra}%
   {vector of algebra}%
\Index%774%V%vector of element of ring
   {vector of element of ring}%
   {vector of ring}%
\Index%782%V%vector of mapping
   {vector of mapping}%
   {vector of mapping}%
\Index%126%V%vector potential
   {vector potential}%
   {vector potential}%
\Index%727%V%vector space over field
   {vector space over field}%
   {vector space over field}%
\Index%485%V%vector space type
   {vector space type}%
   {vector space type}%

\CloseIndex

%% file: Symbol.English.tex
%auto-ignore
\def\indexname{Special Symbols and Notations}
\OpenIndex

\SetIndexSpace%A%0
\Symb%A/0
   {$(^a_b)$\hyph\CR quasideterminant}%
   {a b CR quasideterminant definition}%
\Symb%A/0
   {minor}%
   {A from b a}%
\Symb%A/0
   {minor}%
   {A from columns T}%
\Symb%A/0
   {minor}%
   {A from rows S}%
\Symb%A/0
   {minor}%
   {A without column a}%
\Symb%A/0
   {minor}%
   {A without columns T}%
\Symb%A/0
   {minor}%
   {A without row b}%
\Symb%A/0
   {minor}%
   {A without rows S}%
\Symb%A/0
   {$A\CRcirc$\Hyph linear combination}%
   {ACRcirc linear combination 1}%
\Symb%A/0
   {$A\CRcirc$\Hyph linear combination}%
   {ACRcirc linear combination 2}%
\Symb%A/0
   {active representation of group $G(f)$ in basis manifold $\mathcal B(f)$}%
   {active representation in basis manifold}%
\Symb%A/0
   {active representation of group $G(\Vector f)$ in basis manifold $\mathcal B(\Vector f)$}%
   {active representation in basis manifold, tower of representations}%
\Symb%A/0
   {affine space}%
   {affine space, division ring}%
\Symb%A/0
   {affine space}%
   {An}%
\Symb%A/0
   {associator of $R$\Hyph algebra}%
   {associator of algebra}%
\Symb%A/0
   {\subs row ($c$\hyph row) of matrix}%
   {c row}%
\Symb%A/0
   {commutator of $R$\Hyph algebra}%
   {commutator of algebra}%
\Symb%A/0
   {component of linear map $\Vector{A}$ of $D$\Hyph vector space}%
   {component of linear map, D vector space}%
\Symb%A/0
   {component $p$ of polylinear mapping $\Vector A$}%
   {component of polyadditive map, D vector space}%
\Symb%A/0
   {linear combination of vectors of $A$\Hyph module}%
   {CR linear combination in A module}%
\Symb%A/0
   {\CR power of element $A$ of biring}%
   {cr power}%
\Symb%A/0
   {\CR inverse element of biring}%
   {cr-inverse element}%
\Symb%A/0
   {\CR product of matrices}%
   {cr-product of matrices}%
\Symb%A/0
   {derivative of left shift}%
   {derivative of left shift}%
\Symb%A/0
   {derivative of left shift in $1$\Hyph parameter Lie group}%
   {derivative of left shift, 1-Parameter Group}%
\Symb%A/0
   {derivative of left shift in $1$\Hyph parameter Lie D group}%
   {derivative of left shift, 1-Parameter Group, algebra}%
\Symb%A/0
   {}%
   {derivative of right shift}%
\Symb%A/0
   {}%
   {derivative of right shift}%
\Symb%A/0
   {derivative of right shift in $1$\Hyph parameter Lie group}%
   {derivative of right shift, 1-Parameter Group}%
\Symb%A/0
   {derivative of right shift in $1$\Hyph parameter Lie D group}%
   {derivative of right shift, 1-Parameter Group, algebra}%
\Symb%A/0
   {derivative of left shift}%
   {derivative of Tstar shift}%
\Symb%A/0
   {\drc vector}%
   {drc vector}%
\Symb%A/0
   {hermitian conjugation in division ring}%
   {hermitian conjugation, division ring}%
\Symb%A/0
   {$(ji)$\hyph quasideterminant of matrix $\bfA$}%
   {j i quasideterminant definition}%
\Symb%A/0
   {$(^j_i)$\hyph \RC quasideterminant}%
   {j i RC-quasideterminant definition}%
\Symb%A/0
   {$(^j_i)$\hyph $\RCcirc$\Hyph quasideterminant}%
   {j i RCcirc-quasideterminant definition}%
\Symb%A/0
   {left shift in $D$\Hyph algebra}%
   {left shift, D algebra}%
\Symb%A/0
   {linear combination of vectors of $A$\Hyph module}%
   {linear combination in A module}%
\Symb%A/0
   {transformation of matrix}%
   {matrix, replacing its column}%
\Symb%A/0
   {transformation of matrix}%
   {matrix, replacing its row}%
\Symb%A/0
   {norm of map $\Vector A$ of normed $D$\hyph vector space}%
   {norm of map, D vector space}%
\Symb%A/0
   {opposite algebra to algebra $A$}%
   {opposite algebra}%
\Symb%A/0
   {orbit of linear mapping}%
   {orbit of linear mapping}%
\Symb%A/0
   {derivative}%
   {overline nabla_l, definition 2}%
\Symb%A/0
   {partial linear map of variable $v^i$}%
   {partial linear map of variable}%
\Symb%A/0
   {quasideterminant of matrix $\bfA$}%
   {quasideterminant definition}%
\Symb%A/0
   {\sups row ($r$\hyph row) of matrix}%
   {r row}%
\Symb%A/0
   {\RC power of element $A$ of biring}%
   {rc power}%
\Symb%A/0
   {\RC inverse element of biring}%
   {rc-inverse element}%
\Symb%A/0
   {\RC product of matrices}%
   {rc-product of matrices}%
\Symb%A/0
   {\RC quasideterminant}%
   {RC-quasideterminant definition}%
\Symb%A/0
   {$\RCcirc$\Hyph quasideterminant}%
   {RCcirc-quasideterminant definition}%
\Symb%A/0
   {\rcd vector}%
   {rcd vector}%
\Symb%A/0
   {right shift in $D$\Hyph algebra}%
   {right shift, D algebra}%
\Symb%A/0
   {coordinates of vector $\Vector a$ relative to Schauder basis}%
   {Schauder basis, coordinates}%
\Symb%A/0
   {set of $A$\Hyph linear mappings of module $\Vector V$ into module $\Vector W$}%
   {set A linear maps, module}%
\Symb%A/0
   {set of linear maps of division ring $D_1$ into division ring $D_2$}%
   {set linear maps, division ring}%
\Symb%A/0
   {set of polylinear mappings of rings $R_1$, ..., $R_n$ into module $S$}%
   {set polylinear maps, ring}%
\Symb%A/0
   {skew product of vectors $\Vector a_1$, ..., $\Vector a_m$}%
   {skew product of vectors}%
\Symb%A/0
   {standard component of tensor in tensor product of algebras}%
   {standard component of tensor, algebra}%
\Symb%A/0
   {right shift}%
   {starT shift}%
\Symb%A/0
   {\sT shift}%
   {starT shift, fibered group}%
\Symb%A/0
   {tensor power of algebra $A$}%
   {tensor power of algebra}%
\Symb%A/0
   {tensor product of algebras}%
   {tensor product of algebras}%
\Symb%A/0
   {left shift}%
   {Tstar shift}%
\Symb%A/0
   {\Ts shift}%
   {Tstar shift, fibered group}%
\Symb%A/0
   {anholonomic coordinates of vector}%
   {vector anholonomic coordinates}%
\Symb%A/0
   {holonomic coordinates of vector}%
   {vector holonomic coordinates}%

\SetIndexSpace%B%0
\Symb%B/0
   {basis manifold of \rcd vector space $\Vector V$}%
   {basis manifold of rcd vector space}%
\Symb%B/0
   {basis manifold of vector space}%
   {basis manifold of vector space}%
\Symb%B/0
   {basis manifold of representation $f$}%
   {basis manifold representation F algebra}%
\Symb%B/0
   {basis manifold of tower of representations $\Vector f$}%
   {basis manifold tower of representations}%
\Symb%B/0
   {basis manifold of affine space}%
   {Basis Manifold, Affine Space}%
\Symb%B/0
   {basis manifold of \rcd affine space}%
   {Basis Manifold, rcd Affine Space, division ring}%
\Symb%B/0
   {basis manifold of central affine space}%
   {BCAn}%
\Symb%B/0
   {basis manifold of Euclid space}%
   {BEn}%
\Symb%B/0
   {Cartesian power $\Bundle A$ of bundle $\Bundle B$}%
   {Cartesian power A of bundle B}%
\Symb%B/0
   {Cartesian power $A$ of set $B$}%
   {Cartesian power of set}%
\Symb%B/0
   {basis manifold of central affine space}%
   {FCAn}%
\Symb%B/0
   {basis manifold of Euclid space}%
   {FEn}%
\Symb%B/0
   {lattice of subrepresentations of representation $f$}%
   {lattice of subrepresentations}%
\Symb%B/0
   {lattice of towers of subrepresentations of tower of representations $\Vector f$}%
   {lattice of subrepresentations, tower of representations}%
\Symb%B/0
   {product of objects $B_1$, ..., $B_n$ in category $\mathcal A$}%
   {product of objects in category, 1 n}%
\Symb%B/0
   {structural constants of division ring $D$ over field $F$}%
   {structural constants of division ring over field}%
\Symb%B/0
   {tensor power of representation}%
   {tensor power of representation}%
\Symb%B/0
   {tensor product of representations}%
   {tensor product of representations}%

\SetIndexSpace%C%0
\Symb%C/0
   {central affine space}%
   {CAn}%
\Symb%C/0
   {central affine space}%
   {central affine space}%
\Symb%C/0
   {$j$th column determinant of matrix $\bfA$}%
   {column determinant}%
\Symb%C/0
   {$\CRcirc$\Hyph product of matrices of mappings}%
   {cr product of matrices of mappings}%
\Symb%C/0
   {left structural constant of Lie algebra}%
   {left structural constant of Lie algebra}%
\Symb%C/0
   {right structural constant of Lie algebra}%
   {right structural constant of Lie algebra}%
\Symb%C/0
   {structural constants of algebra $A$ over ring $D$}%
   {structural constants of algebra}%

\SetIndexSpace%D%0
\Symb%D/0
   {basis vector of representation of Lie group}%
   {basis vector of representation of Lie group}%
\Symb%D/0
   {basis vector of representation of Lie group over algebra $A$}%
   {basis vector of representation of Lie group over algebra A}%
\Symb%D/0
   {coordinates of basis vector of representation of Lie group over algebra $A$}%
   {basis vector of representation of Lie group over algebra A, coordinates}%
\Symb%D/0
   {coordinates of basis vector of representation of Lie group}%
   {basis vector of representation of Lie group, coordinates}%
\Symb%D/0
   {component of the G\^ateaux derivative of map $f(x)$ of algebra}%
   {component of Gateaux derivative of map, algebra}%
\Symb%D/0
   {component of the G\^ateaux derivative of map $\Vector f(\Vector x)$}%
   {component of Gateaux derivative of map, D vector space}%
\Symb%D/0
   {component of the G\^ateaux derivative of map $\Vector f(\Vector x)$}%
   {component of Gateaux derivative of map, D vector space, short}%
\Symb%D/0
   {component of the G\^ateaux differential of map $f(x)$}%
   {component of Gateaux derivative of map, division ring}%
\Symb%D/0
   {component of the G\^ateaux derivative of second order of map $f(x)$ of algebra}%
   {component of Gateaux derivative of Second Order, algebra}%
\Symb%D/0
   {component of the G\^ateaux derivative of second order of map $\Vector f(\Vector x)$}%
   {component of Gateaux derivative of Second Order, D vector space}%
\Symb%D/0
   {component of the G\^ateaux derivative of second order of map $f(x)$ of division ring}%
   {component of Gateaux derivative of Second Order, division ring}%
\Symb%D/0
   {conjugation in algebra}%
   {conjugation in algebra}%
\Symb%D/0
   {conjugation in ring}%
   {conjugation in ring}%
\Symb%D/0
   {coordinate \Drc vector bundle}%
   {coordinate drc vector bundle}%
\Symb%D/0
   {coordinate \rcd vector space}%
   {coordinate rcd vector space}%
\Symb%D/0
   {coordinate reference frame}%
   {coordinate reference frame, extensive definition}%
\Symb%D/0
   {diagonal in bundle $\Bundle A$}%
   {diagonal in bundle, 1}%
\Symb%D/0
   {direct product of division rings $D_1$, ..., $D_n$}%
   {direct product of division rings, 1 n}%
\Symb%D/0
   {double determinant of matrix $\bfA$}%
   {double determinant}%
\Symb%D/0
   {the Fr\'echet \Ds derivative of map $f$ of division ring}%
   {Frechet Dstar derivative of map, division ring}%
\Symb%D/0
   {the G\^ateaux \crd derivative of map $\Vector f$ of $D$\hyph vector space $\Vector V$ to $D$\hyph vector space $\Vector W$}%
   {Gateaux crd derivative of map, D vector space}%
\Symb%D/0
   {the G\^ateaux derivative of map $f$ of algebra}%
   {Gateaux derivative of map, algebra}%
\Symb%D/0
   {the G\^ateaux derivative of map $\Vector f$ of normed $D$\Hyph vector space $\Vector{V}$ to normed $D$\Hyph vector space $\Vector{W}$}%
   {Gateaux derivative of map, D vector space}%
\Symb%D/0
   {the G\^ateaux derivative of map $f$}%
   {Gateaux derivative of map, division ring}%
\Symb%D/0
   {the G\^ateaux derivative of map $f$ of algebra}%
   {Gateaux derivative of map, fraction, algebra}%
\Symb%D/0
   {the G\^ateaux derivative of map $f$}%
   {Gateaux derivative of map, fraction, division ring}%
\Symb%D/0
   {the G\^ateaux derivative of order $n$ of map $f$ of algebra}%
   {Gateaux derivative of Order n, algebra}%
\Symb%D/0
   {the G\^ateaux derivative of order $n$ of map $\Vector f$}%
   {Gateaux derivative of Order n, D vector space}%
\Symb%D/0
   {the G\^ateaux derivative of order $n$ of map $f$ of division ring}%
   {Gateaux derivative of Order n, division ring}%
\Symb%D/0
   {the G\^ateaux derivative of order $n$ of map $f$ of algebra}%
   {Gateaux derivative of Order n, fraction, algebra}%
\Symb%D/0
   {the G\^ateaux derivative of order $n$ of map $f$ of division ring}%
   {Gateaux derivative of Order n, fraction, division ring}%
\Symb%D/0
   {the G\^ateaux derivative of second order of mapping $f$ of algebra}%
   {Gateaux derivative of Second Order, algebra}%
\Symb%D/0
   {the G\^ateaux derivative of second order of map $\Vector f$}%
   {Gateaux derivative of Second Order, D vector space}%
\Symb%D/0
   {the G\^ateaux derivative of second order of map $f$ of division ring}%
   {Gateaux derivative of Second Order, division ring}%
\Symb%D/0
   {the G\^ateaux derivative of second order of mapping $f$ of algebra}%
   {Gateaux derivative of Second Order, fraction, algebra}%
\Symb%D/0
   {the G\^ateaux derivative of second order of map $f$ of division ring}%
   {Gateaux derivative of Second Order, fraction, division ring}%
\Symb%D/0
   {the G\^ateaux differential of mapping $f$ of algebra}%
   {Gateaux differential of map, algebra}%
\Symb%D/0
   {the G\^ateaux differential of map $\Vector f$ of normed $D$\Hyph vector space $\Vector{V}$ to normed $D$\Hyph vector space $\Vector{W}$}%
   {Gateaux differential of map, D vector space}%
\Symb%D/0
   {the G\^ateaux differential of map $f$}%
   {Gateaux differential of map, division ring}%
\Symb%D/0
   {the G\^ateaux differential of second order of mapping $f$ of algebra}%
   {Gateaux differential of Second Order, algebra}%
\Symb%D/0
   {the G\^ateaux differential of second order of mapping $\Vector f$}%
   {Gateaux differential of Second Order, D vector space}%
\Symb%D/0
   {the G\^ateaux differential of second order of mapping $f$ of division ring}%
   {Gateaux differential of Second Order, division ring}%
\Symb%D/0
   {the G\^ateaux \drc derivative of map $\Vector f$ of $D$\Hyph vector space $\Vector V$ to $D$\Hyph vector space $\Vector W$}%
   {Gateaux drc derivative of map, D vector space}%
\Symb%D/0
   {the G\^ateaux \Ds derivative of map $f$ of division ring $D$}%
   {Gateaux Dstar derivative of map, division ring}%
\Symb%D/0
   {the G\^ateaux Jacobian of map of $D$\Hyph vector space}%
   {Gateaux Jacobian of map, D vector space}%
\Symb%D/0
   {the G\^ateaux partial \drc derivative of map $f^b$ with respect to variable $v^a$}%
   {Gateaux partial crd derivative of map, 1, D vector space}%
\Symb%D/0
   {the G\^ateaux partial \drc derivative of map $f^b$ with respect to variable $v^a$}%
   {Gateaux partial crd derivative of map, 2, D vector space}%
\Symb%D/0
   {the G\^ateaux partial \drc derivative of map $f^b$ with respect to variable $v^a$}%
   {Gateaux partial crd derivative of map, 3, D vector space}%
\Symb%D/0
   {the G\^ateaux mixed partial derivative of map $f^j$ with respect to variables $v^i$, $v^j$}%
   {Gateaux partial derivative of Second Order, D vector space}%
\Symb%D/0
   {the G\^ateaux partial derivative of map $f^j$ with respect to variable $v^i$}%
   {Gateaux partial derivative, D vector space}%
\Symb%D/0
   {the G\^ateaux partial \drc derivative of map $f^b$ with respect to variable $v^a$}%
   {Gateaux partial drc derivative of map, 1, D vector space}%
\Symb%D/0
   {the G\^ateaux partial \drc derivative of map $f^b$ with respect to variable $v^a$}%
   {Gateaux partial drc derivative of map, 2, D vector space}%
\Symb%D/0
   {the G\^ateaux partial \drc derivative of map $f^b$ with respect to variable $v^a$}%
   {Gateaux partial drc derivative of map, 3, D vector space}%
\Symb%D/0
   {the G\^ateaux \sD derivative of map $f$ of division ring $D$}%
   {Gateaux starD derivative of map, division ring}%
\Symb%D/0
   {matrices vector space}%
   {matrices vector space}%
\Symb%D/0
   {Cartan derivative}%
   {overbrace D}%
\Symb%D/0
   {derivative}%
   {overline D}%
\Symb%D/0
   {derivative $e_{(k)}$}%
   {partial(k)}%
\Symb%D/0
   {\subs rows \drc vector space}%
   {r rows drc vector space}%
\Symb%D/0
   {speed of deviation}%
   {speed of deviation}%
\Symb%D/0
   {standard component of the G\^ateaux derivative of mapping $f$}%
   {standard component of Gateaux derivative, algebra}%
\Symb%D/0
   {standard component of the G\^ateaux differential of map $f$}%
   {standard component of Gateaux differential, division ring}%
\Symb%D/0
   {tensor product of division rings}%
   {tensor product of division rings}%
\Symb%D/0
   {tensor product of rings}%
   {tensor product of rings}%
\Symb%D/0
   {vector space type}%
   {vector space type}%

\SetIndexSpace%E%0
\Symb%E/0
   {$A\CRcirc$\Hyph basis for module}%
   {A CRcirc basis, module over algebra}%
\Symb%E/0
   {affine basis}%
   {Affine Basis}%
\Symb%E/0
   {basis of vector space}%
   {Basis e}%
\Symb%E/0
   {basis in vector space $\Vector V$}%
   {basis in V}%
\Symb%E/0
   {basis of vector space}%
   {basis, vector space}%
\Symb%E/0
   {basis of $(n)$\hyph vector space}%
   {basis,n vector space}%
\Symb%E/0
   {Cartesian power of total spaces}%
   {Cartesian power of total spaces}%
\Symb%E/0
   {Cartesian product of total spaces}%
   {Cartesian product of total spaces, definition 1}%
\Symb%E/0
   {central affine basis}%
   {Central Affine Basis}%
\Symb%E/0
   {basis for \Drc vector bundle}%
   {drc basis, vector bundle}%
\Symb%E/0
   {form of reference frame}%
   {dual forms, reference frame}%
\Symb%E/0
   {Euclid space}%
   {Euclid space}%
\Symb%E/0
   {Euclid space}%
   {Euclid space, division ring}%
\Symb%E/0
   {identical transformation of bundle}%
   {identical transformation of bundle}%
\Symb%E/0
   {linear automorphism of quaternioin algebra}%
   {mapping E, quaternion}%
\Symb%E/0
   {linear automorphism of quaternioin algebra}%
   {mapping E_1, quaternion}%
\Symb%E/0
   {linear automorphism of quaternioin algebra}%
   {mapping E_2, quaternion}%
\Symb%E/0
   {orthonornal basis}%
   {Orthonornal Basis}%
\Symb%E/0
   {pseudo Euclid space}%
   {pseudo Euclid space}%
\Symb%E/0
   {pseudo Euclid space}%
   {pseudo Euclid space, division ring}%
\Symb%E/0
   {quaternion algebra over the field $F$}%
   {quaternion algebra over the field}%
\Symb%E/0
   {quaternion division algebra over the field}%
   {quaternion division algebra over the fieldL}%
\Symb%E/0
   {$\RCcirc A$\Hyph linear combination}%
   {RCcircA linear combination 1}%
\Symb%E/0
   {$\RCcirc A$\Hyph linear combination}%
   {RCcircA linear combination 2}%
\Symb%E/0
   {\rcd affine basis}%
   {rcd affine basis, division ring}%
\Symb%E/0
   {reduced Cartesian product of total spaces}%
   {reduced Cartesian product of total spaces, definition 1}%
\Symb%E/0
   {Schauder basis}%
   {Schauder basis}%
\Symb%E/0
   {set of nonsingular \sT transformations of bundle $\Bundle E$}%
   {set of starT nonsingular transformations of bundle}%
\Symb%E/0
   {set of nonsingular \Ts transformations of bundle $\Bundle E$}%
   {set of Tstar nonsingular transformations of bundle}%
\Symb%E/0
   {standard coordinates of basis}%
   {standard coordinates of basis}%
\Symb%E/0
   {standard coordinates of reference frame}%
   {standard coordinates of reference frame}%
\Symb%E/0
   {vector field of reference frame}%
   {vector field of reference frame}%
\Symb%E/0
   {vector of basis}%
   {vector of basis}%

\SetIndexSpace%F%0
\Symb%F/0
   {coordinates of basis in \sups rows \rcd vector space}%
   {basis coordinates, c rows rcd vector space}%
\Symb%F/0
   {coordinates of basis in \subs rows \drc vector space}%
   {basis coordinates, r rows drc vector space}%
\Symb%F/0
   {basis for \sups rows \rcd vector space}%
   {basis, c rows rcd vector space}%
\Symb%F/0
   {basis for \subs rows \drc vector space}%
   {basis, r rows drc vector space}%
\Symb%F/0
   {central affine basis}%
   {Central Affine Basis, division ring}%
\Symb%F/0
   {component of linear map $f$ of division ring}%
   {component of linear map, division ring}%
\Symb%F/0
   {component of polylinear map into associative algebra}%
   {component of polylinear map, associative algebra}%
\Symb%F/0
   {component of polylinear map of division ring}%
   {component of polylinear map, division ring}%
\Symb%F/0
   {fibered morphism from bundle $\Bundle A$ into $\Bundle B$}%
   {fibered morphism from A into B}%
\Symb%F/0
   {filter $\mathfrak{F}$ converges to set $A$}%
   {filter converges}%
\Symb%F/0
   {homomorphism of fibered universal algebras}%
   {homomorphism of fibered universal algebras}%
\Symb%F/0
   {inverse fibered correspondence}%
   {inverse fibered correspondence, 1}%
\Symb%F/0
   {inverse reduced fibered correspondence}%
   {inverse reduced fibered correspondence, 1}%
\Symb%F/0
   {map to Cartesian product}%
   {map to Cartesian product}%
\Symb%F/0
   {norm of mapping into $D$\Hyph algebra}%
   {norm of map, algebra}%
\Symb%F/0
   {norm of map $f$ of division ring}%
   {norm of map, division ring}%
\Symb%F/0
   {representation orbit of group $G$}%
   {orbit of Tstar representation of group}%
\Symb%F/0
   {orthonornal basis}%
   {Orthonornal Basis, division ring}%
\Symb%F/0
   {quaternion algebra  over field ${\rm {\mathbb{F}}}$}%
   {quaternion algebra F a b}%
\Symb%F/0
   {reference frame}%
   {reference frame}%
\Symb%F/0
   {reference frame, extensive definition}%
   {reference frame, extensive definition}%
\Symb%F/0
   {standard component of biadditive map $f$ over field $F$}%
   {standard component of biadditive map, division ring}%
\Symb%F/0
   {standard $F$\Hyph component of linear mapping $f$}%
   {standard component of linear map, division ring}%
\Symb%F/0
   {standard component of polylinear mapping into associative algebra}%
   {standard component of polylinear map, associative algebra}%
\Symb%F/0
   {standard component of polylinear map $f$ of division ring}%
   {standard component of polylinear map, division ring}%
\Symb%F/0
   {standard component of quadratic map $f$ over field $F$}%
   {standard component of quadratic map, division ring}%
\Symb%F/0
   {standard component of tensor}%
   {standard component of tensor, division ring}%

\SetIndexSpace%G%0
\Symb%G/0
   {affine transformation group}%
   {affine transformation group}%
\Symb%G/0
   {\CR matrix group}%
   {cr-matrix group}%
\Symb%G/0
   {affine transformation group}%
   {drc affine transformation group}%
\Symb%G/0
   {fibered little group of section $h$}%
   {fibered little group}%
\Symb%G/0
   {fibered stability group of section $h$}%
   {fibered stability group}%
\Symb%G/0
   {algebra Lie of group Lie}%
   {g}%
\Symb%G/0
   {left defined algebra Lie of group Lie}%
   {gl}%
\Symb%G/0
   {right defined algebra Lie of group Lie}%
   {gr}%
\Symb%G/0
   {group of automorphisms of representation $f$}%
   {group of automorphisms of representation}%
\Symb%G/0
   {group of homomorphisms of vector space $\Vector V$}%
   {GV}%
\Symb%G/0
   {little group of $x$}%
   {little group}%
\Symb%G/0
   {orbit of effective covariant \sT representation of fibered group}%
   {orbit of effective covariant starT representation of fibered group}%
\Symb%G/0
   {orbit of effective covariant \sT representation of group}%
   {orbit of effective covariant starT representation of group}%
\Symb%G/0
   {orbit of effective covariant \Ts representation of fibered group}%
   {orbit of effective covariant Tstar representation of fibered group}%
\Symb%G/0
   {orbit of effective covariant \Ts representation of group}%
   {orbit of effective covariant Tstar representation of group}%
\Symb%G/0
   {product of groups $G_1$, ..., $G_n$}%
   {product of groups, 1 n}%
\Symb%G/0
   {\RC matrix group}%
   {rc-matrix group}%
\Symb%G/0
   {stability group of $x$}%
   {stability group}%

\SetIndexSpace%H%0
\Symb%H/0
   {Hadamard inverse of matrix}%
   {Hadamard inverse of matrix}%
\Symb%H/0
   {quaternion algebra}%
   {quaternion algebra H a b}%
\Symb%H/0
   {quaternion algebra over real field}%
   {quaternion algebra over real field}%

\SetIndexSpace%I%0
\Symb%I/0
   {infinitesimal generator of representation}%
   {infinitesimal generator I of representation}%
\Symb%I/0
   {infinitesimal generator of representation}%
   {infinitesimal generator i of representation}%
\Symb%I/0
   {Lie group infinitesimal generators}%
   {Lie group infinitesimal generators}%
\Symb%I/0
   {vector module of algebra $A$}%
   {vector module of algebra}%
\Symb%I/0
   {vector module of ring $D$}%
   {vector module of ring}%
\Symb%I/0
   {vector of element $d$ of algebra}%
   {vector of algebra}%
\Symb%I/0
   {vector of mapping $f$}%
   {vector of mapping}%
\Symb%I/0
   {vector of element $d$ of ring}%
   {vector of ring}%

\SetIndexSpace%J%0
\Symb%J/0
   {Jacobian matrix of left shift}%
   {aE, quaternion, Jacobian matrix}%
\Symb%J/0
   {closure operator of representation $f$}%
   {closure operator, representation}%
\Symb%J/0
   {closure operator of tower of representations $\Vector f$}%
   {closure operator, tower of representations}%
\Symb%J/0
   {Jacobian matrix of right shift}%
   {Ea, quaternion, Jacobian matrix}%
\Symb%J/0
   {tower of subrepresentations of tower of representations $\Vector f$ generated by tuple of sets $\VX X$}%
   {subrepresentation generated by tuple of sets}%

\SetIndexSpace%K%0
\Symb%K/0
   {kernel of linear mapping of $D$\Hyph vector space}%
   {kernel of linear map, D vector space}%
\Symb%K/0
   {kernel of linear mapping of division ring}%
   {kernel of linear map, division ring}%

\SetIndexSpace%L%0
\Symb%L/0
   {left $ij$th cofactor of entry of matrix}%
   {left cofactor, matrix}%
\Symb%L/0
   {left double $ij$th cofactor of entry of matrix}%
   {left double cofactor}%
\Symb%L/0
   {left shift}%
   {left shift}%
\Symb%L/0
   {Lie derivative of connection}%
   {Lie derivative of connection}%
\Symb%L/0
   {Lie derivative of metric}%
   {Lie derivative of metric}%
\Symb%L/0
   {limit of correspondence $\Phi$ with respect to the filter $\mathfrak{F}$}%
   {limit of correspondence with respect to the filter}%
\Symb%L/0
   {limit of sequence in normed ring}%
   {limit of sequence, normed ring}%
\Symb%L/0
   {limit of sequence in valued division ring}%
   {limit of sequence, valued division ring}%
\Symb%L/0
   {limit of sequence in valued ring}%
   {limit of sequence, valued ring}%
\Symb%L/0
   {passive transformation}%
   {passive transformation}%
\Symb%L/0
   {set of \Acr linear mappings of module $\Vector V$ into module $\Vector W$}%
   {set Acr linear maps, module}%
\Symb%L/0
   {\rcd vector space of \drc linear maps of \drc vector space $\Vector V$ into \drc vector space $\Vector W$}%
   {set drc linear maps, drc vector space}%
\Symb%L/0
   {set of linear mappings of algebra $A_1$ into algebra $A_2$}%
   {set linear mappings, algebra}%
\Symb%L/0
   {set of linear mappings of $D$\Hyph vector space $\Vector{V}$ into $D$\Hyph vector space $\Vector{W}$}%
   {set linear maps, D vector space}%
\Symb%L/0
   {set of left-side nonsingular transformations of set $M$}%
   {set of left-side nonsingular transformations}%
\Symb%L/0
   {set of $n$\hyph linear mappings of algebra $A$ into module $S$}%
   {set polylinear mappings An, algebra}%
\Symb%L/0
   {set of polylinear maps of algebras $A_1$, ..., $A_n$ into module $S$}%
   {set polylinear mappings, algebra}%
\Symb%L/0
   {set of polylinear mappings}%
   {set polylinear maps, D vector space}%
\Symb%L/0
   {set of polylinear maps of algebras $A_1$, ..., $A_n$ into algebra $A$}%
   {set polylinear maps, finite dimensional algebra}%
\Symb%L/0
   {\drc vector space of \rcd linear maps of \rcd vector space $\Vector{V}$ into \rcd vector space $\Vector{W}$}%
   {set rcd linear maps, rcd vector space}%
\Symb%L/0
   {set of \sT representations of division ring $S$ in additive group of division ring $R$}%
   {set sT representations, division ring}%
\Symb%L/0
   {set of \Ts representations of division ring $S$ in additive group of division ring $R$}%
   {set Ts representations, division ring}%

\SetIndexSpace%M%0
\Symb%M/0
   {set of \sT transformations of set $M$}%
   {set of starT transformations}%
\Symb%M/0
   {set of transformations of set $M$}%
   {set of transformations}%
\Symb%M/0
   {set of \Ts transformations of set $M$}%
   {set of Tstar transformations}%
\Symb%M/0
   {space of orbits of effective \sT covariant representation of the group}%
   {space of orbits of effective sT representation}%
\Symb%M/0
   {space of orbits of effective \Ts covariant representation of the group}%
   {space of orbits of effective Ts representation}%
\Symb%M/0
   {space of orbits of \Ts representation $f$ of group $G$ in set $M$}%
   {space of orbits of Ts representation}%

\SetIndexSpace%N%0
\Symb%N/0
   {norm of quaternion $x$}%
   {norm, quaternion algebra}%
\Symb%N/0
   {nucleus of $R$\Hyph algebra $A$}%
   {nucleus of algebra}%

\SetIndexSpace%O%0
\Symb%O/0
   {geometric object in coordinate representation defined in \rcd vector space}%
   {geometric object, coordinate rcd vector space}%
\Symb%O/0
   {geometric object in coordinate representation}%
   {geometric object, coordinate vector space}%
\Symb%O/0
   {geometric object defined in \rcd vector space}%
   {geometric object, rcd vector space}%
\Symb%O/0
   {octonion algebra}%
   {octonion algebra}%
\Symb%O/0
   {orbit of representation of fibered group $\Bundle G$}%
   {orbit of representation of fibered group}%
\Symb%O/0
   {orbit of representation of the group $G$}%
   {orbit of representation of group}%

\SetIndexSpace%P%0
\Symb%P/0
   {bundle}%
   {bundle}%
\Symb%P/0
   {bundle of level $2$}%
   {bundle of level 2}%
\Symb%P/0
   {bundle of level $n$}%
   {bundle of level n}%
\Symb%P/0
   {Cartesian power $n$ of bundle $\bundle{}{p}{E}{}$}%
   {Cartesian power of bundle}%
\Symb%P/0
   {Cartesian product of bundles}%
   {Cartesian product of bundles, definition 1}%
\Symb%P/0
   {passive representation of group $G(f)$ in basis manifold $\mathcal B(f)$}%
   {passive representation in basis manifold}%
\Symb%P/0
   {passive representation of group $G(\Vector f)$ in basis manifold $\mathcal B(\Vector f)$}%
   {passive representation in basis manifold, tower of representations}%
\Symb%P/0
   {reduced Cartesian product of bundles}%
   {reduced Cartesian product of bundles, definition 1}%
\Symb%P/0
   {set of nonsingular \sT transformations of bundle $\bundle{}pE{}$}%
   {set of starT nonsingular transformations of bundle, projection}%
\Symb%P/0
   {set of nonsingular \Ts transformations of bundle $\bundle{}pE{}$}%
   {set of Tstar nonsingular transformations of bundle, projection}%

\SetIndexSpace%R%0
\Symb%R/0
   {active transformation}%
   {active transformation}%
\Symb%R/0
   {\sups rows \rcd vector space}%
   {c rows rcd vector space}%
\Symb%R/0
   {Cartan curvature}%
   {Cartan curvature}%
\Symb%R/0
   {\CR rank of matrix}%
   {cr-rank of matrix}%
\Symb%R/0
   {diagonal in bundle  $\bundle{}pA{}$}%
   {diagonal in bundle, 2}%
\Symb%R/0
   {diagonal in bundle $\Bundle A$}%
   {diagonal in reduced bundle, 2}%
\Symb%R/0
   {\Ds component of coordinates of vector $\Vector r$}%
   {Dstar component of coordinates of vector, D vector space}%
\Symb%R/0
   {image of $m$ under endomorphism $R$ of effective representation}%
   {endomorphism image, effective representation}%
\Symb%R/0
   {image of tuple $\VX a$ under endomorphism $\VX r$ of tower of effective representations}%
   {endomorphism image, tower of effective representations}%
\Symb%R/0
   {curvature}%
   {GLn curvature_overline}%
\Symb%R/0
   {$\RCcirc$\Hyph product of matrices of mappings}%
   {rc product of matrices of mappings}%
\Symb%R/0
   {\RC rank of matrix}%
   {rc-rank of matrix}%
\Symb%R/0
   {right $ij$th cofactor of entry of matrix}%
   {right cofactor, matrix}%
\Symb%R/0
   {right double $ij$th cofactor of entry of matrix}%
   {right double cofactor}%
\Symb%R/0
   {right shift}%
   {right shift}%
\Symb%R/0
   {$i$th row determinant of matrix $\bfA$}%
   {row determinant}%
\Symb%R/0
   {scalar algebra of algebra $A$}%
   {scalar algebra of algebra}%
\Symb%R/0
   {scalar algebra of ring $D$}%
   {scalar algebra of ring}%
\Symb%R/0
   {scalar of element $d$ of algebra}%
   {scalar of algebra}%
\Symb%R/0
   {scalar of mapping $f$}%
   {scalar of mapping}%
\Symb%R/0
   {scalar of element $d$ of ring}%
   {scalar of ring}%
\Symb%R/0
   {set of right-side nonsingular transformations of set $M$}%
   {set of right-side nonsingular transformations}%
\Symb%R/0
   {\sD component of coordinates of vector $\Vector r$}%
   {starD component of coordinates of vector, D vector space}%

\SetIndexSpace%S%0
\Symb%S/0
   {composition of fibered correspondences}%
   {composition of fibered correspondences}%
\Symb%S/0
   {inverse fibered correspondence}%
   {inverse fibered correspondence, 2}%
\Symb%S/0
   {inverse reduced fibered correspondence}%
   {inverse reduced fibered correspondence, 2}%
\Symb%S/0
   {linear span in vector space}%
   {linear span, vector space}%
\Symb%S/0
   {image of basis $X$ under passive transformation $S$}%
   {passive transformation of basis, representation}%
\Symb%S/0
   {image of basis $\VX  X$ under passive transformation $\VX s$}%
   {passive transformation of basis, tower of representations}%
\Symb%S/0
   {symmetric group}%
   {symmetric group}%

\SetIndexSpace%T%0
\Symb%T/0
   {category of \Ts representations of $\Omega_1$\Hyph algebra $A$}%
   {category of Tstar representations of Omega1 algebra}%
\Symb%T/0
   {category of \Ts representations of $\Omega_1$\Hyph algebra from category $\mathcal A$}%
   {category of Tstar representations of Omega1 algebra from category}%
\Symb%T/0
   {tangent plane to group $G$}%
   {TaG}%
\Symb%T/0
   {trace of quaternion $x$}%
   {trace, quaternion algebra}%

\SetIndexSpace%V%0
\Symb%V/0
   {coordinate vector space}%
   {coordinate vector space}%
\Symb%V/0
   {coordinates in vector space}%
   {coordinates in vector space}%
\Symb%V/0
   {direct product of $\RCstar D_i$\hyph vector spaces $\Vector V_1$, ..., $\Vector V_n$}%
   {direct product, rcd vector space, 1 n}%
\Symb%V/0
   {dual space of \rcd vector space $\Vector V$}%
   {dual space of rcd vector space}%
\Symb%V/0
   {hermitian conjugated vector}%
   {hermitian conjugated vector}%
\Symb%V/0
   {\dcr vector space}%
   {left CR vector space}%
\Symb%V/0
   {\drc vector space}%
   {left RC vector space}%
\Symb%V/0
   {\crd vector space}%
   {right CR vector space}%
\Symb%V/0
   {\rcd vector space}%
   {right RC vector space}%
\Symb%V/0
   {tensor product of $D$\Hyph vector spaces}%
   {tensor product of D vector spaces}%
\Symb%V/0
   {tensor product of \Ds vector spaces}%
   {tensor product of Dstar vector spaces}%
\Symb%V/0
   {vector space}%
   {V}%

\SetIndexSpace%W%0
\Symb%W/0
   {set of coordinates of representation $J(f,X)$}%
   {coordinate set of representation}%
\Symb%W/0
   {set of tuples of coordinates of tower of representations $\Vector J(\Vector f,\VX X)$}%
   {coordinate set of tower of representations}%
\Symb%W/0
   {coordinates of basis $X'$ relative to basis $X$ of representation}%
   {coordinates of basis relative to basis, representation}%
\Symb%W/0
   {coordinates of element $m$ of representation $f$ relative to set $X$}%
   {coordinates of element relative to generating set, representation}%
\Symb%W/0
   {coordinates of element $m$ relative to set $X$}%
   {coordinates of element relative to set, representation}%
\Symb%W/0
   {tuple of coordinates of element $\Vector a*$ relative to tuple of sets $\VX X$}%
   {coordinates of element, tower of representations}%
\Symb%W/0
   {geometric object in coordinate representation defined in $\Omega_2$\Hyph algebra $M$}%
   {geometric object, coordinate representation g}%
\Symb%W/0
   {geometric object in coordinate representation defined in tuple of $\VX\Omega$\Hyph algebras $\VX A$}%
   {geometric object, coordinate tower of representations g}%
\Symb%W/0
   {geometric object defined in $\Omega_2$\Hyph algebra $M$}%
   {geometric object, representation g}%
\Symb%W/0
   {geometric object defined in tuple of $\VX\Omega$\Hyph algebras $\VX A$}%
   {geometric object, tower of representations g}%
\Symb%W/0
   {geometric object in vector space}%
   {geometric object, vector space}%
\Symb%W/0
   {set of coordinates of set $B\subset J(f,X)$}%
   {subset of coordinates of representation}%
\Symb%W/0
   {coordinates of tuple of sets $\VX B$ relative to tuple of sets $\VX X$}%
   {subset of coordinates of tower of representations}%
\Symb%W/0
   {coordinates of set $B_k$ relative to tuple of sets $\VX X$}%
   {subset of coordinates of tower of representations, k}%
\Symb%W/0
   {set of $\Omega_2$\Hyph words representing set $B\subset J(f,X)$}%
   {subset of words of representation}%
\Symb%W/0
   {superposition of coordinates of the representation $f$ and the element $m$}%
   {superposition of coordinates, representation}%
\Symb%W/0
   {superposition of coordinates of the tower of representations $\Vector f$ and the element $\VX a$}%
   {superposition of coordinates, tower of representations}%
\Symb%W/0
   {$\Omega_2$\Hyph word representing element $m\in J(f,X)$}%
   {word of element relative to generating set, representation}%
\Symb%W/0
   {set of $\Omega_2$\Hyph words of representation $J(f,X)$}%
   {word set of representation}%
\Symb%W/0
   {set of tuples of $\VX{\Omega}$\Hyph words of tower of representations $\Vector J(\Vector f,\VX X)$}%
   {word set of tower of representations}%
\Symb%W/0
   {tuple of words of element $\Vector a*$ relative to tuple of sets $\VX X$}%
   {words of element, tower of representations}%

\SetIndexSpace%X%0
\Symb%X/0
   {conjugate of quaternion $x$}%
   {conjugate of quaternion}%
\Symb%X/0
   {local basis of affine space}%
   {local basis of affine space}%
\Symb%X/0
   {anholonomic coordinate}%
   {x(k)}%

\SetIndexSpace%Z%0
\Symb%Z/0
   {center of an $R$\Hyph algebra $A$}%
   {center of algebra}%
\Symb%Z/0
   {center of ring $D$}%
   {center of ring}%

\SetIndexSpace%Delta%1
\Symb%Delta/1
   {deviation of trajectories}%
   {deviation of trajectories}%
\Symb%Delta/1
   {identical transformation}%
   {identical transformation}%
\Symb%Delta/1
   {image of vector $\Vector e_k\in\Basis e$ under isomorphism to coordinate vector space}%
   {image of vector e_k, coordinate vector space}%
\Symb%Delta/1
   {Kronecker symbol}%
   {Kronecker symbol}%

\SetIndexSpace%Gamma%1
\Symb%Gamma/1
   {anholonomic coordinates of connection}%
   {anholonomic coordinates of connection}%
\Symb%Gamma/1
   {Cartan symbol}%
   {Cartan symbol}%
\Symb%Gamma/1
   {connection}%
   {conection overline}%
\Symb%Gamma/1
   {connection coefficients in $D$\Hyph affine space}%
   {connection coefficients, D affine space}%
\Symb%Gamma/1
   {connection in $D$\Hyph affine manifold}%
   {connection, affine manifold}%
\Symb%Gamma/1
   {$D$\Hyph affine connection coefficients on manifold}%
   {D affine connection coefficients, manifold}%
\Symb%Gamma/1
   {holonomic coordinates of connection}%
   {holonomic coordinates of connection}%
\Symb%Gamma/1
   {Cartan connection}%
   {overbrace Gamma i kl}%
\Symb%Gamma/1
   {set of sections of bundle}%
   {set of sections of bundle}%

\SetIndexSpace%Lambda%1
\Symb%Lambda/1
   {inverse operator to operator $\psi_l$}%
   {inverse operator to operator psi l}%
\Symb%Lambda/1
   {inverse operator to operator $\psi_r$}%
   {inverse operator to operator psi r}%

\SetIndexSpace%Omega%1
\Symb%Omega/1
   {anholonomity object}%
   {anholonomity object}%

\SetIndexSpace%Psi%1
\Symb%Psi/1
   {left basic operator of Lie group over algebra $A$}%
   {L basic operator of Lie group over algebra A}%
\Symb%Psi/1
   {left basic operator of group Lie}%
   {Lie Basic Operator L}%
\Symb%Psi/1
   {left basic operator of Lie 1-parameter group}%
   {Lie Basic Operator L, 1-Parameter Group}%
\Symb%Psi/1
   {left basic operator of Lie 1-parameter group over algebra $A$}%
   {Lie Basic Operator L, 1-Parameter Group, algebra}%
\Symb%Psi/1
   {right basic operator of group Lie}%
   {Lie Basic Operator R}%
\Symb%Psi/1
   {right basic operator of Lie 1-parameter group}%
   {Lie Basic Operator R, 1-Parameter Group}%
\Symb%Psi/1
   {right basic operator of Lie 1-parameter group over algebra $A$}%
   {Lie Basic Operator R, 1-Parameter Group, algebra}%
\Symb%Psi/1
   {right basic operator of Lie group over algebra $A$}%
   {R basic operator of Lie group over algebra A}%

\SetIndexSpace%Fi%2
\Symb%Fi/2
   {Lie group composition law}%
   {Lie group composition law}%

\SetIndexSpace%Nabla%2
\Symb%Nabla/2
   {Cartan derivative}%
   {overbrace nabla_l}%
\Symb%Nabla/2
   {derivative}%
   {overline nabla_l, definition 1}%

\SetIndexSpace%Phi%2
\Symb%Phi/2
   {restriction of correspondence $\Phi$ to set $C$}%
   {restriction of correspondence}%

\SetIndexSpace%Pi%2
\Symb%Pi/2
   {Cartesian product of bundles}%
   {Cartesian product of bundles, definition 2}%
\Symb%Pi/2
   {Cartesian product of total spaces}%
   {Cartesian product of total spaces, definition 2}%
\Symb%Pi/2
   {direct product of division rings $D_i$, $i\in I$}%
   {direct product of division rings}%
\Symb%Pi/2
   {direct product of division rings $D_1$, ..., $D_n$}%
   {direct product of division rings, i 1 n}%
\Symb%Pi/2
   {direct product of $\RCstar D_i$\hyph vector spaces $\Vector V_i$, $i\in I$}%
   {direct product, rcd vector space}%
\Symb%Pi/2
   {direct product of $\RCstar D_i$\hyph vector spaces}%
   {direct product, rcd vector space, i 1 n}%
\Symb%Pi/2
   {product of groups $G_i$, $i\in I$}%
   {product of groups}%
\Symb%Pi/2
   {product of groups $G_1$, ..., $G_n$}%
   {product of groups, i 1 n}%
\Symb%Pi/2
   {product of objects $\{B_i,i\in I\}$ in category $\mathcal A$}%
   {product of objects in category}%
\Symb%Pi/2
   {product of objects $B_1$, ..., $B_n$ in category $\mathcal A$}%
   {product of objects in category, i 1 n}%
\Symb%Pi/2
   {reduced Cartesian product of bundles}%
   {reduced Cartesian product of bundles, definition 2}%
\Symb%Pi/2
   {reduced Cartesian product of total spaces}%
   {reduced Cartesian product of total spaces, definition 2}%

\SetIndexSpace%S%2
\Symb%S/2
   {fibered subset}%
   {fibered subset}%
\Symb%S/2
   {subbundle}%
   {subbundle}%

\CloseIndex

%% file: Abstract.1006.2597.English.tex
%auto-ignore
%\documentclass{amsart}
%\author{Aleks Kleyn}
%\begin{document}
%\title{The G\^ateaux Derivative and Integral over Division Ring}
%\maketitle
\def\P{\partial}
\def\S{\sum_{n=0}^{\infty}(n!)^{-1}}
\def\G{G\^ateaux\ }
\def\M{$f:A\rightarrow A\ $}
Let $A$ be algebra over commutative ring $D$. Map
\M
is linear if for any $a,b\in A$ and any $c\in D$
\[
f\circ (a+b)=f\circ a+f\circ b
\]
\[
f\circ(ca)=c\ f\circ a
\]
Map
\M
is called differentiable in the \G sense,if
\[f(x+a)-f(x)=\P f(x)\circ a+o(a)\]
where
the \G derivative
$\P f(x)$
of map $f$
is linear map of increment $a$ and
$o$ is such continuous map that
\[
\lim_{a\rightarrow 0}\frac{|o(a)|}{|a|}=0
\]
For instance
\[
\P (x^2)\circ h=xh+hx
\]
\[
\P (x^{-1})\circ h=-x^{-1}hx^{-1}
\]

Assuming that we defined the \G derivative
$\P^{n-1} f(x)$ of order $n-1$,we define
\[
\P^n f(x)\circ(a_1\otimes...\otimes a_n)
=\P(\P^{n-1} f(x)\circ(a_1\otimes...\otimes a_{n-1}))\circ a_n
\]
the \G derivative of order $n$ of map $f$.
When $h_1=...=h_n=h$, we assume
\[
\P^nf(x)\circ h=\P^nf(x)\circ(h_1\otimes...\otimes h_n)
\]

Function $f(x)$ has Taylor series expansion
\[
f(x)=\S\P^n f(x_0)\circ(x-x_0)
\]

Differential equation over division ring
\[
\P (y)\circ h=hx^2+xhx+x^2h
\]
\[
y(0)=0
\]
has solution
\[
y=x^3
\]

The solution of differential equation
\[
\P(y)\circ h=\frac 12(yh+hy)
\]
\[y(0)=1\]
is exponent
$y=e^x$
that has following Taylor series expansion
\[
e^x=\S x^n
\]
The equation
\[
e^{a+b}=e^ae^b
\]
is true iff
$ab=ba$
%\end{document}